\documentclass[10pt]{article}

\usepackage{authblk}
\usepackage[toc, page]{appendix}
\usepackage{amsthm}
\usepackage{lipsum}
\usepackage{amsfonts}
\usepackage{graphicx}
\usepackage{epstopdf}
\usepackage{amssymb}
\usepackage{mathtools}
\usepackage{bm}
\usepackage{comment}
\usepackage{cases}
\usepackage{ascmac}
\usepackage{algorithm}
\usepackage{algpseudocode}
\usepackage{caption}
\usepackage{subcaption}
\usepackage{listings}
\usepackage{mathtools}
\usepackage{setspace}  

\usepackage{enumitem}
\setlist[enumerate]{leftmargin=.5in}
\setlist[itemize]{leftmargin=.5in}

\newtheorem{remark}{Remark}
\newtheorem{ex}{Example}

\newtheorem{prob}{Problem}
\newtheorem{prop}{Proposition}

\newcommand{\tr}{\mathsf{T}}

\newcommand{\takespace}{\vspace{7mm}}

\usepackage{color}
\definecolor{gr}{RGB}{0, 150, 0}
\usepackage{hyperref}
\hypersetup{
	setpagesize=false,
	bookmarksnumbered=true,
	bookmarksopen=true,
	colorlinks=true,
	linkcolor=gr,
	citecolor=gr,
}

\usepackage[a4paper, total={6in, 8in}]{geometry}

\title{
	\begin{center}
		{Extracting Koopman Operators for Prediction and Control\\ of Non-linear Dynamics Using Two-stage Learning\\ and Oblique Projections}
	\end{center}
}
\author[1]{Daisuke Uchida}
\author[1]{Karthik Duraisamy}
\affil[1]{Department of Aerospace Engineering, University of Michigan}
\date{}

\begin{document}

\maketitle

% REQUIRED
\begin{abstract}
The Koopman operator framework provides a perspective that non-linear dynamics can be described through the lens of linear operators acting on function spaces. As the framework naturally yields linear embedding models, there have been extensive efforts to utilize it for control, where linear controller designs 
%in the embedded space 
can be applied to control possibly nonlinear dynamics. However, it is challenging to successfully deploy this modeling procedure in a wide range of applications. In this work, some of the fundamental limitations of linear embedding models are addressed. We show a necessary condition for a linear embedding model to achieve zero modeling error, highlighting a trade-off relation between the model expressivity and a restriction on the model structure to allow the use of linear systems theories for nonlinear dynamics. 
To achieve good performance despite this trade-off,
neural network-based modeling is proposed based on linear embedding with oblique projection, which is derived from a weak formulation of projection-based linear operator learning. We train the proposed model using a two-stage learning procedure, wherein the features and operators are initialized with orthogonal projection, followed by the main training process in which test functions characterizing the oblique projection are learned from data. The first stage achieves an optimality ensured by the orthogonal projection and the second stage improves the generalizability to various tasks by optimizing the model with the oblique projection. We demonstrate the effectiveness of the proposed method over other data-driven modeling methods by providing comprehensive numerical evaluations where four tasks are considered targeting three different systems.
\end{abstract}

\section{Introduction}
\label{sec. intro}
Data-driven approaches\cite{Data_driven_book} to the modeling and control of dynamical systems have been gaining attention and popularity in recent years.  Among many different data-driven modeling frameworks, the Koopman operator is emerging as a promising tool
that can represent a nonlinear dynamical system as a linear one in a higher (possibly infinite) dimensional embedded space.
In addition to theoretical works concerning the development of computational methods to estimate the Koopman operator\cite{applied_Koopmanism,Koopman_book,Otto_Koopman_review,modern_Koopman_theory}, a number of applications in various fields have been also reported in the literature in recent years\cite{survey_Koopman_vehicle,building_energy_modeling,application_Koopman_power_system,neural_recording_DMD}. 
While the numerical estimation of the spectral properties of the Koopman operator has been one of the central topics (e.g., \cite{Mezic_2005,data_driven_spectral_analysis_of_Koopman, DAS202175, rigorous_data_driven_computation}), there have been also many efforts that extend the Koopman operator framework to systems with control inputs\cite{DMDc,EDMD_actuated_system,Koopman_control_PLOS_ONE}.

After a formal extension of the Koopman operator formalism to general non-autonomous systems was established in \cite{KORDA_Koopman_MPC},
various types of controller design problems
have been addressed.
The most notable feature of this framework is that linear control theories are utilized to control possibly nonlinear dynamics, whose governing equations are even unknown.
Linear Quadratic Regulator (LQR) is one of the simplest designs that may be adopted in this scenario and several studies showed its applicability and effectiveness in both simulations and real experiments targeting robotic systems\cite{Local_Koopman_operator_robot_control,Derivarive_based_Koopman_robot_control,Mezic_group_soft_robot_control}.
Linear Model Predictive Control (MPC) is currently the most popular choice of Koopman-based controller design and a number of papers address linear MPC design problems to control nonlinear dynamics based on this framework\cite{KORDA_Koopman_MPC, Koopman_Lyapunov_based_MPC,MPC_for_PDE,application_of_Koopman_pulp,handling_plant_model_mismatch,hybrid_Koopman_pulp,tube_based_MPC}.
Several other works combine the data-driven Koopman operator framework with different tools, e.g., eigenfunctions of the Koopman operator\cite{Kaiser_data_driven_control_eigen}, control barrier functions\cite{Koopman_CBF}, and switching time optimization\cite{PEITZ_switched_control}.

While many Koopman-based controller designs have been developed thus far, challenges and difficulties persist, especially in accurately representing the non-autonomous dynamics using the Koopman operator.
Extended Dynamic Mode Decomposition (EDMD)\cite{Williams2015}, which is a linear regression type of method to compute a finite-dimensional approximation of the Koopman operator, is widely used in the modeling phase.
However, it is known that its convergence property to the true Koopman operator does not hold for non-autonomous systems\cite{KORDA_Koopman_MPC}.
This implies that collecting a sufficiently large number of observables and data points may not necessarily lead to better model accuracy.
Thus, one may not be able to establish a reliable theoretical basis for the use of the Koopman operator for modeling the unknown, non-autonomous dynamics.
In addition to this convergence issue, 
the nature of data-driven problem settings also renders the occurrence of the modeling error unavoidable.
For instance, model training can yield biased or poor results if the data is obtained only from very limited operating points or the amount of data is too small.

To tackle the incompleteness of EDMD-based modeling, robust controller designs have been integrated into the Koopman framework.
In \cite{tube_based_MPC}, a robust Koopman MPC design is proposed on the basis of EDMD and tube-based MPC that accounts for both modeling error and additive disturbances.
An off-set free MPC framework was combined with Koopman Lyapunov-based MPC, which also takes the modeling error of EDMD approximation into account\cite{handling_plant_model_mismatch}.
Another robust controller design utilizes an $\mathcal{H}_2$ norm characterization of discrete-time linear systems to represent the modeling error of EDMD models using polytope sets\cite{data-driven_Koopman_H2}.
Despite these efforts focused on the controller design aspect, the fact that EDMD often fails to accurately reproduce the behavior of complex nonlinear dynamics may become an issue. 
Indeed, state prediction is one of the most basic yet essential tasks from a modeling perspective.

When EDMD-based methods do not result in satisfactory model accuracy, one may resort to using more expressive schemes such as neural networks to learn observables or feature maps from data.
The use of neural networks in Koopman-based modeling has been shown to be promising in many studies (e.g., \cite{Physics-based_robabilistic_learning,Learning_Koopman_Invariant_Subspaces,Deep_learning_universal_linear_embeddings,Learning_DNN_ACC2019,linearly_recurrent_autoencoder}).
Specifically, neural network-based models can be considerably more expressive than EDMD models and it is expected that they can address more challenging dynamical systems such as those with high degrees of nonlinearities and large dimensions. 
This strategy can be also adopted in control applications and various papers developed data-driven Koopman-based controller designs with neural networks\cite{DeSKO,deep_learning_Koopman_CDC2020,Deep_Koopman_vehicles,control_aware_Koopman}.

Although the use of neural networks may be a more preferable choice that yields expressive and flexible models, issues concerning the modeling error can hinder the successful application of Koopman-based modeling.
Compared to EDMD, learning with neural networks is generally more sensitive to the quality of data as well as other factors in the learning procedure since it is a high-dimensional non-convex optimization with many learnable parameters.
As a result, neural network-based models may be prone to overfitting or poor learning due to biased data-collecting procedures or lack of data.
In \cite{DeSKO}, a probabilistic neural network is adopted to introduce the uncertainty of modeling as an additive noise to the dynamics and a nominal controller is compensated by an additional control input.
Reference \cite{control_aware_Koopman} points out that the modeling error regarding closed-loop dynamics may be more difficult to manage than that of state prediction.
Based on this fact, a model refinement technique was proposed that incorporates data points collected from a closed-loop dynamics into neural network training to directly mitigate the modeling error in the closed-loop environment.

Given the large number of Koopman operator learning approaches in the literature, it remains a challenge to ensure that the learned model is not just capable of predicting the behavior of a nonlinear dynamics accurately, but is also amenable to many types of tasks such as controller design and decision-making.
The present work is focused on the use of oblique projection in a Hilbert space to develop a new linear operator-based data-driven modeling method.
The concept of projection has been used and proven to be useful in other disciplines of dynamical systems modeling\cite{survey_projection_based_ROM, Galerkin_v_Petrov_Galerkin}.
It also appears in the formulations of the data-driven Koopman operator-based modeling.
EDMD can be considered as an algorithm to compute the orthogonal projection of the action of the Koopman operator in the $L_2$ space endowed with the empirical measure supported on the given data points\cite{On_convergence_Klus,on_convergence_of_EDMD}.
This is also termed the Galerkin projection or Galerkin approximation of the Koopman operator in the literature\cite{Otto_Koopman_review,gEDMD}.
While this perspective is mostly adopted in theoretical analyses of EDMD or when the EDMD procedure is introduced, there are few works that utilize the concept of projection directly to develop Koopman operator-based data-driven methods.

This paper develops a new data-driven modeling method by reformulating a problem of dynamical systems modeling as projection-based linear operator-learning in a Hilbert space.
A notable feature of the proposed method is that its model structure is considered an extension of EDMD.
Consider a discrete-time, non-autonomous system $\chi^+=F(\chi,u)$, where $\chi\in \mathcal{X}\subseteq\mathbb{R}^{n}$, $u\in\mathcal{U}\subseteq \mathbb{R}^p$, $\chi^+\in \mathcal{X}$, and $F:\mathcal{X}\rightarrow \mathcal{X}$ are the state, the input, their corresponding successor, and a nonlinear mapping describing the dynamics, respectively.
Given feature maps or observables $g_i:\mathcal{X}\rightarrow \mathbb{R}$ ($i=1,\cdots,N_x$) and a data set $\{ (\chi_i,u_i,y_i)\mid y_i=F(\chi_i,u_i),\chi_i\in \mathcal{X},u_i\in \mathcal{U},i=1,\cdots,M \}$, we consider the following dynamic model that fits the data:
\begin{align}
	\left[
	\begin{array}{c}
		g_1(y_i)
		\\
		\vdots 
		\\
		g_{N_x}(y_i)
	\end{array}
	\right]
	\approx 
	\bm{A}
	\left[
	\begin{array}{c}
		g_1(\chi_i)
		\\
		\vdots 
		\\
		g_{N_x}(\chi_i)
	\end{array}
	\right]
	+
	\bm{B}
	u_i,\ \ i=1,\cdots, M,
	\nonumber
\end{align}
which is formally introduced as a linear embedding model in Section \ref{sec. defining the problem of data-driven modeling} (Eq. \eqref{eq. approx equation of the model}).
Note that the dynamics of this model is represented in the embedded state $[g_1(\chi),\cdots,g_{N_x}(\chi)]^\tr\in \mathbb{R}^{N_x}$, rather than the original state $\chi\in \mathbb{R}^n$.
The standard EDMD solves a linear regression problem that minimizes the sum of squared errors over the given data points.
This yields the model parameters $[\bm{A}\ \bm{B}]$
as its unique least square solution, which depends on the design of the feature maps $g_i$.
In the proposed method, we extend this structure of $[\bm{A}\ \bm{B}]$ in a way that they are not only dependent on the feature maps $g_i$, but also test functions $\varphi_i:\mathcal{X}\times \mathcal{U}\rightarrow \mathbb{R}$, ($i=1,\cdots,\hat{N}$), which are introduced by the idea of oblique projection in the context of linear operator-learning in a Hilbert space.

\begin{figure}[t]
	\centering
	\includegraphics[width=1.0\linewidth]{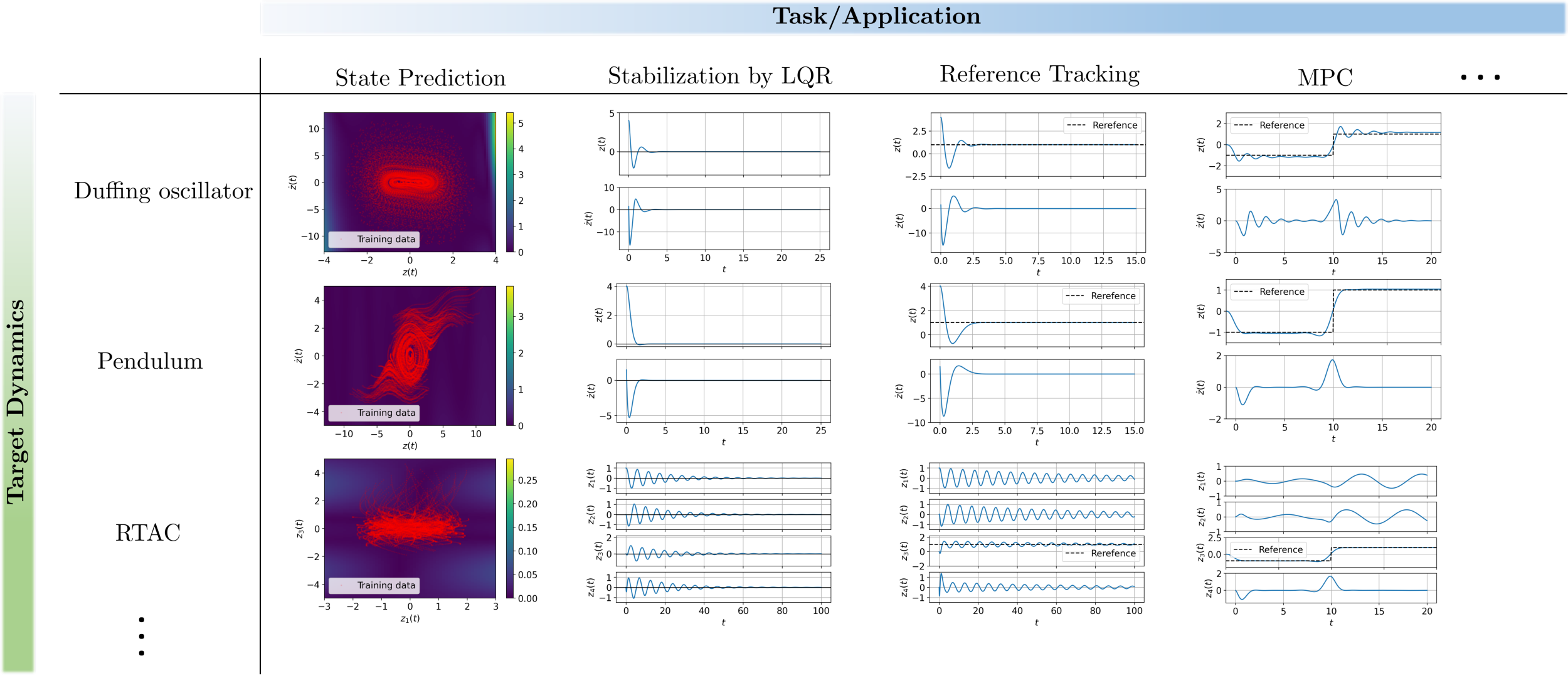}
	\caption{Chart showing the matrix of examples and tasks to assess generalizability of the present approach.}
	\label{fig. matrix chart}
\end{figure}

A major difference between the proposed method and EDMD, in addition to the different structures of $[\bm{A}\ \bm{B}]$, is that the feature maps $g_i$ and the test functions $\varphi_i$ are learned from data in the proposed method, whereas EDMD  seeks $[\bm{A}\ \bm{B}]$, with user-specified feature maps $g_i$.
%In this paper, we also adopt the same strategy with the test functions $\varphi_i$ be also characterized by neural networks.
Since the new structure of $[\bm{A}\ \bm{B}]$ includes the test functions $\varphi_i$ as additional tunable parameters, it allows for flexibility of the model, which could, in turn, lead to better accuracy and generalizability of the model.
Optimizing oblique projections has been also found to yield better performance in the context of reduced-order modeling\cite{Optimizing_oblique,huang2022model}.
Whereas \cite{Optimizing_oblique,huang2022model} seek to find optimal subspaces that define the oblique projection, we directly parameterize the feature maps and the test functions and implement the model learning as neural network training.
As shown in Fig. \ref{fig. matrix chart}, the generalizability of our modeling method is evaluated by two factors in this paper: different dynamical systems; and the range of tasks the model can be applied to.
We provide comprehensive numerical evaluations that compare the generalizability of different data-driven models using four tasks: state prediction, stabilization by LQR, reference tracking, and linear MPC, each of which targets three nonlinear dynamical systems: the Duffing oscillator, the simple pendulum, and the Rotational/Translational Actuator (RTAC).

In Section \ref{sec. ch 2}, we formally introduce a problem of data-driven modeling for discrete-time, non-autonomous systems and show that while the linear embedding model leads to an advantage that linear controller designs in the embedded space can be applied to control possibly nonlinear dynamics, it may suffer from fundamental limitations to realize accurate data-driven models.
In Section \ref{sec. Linear embedding models with oblique projection}, the same problem is reformulated in the context of projection-based linear operator learning in a Hilbert space and the EDMD operator is generalized to the one that is derived based on a weak formulation obtained by oblique projection.
In Section \ref{sec. learning procedures}, learning problems are formulated to train the proposed model, where neural networks are used to characterize both the feature maps and the test functions.
Numerical analyses are provided in Section \ref{sec. numerical evaluations} to show the capability and generalizability of the proposed method. 
In addition to the predictive accuracy of the model, three different control applications: stabilization by LQR, reference tracking, and linear MPC are considered and the proposed modeling method is compared to other Koopman-based data-driven models. 

\section{Linear Embedding Models for Nonlinear Dynamics}
\label{sec. ch 2}
\subsection{Defining the Modeling Problem}
\label{sec. defining the problem of data-driven modeling}
This section introduces a problem of data-driven modeling for a discrete-time dynamical system:
\begin{align}
	\chi^+=F(\chi,u),\ \ 
	\chi\in \mathcal{X}\subseteq \mathbb{R}^n,\  
	u\in \mathcal{U}\subseteq \mathbb{R}^p, \
	F:\mathcal{X}\times \mathcal{U}\rightarrow \mathcal{X},
	\label{eq. discrete time dynamical system}
\end{align}
where $\chi$ is an $n$ dimensional state in the state space $\mathcal{X}$, $u$ is an input from the admissible set $\mathcal{U}$ of inputs to the system, $\chi^+$ is the successor corresponding to $\chi$ and $u$, and $F$ is a possibly nonlinear map describing the dynamics, respectively.
It is assumed throughout the paper that the explicit knowledge of $F$ is not available and a model that reproduces the dynamics \eqref{eq. discrete time dynamical system} is obtained only from data of the form $\{(\chi^+, \chi, u)\mid \chi^+, \chi, \text{ and } u\text{ satisfy }\eqref{eq. discrete time dynamical system} \}$.

We introduce $N_x$ feature maps $g_i:\mathcal{X}\rightarrow \mathbb{R}$ ($i=1,\cdots, N_x$), which are functions from the state space $\mathcal{X}$ to $\mathbb{R}$, and consider approximating the unknown dynamics \eqref{eq. discrete time dynamical system} by a linear operator $[\bm{A}\ \bm{B}]$ s.t.
\begin{align}
\left\{
\begin{array}{l}
	\left[
		\begin{array}{c}
			g_1(\chi^+)
		\\
			\vdots 
		\\
			g_{N_x}(\chi^+)
		\end{array}
	\right]
	\approx 
	[\bm{A}\ \bm{B}]
	\left[
	\begin{array}{c}
		g_1(\chi)
		\\
		\vdots 
		\\
		g_{N_x}(\chi)
		\\
		u
	\end{array}
	\right]
	=
	\bm{A}
	\left[
	\begin{array}{c}
		g_1(\chi)
		\\
		\vdots 
		\\
		g_{N_x}(\chi)
	\end{array}
	\right]
	+
	\bm{B}
	u,
\\
	\chi^+\approx 
	\omega
	\left(
		\left[
		\begin{array}{c}
			g_1(\chi^+)
			\\
			\vdots 
			\\
			g_{N_x}(\chi^+)
		\end{array}
		\right]
	\right),
\end{array}
\right. 
\label{eq. approx equation of the model}
\end{align}
where $\bm{A}\in \mathbb{R}^{N_x\times N_x}$ and $\bm{B}\in \mathbb{R}^{N_x\times p}$ are matrices and $\omega:\mathbb{R}^{N_x}\rightarrow \mathbb{R}^n$ denotes a decoder from the embedded state $[g_1(\chi)\cdots g_{N_x}(\chi)]^\tr$, which is characterized by the feature maps $g_i$, to the original state $\chi$.
We call \eqref{eq. approx equation of the model} a linear embedding model in this paper. 
It is a common control model architecture in the Koopman literature, which has several unique features:
\vspace{5mm}
\begin{enumerate}
	\item 
	The model is represented in an $N_x$ dimensional embedded space, which is characterized by the nonlinear feature maps $g_i$.
	\item 
	The model dynamics in the embedded space is linear and is governed by
	the linear operator $[\bm{A}\ \bm{B}]$.
	\item 
	While nonlinearity w.r.t. the state $\chi$ can be introduced through the feature maps $g_i$, the model is strictly linear w.r.t. the input $u$ (Being strictly linear w.r.t. $u$ implies that the term involving $u$ is given in the form $\bm{B}u$, where $\bm{B}$ is a constant matrix). 
\end{enumerate}
\vspace{5mm}

The linear dynamics structure in the embedded space is especially advantageous when one uses the data-driven modeling method for control.
For simplicity, if the first $n$ feature maps are defined as the state $\chi\in \mathbb{R}^n$ itself so that $[g_1(\chi)\cdots g_n(\chi)]^\tr:=\chi$, the above equations can be written as
\begin{align}
	\left[
		\begin{array}{c}
			\chi^+
		\\
			g_{n+1}(\chi^+)
		\\
			\vdots 
		\\
			g_{N_x}(\chi^+)
		\end{array}
	\right]
	\approx 
	\bm{A}
	\left[
	\begin{array}{c}
		\chi 
		\\
		g_{n+1}(\chi)
		\\
		\vdots 
		\\
		g_{N_x}(\chi)
	\end{array}
	\right]
	+
	\bm{B}u,
	\label{eq. approx equation of the model with state obs included}
\end{align}
where the decoder $\omega$ need not be explicitly defined since the first $n$ components of the embedded state correspond to the prediction of the state $\chi^+$.
The above equation implies that linear controller designs can be utilized to control possibly nonlinear dynamics \eqref{eq. discrete time dynamical system} by viewing \eqref{eq. approx equation of the model with state obs included} as a Linear Time-Invariant (LTI) system in a new coordinate $[\chi^\tr\ g_{n+1}(\chi)\cdots\ g_{N_x}(\chi)]^\tr$. 
For instance, stabilization of the original dynamics \eqref{eq. discrete time dynamical system} may be achieved by designing an LQR that stabilizes the LTI system \eqref{eq. approx equation of the model with state obs included} since $[\chi^\tr\ g_{n+1}(\chi)\cdots\ g_{N_x}(\chi)]^\tr\rightarrow 0$ implies $\chi\rightarrow 0$.

On the other hand, the validity of this model w.r.t. the target dynamics \eqref{eq. discrete time dynamical system} may need to be established in practice since its construction is based on specific model structures as stated above, which in general limit the performance of the model for a wide range of dynamical systems.
Specifically, the modeling error, which is not yet stated in \eqref{eq. approx equation of the model} or \eqref{eq. approx equation of the model with state obs included} explicitly since the approximate equal sign $\approx$ is used instead, may become so dominant that the model will be no longer valid or reasonable to approximate the unknown dynamics \eqref{eq. discrete time dynamical system}.

\subsection{Koopman Operator}
\label{sec. introduce Koopman operator in the autonomous setting}
The Koopman operator\cite{Koopman_1931}, which is a type of linear operator acting on a function space, is well-known as a promising mathematical tool that may be able to justify the use of linear embedding models introduced in the previous section as well as serve as a basis to develop methods to learn such models from data.
For an autonomous system s.t.
\begin{align}
	\chi^+ = f(\chi),\ \ f:\mathcal{X}\rightarrow \mathcal{X},
	\label{eq. autonomous system}
\end{align}
the Koopman operator $\mathcal{K}:\mathcal{F}\rightarrow \mathcal{F}$ is defined as a linear operator s.t.
\begin{align}
	\mathcal{K}g = g\circ f
	\ \ \Leftrightarrow\ \ 
	\left(
	 \mathcal{K}g
	\right)(\chi)
	=
	g(f(\chi)),
	\ 
	g\in \mathcal{F},\ \chi\in \mathcal{X},
	\label{eq. def of Koopman operator for autonomous setting}
\end{align}
where $\mathcal{F}$ is a space of functions and its elements $g:\mathcal{X}\rightarrow \mathbb{R}$ may be considered as feature maps introduced in the previous section.
In the Koopman literature, $g$ are also referred to as observables.
From \eqref{eq. def of Koopman operator for autonomous setting}, it is easily seen that the Koopman operator $\mathcal{K}$ plays a role in describing the r.h.s. of autonomous dynamics \eqref{eq. autonomous system} in the function space $\mathcal{F}$ instead of the state space $\mathcal{X}$ as
\begin{align}
	g(\chi^+)=g(f(\chi))=(\mathcal{K}g)(\chi).
	\label{eq. composition relation}
\end{align}

Moreover, if we consider the same type of model representation as \eqref{eq. approx equation of the model} but ignoring the input effect:
\begin{align}
		\left[
		\begin{array}{c}
			g_1(\chi^+)
			\\
			\vdots 
			\\
			g_{N_x}(\chi^+)
		\end{array}
		\right]
		\approx 
		\bm{A}
		\left[
		\begin{array}{c}
			g_1(\chi)
			\\
			\vdots 
			\\
			g_{N_x}(\chi)
		\end{array}
		\right],
		\label{eq. autonomous model ignoring inputs}
\end{align}
this model can be naturally derived from a finite-dimensional approximation of the Koopman operator $\mathcal{K}$ by the following proposition.

%\takespace
\begin{screen}
\begin{prop}
\rm{} 
Suppose $g_i\in \mathcal{F}$ ($i=1,\cdots,D$, $D\in \mathbb{N}$). There exists $\bm{K}\in \mathbb{R}^{D\times D}$ s.t.
\begin{align}
	\left[
		\begin{array}{c}
			\mathcal{K}g_1
		\\
			\vdots 
		\\
			\mathcal{K}g_{D}
		\end{array}
	\right]
	=
	\bm{K}
	\left[
	\begin{array}{c}
		g_1
		\\
		\vdots 
		\\
		g_{D}
	\end{array}
	\right],
	\label{eq. prop finite dimensional approx of K}
\end{align}
if and only if $\text{span}(g_1,\cdots,g_D)$ is an invariant subspace under the action of the Koopman operator $\mathcal{K}$.
\label{prop. finite dimensional approx of Koopman}
\end{prop}
\end{screen}
\begin{proof}
	For instance, see \cite{control_aware_Koopman}, in which the non-autonomous setting is adopted but the same argument obviously holds for autonomous systems as well.
\end{proof}
%\takespace

If \eqref{eq. prop finite dimensional approx of K} holds, 
\begin{align}
	\left[
	\begin{array}{c}
		g_1(\chi^+)
		\\
		\vdots 
		\\
		g_{N_x}(\chi^+)
	\end{array}
	\right]=
	\left[
	\begin{array}{c}
		(\mathcal{K}g_1)(\chi)
		\\
		\vdots 
		\\
		(\mathcal{K}g_{N_x})(\chi)
	\end{array}
	\right]
	=
	\bm{K}
	\left[
		\begin{array}{c}
			g_1(\chi)
		\\
			\vdots
		\\
			g_{N_x}(\chi)
		\end{array}
	\right],
\end{align}
which implies that we can use $\bm{K}$ as the parameter $\bm{A}$ of model \eqref{eq. autonomous model ignoring inputs} on the assumption that the feature maps $g_i$ are appropriately designed so that \eqref{eq. prop finite dimensional approx of K} holds and $\bm{K}$ can be computed in some way.
Indeed, with the design of $g_i$ fixed, $\bm{K}$ can be easily estimated by EDMD\cite{Williams2015}.

Furthermore, it has been shown that the approximation $\bm{K}$ obtained by EDMD converges to the true Koopman operator in the Strong Operator Topology (SOT)\cite{on_convergence_of_EDMD}.
On several mathematical assumptions, as the number of feature maps $g_i$ and the number of data points used to estimate $\bm{K}$ tend to infinity, the EDMD operator converges to $\mathcal{K}$ in SOT.
Therefore, it may be judicious in the autonomous setting to add as many feature maps and data points as possible. 
This will eventually render the EDMD approximation converged to the true Koopman operator, which is sufficient to guarantee an accurate model for most modeling purposes.

\subsection{Extending the Koopman Operator Framework to Non-Autonomous Systems}
\label{sec. introduce Koopman operator in the non-autonomous setting}
While the Koopman operator framework can be used to form a rigorous basis for linear embedding models in the autonomous setting as seen in the previous section, the non-autonomous setting obstructs some of its useful properties from being true, which brings about several difficulties and challenges to the problem of interest.
In this section, we review the Koopman operator formalism for non-autonomous systems.
It is then followed by theoretical analyses in the next section, where we show some fundamental limitations of linear embedding models.
%to motivate the use of the proposed model architecture that combines the concepts of learning linear operators in a Hilbert space and oblique projection. 

To define the Koopman operator for the non-autonomous system \eqref{eq. discrete time dynamical system}, 
we first need to define a set $l(\mathcal{U})$ of infinite sequences of admissible inputs, i.e.,
\begin{align}
	l(\mathcal{U}):=
	\left\{
	\bm{U}:\mathbb{Z}_{\geq 0}\rightarrow \mathcal{U}:k\mapsto u_k
	\right\}
	=
	\{
	U:=
	(u_0,u_1,\cdots)\mid u_k\in \mathcal{U}, k\in \mathbb{Z}_{\geq 0}
	\}.
	\label{eq. def of space of infinite sequences of inputs}
\end{align}

Also, let $\hat{F}$ be defined as:
\begin{align}
	\hat{F}:\mathcal{X}\times l(\mathcal{U})\rightarrow \mathcal{X}\times l(\mathcal{U}):
	(\chi, U)\mapsto (F(\chi,U(0)), \mathcal{S}U),
	\label{eq. def of Fhat}
\end{align}
where $\mathcal{S}$ denotes the shift operator s.t.
\begin{align}
	\mathcal{S}U=
	\mathcal{S}(u_0,u_1,\cdots)
	:=
	(u_1,u_2,\cdots),
\end{align}
and $U(0)$ in \eqref{eq. def of Fhat} refers to the first element of $U=(u_0,u_1,\cdots)$.

Then, the Koopman operator $\hat{\mathcal{K}}$ for the non-autonomous system \eqref{eq. discrete time dynamical system} is defined as
\begin{align}
	\hat{\mathcal{K}}h = h\circ \hat{F}
	\ \ \Leftrightarrow\ \ 
	\left(
		\hat{\mathcal{K}}h
	\right)(\chi,U)
	=
	h(\hat{F}(\chi,U)),
	\
	h\in \mathcal{H},
	\label{eq. def of Koopman operator for non autonomous setting}
\end{align}
where the feature map $h:\mathcal{X}\times l(\mathcal{U})\rightarrow \mathbb{R}$ is now a function from the direct product $\mathcal{X}\times l(\mathcal{U})$ to the real numbers $\mathbb{R}$ belonging to some function space $\mathcal{H}$.
Note that $l(\mathcal{U})$ needs to be introduced to properly define the domain of $h$, not $\mathcal{U}$ itself, 
since if its domain were 
$\mathcal{X}\times \mathcal{U}$ 
so that
$h:\mathcal{X}\times \mathcal{U}\rightarrow \mathbb{R}$, 
the definition corresponding to \eqref{eq. def of Koopman operator for non autonomous setting} would become
$(\hat{\mathcal{K}}h)(\chi,u) = h(F(x,u), u^+)=h(x^+,u^+)$, 
which is clearly not well-defined since there exists no successor $u^+$ of input defined in the dynamics \eqref{eq. discrete time dynamical system}.
To avoid this issue, we assume an infinite sequence $U=(u_0,u_1,\cdots)$ of inputs so that what $\hat{\mathcal{K}}$ returns after composing two maps can be still well-defined. 
This also corresponds to the fact that an input \textit{signal} $\bm{U}:\mathbb{Z}_{\geq 0}\rightarrow \mathcal{U}$ in \eqref{eq. def of space of infinite sequences of inputs}, which is herein viewed as a function, is assumed to be specified when we introduce the Koopman operator along with the governing equation \eqref{eq. discrete time dynamical system}.
Note that only specifying a relation between the consecutive two states in time by the difference equation \eqref{eq. discrete time dynamical system} is sufficient to define the discrete-time, non-autonomous dynamics itself.
This formal extension to the non-autonomous setting was first developed in \cite{KORDA_Koopman_MPC}.

Similar to the previous section, the linear embedding model \eqref{eq. approx equation of the model} can be obtained by a finite-dimensional approximation of the Koopman operator $\hat{\mathcal{K}}$ for non-autonomous systems. 
Consider the following feature maps $h_i:\mathcal{X}\times l(\mathcal{U})\rightarrow \mathbb{R}$:
\begin{align}
	\left[
		\begin{array}{c}
			h_1(\chi, U)
		\\
			\vdots 
		\\	
			h_{N_x}(\chi,U)
		\\
			h_{N_x+1}(\chi, U)
		\\
			\vdots 
		\\
			h_{N_x+p}(\chi, U)
		\end{array}
	\right]
	=
	\left[
		\begin{array}{c}
			g_1(\chi)
		\\
			\vdots 
		\\
			g_{N_x}(\chi)
		\\
			U(0)
		\end{array}
	\right],
	\label{eq. def of each h_i for non autonomous setting Koopman}
\end{align}
where $g_i:\mathcal{X}\rightarrow \mathbb{R}$ are some functions from the state space $\mathcal{X}$ into $\mathbb{R}$.

Proposition \ref{prop. finite dimensional approx of Koopman} also holds in the non-autonomous setting, i.e., it is also true if $g_i\in \mathcal{F}$ and $\mathcal{K}$ are replaced with $h_i\in \mathcal{H}$ and $\hat{\mathcal{K}}$, respectively, as stated in its proof.
Suppose $\text{span}(h_1,\cdots,h_{N_x+p})$ is an invariant subspace under the action of $\hat{\mathcal{K}}$ so that there exists $\hat{\bm{K}}\in \mathbb{R}^{(N_x+p)\times (N_x+p)}$ s.t.
\begin{align}
	\left[
		\begin{array}{c}
			\hat{\mathcal{K}}h_1
		\\
			\vdots 
		\\
			\hat{\mathcal{K}}h_{N_x+p}
		\end{array}
	\right]
	=
	\hat{\bm{K}}
	\left[
	\begin{array}{c}
		h_1
		\\
		\vdots 
		\\
		h_{N_x+p}
	\end{array}
	\right].
	\label{eq. exact relation of finite approximation of K for non autonomous systems}
\end{align}

Also, let $\chi^+$ be the successor that corresponds to $\chi\in \mathcal{X}$ and $U(0)\in \mathcal{U}$. 
Then, we have
\begin{align}
	\left[
	\begin{array}{c}
		g_1(\chi^+)
		\\
		\vdots 
		\\
		g_{N_x}(\chi^+)
		\\
		U(1)
	\end{array}
	\right]
	=&
	\left[
	\begin{array}{c}
		g_1(F(\chi, U(0)))
		\\
		\vdots 
		\\
		g_{N_x}(F(\chi, U(0)))
		\\
		U(1)
	\end{array}
	\right]
%	&\nonumber
%\\
	=
	\left[
		\begin{array}{c}
			h_1( \hat{F}(\chi, U) )
		\\
			\vdots 
		\\ 
			h_{N_x+p}( \hat{F}(\chi, U) )
		\end{array}
	\right]
%	\ \ 
%	( \because \eqref{eq. def of each h_i for non autonomous setting Koopman},\
%	 \eqref{eq. def of Fhat})
%	&\nonumber
%\\
	=
	\left[
	\begin{array}{c}
		(\hat{\mathcal{K}}h_1)(\chi, U)
		\\
		\vdots 
		\\
		(\hat{\mathcal{K}}h_{N_x+p})(\chi, U)
	\end{array}
	\right]
	\ \ 
%	(\because \eqref{eq. def of Koopman operator for non autonomous setting})
	&\nonumber
\\
	=&
	\hat{\bm{K}}
	\left[
		\begin{array}{c}
			h_1(\chi, U)
		\\
			\vdots 
		\\
			h_{N_x+p}(\chi, U)
		\end{array}
	\right]
	\ \ 
%	(\because \eqref{eq. exact relation of finite approximation of K for non autonomous systems})
%	&\nonumber
%\\
	=
	\hat{\bm{K}}
	\left[
	\begin{array}{c}
		g_1(\chi)
		\\
		\vdots 
		\\
		g_{N_x}(\chi)
		\\
		U(0)
	\end{array}
	\right].
	\ \ 
%	( \because \eqref{eq. def of each h_i for non autonomous setting Koopman})
	&
 \label{eq. linear embedding model from Koopman middle eq}
\end{align}

If we define $[\bm{A}\ \bm{B}]\in \mathbb{R}^{N_x\times (N_x+p)}$ as the upper $N_x$ rows of $\hat{\bm{K}}$, the first $N_x$ equations of \eqref{eq. linear embedding model from Koopman middle eq} read:
\begin{align}
		\left[
		\begin{array}{c}
			g_1(\chi^+)
			\\
			\vdots 
			\\
			g_{N_x}(\chi^+)
		\end{array}
		\right]
		=
		[\bm{A}\ \bm{B}]
		\left[
		\begin{array}{c}
			g_1(\chi)
			\\
			\vdots 
			\\
			g_{N_x}(\chi)
			\\
			U(0)
		\end{array}
		\right]
		=
		\bm{A}
		\left[
		\begin{array}{c}
			g_1(\chi)
			\\
			\vdots 
			\\
			g_{N_x}(\chi)
		\end{array}
		\right]
		+
		\bm{B}
		U(0),
\end{align}
which justifies the use of model \eqref{eq. approx equation of the model} for the nonlinear dynamics \eqref{eq. discrete time dynamical system} if the assumptions made so far are satisfied.
To estimate $[\bm{A}\ \bm{B}]$, one can also use EDMD with pre-specified feature maps $g_i$ similar to the autonomous case\cite{KORDA_Koopman_MPC, DMDc}.
Specifically, given a data set $\{ (\chi_i,u_i,y_i)\mid y_i=F(\chi_i,u_i),\chi_i\in \mathcal{X},u_i\in \mathcal{U},i=1,\cdots,M \}$, EDMD solves a linear regression problem that minimizes the sum of squared errors over the data points:
\begin{align}
	\underset{[\bm{A}\ \bm{B}]}{\text{min\ }}
	\sum_{i=1}^{M}
	\left\|
	\left[
	\begin{array}{c}
		g_1(y_i)
		\\
		\vdots 
		\\
		g_{N_x}(y_i)
	\end{array}
	\right]
	-
	\left(
	\bm{A}
	\left[
	\begin{array}{c}
		g_1(\chi_i)
		\\
		\vdots 
		\\
		g_{N_x}(\chi_i)
	\end{array}
	\right]
	+
	\bm{B}
	u_i
	\right)
	\right\|_2^2,
\end{align}
which results in the following least square solution:
\begin{align}
	[\bm{A}\ \bm{B}]=&
	\left[
		\begin{array}{ccc}
			g_1(y_1) & \cdots & g_1(y_M)
		\\
			\vdots & \ddots & \vdots
		\\
			g_{N_x}(y_1) & \cdots & g_{N_x}(y_M) 
		\end{array}
	\right]
	\left[
	\begin{array}{ccc}
		g_1(\chi_1) & \cdots & g_1(\chi_M)
		\\
		\vdots & \ddots & \vdots
		\\
		g_{N_x}(\chi_1) & \cdots & g_{N_x}(\chi_M)
		\\
		u_1 & \cdots & u_M 
	\end{array}
	\right]^\tr
	&\nonumber
\\
	&
	\times  
	\left(
		\left[
		\begin{array}{ccc}
			g_1(\chi_1) & \cdots & g_1(\chi_M)
			\\
			\vdots & \ddots & \vdots
			\\
			g_{N_x}(\chi_1) & \cdots & g_{N_x}(\chi_M)
			\\
			u_1 & \cdots & u_M 
		\end{array}
		\right]
		\left[
		\begin{array}{ccc}
			g_1(\chi_1) & \cdots & g_1(\chi_M)
			\\
			\vdots & \ddots & \vdots
			\\
			g_{N_x}(\chi_1) & \cdots & g_{N_x}(\chi_M)
			\\
			u_1 & \cdots & u_M 
		\end{array}
		\right]^\tr 
	\right)^\dagger,
	\label{eq. intro [A B] as in EDMD}
\end{align}
where $\dagger$ denotes the Moore-Penrose pseudo inverse.

\subsection{Fundamental Limitations of Linear Embedding Models for Non-Autonomous Dynamics}
\label{sec. fundamental limitations}
\subsubsection{On the Convergence Property of EDMD for Non-Autonomous Systems}
\label{sec. on no convergence property of EDMD}
While linear embedding models for both autonomous and non-autonomous settings can be associated with the Koopman operator formalism as seen in Sections \ref{sec. introduce Koopman operator in the autonomous setting} and \ref{sec. introduce Koopman operator in the non-autonomous setting}, it is known that the convergence property of EDMD that holds for autonomous systems is no longer true in the non-autonomous setting\cite{KORDA_Koopman_MPC}.
One of the assumptions required for the convergence property to be true is that the feature maps span an orthonormal basis of the $L_2$ space as the number of them tends to infinity.
However, it is not satisfied in the non-autonomous setting we adopt in this paper, where the feature maps $h_i:\mathcal{X}\times l(\mathcal{U})\rightarrow \mathbb{R}$ are defined as in \eqref{eq. def of each h_i for non autonomous setting Koopman}.
Indeed, only the first component of the infinite sequence $U$ appears in \eqref{eq. def of each h_i for non autonomous setting Koopman} and they obviously cannot form any basis of $L_2(\mathcal{X}\times l(\mathcal{U}))$ even if we increase the number of the first $N_x$ feature maps $g_i$.

A possible way to circumvent this issue is relaxing the specific model structure that the model dynamics is strictly linear w.r.t. $u$ so that the term $\bm{B}u$ or $\bm{B}U(0)$ in the model representations will be replaced with a general nonlinear function from $l(\mathcal{U})$ to $\mathbb{R}^{N_x}$.
In fact, if the unknown dynamics can be assumed to be a control-affine Ordinary Differential Equation (ODE), it may be beneficial to consider a bilinear model representation by making use of the Koopman canonical transform\cite{KCT_control_affine} and an associated bilinearization of the original control-affine dynamics using Koopman generators\cite{Koopman_bilinearization}.

However, retaining the strictly linear model dynamics w.r.t. the input $u$ is obviously a quite advantageous factor that allows for the use of well-developed linear systems theories to control unknown nonlinear dynamics as seen in Section \ref{sec. defining the problem of data-driven modeling}.
Therefore, while preserving this linear structure to realize linear controller designs for possibly nonlinear dynamics, we aim to develop a data-driven modeling method that can achieve reasonable approximations of the true dynamics even without the existence of a convergence property.
Note that, in contrast to our main focus on the data-driven modeling aspect, several works also deal with the inevitable modeling error due to this no convergence property and the data-driven nature of EDMD from the controller design perspective\cite{tube_based_MPC,handling_plant_model_mismatch,data-driven_Koopman_H2} (see Section \ref{sec. intro} for details).

Neural networks are adopted in the proposed method to learn the feature maps $g_i$ from data, which has been proven to be a promising strategy for Koopman data-driven methods among many other modeling frameworks\cite{Physics-based_robabilistic_learning,Learning_Koopman_Invariant_Subspaces,Deep_learning_universal_linear_embeddings,Learning_DNN_ACC2019,linearly_recurrent_autoencoder,DeSKO,deep_learning_Koopman_CDC2020,Deep_Koopman_vehicles,control_aware_Koopman}.
Compared to EDMD, in which the feature maps $g_i$ need to be specified by the user, learning a model with $g_i$ characterized by neural networks allows for greater expressivity of the models.
On the other hand, it is noted that the neural network training is reduced to a high dimensional non-convex optimization and one may have difficulty training the model properly, e.g., the learning result can easily suffer from overfitting if data is not collected in an appropriate way or the optimization is terminated at an undesired local minimum before the loss function reaches a sufficiently low value.
Therefore, the use of neural networks does not necessarily resolve the issues about the modeling error on its own.

\subsubsection{Necessary Condition for the Model to Achieve Zero Modeling Error}
As another fundamental limitation of the linear embedding models for non-autonomous systems, we provide a necessary condition for the model to achieve precisely zero modeling error.
Henceforth, we assume that the model is of the form \eqref{eq. approx equation of the model with state obs included}, which includes the state $\chi\in \mathbb{R}^n$ itself as its first $n$ feature maps. 
Note that this is primarily for the sake of model simplicity and most statements given in this paper can be also true for its generalized representation \eqref{eq. approx equation of the model}.
We define the modeling error $\mathcal{E}(\chi,u)$ at $\chi\in \mathcal{X}$ and $u\in \mathcal{U}$ by
\begin{align}
	\mathcal{E}(\chi,u):=
	\left\| 
		\left[
	\begin{array}{c}
		\chi^+
		\\
		g_{n+1}(\chi^+)
		\\
		\vdots 
		\\
		g_{N_x}(\chi^+)
	\end{array}
	\right]
	-
	\left(
	\bm{A}
	\left[
	\begin{array}{c}
		\chi 
		\\
		g_{n+1}(\chi)
		\\
		\vdots 
		\\
		g_{N_x}(\chi)
	\end{array}
	\right]
	+
	\bm{B}u
	\right)
	\right\|_2
        \ \ 
        (\chi^+=F(\chi,u)),
	\label{eq. def of modeling error}
\end{align}
which is a point-wise error at $(\chi,u)$ measured in the embedded space.
This represents the accuracy of the model as a finite-dimensional approximation of the Koopman operator.
Specifically, we can have $\mathcal{E}(\chi,u)=0$ for $\forall \chi,u$ if and only if the feature maps $g_i$ span an invariant subspace under the action of the Koopman operator.
Also, it reflects the model performance when applied to control applications since $\mathcal{E}(\chi,u)=0$ is implicitly assumed when controllers are designed for the virtual LTI system \eqref{eq. approx equation of the model with state obs included}.
\\ \ \\
On the other hand, the state predictive accuracy is dictated by the state prediction error $\mathcal{E}_\text{state}(\chi,u)$, which is defined as follows:

\begin{align}
    \mathcal{E}_\text{state}(\chi,u):=
    \left\|
        \chi^+ - 
        \underset{\text{decoder}}{\underbrace{
        [I_n\ 0]}}
        \left(
	\bm{A}
	\left[
	\begin{array}{c}
		\chi 
		\\
		g_{n+1}(\chi)
		\\
		\vdots 
		\\
		g_{N_x}(\chi)
	\end{array}
	\right]
	+
	\bm{B}u
	\right)
    \right\|_2.
    \label{eq. def of state prediction error}
\end{align}

Since $\mathcal{E}(\chi,u)=0$ implies $\mathcal{E}_\text{state}(\chi,u)=0$, we focus on conditions for $\mathcal{E}(\chi,u)=0$ in this section.
An equivalent condition for $\mathcal{E}_\text{state}(\chi,u)=0$ to be true is provided in \cite{control_aware_Koopman} along with related discussions.

Noticing that individual components of the embedded state $[\chi^\tr\ g_{n+1}(\chi)\ \cdots\ g_{N_x}(\chi)]^\tr$ are mutually dependent through the state $\chi$,
a necessary condition to eliminate the modeling error $\mathcal{E}(\chi,u)$ is stated as follows.
%\takespace
\begin{screen}
	\begin{prop}
		\label{prop. necessary condition manifold}
		\rm{}
		For arbitrary $\chi\in \mathcal{X}$ and $u\in \mathcal{U}$, the following is a necessary condition for $\mathcal{E}(\chi,u)=0$ to hold:
		\begin{align}
			\bm{A}
			\left[
			\begin{array}{c}
				\chi
				\\
				g_{n+1}(\chi)
				\\
				\vdots 
				\\
				g_{N_x}(\chi)
			\end{array}
			\right]
			+
			\bm{B}
			u
			\in 
			\mathcal{M},
			\label{eq. necessary condition for zero error}
		\end{align}
		where $\mathcal{M}$ is a set on which the embedded state is defined, i.e.,
		\begin{align}
			\mathcal{M}:=
			\left\{
			v=[v_1\cdots v_{N_x}]^\tr \in \mathbb{R}^{N_x}
			\mid 
			\chi\in \mathcal{X},\ 
			\left\{
			\begin{array}{l}
				[v_1\cdots v_n]^\tr = \chi
				\\	
				v_i = g_{i}(\chi)
				(i=n+1,\cdots, N_x)
			\end{array}
			\right.
			\right\}.
		\end{align}
	\end{prop}
\end{screen} 
\begin{comment}
\begin{proof}
	If $\mathcal{E}(\chi,i)=0$, \eqref{eq. approx equation of the model with state obs included} is an exact relation and can be written as
	\begin{align}
		\left[
		\begin{array}{c}
			F(\chi, u)
			\\
			g_{n+1}(F(\chi, u))
			\\
			\vdots 
			\\
			g_{N_x}(F(\chi, u))
		\end{array}
		\right]
		=
		\bm{A}
		\left[
		\begin{array}{c}
			\chi
			\\
			g_{n+1}(\chi)
			\\
			\vdots 
			\\
			g_{N_x}(\chi)
		\end{array}
		\right]
		+
		\bm{B}
		u,
		\label{eq. used in proof prop 2}
	\end{align}
	which implies \eqref{eq. necessary condition for zero error} since $F(\chi,u)\in \mathcal{X}$.
\end{proof}
\end{comment}
%\takespace

Proposition \ref{prop. necessary condition manifold} states
that for the model \eqref{eq. approx equation of the model with state obs included} to achieve $\mathcal{E}(\chi,u)=0$, it is necessary that its trajectory is constrained on the set $\mathcal{M}$.
Note that the output of the model, or the r.h.s. of \eqref{eq. approx equation of the model with state obs included}, may not be on $\mathcal{M}$ in practice due to the incompleteness of the model.

This necessary condition can be written in another form, which suggests the difficulty of achieving $\mathcal{E}(\chi,u)=0$ for various input values, as shown in the following proposition.

\begin{screen}
%\vspace{7mm}
	\begin{prop}
		\label{prop. on the necessary condition}
		\rm{}
		Suppose that there exist a model of the form \eqref{eq. approx equation of the model with state obs included}, $\chi_0\in \mathcal{X}$, and $u_0\in \mathcal{U}$ s.t. \eqref{eq. necessary condition for zero error} holds at $\chi_{0}$ and $u_0$, i.e., for $\chi_0$ and $u_0$, there exists $\hat{\chi}_0\in \mathcal{X}$ s.t. 
		\begin{align}
			\bm{A}
			\left[
			\begin{array}{c}
				\chi_0
				\\
				g_{n+1}(\chi_0)
				\\
				\vdots 
				\\
				g_{N_x}(\chi_0)
			\end{array}
			\right]
			+
			\bm{B}
			u_0
			=
			\left[
			\begin{array}{c}
				\hat{\chi}_0
				\\
				g_{n+1}(\hat{\chi}_0)
				\\
				\vdots 
				\\
				g_{N_x}(\hat{\chi}_0)
			\end{array}
			\right].
			\label{eq. used in proof prop 3}
		\end{align}
		The following is an equivalent condition for \eqref{eq. necessary condition for zero error} to be also true at $\chi_0$ and any other $u\in \mathcal{U}$:\\
		For any $u\in \mathcal{U}$, there exists $\chi\in \mathcal{X}$ s.t.
		\begin{align}
			\bm{B}\Delta u
			=
			\left[
			\begin{array}{c}
				\chi - \hat{\chi}_0
				\\	
				g_{n+1}(\chi) - g_{n+1}(\hat{\chi}_0) 
				\\
				\vdots 
				\\
				g_{N_x}(\chi) - g_{N_x}(\hat{\chi}_0) 
			\end{array}
			\right],
			\label{eq. condition for the necessity with any input}
		\end{align}
		where $\Delta u:=u-u_0$.
	\end{prop}
\end{screen}
\begin{proof}
	Suppose that \eqref{eq. used in proof prop 3} holds for some $\hat{\chi}_0,\chi_0\in \mathcal{X}$, and $u_0\in \mathcal{U}$.
	If \eqref{eq. necessary condition for zero error} is true for $\chi_0$ and an arbitrary $u\in \mathcal{U}$, there exists $\chi\in\mathcal{X}$ s.t.
	\begin{align}
		\left[
		\begin{array}{c}
			\chi
			\\
			g_{n+1}(\chi)
			\\
			\vdots 
			\\
			g_{N_x}(\chi)
		\end{array}
		\right]
		=
		\bm{A}
		\left[
		\begin{array}{c}
			\chi_0
			\\
			g_{n+1}(\chi_0)
			\\
			\vdots 
			\\
			g_{N_x}(\chi_0)
		\end{array}
		\right]
		+
		\bm{B}
		u_0
		+
		\bm{B}
		\Delta u,
	\end{align}
	where $\Delta u:=u-u_0$.
	Subtracting \eqref{eq. used in proof prop 3} from the above equation yields \eqref{eq. condition for the necessity with any input}.

	Conversely, if there exists $\chi\in \mathcal{X}$ s.t. \eqref{eq. condition for the necessity with any input} is true for an arbitrary $u\in \mathcal{U}$, we have
	\begin{align}
		&\eqref{eq. condition for the necessity with any input}
		\ \ 
            \Leftrightarrow
		\ \
		\bm{B}u - \bm{B}u_0
		+
		\left[
		\begin{array}{c}
			\hat{\chi}_0
			\\
			g_{n+1}(\hat{\chi}_0)
			\\
			\vdots 
			\\
			g_{N_x}(\hat{\chi}_0)
		\end{array}
		\right]
		=
		\left[
		\begin{array}{c}
			\chi
			\\
			g_{n+1}(\chi)
			\\
			\vdots 
			\\
			g_{N_x}(\chi)
		\end{array}
		\right],
		&\nonumber
        \end{align}%\\
        \begin{align}
		\ \
		\Leftrightarrow
		\ \
		&
		\bm{A}
			\left[
		\begin{array}{c}
			\chi_0
			\\
			g_{n+1}(\chi_0)
			\\
			\vdots 
			\\
			g_{N_x}(\chi_0)
		\end{array}
		\right]
		\hspace{-1mm}+\hspace{-1mm}
		\bm{B}u
		\hspace{-1mm}-\hspace{-1mm}
		\left(
			\bm{A}
			\left[
			\begin{array}{c}
				\chi_0
				\\
				g_{n+1}(\chi_0)
				\\
				\vdots 
				\\
				g_{N_x}(\chi_0)
			\end{array}
			\right]
			+
			\bm{B}u_0
		\right)
		\hspace{-1mm}+\hspace{-1mm}
		\left[
		\begin{array}{c}
			\hat{\chi}_0
			\\
			g_{n+1}(\hat{\chi}_0)
			\\
			\vdots 
			\\
			g_{N_x}(\hat{\chi}_0)
		\end{array}
		\right]
		\hspace{-1mm}=\hspace{-1mm}
		\left[
		\begin{array}{c}
			\chi
			\\
			g_{n+1}(\chi)
			\\
			\vdots 
			\\
			g_{N_x}(\chi)
		\end{array}
		\right],
		&\nonumber
	\\
		\ \ 
		\Leftrightarrow 
		\ \ 
		&
		\bm{A}
		\left[
		\begin{array}{c}
			\chi_0
			\\
			g_{n+1}(\chi_0)
			\\
			\vdots 
			\\
			g_{N_x}(\chi_0)
		\end{array}
		\right]
		+
		\bm{B}u
		=
		\left[
		\begin{array}{c}
			\chi
			\\
			g_{n+1}(\chi)
			\\
			\vdots 
			\\
			g_{N_x}(\chi)
		\end{array}
		\right],
		\ \ 
		(\because \eqref{eq. used in proof prop 3})
		&\nonumber
	\end{align}
	which implies that \eqref{eq. necessary condition for zero error} holds for $\chi_0$ and $u$ since $\chi\in \mathcal{X}$.
\end{proof}
%\vspace{7mm}

	In Proposition \ref{prop. on the necessary condition},
	suppose $\bm{B}\neq 0$ and consider the model dynamics at some point $\chi_0\in \mathcal{X}$.
	Even if there exists a model for which the necessary condition \eqref{eq. necessary condition for zero error} holds at $\chi_0$ and some input $u_0$, it is in general difficult to be true at $\chi_0$ and a different input $u_1$ s.t. $\bm{B}u_0\neq \bm{B}u_1$.
	Indeed, satisfying \eqref{eq. condition for the necessity with any input} for all $u\in\mathcal{U}$ is equivalent to the linear subspace range($\bm{B}$) being a tangent space at every point of $\mathcal{M}$.
	%\eqref{eq. condition for the necessity with any input} implies that the set $\mathcal{M}$, which is characterized by the feature maps $g_i$, may be represented by a linear subspace $\text{range}(\bm{B})$.
	Considering that $g_i$ needs to be nonlinear to capture the nonlinearity of the original dynamics \eqref{eq. discrete time dynamical system}, satisfying \eqref{eq. condition for the necessity with any input} for various inputs while including sufficient complexity or nonlinearity in the model through the design of $g_i$ is not likely to be possible in practice, which can be even seen in the following simple example.

%\takespace
\begin{ex}
	\rm{} 
	\label{ex. necessary condition}
	For one-dimensional nonlinear dynamics with a scalar input given by:
	\begin{equation}
		\chi^+ = \sin \chi + u,
		\ \ 
            \chi\in \mathbb{R},
            \ \ 
		u\in [0,\pi],
	\end{equation}
	consider the following model of the form \eqref{eq. approx equation of the model with state obs included}:
	\begin{align}
		\left[
			\begin{array}{c}
				\chi^+
			\\
				\sin \chi^+
			\end{array}
		\right]
		\approx 
		\bm{A}
		\left[
			\begin{array}{c}
				\chi
			\\
				\sin \chi 
			\end{array}
		\right]
		+
		\underset{=\bm{B}}{\underbrace{
		\left[
			\begin{array}{c}
				1
			\\
				0
			\end{array}
		\right]}}
		u,
	\label{eq. model for ex}		
	\end{align}
	where $\bm{A}\in \mathbb{R}^{2\times 2}$ can be arbitrary for the purpose of this example. 
	Note that this model includes a single feature map $g_2(\chi)=\sin \chi$, which can realize the perfect state prediction, i.e., if the first row of $\bm{A}$ is $[0\ 1]$, the first component of the model's output, which is equivalent to its state prediction, becomes $\sin \chi + u$.

	At the origin $\chi_0:=0$ of the state space, this model satisfies the necessary condition given by Proposition \ref{prop. necessary condition manifold}, or equivalently \eqref{eq. used in proof prop 3}, with the input $u_0:=0$ since there exists a state $\hat{\chi}_0:=0$, which satisfies the following:
	\begin{align}
		\left(
		\bm{A}
		\left[
		\begin{array}{c}
			\chi_0
			\\
			\sin \chi_0
		\end{array}
		\right]
		+
		\bm{B}u_0
		=
		\right)
		\bm{A}
		\left[
		\begin{array}{c}
			0
			\\
			0
		\end{array}
		\right]
		+
		\bm{B}\cdot 0
		=
		\left[
		\begin{array}{c}
			0
			\\
			0
		\end{array}
		\right]
		\left(
		=
		\left[
		\begin{array}{c}
			\hat{\chi}_0
			\\	
			\sin \hat{\chi}_0
		\end{array}
		\right]
		\right).
	\end{align}

	Then, by inspecting Fig. \ref{fig. example for proposition about the necessary condition} or finding roots of the following equation:
	\begin{align}
		\sin\chi = 0,
	\end{align}
	we obtain $\chi_1:=\pi$, which satisfies \eqref{eq. condition for the necessity with any input} since there exists $u_1:=\pi$ s.t. $\Delta u_1:=u_1 - u_0$ and
	\begin{align}
		\left(
		\bm{B}
		\Delta u_1
		=
		\right)
		\left[
		\begin{array}{c}
			1
			\\
			0
		\end{array}
		\right]
		\pi
		=
		\left[
		\begin{array}{c}
			\pi - 0
			\\	
			0 - 0
		\end{array}
		\right]
		\left(
		=
		\left[
		\begin{array}{c}
			\chi_1 - \hat{\chi}_0
			\\
			\sin \chi_1 - \sin \hat{\chi}_0
		\end{array}
		\right]
		\right).
	\end{align}

	Obviously, $u_0=0$ and $u_1=\pi$ are the only inputs $u$ s.t. $u\in [0,\pi]$ and they satisfy \eqref{eq. condition for the necessity with any input}. 
	Therefore, at the reference state $\chi_0$, the model \eqref{eq. model for ex} can satisfy the necessary condition \eqref{eq. necessary condition for zero error} for $\mathcal{E}(\chi_0,u)=0$ to be true only if the input $u$ is either $u=u_0$ or $u=u_1$.

		\begin{figure}[t]
		\centering 
		\includegraphics[width=0.6\columnwidth]{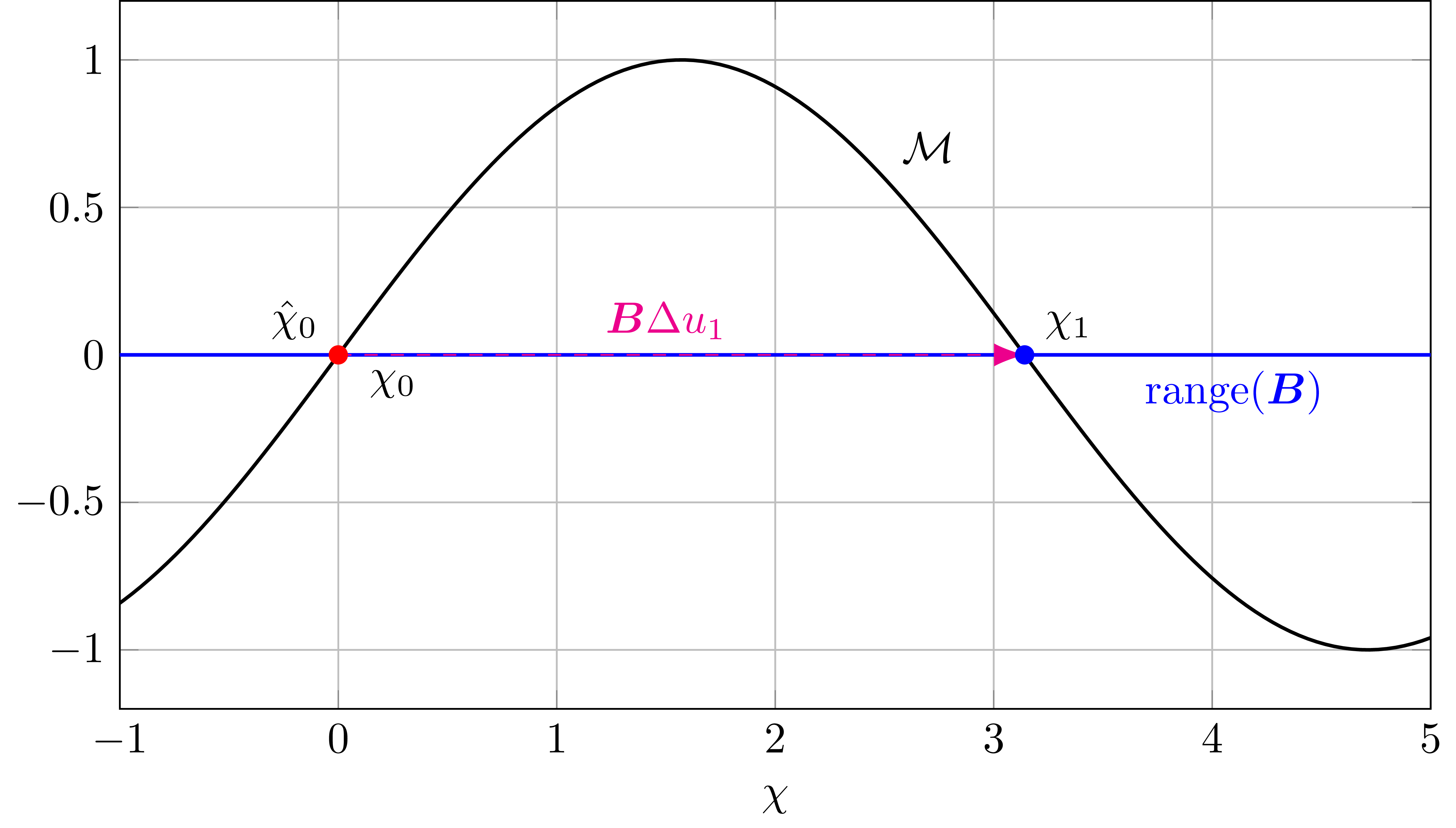}
		\caption{Embedded space $\mathbb{R}^2$ associated with the model \eqref{eq. model for ex}.}
		\label{fig. example for proposition about the necessary condition}
	\end{figure}
	\begin{table}[]
		\centering
		\caption{Modeling errors $\mathcal{E}(\chi_0,u)$ of model \eqref{eq. model for ex}.}
		\label{table ex necessary condition}
		\begin{tabular}{c|ccccccc}
			$u$ & 0(=$u_0$) & 0.5 & 1.0 & 1.5 & 2.0 & 2.5 & $\pi$(=$u_1$) 
			\\ \hline 
			$\mathcal{E}(\chi_0,u)$ & 0 & 0.479
			& 0.841
			& 0.997
			& 0.909
			& 0.598
			& 0
		\end{tabular}
	\end{table}

	As a numerical evaluation, modeling errors at $\chi_0$ and various inputs are listed in Table \ref{table ex necessary condition}.
	As expected from Proposition \ref{prop. on the necessary condition}, the cases with $u=u_0$ and $u=u_1$ achieve zero modeling errors whereas other input values result in non-zero errors.
	These errors can be also computed analytically in this example as
	\begin{align}
            \mathcal{E}(\chi_0,u)=
		\left\|
		\left[
		\begin{array}{c}
			\chi^+
			\\
			\sin \chi^+
		\end{array}
		\right]
		-
	\bm{A}
		\left[
		\begin{array}{c}
			\chi_0
			\\
			\sin \chi_0
		\end{array}
		\right]
		-
		\left[
		\begin{array}{c}
			1
			\\
			0
		\end{array}
		\right]u
		\right\|_2
		=&
		\left\| 
		\left[
		\begin{array}{c}
			u
			\\
			\sin u
		\end{array}
		\right]
		-
		\left[
		\begin{array}{c}
			u
			\\
			0
		\end{array}
		\right]
		\right\|_2
		=
		| \sin u |,
		&
		\label{eq. ex analytical expression of error}
	\end{align} 
	which becomes 0 only if we choose $u=0,\pi$, within the admissible set $[0,\pi]$.

	This example implies difficulty of deploying the linear embedding model \eqref{eq. approx equation of the model with state obs included} such that it remains accurate across a wide range of operating regimes of dynamics.
	For instance, if the model \eqref{eq. model for ex} is to be used to design a feedback controller represented by:
	\begin{align}
		u=\bm{Q}
		\left[
		\begin{array}{c}
			\chi 
			\\
			\sin \chi
		\end{array}
		\right],
		\ \ \bm{Q}\in \mathbb{R}^{1\times 2},
	\end{align}
	the input $u$ will vary continuously w.r.t. $\chi$ and the controller performance may be deteriorated due to the modeling error whenever $u\neq u_0,u_1$, as seen in Table \ref{table ex necessary condition} or \eqref{eq. ex analytical expression of error}.	
	Note that this analysis based on Proposition \ref{prop. on the necessary condition} is not even related to the Koopman operator formalism and is only concerned with the model structure itself.
\end{ex}
\takespace

	The discussions provided in this section imply that even if the model is quite expressive, e.g., with the use of neural networks, and its optimal model parameters minimizing the loss function are obtained in the actual training process, the performance of the model can still suffer from the fundamental limitation as shown in Proposition \ref{prop. on the necessary condition}.
    Indeed, it may be
    possible that the modeling error $\mathcal{E}(\chi,u)$ is not negligible for a wide range of state/input values due to this limitation.

 This observation, combined with the model's lack of convergence discussed in Section \ref{sec. on no convergence property of EDMD}, suggests that learning linear embedding models for non-autonomous dynamics is challenging.
		In this paper, we propose extending the EDMD-type structure of the model parameters $[\bm{A}\ \bm{B}]$ by introducing the oblique projection in the context of linear operator-learning in a Hilbert space.
		This effectively imposes a new type of constraint on $[\bm{A}\ \bm{B}]$ through the use of test functions that can be viewed as additional degrees of freedom in  the training process.
		As a result, both the accuracy and the generalizability of the model are expected to be improved even in the presence of the fundamental limitations of linear embedding models discussed in this section.
	%This observation, along with the no-convergence property of the model as seen in Section \ref{sec. on no convergence property of EDMD}, motivates a modification to the existing model structure so that both the accuracy and the generalizability of the model can be further improved under the fundamental limitations of linear embedding models. 
	%Specifically, we propose extending the EDMD-type structure of the model parameters $[\bm{A}\ \bm{B}]$ by introducing the oblique projection in the context of linear operator-learning in a Hilbert space.
    %This effectively imposes a new type of constraint on $[\bm{A}\ \bm{B}]$ through the use of test functions that can be viewed as additional tuning knobs of the training process.

\section{Linear Embedding Models with Oblique Projection}
\label{sec. Linear embedding models with oblique projection}
In this section, we derive the proposed model architecture, starting from a weak formulation obtained by characterization of the action of a general linear operator w.r.t. a finite dimensional subspace of a Hilbert space. 
We focus on a finite dimensional linear operator that appears in the weak formulation, which will be utilized to obtain the operator that governs the dynamics of the proposed linear embedding model.
It is numerically estimated through a finite data approximation and used as the parameters $[\bm{A}\ \bm{B}]$.

\subsection{Weak Formulation for Characterizing a Linear Operator}
\label{sec. 3}
Let $\mathcal{L}:H\rightarrow H$ be a linear operator acting on a Hilbert space $H$ and consider the problem of finding $\mathcal{L}$ 
given data obtained from simulations or experiments.
If $H$ is a space of functions defined on a set $\mathbb{X}$, the data takes the form $\{ g(x_i),(\mathcal{L}g)(x_i)\mid x_i\in \mathbb{X}, g:\mathbb{X}\rightarrow \mathbb{R}\in H \}$.

Since it is quite challenging to find a possibly infinite-dimensional $\mathcal{L}$ directly, we instead consider its restriction to some finite-dimensional subspace $W$, i.e., $\mathcal{L}_{|W}:W\rightarrow H:g\mapsto \mathcal{L}g$.
If $W$ is invariant under the action of $\mathcal{L}$, there exists a unique matrix representation of $\mathcal{L}_{|W}$ given a basis of $W$.
Therefore, we aim to find a subspace $W$, or equivalently its basis, such that it spans an invariant subspace.
While this strategy itself does not provide a solution to the problem of finding $\mathcal{L}$, only obtaining a matrix representation of its restriction to a finite-dimensional invariant subspace is in fact enough for many engineering applications that utilize the linear operator-learning framework to realize prediction, analyses, and control of dynamical systems.

Given a finite-dimensional subspace $W\subset H$, consider a direct decomposition of $H$ s.t.
\begin{align}
	H = W \oplus E,
	\label{eq. direct decomposition in a general form}
\end{align}
where $E$ is a complementary subspace.
This decomposition is well-defined since we can at least take the orthogonal direct decomposition: $H=W\oplus W^{\perp}$, whose existence is ensured for a general Hilbert space $H$.
Note that 
for $\forall g\in W$,
$\mathcal{L}_{|W}g\in H$ may not stay in $W$ ($\mathcal{L}_{|W}g\in W$ is true if and only if $W$ is an invariant subspace), and
it can be decomposed as follows:
\begin{align}
	\mathcal{L}_{|W}g = Lg + e,
	\label{eq. direct decomposition each element}
\end{align}
where $L:W\rightarrow W$ is an operator that assigns a unique element of $W$ when $\mathcal{L}g$ is decomposed according to \eqref{eq. direct decomposition in a general form} and $e\in E$ is the other corresponding unique element of $E$, respectively.

The second term $e$ may be interpreted as the projection error by considering a (possibly oblique) projection operator $P_{W,E}:H\rightarrow W$, which is defined as:
\begin{align}
	P_{W,E}h:=
        w\ \  
        \text{for }\forall h\in H,\ \text{where}\ 
        h=w+e,\ 
        w\in W,\ 
        e\in E
	.
        \label{eq. def of projection operator}
\end{align}

It is easily confirmed that $P_{W,E}$ is a linear operator.
By definition, the operator $L$ can be represented as
\begin{align}
	L = P_{W,E}\mathcal{L}_{|W}.
	\label{eq. normal L is Pcurl_L}
\end{align}

\begin{comment}
Noticing that $g$ can be represented as 
$g=\sum_{j=1}^{N}\alpha_j \psi_j$ for some $\alpha_j\in \mathbb{R}$, $j=1,\cdots, N$,
\begin{align}
	\eqref{eq. direct decomposition each element after applying projection}
	\ \ 
	\Leftrightarrow
	\ \
	& 
	\sum_{j=1}^{N}
	\alpha_j
	\left\{
		\left(
		P_{W,E}\mathcal{L}
		\right)\psi_j
		-
		L\psi_j
	\right\}
	=
	0,&
	\nonumber
\\
	\Leftrightarrow\ \ 
	&
	\left(
	P_{W,E}\mathcal{L}
	\right)\psi_j
	-
	L\psi_j
	=0,
	\ \ 
	j=1,\cdots,N.
	&
\end{align}
\end{comment}

Herein, we aim to find this finite dimensional operator $L$ since
$\mathcal{L}_{|W}=L$
if and only if $W$ is an invariant subspace.
Let $\{ \psi_i \}_{i=1}^N$ be a basis of $W$, where $N=\text{dim}(W)$.
If we take $g=\psi_i$ in \eqref{eq. direct decomposition each element},
there exists a unique element $e_i\in E$ s.t.
\begin{align}
	\mathcal{L}_{|W}\psi_i=L\psi_i+e_i.
\end{align}

Furthermore, $L\psi_i\in W$ can be represented as a linear combination of the basis $\{\psi_i\}_{i=1}^N$, that is, for each 
$\psi_i\in W$, there exist unique coefficients $l_{ij}\in \mathbb{R}$, $j=1,\cdots,N$, and a unique element $e_i\in E$ s.t.
\begin{align}
	\mathcal{L}_{|W}\psi_i=
	\sum_{j=1}^{N} l_{ij} \psi_j
	+
	e_i,
	\ \ 
	i=1,\cdots,N.
	\label{eq. direct decomposition}
\end{align}

Note that $\bm{L}^\tr$ denotes a matrix representation of $L$, where $\bm{L}:=(l_{ij})\in \mathbb{R}^{N\times N}$.
\begin{comment}
Specifically, 
\begin{align}
	\eqref{eq. direct decomposition}
	\Rightarrow
	L\psi_i = \sum_{j=1}^{N} l_{ij} \psi_j,
	\ \ 
	(\because \eqref{eq. normal L is Pcurl_L})
    \label{eq. equation for showing matrix rep}
\end{align}
and
for $\forall g=\sum_{i=1}^{N}\alpha_i \psi_i \in W$ ($\alpha_i\in \mathbb{R}$),
\begin{align}
	Lg&=\sum_{i=1}^{N}\alpha_i (L\psi_i)
	\ \ (\because \text{ Linearity of }L)
	&\nonumber
\\
	&=
	\sum_{i=1}^{N}\alpha_i \sum_{j=1}^{N} l_{ij} \psi_j
        \ \ (\because \eqref{eq. equation for showing matrix rep})
	&\nonumber
\\
	&=
	\sum_{j=1}^{N} 
	\left(
		\sum_{i=1}^{N}\alpha_i l_{ij} 
	\right)\psi_j,&
\end{align}
whereas there also exist coefficients $\hat{\alpha}_j\in \mathbb{R}$, $j=1,\cdots,N$, s.t.
\begin{align}
	Lg=\sum_{j=1}^{N} \hat{\alpha}_j \psi_j.
\end{align}

Thus, both coefficients $\alpha_i$ and $\hat{\alpha}_i$ are related by:
\begin{align}
	\left[
		\begin{array}{c}
			\hat{\alpha}_1
		\\
			\vdots 
		\\
			\hat{\alpha}_N
		\end{array}
	\right]
	=
	\left[
		\begin{array}{c}
			\sum_{i=1}^{N}\alpha_i l_{i1}
		\\
			\vdots 
		\\
			\sum_{i=1}^{N}\alpha_i l_{iN}
		\end{array}
	\right]
	=
	\left[
		\begin{array}{ccc}
			l_{11} & \cdots & l_{N1}
		\\
			\vdots & & \vdots 
		\\
		 	l_{1N} & \cdots & l_{NN}
		\end{array}
	\right]
	\left[
		\begin{array}{c}
			\alpha_1
		\\
			\vdots 
		\\
			\alpha_N
		\end{array}
	\right]
	=
	\bm{L}^\tr 
	\left[
	\begin{array}{c}
		\alpha_1
		\\
		\vdots 
		\\
		\alpha_N
	\end{array}
	\right].
\end{align}
\end{comment}
A solution to the relaxed problem of finding $\mathcal{L}_{|W}$ can be reduced to seeking $W$ such that it is invariant under $\mathcal{L}$ and computing a matrix representation $\bm{L}^\tr$ of $L$.

To obtain an explicit representation of $\bm{L}^\tr$, 
we introduce a weak formulation: 
let $\{ \varphi_l \}_{l=1}^{\hat{N}}$ be a basis of another finite dimensional subspace $\hat{W}\subset H$ s.t. $\hat{W}=E^\perp$ and
take inner products of both sides of \eqref{eq. direct decomposition} with $\varphi_l$, $l=1,\cdots,\hat{N}$, so that
\begin{align}
        \eqref{eq. direct decomposition} 
        \Rightarrow
        &
	\langle \mathcal{L}_{|W}\psi_i, \varphi_l \rangle=
	\sum_{j=1}^{N} l_{ij} \langle \psi_j, \varphi_l \rangle
	+
	\langle e_i, \varphi_l \rangle
	=
	\left[
		l_{i1}\cdots l_{iN}
	\right]
	\left[
		\begin{array}{c}
			\langle \psi_1, \varphi_l \rangle
		\\
			\vdots 
		\\
			\langle \psi_N, \varphi_l \rangle
		\end{array}
	\right]
	+
	\langle e_i, \varphi_l \rangle,
	&\nonumber
\\
        &\hspace{10cm}
        i=1,\cdots,N,
	\ 
	l=1,\cdots,\hat{N},
        &\nonumber 
%\\ 
%	\Leftrightarrow&
%	\left[
%		\begin{array}{c}
%			\langle \mathcal{L}_{|W}\psi_1, \varphi_l \rangle
%		\\
%			\vdots 
%		\\
%			\langle \mathcal{L}_{|W}\psi_N, \varphi_l \rangle
%		\end{array}
%	\right]
%	=
%	\left[
%		\begin{array}{ccc}
%			l_{11}&\cdots & l_{1N}
%		\\
%			\vdots &  & \vdots 
%		\\	
%			l_{N1}&\cdots & l_{NN}	
%		\end{array}
%	\right]
%	\left[
%	\begin{array}{c}
%		\langle \psi_1, \varphi_l \rangle
%		\\
%		\vdots 
%		\\
%		\langle \psi_N, \varphi_l \rangle
%	\end{array}
%	\right]
%	,
%	\ \ 
%	l=1,\cdots,\hat{N},
%	&\nonumber
%\\
%	\Leftrightarrow&
%	\left[
%	\begin{array}{ccc}
%		\langle \mathcal{L}_{|W}\psi_1, \varphi_1 \rangle & \cdots & \langle \mathcal{L}_{|W}\psi_1, \varphi_{\hat{N}} \rangle
%	\\
%		\vdots & & \vdots  
%	\\
%		\langle \mathcal{L}_{|W}\psi_N, \varphi_1 \rangle & \cdots & \langle \mathcal{L}_{|W}\psi_N, \varphi_{\hat{N}} \rangle
%	\end{array}
%	\right]
%	\hspace{-1.5mm}=\hspace{-1.5mm}
%	\left[
%	\begin{array}{ccc}
%		l_{11}&\cdots & l_{1N}
%		\\
%		\vdots &  & \vdots 
%		\\	
%		l_{N1}&\cdots & l_{NN}	
%	\end{array}
%	\right]
%        \hspace{-2mm}
%	\left[
%	\begin{array}{ccc}
%		\hspace{-1mm}\langle \psi_1, \varphi_1 \rangle & \cdots & \langle \psi_1, \varphi_{\hat{N}} \rangle\hspace{-1mm}
%		\\
%		\vdots & & \vdots  
%		\\
%		\hspace{-1mm}\langle \psi_N, \varphi_1 \rangle & \cdots & \langle \psi_N, \varphi_{\hat{N}} \rangle\hspace{-1mm}
%	\end{array}
%	\right],
%	&\nonumber
\\
	\Leftrightarrow
	&
        \left(
		\langle \mathcal{L}_{|W} \psi_i, \varphi_j \rangle
	\right)
	=
	\bm{L}
	\left(
		\langle \psi_i, \varphi_j \rangle 
	\right).&
	\label{eq. weak formulation}
\end{align}

Note that the weak formulation \eqref{eq. weak formulation} is a necessary condition for \eqref{eq. direct decomposition} to hold.

\begin{comment}
\begin{figure}[b]
	\centering 
	\includegraphics[width=12cm]{figs/projection.png}
	\caption{Finite-dimensional subspace and projections onto it.}
\end{figure}
\end{comment}

\begin{remark}
\label{remark orthogonal projection}
\rm{} 
If we choose $E=W^\perp$, i.e., adopt the orthogonal direct decomposition, we have $\hat{W}=W$ and $\varphi_i=\psi_i$ for $\forall i=1,\cdots,N$.
Also, if $(\langle \psi_i, \psi_j \rangle)$ is invertible, 
$\hat{\bm{L}}:=(\langle \mathcal{L}_{|W} \psi_i, \psi_j \rangle)(\langle \psi_i, \psi_j \rangle)^{-1}$
gives a necessary condition for $\hat{\bm{L}}$ to be
a matrix representation of the operator $L=P_{W,W^{\perp}}\mathcal{L}_{|W}$.
Furthermore, $L$ is in this case an optimal approximation of $\mathcal{L}_{|W}$ in the following sense:
\begin{align}
	Lg = \underset{h\in W}{\text{argmin}}\ 
	\| \mathcal{L}_{|W}g - h  \|,
	\ \ 
	\forall g\in W.
	\label{eq. optimality in a Hilbert space}
\end{align}

\end{remark}
\begin{comment}
\takespace
%\begin{screen}
\begin{prop}
    \label{prop. optimality of orthogonal projection in a Hilbert space}
    \rm{}
    Suppose that the orthogonal direct decomposition is adopted in \eqref{eq. direct decomposition in a general form} so that $E=W^\perp$.
    Then, $L$ gives an optimal approximation of $\mathcal{L}_{|W}$ in the following sense:
    \begin{align}
        Lg = \underset{h\in W}{\text{argmin}}\ 
	\| \mathcal{L}_{|W}g - h  \|,
	\ \ 
	\forall g\in W.
        \label{eq. optimality in a Hilbert space}
    \end{align}
\end{prop}
%\end{screen}
\begin{proof}
    See Appendix \ref{appendix proof of optimality of orthogonal projection in a Hilbert space}.
\end{proof}
\end{comment}

\subsection{Finite Data Approximation of $L$}
\label{sec. finite data approximation}
In practice, it may not be possible to compute the exact values of $\langle \psi_i, \varphi_j \rangle$ or $\langle \mathcal{L}_{|W}\psi_i, \varphi_j \rangle$ in the weak formulation \eqref{eq. weak formulation} but one may be able to compute their approximations, which yield a finite data approximation of $L$, or equivalently $\bm{L}$.

In this section, it is assumed that $H=L_{2}(\mathbb{X})$, i.e., an $L_2$ space of functions defined on some set $\mathbb{X}$.
Given a basis $\{ \psi_i \}_{i=1}^N$ of the finite dimensional subspace $W$ of $H$, suppose that the following set of data points is available from either experiments or simulations:
\begin{align}
	\{ 
		\bm{\psi}(x_l), (\mathcal{L}_{|W}\bm{\psi})(x_l),
		\bm{\varphi}(x_l)
		\mid x_l\in \mathbb{X},\ l=1,\cdots,M 
	\},
	\label{eq. data set}
\end{align}
where
\begin{align}
	\bm{\psi}(x_l)&:=
	\left[
		\begin{array}{ccc}
			\psi_1(x_l) & \cdots & \psi_N(x_l)
		\end{array}
	\right]^\tr\in \mathbb{R}^{N},
	&
	\label{eq. def of data vector psi}
\\
	(\mathcal{L}_{|W}\bm{\psi})(x_l)&:=
	\left[
	\begin{array}{ccc}
		(\mathcal{L}_{|W}\psi_1)(x_l) & \cdots & (\mathcal{L}_{|W}\psi_N)(x_l)
	\end{array}
	\right]^\tr\in \mathbb{R}^{N},
	&
	\label{eq. def of data vector L psi}
\\
	\bm{\varphi}(x_l)&:=
	\left[
	\begin{array}{ccc}
		\varphi_1(x_l) & \cdots & \varphi_{\hat{N}}(x_l)
	\end{array}
	\right]^\tr\in \mathbb{R}^{\hat{N}}.
	&
	\label{eq. def of data vector varphi}
\end{align}

If we adopt the empirical measure defined as
\begin{align}
	\mu_{M}:=
	\frac{1}{M} \sum_{i=1}^{M} \delta_{x_i,}
	\label{eq. def. empirical measure}
\end{align}
where $\delta_{x_i}$ denotes the Dirac measure at point $x_i$, we have
\begin{align}
	\langle \psi_i,\varphi_j \rangle_{L_2(\mu_{M})}
	&=
	\int_{\mathbb{X}} \psi_i\varphi_j d\mu_{M}
	=
	\frac{1}{M} \sum_{l=1}^{M} \psi_i(x_l)\varphi_j(x_l),
	&
\\
	\langle \mathcal{L}_{|W}\psi_i,\varphi_j \rangle_{L_2(\mu_{M})} 
	&=
	\int_{\mathbb{X}} (\mathcal{L}_{|W}\psi_i)\varphi_j d\mu_{M}
	=
	\frac{1}{M} \sum_{l=1}^{M} (\mathcal{L}_{|W}\psi_i)(x_l)\varphi_j(x_l),
	&
\end{align}
for $i=1,\cdots,N$ and $j=1,\cdots,\hat{N}$, where the notation $\langle \cdot,\cdot \rangle_{L_2(\mu_M)}$ is used to explicitly indicate that the specific measure \eqref{eq. def. empirical measure} is taken.
Therefore, the matrices appeared in \eqref{eq. weak formulation} are in this case reduced to:
\begin{align}
	(\langle \psi_i,\varphi_j \rangle_{L_2(\mu_{M})})
%	&=
%	\frac{1}{M} \sum_{l=1}^{M}
%	\left[
%		\begin{array}{ccc}
%			\psi_1(x_l)\varphi_1(x_l) & \cdots & \psi_1(x_l)\varphi_{\hat{N}}(x_l)
%		\\
%			\vdots && \vdots 
%		\\
%			\psi_N(x_l)\varphi_1(x_l) & \cdots & \psi_N(x_l)\varphi_{\hat{N}}(x_l)
%		\end{array}
%	\right]
%	&\nonumber
%\\
%	&=
%	\frac{1}{M} \sum_{l=1}^{M}
%	\left[
%		\begin{array}{c}
%			\psi_1(x_l)
%		\\
%			\vdots 
%		\\
%			\psi_N(x_l)
%		\end{array}
%	\right]
%	\left[
%		\begin{array}{ccc}
%			\varphi_1(x_l) &\cdots & \varphi_{\hat{N}}(x_l)
%		\end{array}
%	\right]
%	&\nonumber
%\\
	&=
	\frac{1}{M} \sum_{l=1}^{M}
	\bm{\psi}(x_l) \bm{\varphi}(x_l)^\tr,
	&
\\
	( \langle \mathcal{L}_{|W}\psi_i, \varphi_j \rangle_{L_2(\mu_{M})} )
	&=
	\frac{1}{M} \sum_{l=1}^{M}
	(\mathcal{L}_{|W}\bm{\psi})(x_l) \bm{\varphi}(x_l)^\tr,
	&
\end{align}
which can be computed from the data set \eqref{eq. data set}.
In the sequel, the notation $\bm{L}_M$ is used to refer to the approximated matrix representation of $L$ for which the empirical measure \eqref{eq. def. empirical measure} is adopted, i.e.,
\begin{align}
	\bm{L}_M
	:=&
	(\langle \mathcal{L}_{|W}\psi_i, \varphi_j\rangle_{L_2(\mu_{M})}) (\langle \psi_i, \varphi_j\rangle_{L_2(\mu_{M})})^{\dagger}
	&\nonumber
\\
	=&
	\frac{1}{M} \sum_{l=1}^{M}
	(\mathcal{L}_{|W}\bm{\psi})(x_l) \bm{\varphi}(x_l)^\tr
	\left\{
		\frac{1}{M} \sum_{l=1}^{M}
		\bm{\psi}(x_l) \bm{\varphi}(x_l)^\tr
	\right\}^{\dagger},
	&
	\label{eq. def of L_m}
\end{align}
where $\dagger$ denotes the Moore-Penrose pseudo inverse.
When $N=\hat{N}$ and the matrix\\
$(\langle \psi_i, \varphi_j\rangle_{L_2(\mu_{M})})$
is invertible, the pseudo inverse can be replaced with the matrix inverse; otherwise, the pseudo inverse is introduced
as the least-square or minimum-norm approximate solution to the problem of finding $\bm{L}$ in \eqref{eq. weak formulation}.

Note that the finite-data approximation $\bm{L}_M$ may need to be compared with $\bm{L}$ derived with the original measure $\mu$, which is supposed to be consistent with actual scenarios in experiments or simulations, e.g., a uniform distribution supported on an interval corresponding to the range of values that can be measured by experimental devices. 
In general, increasing the number $M$ of data points is advisable to mitigate the discrepancy originating from this finite-data approximation, as seen in the following remark.

\begin{remark}
\rm{}
Let $\psi_i:\mathbb{X}\rightarrow \mathbb{R}$ and $\varphi_i:\mathbb{X}\rightarrow \mathbb{R}$ be measurable functions for $\forall i$.
If $\{ x_l \}_{l=1}^M$ in the data set \eqref{eq. data set} are realizations of independent and identically distributed random variables $X_l$ ($l=1,\cdots,M$), $\psi_i(X_l)$ and $\varphi_i(X_l)$ are also random variables for $\forall l,i$.
By the strong law of large numbers, the following holds with probability one:
\begin{align}
	\langle \psi_i,\varphi_j \rangle_{L_2(\mu_{M})}
	&\rightarrow 
	\langle \psi_i,\varphi_j \rangle_{L_2(\mu)}
	\ 
	(M\rightarrow \infty),
	&
\end{align}
where $\mu$ denotes a probability measure that corresponds to the probability distribution from which $X_l$ are sampled.
\end{remark}

\begin{ex}
\label{example EDMD}
\rm{}
The method to compute $\bm{L}_M$ described in this section is equivalent to EDMD when $\mathcal{L}$ is a Koopman operator, $\hat{N}=N$, and $\varphi_i=\psi_i$ for $\forall i=1,\cdots,N$.
This fact is easily inferred from the following proposition.
\end{ex}

%\takespace
\begin{screen}
\begin{prop}
\label{prop. optimality of orthogonal projection in the finite data setting}
\rm{}
Suppose that $N<M$ and
$\left[
	\begin{array}{ccc}
		\bm{\psi} (x_1)&\cdots&\bm{\psi} (x_M)
	\end{array}
	\right]\in \mathbb{R}^{N\times M}$
has full rank.
If the orthogonal direct decomposition is adopted in \eqref{eq. direct decomposition in a general form} so that $E=W^\perp$, $\hat{N}=N$, and $\varphi_i=\psi_i$ for $\forall i=1,\cdots,N$, the finite data approximation $\bm{L}_M$ in \eqref{eq. def of L_m} is the least-square solution to the problem:
\begin{align}
	\underset{\hat{\bm{L}}\in \mathbb{R}^{N\times N}}{\text{min}}
	\sum_{l=1}^{M}
	\| (\mathcal{L}_{|W}\bm{\psi})(x_l) - \hat{\bm{L}}\bm{\psi} (x_l)\|_2^2.
	\label{eq. least square problem}
\end{align}
\end{prop}
\end{screen} 
\begin{proof}
See Appendix \ref{appendix proof of prop on optimality of orthogonal projection}.
\end{proof}

\subsection{Linear Embedding Model Revisited}
In this section, we revisit the linear embedding model \eqref{eq. approx equation of the model with state obs included} and derive the same type of model representation from the projection-based linear operator-learning formalism developed in previous sections.
This new perspective enables us to extend the model structure of EDMD by introducing the test functions $\varphi_i$ into the model parameters $[\bm{A}\ \bm{B}]$, which are expected to increase the expressivity of the model as well as improve its accuracy and generalizability.

Recall that the unknown non-autonomous dynamics to be modeled was given by
\begin{align}
	\chi^+=F(\chi,u),\ \ 
	\chi\in \mathcal{X}\subseteq \mathbb{R}^n,\  
	u\in \mathcal{U}\subseteq \mathbb{R}^p, \
	F:\mathcal{X}\times \mathcal{U}\rightarrow \mathcal{X},
	\tag{\ref{eq. discrete time dynamical system}}
\end{align}
and the space of infinite sequences of admissible inputs was also introduced to define the Koopman operator for this system s.t.
\begin{align}
	l(\mathcal{U}):=
	\left\{
	U:\mathbb{Z}_{\geq 0}\rightarrow \mathcal{U}:k\mapsto u_k
	\right\}
	=
	\{
	(u_0,u_1,\cdots)\mid u_k\in \mathcal{U}, k\in \mathbb{Z}_{\geq 0}
	\}.
	\tag{\ref{eq. def of space of infinite sequences of inputs}}
\end{align}

From this point on, we assume $\mathbb{X}=\mathcal{X}\times l(\mathcal{U})$ so that
the Hilbert space of interest is defined as $H=L_2(\mathcal{X}\times l(\mathcal{U}))$ and the linear operator $\mathcal{L}$ is the Koopman operator $\hat{\mathcal{K}}$, which was defined in \eqref{eq. def of Koopman operator for non autonomous setting}.
\begin{comment}
From this point on, we assume $\mathbb{X}=\mathcal{X}\times l(\mathcal{U})$, where $\mathcal{X}\subseteq \mathbb{R}^n$ is a state-space in which a discrete-time dynamical system with control inputs:
\begin{align}
	\chi^+=F(\chi,u),\ \ 
	\chi\in \mathcal{X}\subseteq \mathbb{R}^n,\  
	u\in \mathcal{U}\subseteq \mathbb{R}^p, \
	F:\mathcal{X}\times \mathcal{U}\rightarrow \mathcal{X},
	\label{eq. discrete time dynamical system}
\end{align}
is defined. The set $\mathcal{U}$ is a collection of possible inputs $u$ for the system and $l(\mathcal{U})$ denotes a set of possible input \textit{signals}, i.e.,
\begin{align}
	l(\mathcal{U}):=
	\left\{
		U:\mathbb{Z}_{\geq 0}\rightarrow \mathcal{U}:k\mapsto u_k
	\right\}
	=
	\{
		(u_0,u_1,\cdots)\mid u_k\in \mathcal{U}, k\in \mathbb{Z}_{\geq 0}
	\}.
\end{align}

Also, we take $\mathcal{L}$ as the Koopman operator associated with a map $\hat{F}$, i.e., $\mathcal{L}$ satisfies
\begin{align}
	\mathcal{L}h = h\circ \hat{F},\ \ 
	\forall h\in L_2(\mathcal{X}\times l(\mathcal{U})),
\end{align}
where $\hat{F}$ is defined by
\begin{align}
	\hat{F}:\mathcal{X}\times l(\mathcal{U})\rightarrow \mathcal{X}\times l(\mathcal{U}):
	(\chi, U)\mapsto (F(\chi,U(0)), \mathcal{S}U),
	\label{eq. def of Fhat}
\end{align}
and $\mathcal{S}$ denotes the shift operator s.t.
\begin{align}
	\mathcal{S}U=
	\mathcal{S}(u_0,u_1,\cdots)
	=
	(u_1,u_2,\cdots).
\end{align}

In \eqref{eq. def of Fhat}, the notation $U(0)$ refers to the first element of $U$.
\end{comment}

Similar to \eqref{eq. def of each h_i for non autonomous setting Koopman}, if we impose a special structure on the basis functions $\psi_i\in H$ s.t.
\begin{align}
	\bm{\psi}(\chi, U)=
	\left[
		\begin{array}{c}
			\psi_1(\chi, U)
		\\
			\vdots 
		\\
			\psi_{N_x}(\chi, U)
		\\
			\psi_{N_x+1}(\chi, U)
		\\
			\vdots 
		\\
			\psi_{N_x+p}(\chi, U)
		\end{array}
	\right]
	=
	\left[
	\begin{array}{c}
		g_1(\chi)
		\\
		\vdots 
		\\
		g_{N_x}(\chi)
		\\
		U(0)
	\end{array}
	\right]
	\in \mathbb{R}^{N_x + p},
	\label{eq. special psi}
\end{align}
where $g_i:\mathcal{X}\rightarrow \mathbb{R}$ are functions from the state-space $\mathcal{X}$ into $\mathbb{R}$,
the action of the Koopman operator $\mathcal{L}$ is characterized as follows:
\begin{align}
	(\mathcal{L}\bm{\psi})(\chi, U)
	=
	\left[
	\begin{array}{c}
		(\mathcal{L}\psi_1)(\chi, U)
		\\
		\vdots 
		\\
		(\mathcal{L}\psi_{N_x})(\chi, U)
		\\
		(\mathcal{L}\psi_{N_x+1})(\chi, U)
		\\
		\vdots 
		\\
		(\mathcal{L}\psi_{N_x+p})(\chi, U)
	\end{array}
	\right]
	=
	\left[
	\begin{array}{c}
		(\psi_1\circ \hat{F})(\chi, U)
		\\
		\vdots 
		\\
		(\psi_{N_x}\circ \hat{F})(\chi, U)
		\\
		(\psi_{N_x+1}\circ \hat{F})(\chi, U)
		\\
		\vdots 
		\\
		(\psi_{N_x+p}\circ \hat{F})(\chi, U)
	\end{array}
	\right]
	=
	\left[
		\begin{array}{c}
			g_1(\chi^+)
		\\
			\vdots 
		\\
			g_{N_x}(\chi^+)
		\\
			U(1)
		\end{array}
	\right]
	\ (\because \eqref{eq. def of Koopman operator for non autonomous setting}),
	\label{eq. action of L to h}
\end{align}
where $\chi^+$ denotes the successor of $\chi$ and $U(0)$ so that $\chi^+=F(\chi, U(0))$.
The direct decomposition \eqref{eq. direct decomposition} reads
\begin{align}
	&\left[
		\begin{array}{c}
			(\mathcal{L}\psi_1)(\chi, U)
		\\
			\vdots 
		\\
			(\mathcal{L}\psi_{N_x})(\chi, U)
		\\
			(\mathcal{L}\psi_{N_x+1})(\chi, U)
		\\
			\vdots 
		\\
			(\mathcal{L}\psi_{N_x+p})(\chi, U)
		\end{array}
	\right]
	=
	\bm{L}
	\left[
		\begin{array}{c}
			\psi_1(\chi, U)
		\\
			\vdots 
		\\
			\psi_{N_x}(\chi, U)
		\\	
			\psi_{N_x+1}(\chi, U)
		\\
			\vdots 
		\\
			\psi_{N_x+p}(\chi, U)
		\end{array}
	\right]
	+
	\left[
		\begin{array}{c}
			e_1(\chi, U)
		\\
			\vdots 
		\\
			e_{N_x+p}(\chi, U)
		\end{array}
	\right],&\nonumber
\\
	\Leftrightarrow\ &
	\left[
	\begin{array}{c}
		g_1(\chi^+)
		\\
		\vdots 
		\\
		g_{N_x}(\chi^+)
		\\
		U(1)
	\end{array}
	\right]
	=
	\bm{L}
	\left[
	\begin{array}{c}
		g_1(\chi)
		\\
		\vdots 
		\\
		g_{N_x}(\chi)
		\\
		U(0)
	\end{array}
	\right]
	+
	\left[
	\begin{array}{c}
		e_1(\chi, U)
		\\
		\vdots 
		\\
		e_{N_x+p}(\chi, U)
	\end{array}
	\right].
	\ \ 
	(\because \eqref{eq. special psi}, \eqref{eq. action of L to h})
	&
\end{align}

Letting $[\bm{A}\ \bm{B}]\in \mathbb{R}^{N_x\times (N_x+p)}$ be the first $N_x$ rows of $\bm{L}$, we have
\begin{align}
	\left[
	\begin{array}{c}
		g_1(\chi^+)
		\\
		\vdots 
		\\
		g_{N_x}(\chi^+)
	\end{array}
	\right]
	=
	\bm{A}
	\left[
		\begin{array}{c}
			g_1(\chi)
		\\
			\vdots 
		\\
			g_{N_x}(\chi)
		\end{array}
	\right]
	+
	\bm{B}
	u
	+
	\left[
	\begin{array}{c}
		e_1(\chi, u)
		\\
		\vdots 
		\\
		e_{N_x}(\chi, u)
	\end{array}
	\right],
	\label{eq. discrete time dynamics in a new coordinate}
\end{align} 
where $u:=U(0)$.
Note that we have $e_i(\chi,U)=e_i(\chi,u)$ in the above equation since only the first element of $U$ depends on the definitions of $\psi_i$ as in \eqref{eq. special psi}.
If we specify\\ $[g_1(\chi)\cdots g_n(\chi)]^\tr := \chi$, this equation is written as
\begin{align}
	\left[
	\begin{array}{c}
		\chi^+
		\\
		g_{n+1}(\chi^+)
		\\
		\vdots 
		\\
		g_{N_x}(\chi^+)
	\end{array}
	\right]
	=
	\bm{A}
	\left[
	\begin{array}{c}
		\chi
		\\
		g_{n+1}(\chi)
		\\
		\vdots 
		\\
		g_{N_x}(\chi)
	\end{array}
	\right]
	+
	\bm{B}
	u
	+
	\left[
	\begin{array}{c}
		e_1(\chi, u)
		\\
		\vdots 
		\\
		e_{N_x}(\chi, u)
	\end{array}
	\right],
	\label{eq. linear embedding model with errors from oblique projection}
\end{align}
which corresponds to the linear embedding model \eqref{eq. approx equation of the model with state obs included}.
	Note that while this paper focuses on linear embedding models so that linear control theories can be readily utilized, other approaches such as bilinear models can be also derived by imposing appropriate structure on the basis functions $\psi_i$.

The error terms are explicitly included in \eqref{eq. linear embedding model with errors from oblique projection} and the model parameters $[\bm{A}\ \bm{B}]$ have a specific interpretation in this case, i.e., they are realized as (the upper rows of) the transpose $\bm{L}$ of a matrix representation of the linear operator $L$ defined in \eqref{eq. normal L is Pcurl_L}.
It is also noted that the error terms $e_i(\chi,u)$ are elements of the complementary subspace $E$ introduced in \eqref{eq. direct decomposition in a general form}.

In the proposed method, the finite data approximation $\bm{L}_M$ in \eqref{eq. def of L_m} is computed to estimate $[\bm{A}\ \bm{B}]$ in the above equation. 
Specifically, 
letting $[\bm{A}\ \bm{B}]$ denote the first $N_x$ rows of $\bm{L}_M$, we have
\begin{align}
	[\bm{A}\ \bm{B}]=&
	\frac{1}{M} \sum_{l=1}^{M}
	\left[
		\begin{array}{c}
			g_{1}(y_l)
		\\
			\vdots
		\\
			g_{N_x}(y_l)
		\end{array}
	\right]
	\left[
	\begin{array}{c}
		\varphi_1(\chi_l,u_l)
		\\
		\vdots 
		\\
		\varphi_{\hat{N}}(\chi_l,u_l)
	\end{array}
	\right]^\tr
	\left\{
	\frac{1}{M} \sum_{l=1}^{M}
	\left[
	\begin{array}{c}
		g_{1}(\chi_l)
		\\
		\vdots 
		\\
		g_{N_x}(\chi_l)
		\\
		u_l
	\end{array}
	\right] 
	\left[
	\begin{array}{c}
		\varphi_1(\chi_l,u_l)
	\\
		\vdots 
	\\
		\varphi_{\hat{N}}(\chi_l,u_l)
	\end{array}
	\right]^\tr
	\right\}^{\dagger}
	&\nonumber
%\\
%	=&
%	\frac{1}{M} \sum_{l=1}^{M}
%	\left[
%		\begin{array}{ccc}
%			g_1(y_l)\varphi_1(\chi_l,u_l) & \cdots & g_1(y_l)\varphi_{\hat{N}}(\chi_l,u_l)
%		\\
%			\vdots & \ddots & \vdots 
%		\\
%			g_{N_x}(y_l)\varphi_1(\chi_l,u_l) & \cdots & g_{N_x}(y_l)\varphi_{\hat{N}}(\chi_l,u_l)
%		\end{array}
%	\right]
%	&\nonumber
%\\
%	&\hspace{25mm}
%	\times 
%	\left\{
%		\frac{1}{M} \sum_{l=1}^{M}
%		\left[
%		\begin{array}{ccc}
%			g_1(\chi_l)\varphi_1(\chi_l,u_l) & \cdots & g_1(\chi_l)\varphi_{\hat{N}}(\chi_l,u_l)
%			\\
%			\vdots & \ddots & \vdots 
%			\\
%			g_{N_x}(\chi_l)\varphi_1(\chi_l,u_l) & \cdots & g_{N_x}(\chi_l)\varphi_{\hat{N}}(\chi_l,u_l)
%			\\
%			\varphi_1(\chi_l,u_l) u_l & \cdots & \varphi_{\hat{N}}(\chi_l,u_l) u_l
%		\end{array}
%		\right]
%	\right\}^\dagger 
%	&\nonumber
\\
	=&
	\left[
		\begin{array}{ccc}
			g_1(y_1) & \cdots & g_1(y_M)
		\\	
			\vdots & \ddots & \vdots 
		\\
			g_{N_x}(y_1) & \cdots & g_{N_x}(y_M)
		\end{array}
	\right]
	\left[
		\begin{array}{ccc}
			\varphi_1(\chi_1,u_1) & \cdots & \varphi_1(\chi_M,u_M)
		\\
			\vdots & \ddots & \vdots 
		\\
			\varphi_{\hat{N}}(\chi_1,u_1) & \cdots & \varphi_{\hat{N}}(\chi_M,u_M)
		\end{array}
	\right]^\tr 
	&\nonumber
\\	
	&
	\hspace{10mm}
	\times 
	\left(
	\left[
		\begin{array}{ccc}
			g_1(\chi_1) & \cdots & g_1(\chi_M)
		\\
					\vdots & \ddots & \vdots 
		\\	
			g_{N_x}(\chi_1) & \cdots & g_{N_x}(\chi_M)
		\\
			u_1 & \cdots & u_M
		\end{array}
	\right]
	\left[
	\begin{array}{ccc}
		\varphi_1(\chi_1,u_1) & \cdots & \varphi_1(\chi_M,u_M)
		\\
		\vdots & \ddots & \vdots 
		\\
		\varphi_{\hat{N}}(\chi_1,u_1) & \cdots & \varphi_{\hat{N}}(\chi_M,u_M)
	\end{array}
	\right]^\tr 
	\right)^\dagger.
	&
	\label{eq. oblique EDMD}
\end{align}

\begin{figure}[]
    \centering
    \includegraphics[width=0.9\linewidth]{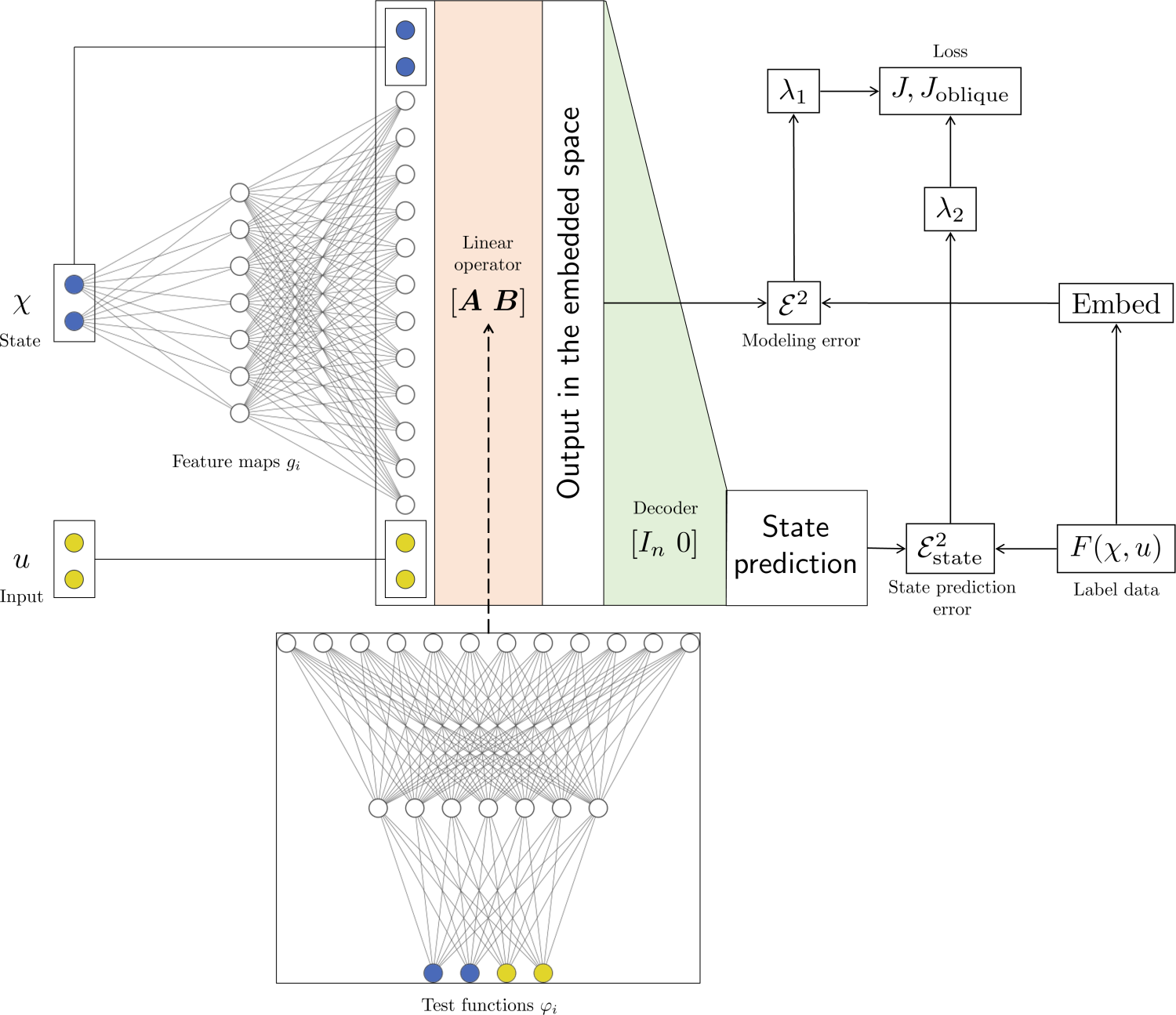}
    \caption{Architecture of the proposed linear embedding model with oblique projection.}
    \label{fig. proposed model architecture}
\end{figure}

This can be viewed as an extension of the model parameters $[\bm{A}\ \bm{B}]$ given by EDMD since substituting $\varphi_i=g_i$ for $\forall i$ with $\hat{N}=N_x$ recovers the EDMD solution \eqref{eq. intro [A B] as in EDMD}.
Figure \ref{fig. proposed model architecture} shows a schematic of the proposed linear embedding model with oblique projection.

\section{Learning Procedures}
\label{sec. learning procedures}
In this section, problems of learning model parameters are formulated.
For the linear embedding model \eqref{eq. linear embedding model with errors from oblique projection}, suppose that a data set is given as $\{ (\chi_i,u_i,y_i)\mid y_i=F(\chi_i,u_i),i=1,\cdots,M \}$. 
The data set of the form \eqref{eq. data set}, which was introduced in the context of linear operator learning in Section \ref{sec. finite data approximation}, is related to this model as follows:
\begin{align}
	&\{ 
	\bm{\psi}(x_l), (\mathcal{L}_{|W}\bm{\psi})(x_l),
	\bm{\varphi}(x_l)
	\mid x_l\in \mathbb{X},\ l=1,\cdots,M 
	\}
	&\nonumber
\\
	=&
	\left\{
	[\chi_l^\tr\ g_{n+1}(\chi_l)\cdots g_{N_x}(\chi_l)\ u_l^\tr]^\tr,\ 
	[y_l^\tr\ g_{n+1}(y_l)\cdots g_{N_x}(y_l)]^\tr,\ 
	\bm{\varphi}(\chi_l, u_l)
        \right. 
        &\nonumber
\\
	&\hspace{8cm}
        \left. 
        \mid 
        y_l=F(\chi_l,u_l),l=1,\cdots,M 
	\right\},&\nonumber
\end{align}
where the basis functions $\psi_i$ are given by \eqref{eq. special psi} along with a condition $[g_1(\chi)\cdots g_n(\chi)]^\tr := \chi$.

Note that the linear operator-learning formulation developed in the previous sections assumes that the functions $g_i$ and $\varphi_i$ are given.
However, the proposed method treats these functions as model parameters to be learned from data.
Therefore, computing the approximation $\bm{L}_M$ in \eqref{eq. def of L_m}, which corresponds to the model parameters $[\bm{A}\ \bm{B}]$, needs to be incorporated into a generalized learning problem where $g_i$ and $\varphi_i$ are trained on data.

First, noticing the similarity to EDMD as explained in Example \ref{example EDMD}, we call the procedure to compute (the upper rows of) $\bm{L}_M$ the oblique EDMD as follows.
%\takespace
\begin{screen}
	\begin{prob}
		\label{problem oblique EDMD}
		\rm{} 
		(Oblique EDMD)\\
		Given feature maps $g_i:\mathcal{X}\rightarrow \mathbb{R}$ ($i=1,\cdots,N_x$) and  test functions $\varphi_i:\mathcal{X}\times l(\mathcal{U})\rightarrow \mathbb{R}$ ($i=1,\cdots,\hat{N}$), compute $[\bm{A}\ \bm{B}]$ in \eqref{eq. oblique EDMD}, which are then used as the model dynamics parameters of the linear embedding model \eqref{eq. linear embedding model with errors from oblique projection}.
	\end{prob}
\end{screen}
%\takespace

The model parameters are initialized by solving the following problem, in which the orthogonal projection is adopted and only the feature maps $g_i$ are the learnable parameters.
In the following problem statements, $\Theta_{g}$ denotes neural network parameters that parameterize functions $g_i$.
%\takespace
\begin{screen}
\begin{prob}
	\label{prob. initialization of model}
	\rm{} 
	(Initialization of model)\\ 
	Given a data set of the form $\{ (\chi_i,u_i,y_i)\mid y_i=F(\chi_i,u_i),i=1,\cdots,M \}$,
	find $\Theta_g^0$ s.t.
 $g_i^0(\cdot;\Theta_g^0)$ are a neural network and
	\begin{align}
		\Theta_g^0
		=
		\underset{\Theta_g^0\in \text{Neural network parameters}}{\text{argmin}}\
		J(\Theta_g^0),
	\end{align}
	where 
	\begin{align}
		J(\Theta_g^0):=&
		\sum_{i=1}^{M}
		\left(
			\lambda_1 
			\left\|
			\left[
			\begin{array}{c}
				y_i
			\\
				g_{n+1}^0(y_i;\Theta_g^0)
				\\
				\vdots 
				\\
				g_{N_x}^0(y_i;\Theta_g^0)
			\end{array}
			\right]
			-
			\bm{A}_0(\Theta_g^0)
			\left[
			\begin{array}{c}
				\chi_i
				\\
				g_{n+1}^0(\chi_i;\Theta_g^0)
				\\
				\vdots 
				\\	
				g_{N_x}^0(\chi_i;\Theta_g^0)
			\end{array}
			\right]
			-
			\bm{B}_0(\Theta_g^0)
			u_i
			\right\|_2^2
			\right.
			&\nonumber
		\\
			&
			\hspace{7mm}
			\left. 
			+
			\lambda_2 
			\left\|
					y_i
					-
					[I\ 0]
					\left(
						\bm{A}_0(\Theta_g^0)
						\left[
						\begin{array}{c}
							\chi_i
							\\
							g_{n+1}^0(\chi_i;\Theta_g^0)
							\\
							\vdots 
							\\	
							g_{N_x}^0(\chi_i;\Theta_g^0)
						\end{array}
						\right]
						+
						\bm{B}_0(\Theta_g^0)
						u_i
					\right)
			\right\|_2^2
		\right)&\nonumber 
        \\
            =&
            \sum_{i=1}^M \left\{
                \lambda_1 \mathcal{E}^2(\chi_i,u_i)
                +
                \lambda_2 \mathcal{E}_\text{state}^2(\chi_i,u_i)
            \right\},
            &
	\end{align}
	and
	$[\bm{A}_0(\Theta_g^0)\ \bm{B}_0(\Theta_g^0)]$ are computed by Problem \ref{problem oblique EDMD} with $\varphi_i=g_i^0$ for $\forall i$.
\end{prob}
\end{screen}
%\takespace

Note that the loss function $J$ is composed of the modeling error and the state prediction error, which were defined in \eqref{eq. def of modeling error} and \eqref{eq. def of state prediction error}, respectively. 
\begin{remark}
\label{remark Initialization of model is same as LKIS-DMD}
\rm{}
Problem \ref{prob. initialization of model} is mathematically equivalent to the implementation of LKIS-DMD \cite{Learning_Koopman_Invariant_Subspaces}. 
Adding the term involving the input $u$ to the model proposed in \cite{Learning_Koopman_Invariant_Subspaces} and setting $\lambda_1=1$, $\lambda_2=0$ in Problem \ref{prob. initialization of model} yield the same problem.
\end{remark}

After the model is initialized by Problem \ref{prob. initialization of model}, test functions $\varphi_i$ are also set as learnable parameters and the oblique projection is optimized on different data points.

%\takespace
\begin{screen}
	\begin{prob}
		\label{prob. proposed method: second training}
		\rm{}
		(Optimizing Oblique Projection)\\
		Suppose a data set $\{ (\chi_i,u_i,y_i)\mid y_i=F(\chi_i,u_i),i=1,\cdots,M \}$ is given, which is either a different one from that of Problem \ref{prob. initialization of model} or the same one with additional new data points added.
		Find $\Theta_g$ and $\Theta_\varphi$ s.t.
  $g_i(\cdot;\Theta_g)$ and
		$\varphi_i(\cdot;\Theta_\varphi)$ are neural networks and
		\begin{align}
			\{ \Theta_g, \Theta_\varphi \} 
			=
			\underset{\{ \Theta_g, \Theta_\varphi \}\in \text{Neural network parameters}}{\text{argmin}}\
			J_\text{oblique}(\Theta_g, \Theta_\varphi),
		\end{align}
		where 
		%$g_{n+1:N_x}:=[g_{n+1}\cdots g_{N_x}]$, $\varphi_{1:\hat{N}}:=[\varphi_1\cdots \varphi_{\hat{N}}]$,
		\begin{align}
			J_\text{oblique}(\Theta_g, \Theta_\varphi):=&
			\sum_{i=1}^{M}
			\left(
				\lambda_1
				\left\|
				\left[
				\begin{array}{c}
					y_i
					\\
					g_{n+1}(y_i;\Theta_g)
					\\
					\vdots 
					\\
					g_{N_x}(y_i;\Theta_g)
				\end{array}
				\right]
				-
				\bm{A}(\Theta_g, \Theta_\varphi)
				\left[
				\begin{array}{c}
					\chi_i
					\\ 
					g_{n+1}(\chi_i;\Theta_g)
					\\
					\vdots 
					\\	
					g_{N_x}(\chi_i;\Theta_g)
				\end{array}
				\right]
				-
				\bm{B}(\Theta_g, \Theta_\varphi)
				u_i
				\right\|_2^2
				\right. 
				&\nonumber
			\\
				&
				\hspace{7mm}
				\left. 
				+
				\lambda_2 
				\left\|
					y_i
					-
					[I\ 0]
					\left(
						\bm{A}(\Theta_g, \Theta_\varphi)
						\left[
						\begin{array}{c}
							\chi_i
							\\ 
							g_{n+1}(\chi_i;\Theta_g)
							\\
							\vdots 
							\\	
							g_{N_x}(\chi_i;\Theta_g)
						\end{array}
						\right]
						+
						\bm{B}(\Theta_g, \Theta_\varphi)
						u_i
					\right)
				\right\|_2^2
			\right)
			&\nonumber 
            \\
            =&
            \sum_{i=1}^M \left\{
                \lambda_1 \mathcal{E}^2(\chi_i,u_i)
                +
                \lambda_2 \mathcal{E}_\text{state}^2(\chi_i,u_i)
            \right\},
            &
		\end{align}
		and
		$[\bm{A}(\Theta_g, \Theta_\varphi)\ \bm{B}(\Theta_g, \Theta_\varphi)]$ are computed by Problem \ref{problem oblique EDMD}.
		In the training process, the initial values of the model parameters are set to the result of Problem \ref{prob. initialization of model}, i.e.,
		$g_i=g_i^0$, $\varphi_i=g_i^0$, $\bm{A}=\bm{A}_0$, and $\bm{B}=\bm{B}_0$.
	\end{prob}
\end{screen}
%\takespace

The intent of the two-staged learning procedure consisting of Problem \ref{prob. initialization of model} (Initialization by the orthogonal projection) and Problem \ref{prob. proposed method: second training} (Optimization with the oblique projection) is that the model generalizability will be increased by letting the training process have multiple chances to optimize its parameters 
on different data sets.
As seen in Example \ref{example EDMD} and Proposition \ref{prop. optimality of orthogonal projection in the finite data setting}, the model derived with the orthogonal projection, which corresponds to the same model structure as EDMD, is already \textit{optimized} being the solution to a linear regression problem.
However, it is optimal w.r.t. the summation of 2-norm errors over the given data points and its optimality is not necessarily ensured 
	when different data points are given or other characteristics are adopted as the losses. 
Thus, the orthogonal projection may not be solely the best choice of model structure, especially from the model's generalizability perspective.

For instance, collecting a sufficient amount of data in a completely unbiased way is difficult in most cases. As a result, the optimized parameters can be varied to a great extent by replacing the data points with different ones.
In such a situation, the model may not be generalizable to other regimes of operating points that are not included in the training data.

%For instance, a worst-case error can be a more important factor than averaged types of characteristics such as the sum of norms depending on the dynamics or task, in which case the orthogonal projection model may not show satisfactory performance even though it is yet optimal in the least-square error sense.
%Moreover, even if the characteristic is appropriate to evaluate the model, it can still lead to an optimality that is not consistent with the generalizability.
%Collecting a sufficient amount of data across the entire region of operating points in a completely unbiased way is not likely to be possible in most cases and the optimized parameters can be even varied to a great extent by replacing the data points with different ones in practice.
%In such a situation, the model may not be generalizable to other regimes of operating points that are not included in the training data.

Based upon these considerations, the proposed method first initializes the model by the orthogonal projection (Problem \ref{prob. initialization of model}), which ensures the least square type of optimality over the given data, and then optimizes the oblique projection (Problem \ref{prob. proposed method: second training}) to allow the training process to have a possible space to increase the generalizability by seeking an optimality on different data points.
	It is expected that the properties of the initialized model will be altered by introducing additional parameters of the test functions $\varphi_i$ in Problem \ref{prob. proposed method: second training} and the resulting oblique projection model acquires better performance for a wide range of applications.
	Similar to optimizations that include regularization terms in the loss function, the proposed modeling method aims to improve the accuracy of the model for an unseen regime of dynamics through the two-stage learning procedure with oblique projections, which may be considered a structural type of constraint imposed on the optimization.

While it is also possible to include other types of errors in the loss function to foster the generalizability further, only the 2-norm errors, $\mathcal{E}(\chi,u)$ and $\mathcal{E}_\text{state}(\chi,u)$, are considered in this paper.
Determining appropriate characteristics for the loss function is left for future research.
Note that the proposed method can also allow the orthogonal projection model by terminating the optimization with $\varphi_i=g^0_i,$ in Problem \ref{prob. proposed method: second training}.
The entire learning procedure is summarized in Algorithm \ref{algorithm proposed method}.
	Since we aim to optimize the model on different data points in Problems \ref{prob. initialization of model} and \ref{prob. proposed method: second training}, the given data set is either divided into two sub data sets or used only partially in Problem \ref{prob. initialization of model}.

\begin{algorithm}
\caption{Proposed Model Learning}
\renewcommand{\algorithmicrequire}{\textbf{Input:}}
\renewcommand{\algorithmicensure}{\textbf{Output:}}
\begin{algorithmic}[1]
\Require Data set $\mathcal{D}=\{(\chi^+, \chi, u)\mid \chi^+ = F(\chi,u)\}$
\Ensure Linear embedding model parameters: a linear operator $[\bm{A}\ \bm{B}]$ and feature maps $g_i$
\State Split $\mathcal{D}$ into sub data sets $\mathcal{D}_1$ and $\mathcal{D}_2$, or extract a part of $\mathcal{D}$ so that $\mathcal{D}_1\subset \mathcal{D}$ and define $\mathcal{D}_2:=\mathcal{D}$
\State Initialize the model on $\mathcal{D}_1$ by solving Problem \ref{prob. initialization of model}
\State Train the model on $\mathcal{D}_2$ by solving Problem \ref{prob. proposed method: second training}
\end{algorithmic}
\label{algorithm proposed method}
\end{algorithm}

\begin{figure}
	\centering 
	\begin{subfigure}{0.45\linewidth}
		\centering 
		\includegraphics[width=0.95\linewidth]{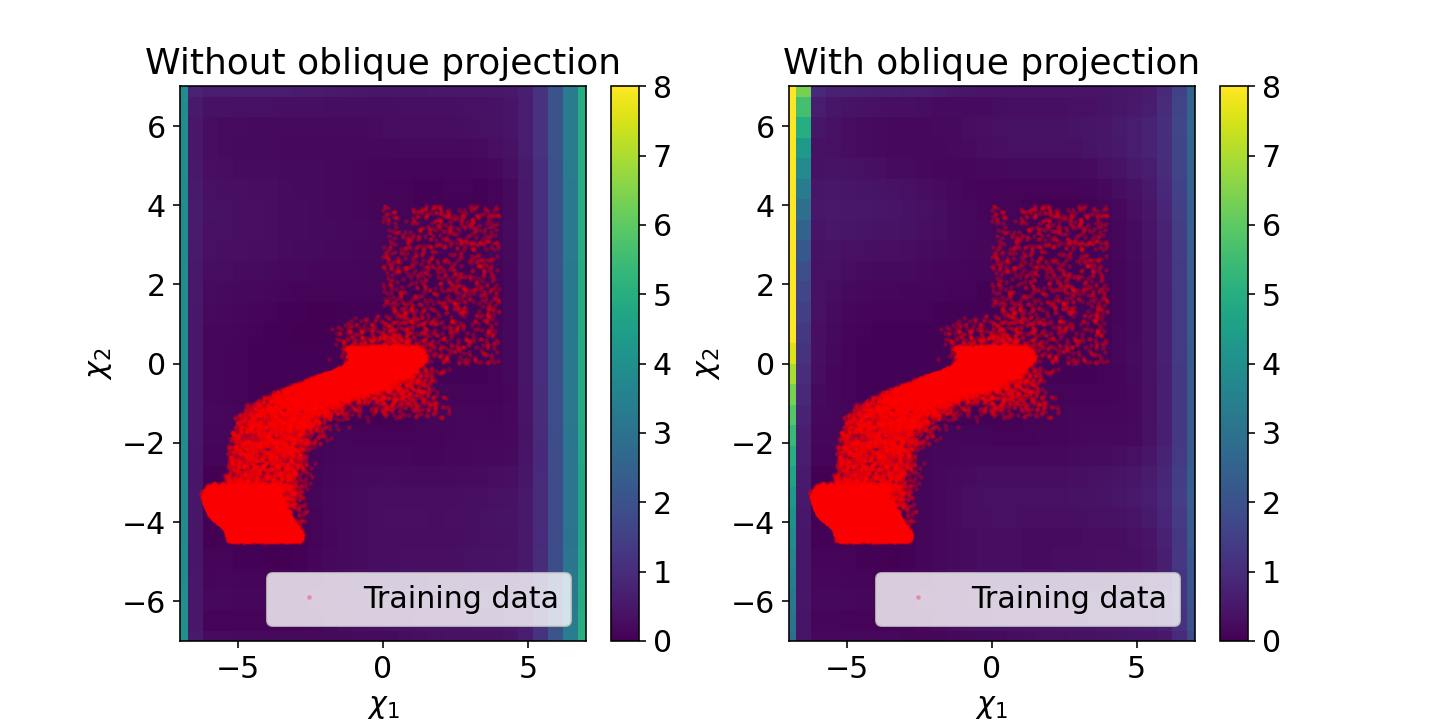}
		%\caption{Error contour ($a=b=0.3$).}
		\caption{$a=b=0.3$.}
		\label{subfig. example normal error contour}
	\end{subfigure}
	\begin{subfigure}{0.45\linewidth}
		\centering 
		\includegraphics[width=0.95\linewidth]{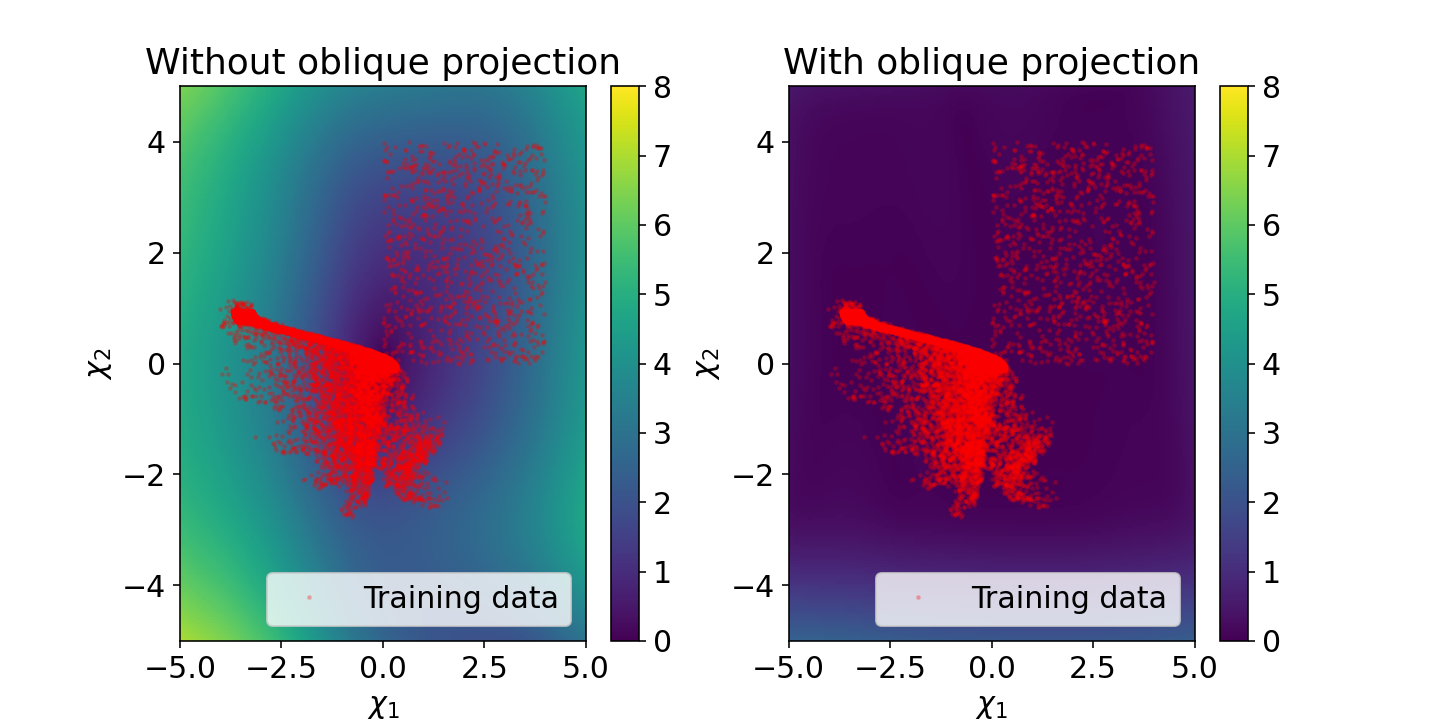}
		%\caption{Error contour ($a=0.3,b=-0.3$).}
		\caption{$a=0.3,b=-0.3$.}
		\label{subfig. example non-normal error contour}
	\end{subfigure}
%	\begin{subfigure}{0.45\linewidth}
%		\centering 
%		\includegraphics[width=0.95\linewidth]{figs/ex_normal_nonnormal/discrete_2_state_prediction_normal.png}
%		\caption{State prediction ($a=b=0.3$).}
%		\label{subfig. example normal state prediction}
%	\end{subfigure}
%	\begin{subfigure}{0.45\linewidth}
%		\centering 
%		\includegraphics[width=0.95\linewidth]{figs/ex_normal_nonnormal/discrete_2_state_prediction_non_normal.png}
%		\caption{State prediction ($a=0.3,b=-0.3$).}
%		\label{subfig. example non-normal state prediction}
%	\end{subfigure}
	\caption{Comparison of the modeling error $\mathcal{E}(\chi,u\equiv0)$ between linear embedding models with and without oblique projections. 
%		(a),(b): Error contour plots of the modeling error $\mathcal{E}(\chi,u\equiv0)$. 
%		(c),(d): Results of the state prediction task in Section \ref{sec. numerical evaluations}.
	}
	\label{fig. example normal non-normal}
\end{figure}

The use of oblique projections in Problem \ref{prob. proposed method: second training} is expected to alter the properties of the initial model obtained by Problem \ref{prob. initialization of model}, in a manner comparable to that observed in more straightforward scenarios, e.g., the stability of oblique projection-based models for linear dynamics can be investigated (therefore, can be altered in the model development process) in the context of intrusive reduced-order modeling\cite{oblique_projection_model_rom}. 
The following example showcases a situation where the use of oblique projections plays an important role when the target dynamics possesses the property of non-normality.
\begin{ex}
\rm{}
\label{ex. normal and non-normal}
Consider a two-dimensional nonlinear dynamics $\chi^+ = F(\chi,u)$ given by:
\begin{align}
\left[
	\begin{array}{c}
		\chi_{1}^+
		\\
		\chi_{2}^+
	\end{array}
	\right]&=
			\left[
			\begin{array}{cc}
				0.5 & a 
				\\
				b & 0.2
			\end{array}
			\right]
	\left[
	\begin{array}{c}
		\chi_{1}
		\\
		\chi_{2}
	\end{array}
	\right]
	+
			\left[
			\begin{array}{c}
				\cos \chi_{1} - 1
				\\	
				\cos \chi_{2} - 1
			\end{array}
			\right]
	+
	\left[
	\begin{array}{c}
		\frac{1}{2}
		\\
		0
	\end{array}
	\right]u, 
	&
	\label{eq. example normal non-normal target dynamics}
\end{align}
where $a,b\in \mathbb{R}$ are constant values and $\chi=[\chi_1\ \chi_2]^\tr\in \mathbb{R}^2$.
With no control $u\equiv 0$, the origin $\chi=[0\ 0]^\tr$ is a fixed point and the corresponding Jacobian:
\begin{align}
	\cfrac{\partial F}{\partial \chi}([0\ 0]^\tr)
	=
	\left[
	\begin{array}{cc}
		0.5 & a 
		\\
		b & 0.2
	\end{array}
	\right],
\end{align}
is normal if and only if $a=b$.
Figure \ref{fig. example normal non-normal} shows error contour plots of the modeling error $\mathcal{E}(\chi,u\equiv 0)$ defined in \eqref{eq. def of modeling error}, where two cases are compared: 1) linear embedding models are trained without oblique projections (model parameters are learned by only Problem \ref{prob. initialization of model}) and 2) models are trained with oblique projections (model parameters are learned by the proposed Algorithm \ref{algorithm proposed method}). 
Two feature maps are included so that the embedded state of the model becomes four dimensional and four test functions are used for the oblique projection.
These two types of models are tested on the target systems \eqref{eq. example normal non-normal target dynamics} with two different sets of parameters: 1) $a=b=0.3$ (Jacobian is normal) and 2) $a=0.3,b=-0.3$ (Jacobian is non-normal).
\\ \ \\
When $a=b=0.3$, in which case the Jacobian is normal, there is no significant difference irrespective of whether the oblique projections are employed (Fig. \ref{subfig. example normal error contour}). 
On the other hand, when the Jacobian is non-normal, the model trained without oblique projections results in high modeling errors apart from the training data, whereas the proposed method with oblique projections achieves a quite lower error profile (Fig. \ref{subfig. example non-normal error contour}). These results imply that the use of oblique projections is necessary in the case of non-normal dynamics of this example to obtain an accurate model. 
Note that the same training data is used in both training cases with and without oblique projections. 
\end{ex}

\section{Numerical Evaluations}
\label{sec. numerical evaluations}
In this section, numerical evaluations are provided to show the effectiveness of the proposed method.
In terms of generalizability,
a \textit{good} model should be capable of not just capturing the dynamics of the target system accurately, but also be applicable to many different control and decision-making tasks.
For instance, MPC is widely employed in the Koopman literature addressing control problems.
Long-term predictive accuracy of the model may not be an important factor in this situation since the optimization solved in the MPC procedure is only concerned with the predictions of the model up to a specified finite time horizon and the control objectives may be still achieved even if the model is not quite long-term accurate.
On the other hand, long-term predictive accuracy will be essential if the model is to be used for forecasting the system's behavior for a long time.
While it may be a practical strategy that one creates a model that is only for a specific task, developing a model that is generalizable to different problems is also of great importance and interest from the modeling perspective.
In this paper, we consider the following four tasks to compare the performance and generalizability of different linear embedding models.
\subsubsection*{State prediction}
As one of the most basic tasks of dynamical systems modeling, we consider state prediction.
Given an initial condition $\chi_0\in \mathcal{X}$ and a sequence of inputs $U=(u_0,u_1,\cdots)$, $u_i\in\mathcal{U}$, state prediction of a linear embedding model \eqref{eq. approx equation of the model with state obs included} or \eqref{eq. linear embedding model with errors from oblique projection} is implemented according to the following equation:
\begin{align}
	\chi_{0}^{\text{est}}=\chi_0,
	\ \ \ 
	\chi_{k+1}^{\text{est}}=
	[I\ 0]
	\left(
		\bm{A}
		\left[
		\begin{array}{c}
			\chi_{k}^{\text{est}}
			\\ 
			g_{n+1}(\chi_{k}^{\text{est}})
			\\
			\vdots 
			\\	
			g_{N_x}(\chi_{k}^{\text{est}})
		\end{array}
		\right]
		+
		\bm{B}
		u_k
	\right),
	\ \ 
	k=0,1,\cdots,
	\label{eq. state prediction}
\end{align}
where $\{ \chi_{0}^{\text{est}}, \chi_{1}^{\text{est}},\cdots \}$ are the predictions of states.
Note that $[I\ 0]\in \mathbb{R}^{n\times N_x}$ is the decoder $\omega$ of the model.
In this paper, only state prediction with $u_k\equiv 0$ is considered for simplicity.
Refer to \cite{control_aware_Koopman} for properties of the state prediction \eqref{eq. state prediction}, including an equivalent condition to achieve zero state prediction error.

\subsubsection*{Stabilization by LQR}
The second task is the stabilization of the unknown dynamics \eqref{eq. discrete time dynamical system} by an LQR designed for the obtained data-driven model.
The LQR gain $\bm{Q}\in \mathbb{R}^{p\times N_x}$ is computed s.t. the control input is given by $u_k=\bm{Q}\xi_k$ and the resulting closed-loop dynamics in $\mathbb{R}^{N_x}$: 
\begin{align}
	\xi_{k+1} = \bm{A}\xi_k + \bm{B}u_k
	=
	(\bm{A}+\bm{B}\bm{Q})\xi_k,\ \ 
	\xi_k\in \mathbb{R}^{N_x},
\end{align}
is stabilized at the origin while
minimizing the cost function $\sum_{k=0}^{\infty} \xi_k^\tr Q_w \xi_k + u_k^\tr R_w u_k$, where $Q_w$ and $R_w$ are the weight matrices.
Note that the matrices $\bm{A}$ and $\bm{B}$ in the above equation are given by the data-driven linear embedding model. For instance, they are the solution $[
\bm{A}\ \bm{B}]$ to Problem \ref{prob. proposed method: second training} if the model is obtained by the proposed method.
Python Control Systems Library\cite{python-control2021} is used to compute the LQR gain.

\subsubsection*{Reference tracking}
We consider reference tracking as the third task.
Specifically, a state-feedback controller for a non-zero, steady reference $r$ is designed for a linear embedding model, which is formulated as the following problem\cite{discrete_control_systems_ogata}:
\begin{align}
	\text{Find\ }\bm{Q}_s, \bm{Q}_I\text{\ \ s.t.\ \ }
	\left\{
	\begin{array}{l}
		u_k=-\bm{Q}_s
		\xi_k
		+\bm{Q}_I\nu_k
		\\
		\nu_k=\nu_{k-1} + r - y_k
		\\
		y_k = \bm{C}
		\xi_k
		\\
		y_k\rightarrow r\ (k\rightarrow \infty)
	\end{array}
	\right.,
	\label{eq. controller ref tracking}
\end{align}
where the state $\xi_k\in \mathbb{R}^{N_x}$ is subject to an LTI system
$\xi_{k+1} =  \bm{A}\xi_k + \bm{B}u_k$ ($\bm{A}$ and $\bm{B}$ are given by the data-driven linear embedding model), $\nu_k$ denotes a slack variable that corresponds to an integrator of the controller, and $\bm{C}$ is an arbitrary matrix with appropriate dimensions that specifies $y_k$ which is expected to follow the reference $r$, respectively.
For more details of the controller design, see \cite{discrete_control_systems_ogata}.

\subsubsection*{Linear MPC}
As the fourth problem, we use linear embedding models for linear MPC, in which
we solve the following optimization problem:
\begin{align}
&
\underset{u_0,\cdots, u_{N_h}}{\text{min}}
\sum_{k=0}^{N_h + 1} \left\{
l(\xi_k,u_k,k)
+
(u_k-u_{k-1})^\tr R (u_k-u_{k-1})
\right\}
&\nonumber
\\
&\hspace{15mm}
\text{\ \ subject\ to:\ }
\left\{
\begin{array}{l}
	\xi_{k+1}=\bm{A}\xi_k + \bm{B}u_k
	\\
	\text{other constraints on }\xi_k,u_k
\end{array}
\right.,
&
\label{eq. MPC formulation}
\end{align}
and implement only the first element $u_0$ at each time.
The parameter $N_h$ specifies the finite horizon of the optimization problem and $\bm{A}$ and $\bm{B}$ in the constraints are given by a data-driven linear embedding model.
We set the objective function $l(\xi_k,u_k,k)$ s.t. the specified state will follow a time-varying reference signal.
Model predictive control python toolbox\cite{do-mpc} is used to produce the results of this paper.

In addition to the generalizability of the model and its performance for individual tasks, the sensitivity of the learning procedure is also of relevance to the numerical evaluation of the method. The sensitivity analysis of the modeling methods considered in this paper is provided in Section 1 of the Supplementary Material.
Also, a discussion on the modeling error $\mathcal{E}(\chi,u)$, or the invariance proximity of the model, is provided in Section 2 of the Supplementary Material.

\subsection{Duffing Oscillator}
As the first example, we consider the Duffing oscillator, which is given by the following ODE:
\begin{align}
	\ddot{z}(t)
	=
	-0.5\dot{z}(t) + z(t) - 4 z^3(t) + u(t),
	\label{eq. duffing oscillator}
\end{align}
where the state $z(t)$ and the input $u(t)$ are continuous variables.
With no control $u(t)\equiv 0$, this system has an unstable fixed point $z=0$ and two stable fixed points $z=\pm 1/2$.
We can conceptually consider a first-order time discretization of \eqref{eq. duffing oscillator}, which is assumed to be the unknown dynamics \eqref{eq. discrete time dynamical system} in the formulations of this paper.
The variables $z(t)$ and $\dot{z}(t)$ correspond to the first and second components of $\chi\in \mathbb{R}^2$, respectively.
%The variables $x(t)$ and $\dot{x}(t)$ are related to the data set $\{ (\chi_i,u_i,y_i)\mid y_i=F(\chi_i,u_i) \}$ as $\chi_i=[]$

In the modeling phase, 600 trajectories of states were generated by \eqref{eq. duffing oscillator}, each of which consists of a single trajectory of 50 discrete-time states sampled at every 0.05 time units, starting from an initial condition $[z(0)\ \dot{z}(0)]^\tr$ drawn from the uniform distribution over $[-3,3]^2\subset \mathbb{R}^2$.
The Runge-Kutta method was used to solve the ODE \eqref{eq. duffing oscillator} with the step size of 0.01.
Following a data generating strategy suggested in \cite{control_aware_Koopman}, the input data was generated according to $u_k=\cos(\omega_i k \Delta t)$, $\omega_i:=20i$, $i=0,1,\cdots,5$, where $\Delta t=0.05$ denotes the sampling period of data.

In the model learning phase, we trained a proposed model with two feature maps, which yields the embedded state of the form $[\chi^\tr\ g_3(\chi)\ g_4(\chi)]^\tr$, and four test functions\\ $[\varphi_1(\chi,u)\cdots \varphi_4(\chi,u)]^\tr$.
Both $[g_3(\chi)\ g_4(\chi)]^\tr$ and $[\varphi_1(\chi,u)\cdots \varphi_4(\chi,u)]^\tr$ are characterized by a fully connected feed-forward neural network with a single hidden layer consisting of 10 neurons, respectively.
The swish function was used as the activation.
Problem \ref{prob. initialization of model} was solved with half of the data set to initialize the model, which was then followed by solving Problem \ref{prob. proposed method: second training} with the rest of data points added. 
Model training was implemented in TensorFlow.

\begin{figure}[]
	\centering
	\includegraphics[width=0.9\linewidth]{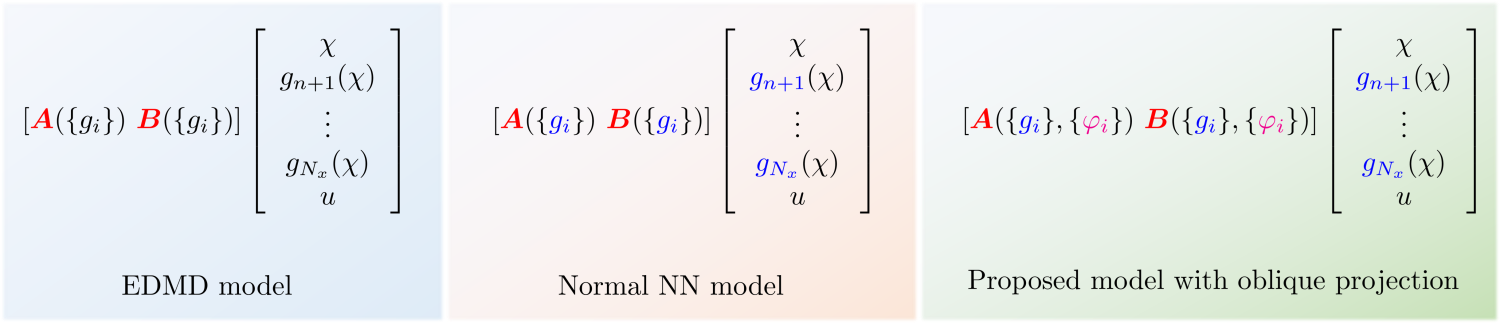}
	\caption{Comparison of different linear embedding models. Each equation represents the output of the model given the state $\chi$ and the input $u$. Colored variables indicate dependency on data: derivation of $\bm{A}$ and $\bm{B}$ depends on data in the EDMD model; $\bm{A}$, $\bm{B}$, and $\{g_i\}$ in the normal NN model; $\bm{A}$, $\bm{B}$, $\{g_i\}$, and $\{\varphi_i\}$ in the proposed model. The notation $\bm{A}(\{g_i\},\{\varphi_i\})$ represents the equation defining the matrix $\bm{A}$ depends on designs of the feature maps $\{g_i\}$ and the test functions $\{\varphi_i\}$.}
	\label{fig. comparison model structures}
\end{figure}

For comparison with the proposed method, we also trained two other Koopman-based data-driven models: an EDMD model and a normal neural network model.
The EDMD model was obtained by \eqref{eq. intro [A B] as in EDMD} with monomial features up to the third order, which yields the following nine-dimensional embedded state of the model \eqref{eq. approx equation of the model with state obs included}:
\begin{align}
	[g_1(\chi)\cdots g_9(\chi)]^\tr=
	[\chi_{(1)}\ \ \chi_{(2)}\ \ \chi_{(1)}^2\ \ \chi_{(2)}^2\ \ \chi_{(1)}\chi_{(2)}\ \ 
	\chi_{(1)}^3\ \ \chi_{(2)}^3\ \ \chi_{(1)}^2\chi_{(2)}\ \ \chi_{(1)}\chi_{(2)}^2  
	]^\tr,
	\label{eq. EDMD model embedded state in numerical examples}
\end{align} 
where we used the notation $[\chi_{(1)}\ \chi_{(2)}]^\tr:=\chi\in \mathbb{R}^2$.
As the normal neural network model, we trained a model by solving Problem \ref{prob. initialization of model} only, which is the non-autonomous version of LKIS-DMD\cite{Learning_Koopman_Invariant_Subspaces} as mentioned in Remark \ref{remark Initialization of model is same as LKIS-DMD}.
It has the same embedding architecture as the proposed model, i.e., two feature maps parameterized by a neural network with the same structure.
It is referred to as the normal NN model throughout the rest of the paper.
Note that both EDMD and normal NN models were also trained on exactly the same data points as the proposed method.
All three models considered here are linear embedding models, but they differ from each other w.r.t. dependency on data or model structures as summarized in Fig. \ref{fig. comparison model structures}.

In the LQR design, we set the weight matrices as
\begin{align}
    Q_w=
    \left[
        \begin{array}{cc}
             Q_\text{state} & 0  \\
             0 & 0 
        \end{array}
    \right],
    \ \ 
    Q_\text{state}:=
    \left[
        \begin{array}{cc}
             100 & 0  \\
             0 & 1
        \end{array}
    \right],
    \ \ 
    R_w = 1.
    \label{eq. LQR weights duffing}
\end{align}

In the reference tracking problem, we defined $\bm{C}=[1\ 0]$ and $r=1$.
In the MPC design, we set the objective function in \eqref{eq. MPC formulation} as follows:
\vspace{-2mm}
\begin{align}
    l(\xi_k,u_k,k)=
    \left(
        \xi_{k,(1)} - r_\text{MPC}(k)
    \right)^2,
    \ \ 
    R = 1,
    \label{eq. MPC objectives duffing}
\end{align}
where $\xi_{k,(1)}$ denotes the first component of $\xi_k$ and $r_\text{MPC}(k)$ is a time-varying reference signal s.t.
\vspace{-2mm}
\begin{align}
    r_\text{MPC}(k):=
    \left\{
        \begin{array}{r}
             -1\ (k\leq k_\text{c})
        \\
            1\ (k> k_\text{c})
        \end{array}
    \right..
\end{align}

We set $k_\text{c}$ s.t. it corresponds to $t=10$ in the simulation environment.
No additional constraints on $\xi_k,u_k$ were imposed.

\begin{figure}[]
	\centering 
	\begin{subfigure}{0.3\linewidth}
		\centering
		\includegraphics[width=0.95\linewidth]{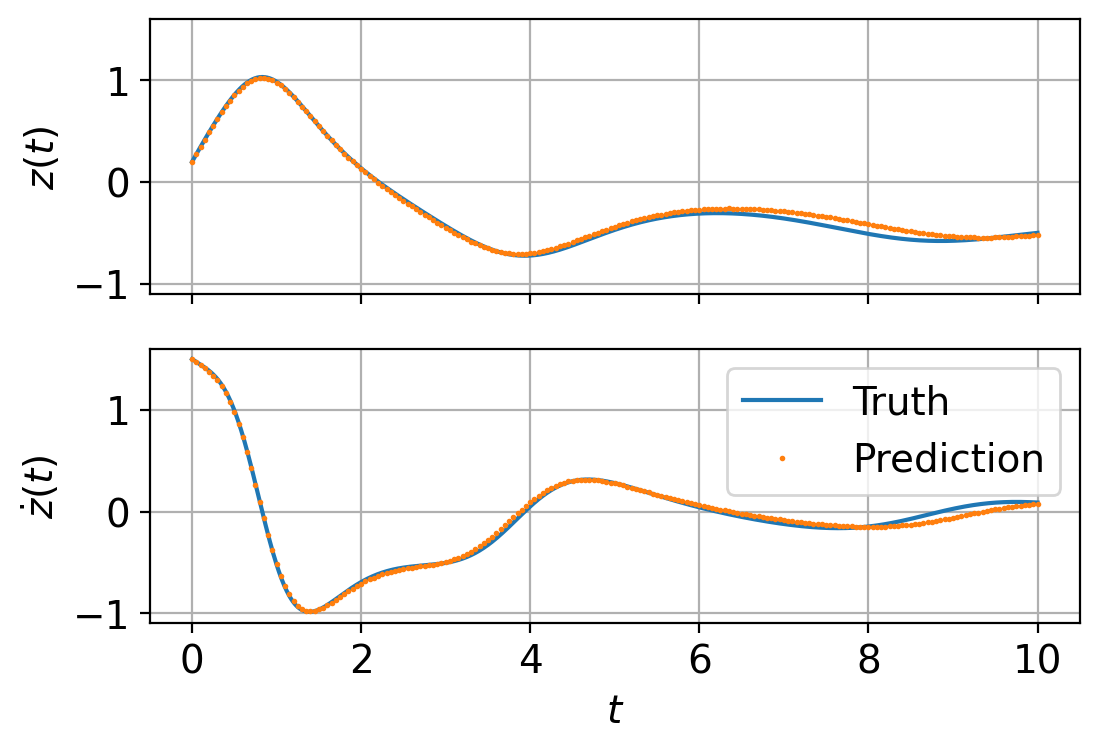}
		\caption{State prediction (EDMD).}
		\label{subfig. duffing state prediction EDMD}
	\end{subfigure}
	\begin{subfigure}{0.3\linewidth}
		\centering
		\includegraphics[width=0.95\linewidth]{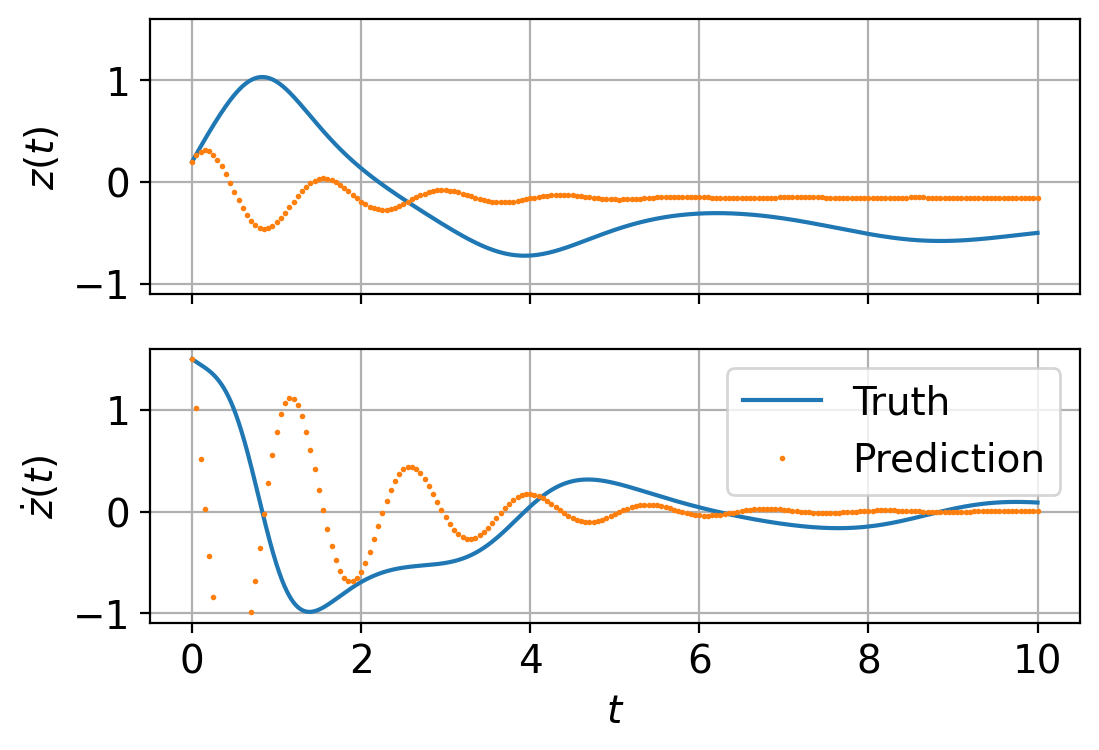}
		\caption{\scriptsize{State prediction (normal NN).}}
		\label{subfig. duffing state prediction normal}
	\end{subfigure}
	\begin{subfigure}{0.3\linewidth}
		\centering
		\includegraphics[width=0.95\linewidth]{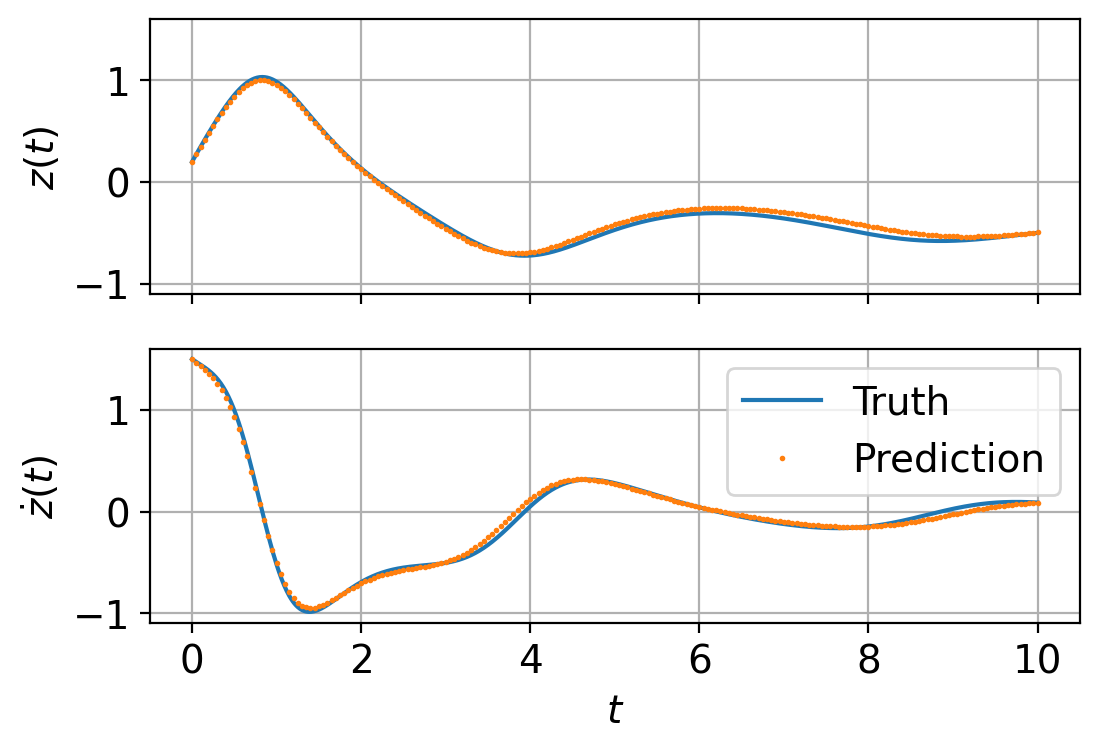}
		\caption{\footnotesize{State prediction (proposed).}}
		\label{subfig. duffing state prediction oblique}
	\end{subfigure}
	\begin{subfigure}{0.3\linewidth}
		\centering
		\includegraphics[width=0.95\linewidth]{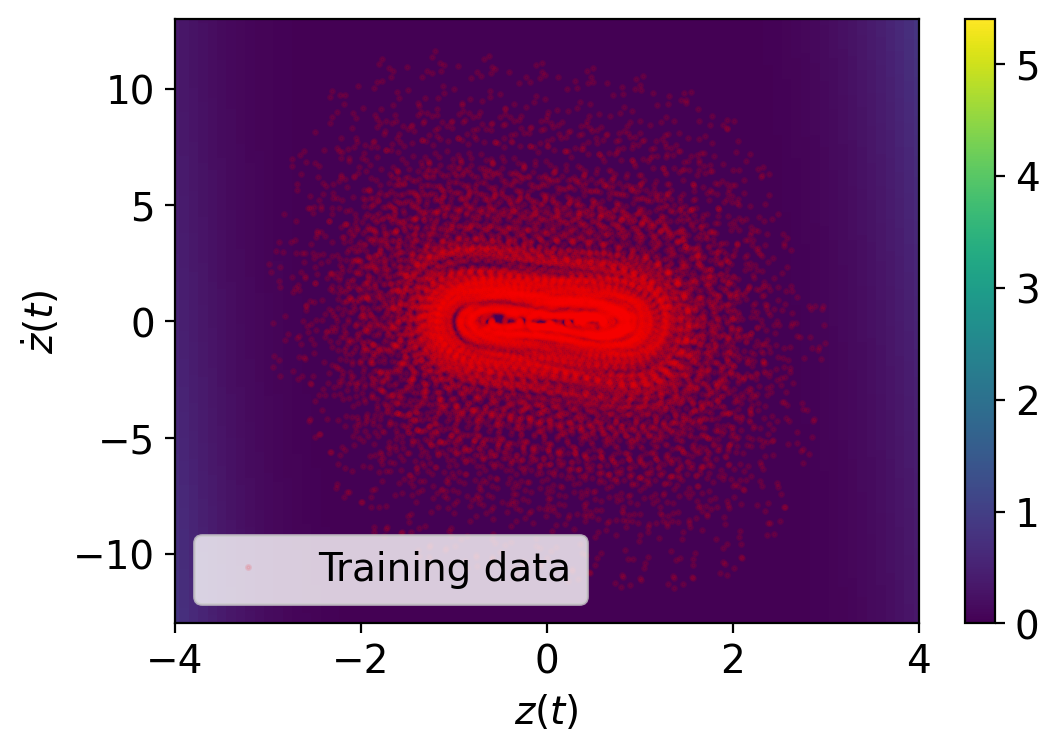}
		\caption{Error contour (EDMD).}
		\label{subfig. duffing error contour EDMD}
	\end{subfigure}
	\begin{subfigure}{0.3\linewidth}
		\centering
		\includegraphics[width=0.95\linewidth]{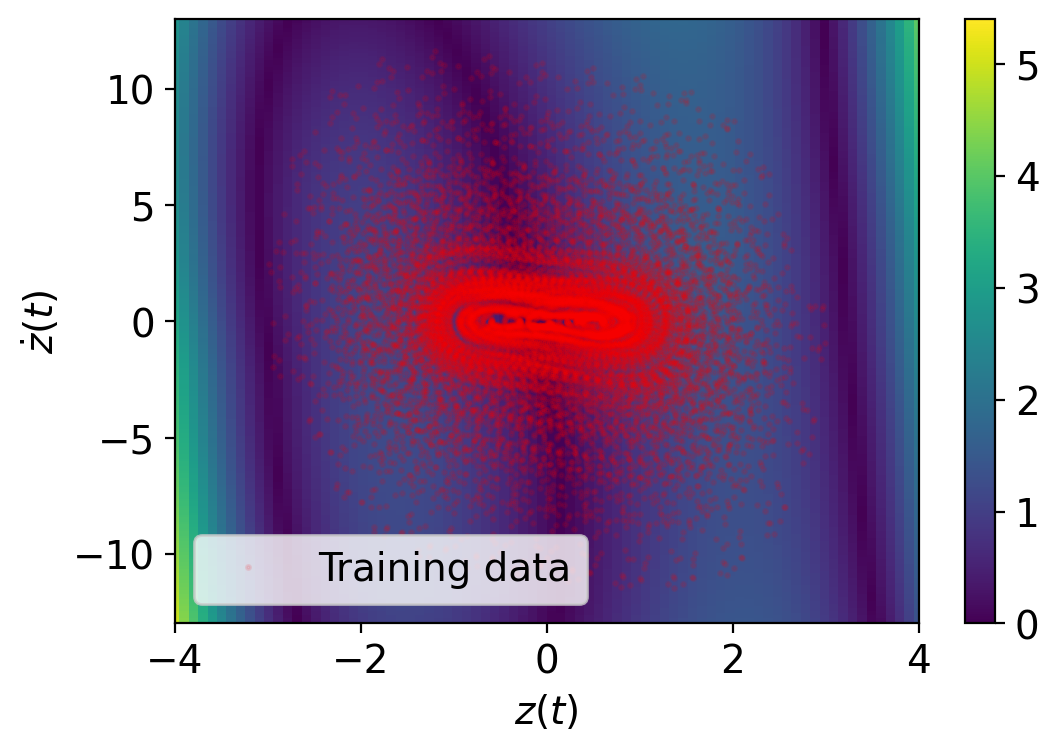}
		\caption{\footnotesize{Error contour (normal NN).}}
		\label{subfig. duffing error contour normal}
	\end{subfigure}
	\begin{subfigure}{0.3\linewidth}
		\centering
		\includegraphics[width=0.95\linewidth]{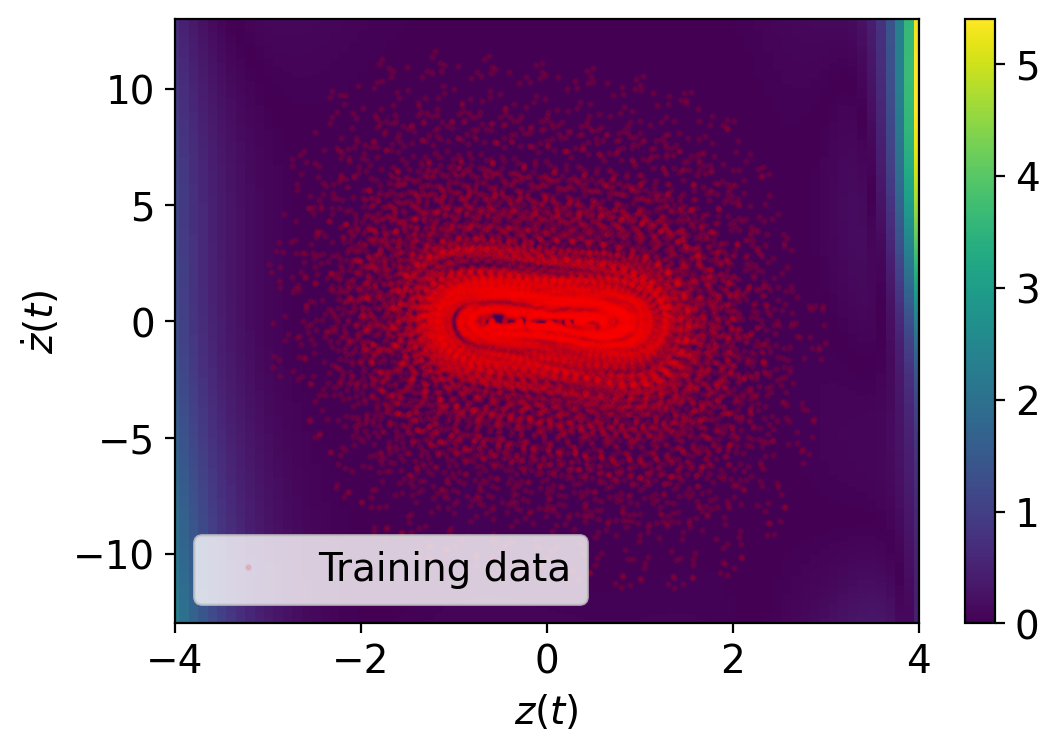}
		\caption{Error contour (proposed).}
		\label{subfig. duffing error contour oblique}
	\end{subfigure}
	\caption{Results of the Duffing oscillator (state prediction).}
	\label{fig. duffing state prediction}
\end{figure}

\begin{figure}[H]
	\centering 
	\begin{subfigure}{0.3\linewidth}
		\centering
		\includegraphics[width=0.95\linewidth]{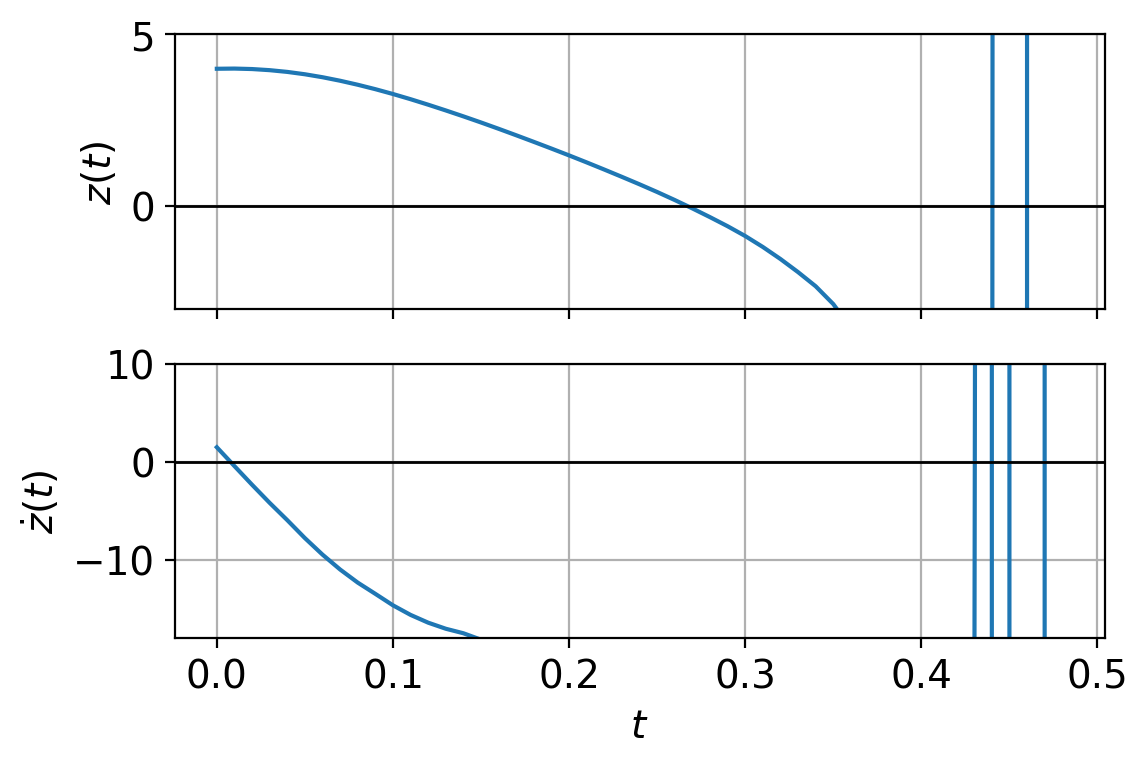}
		\caption{Stabilization by LQR (EDMD).}
		\label{subfig. duffing cl EDMD}
	\end{subfigure}
	\begin{subfigure}{0.3\linewidth}
		\centering
		\includegraphics[width=0.95\linewidth]{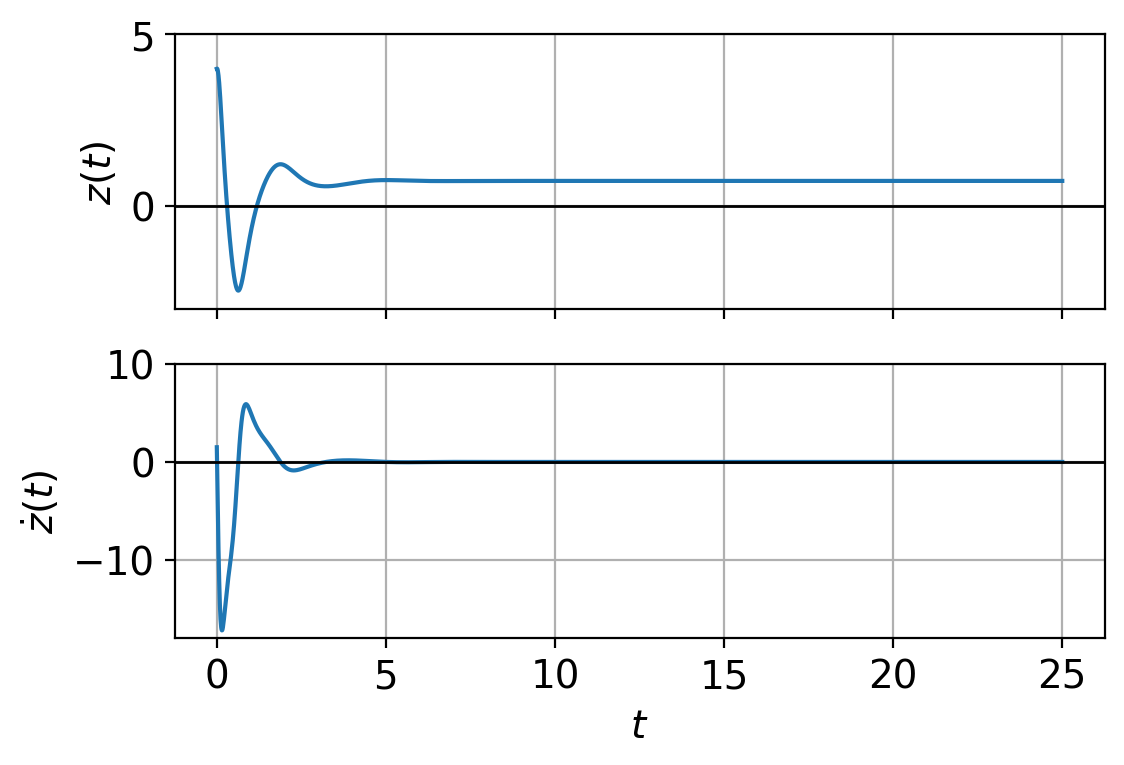}
		\caption{\footnotesize{Stabilization by LQR (normal NN).}}
		\label{subfig. duffing cl normal}
	\end{subfigure}
	\begin{subfigure}{0.3\linewidth}
		\centering
		\includegraphics[width=0.95\linewidth]{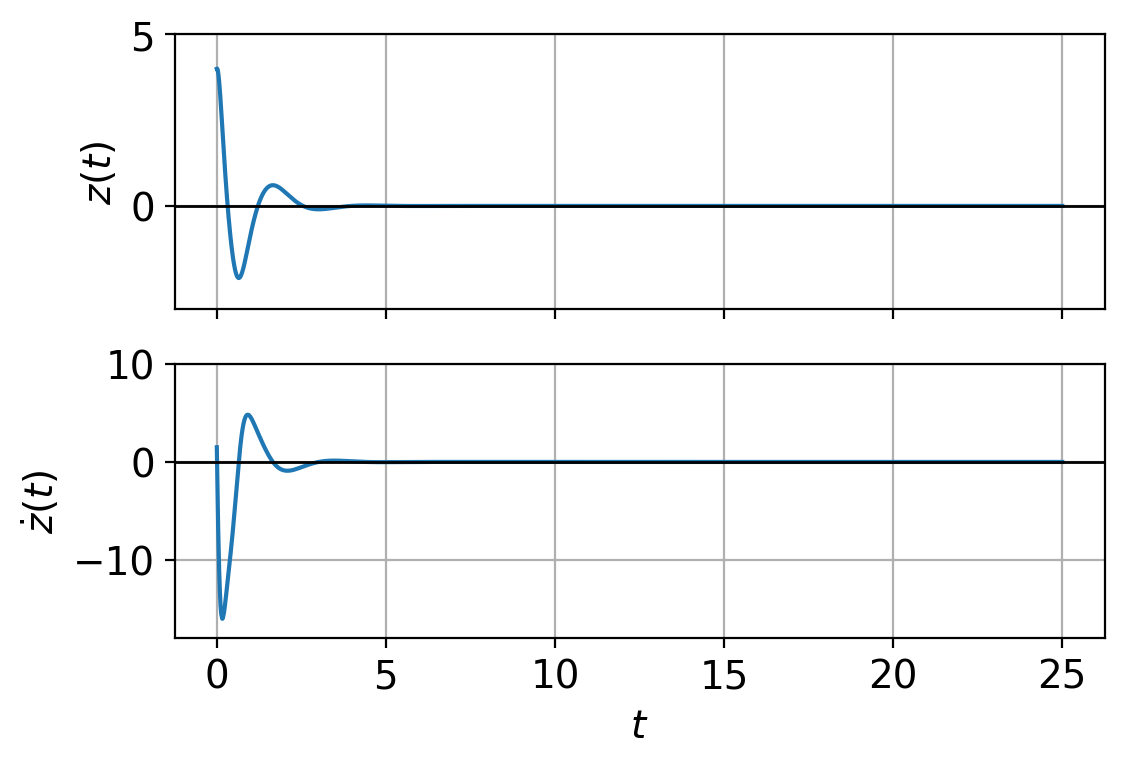}
		\caption{Stabilization by LQR (proposed).}
		\label{subfig. duffing cl oblique}
	\end{subfigure}
	\begin{subfigure}{0.3\linewidth}
		\centering
		\includegraphics[width=0.95\linewidth]{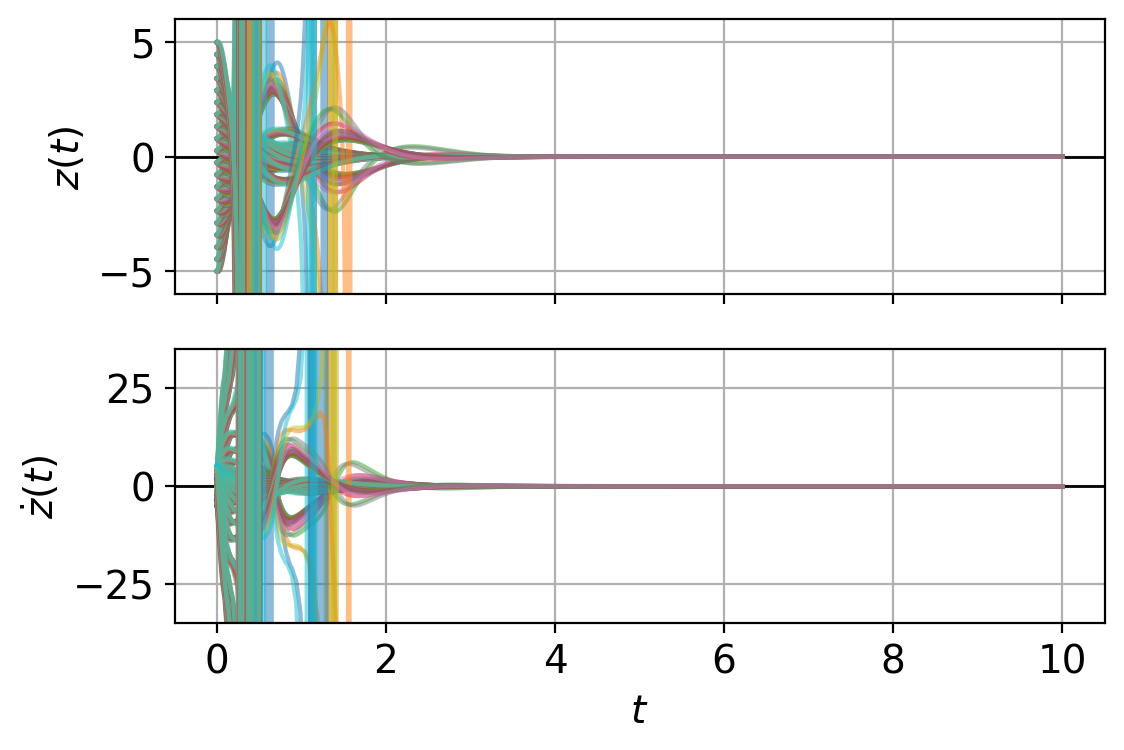}
		\caption{Estimation of basin of attraction (EDMD).}
		\label{subfig. duffing basin of attraction EDMD}
	\end{subfigure}
	\begin{subfigure}{0.3\linewidth}
		\centering
		\includegraphics[width=0.95\linewidth]{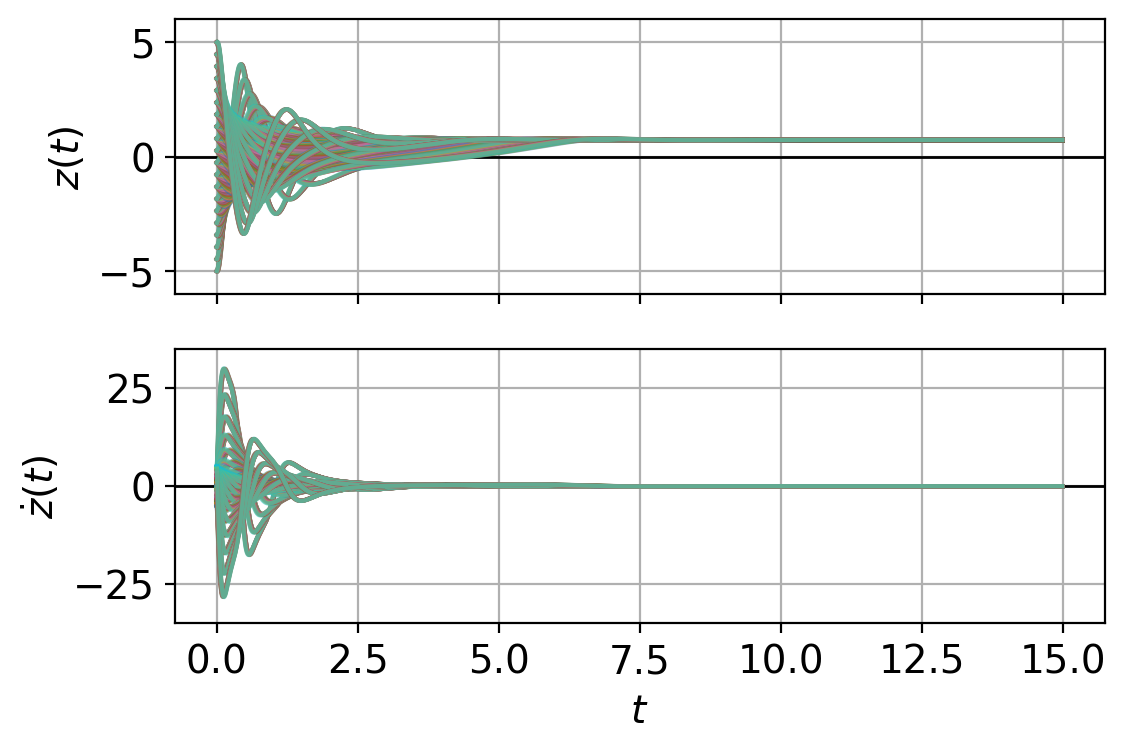}
		\caption{Estimation of basin of attraction (normal NN).}
		\label{subfig. duffing basin of attraction normal}
	\end{subfigure}
	\begin{subfigure}{0.3\linewidth}
		\centering
		\includegraphics[width=0.95\linewidth]{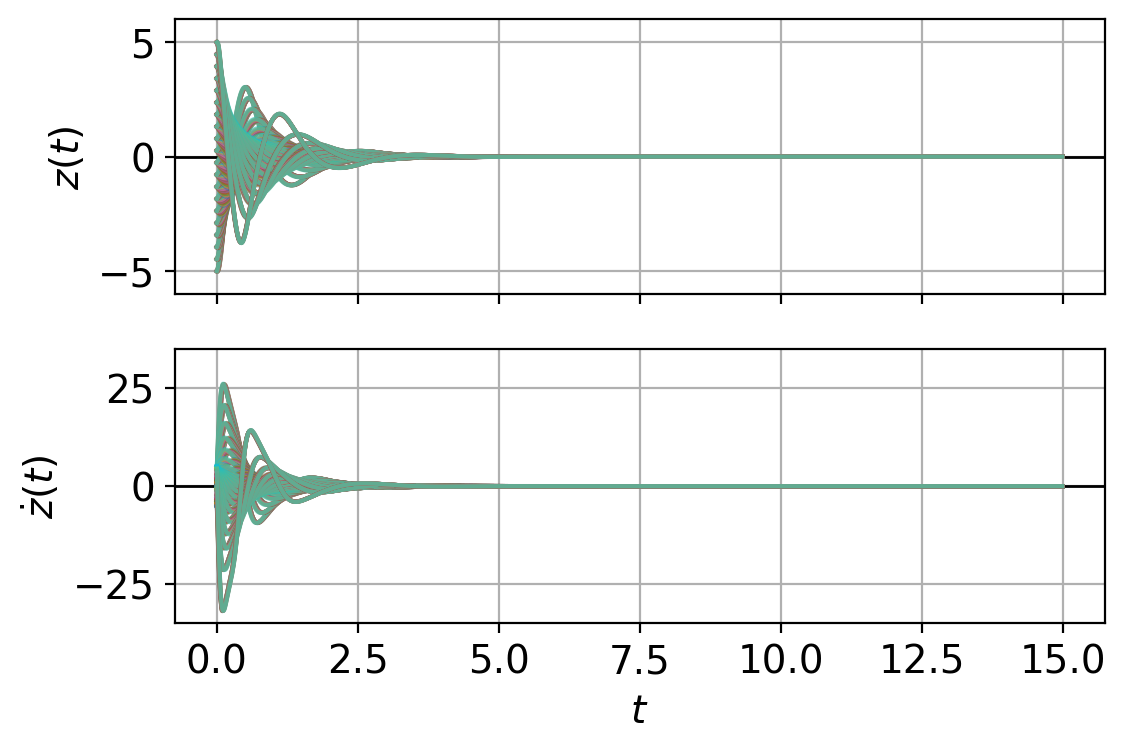}
		\caption{Estimation of basin of attraction (proposed).}
		\label{subfig. duffing basin of attraction oblique}
	\end{subfigure}
	\begin{subfigure}{0.3\linewidth}
		\centering
		\includegraphics[width=0.95\linewidth]{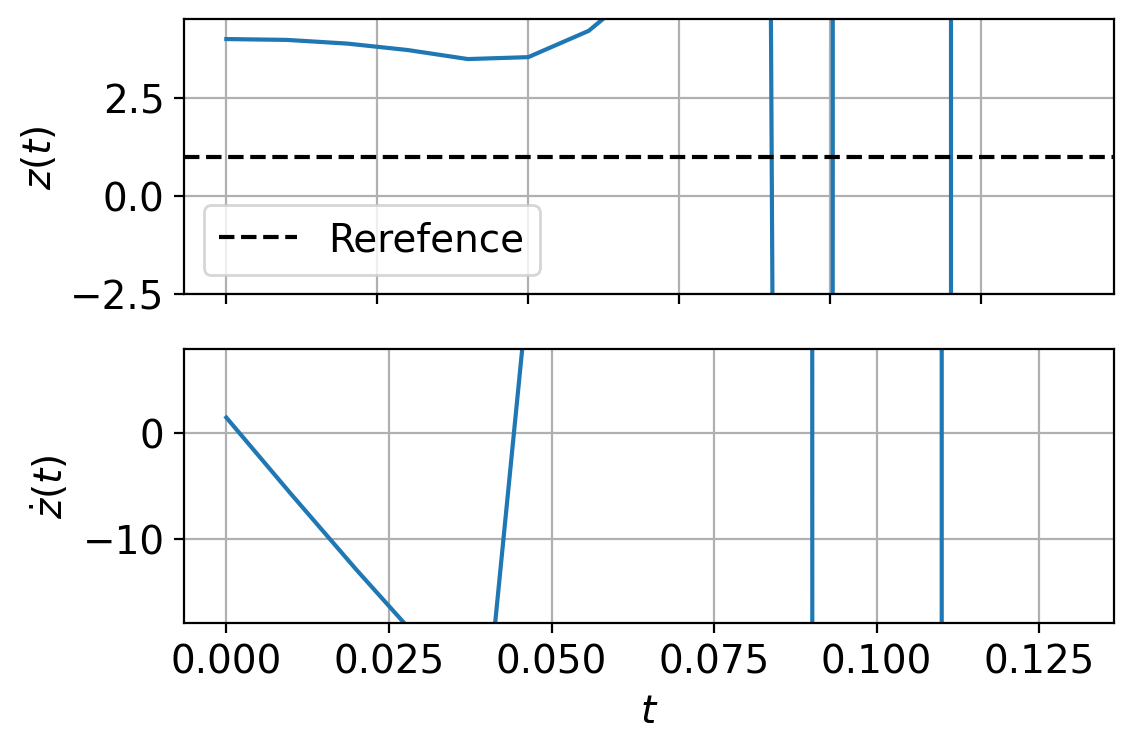}
		\caption{Reference tracking (EDMD).}
		\label{subfig. duffing cl servo EDMD}
	\end{subfigure}
	\begin{subfigure}{0.3\linewidth}
		\centering
		\includegraphics[width=0.95\linewidth]{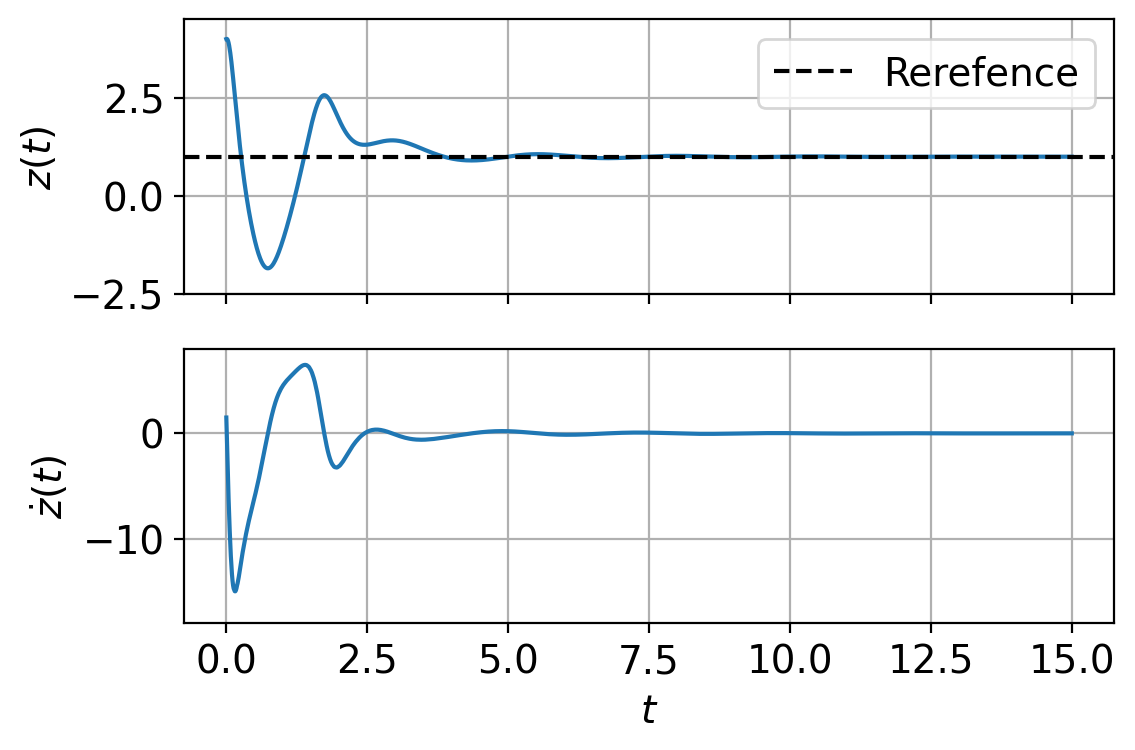}
		\caption{Reference tracking (normal NN).}
		\label{subfig. duffing cl servo normal}
	\end{subfigure}
	\begin{subfigure}{0.3\linewidth}
		\centering
		\includegraphics[width=0.95\linewidth]{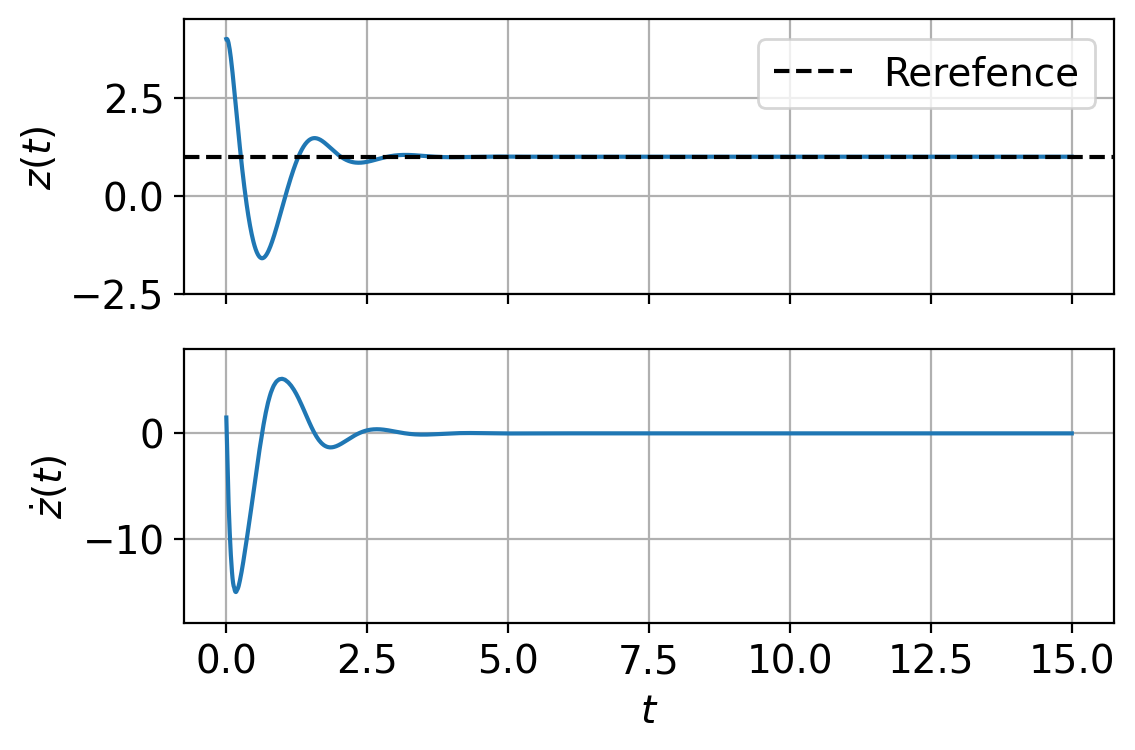}
		\caption{Reference tracking (proposed).}
		\label{subfig. duffing cl servo oblique}
	\end{subfigure}
	\begin{subfigure}{0.3\linewidth}
		\centering
		\includegraphics[width=0.95\linewidth]{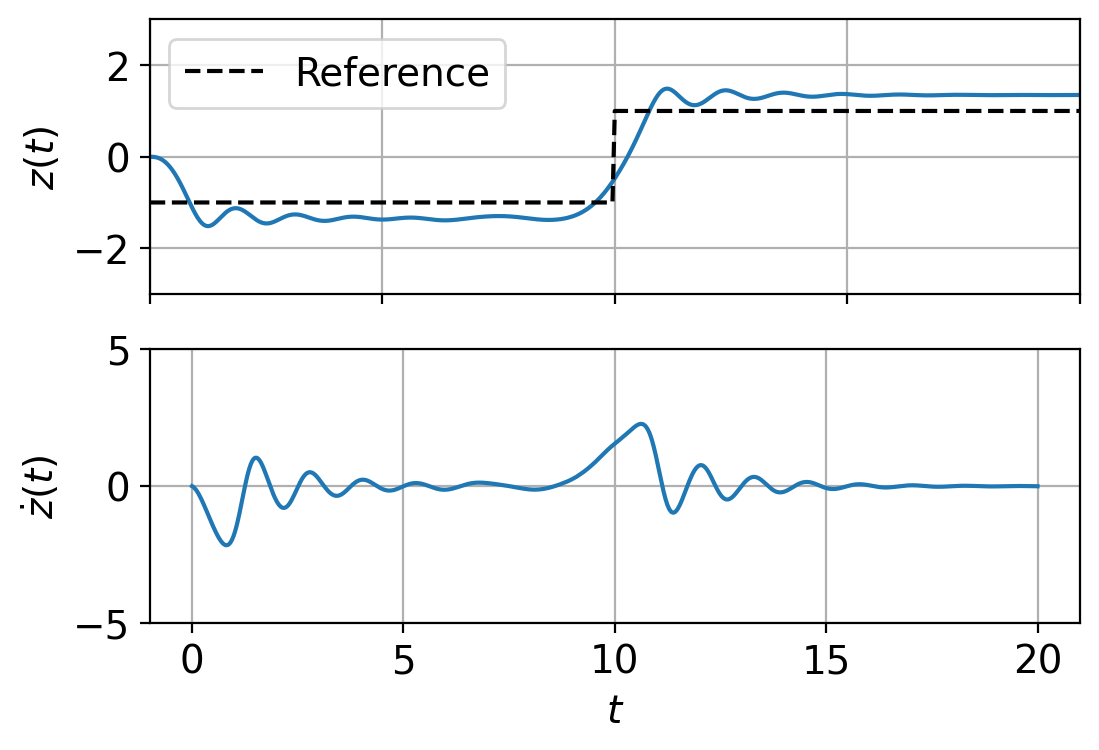}
		\caption{MPC (EDMD).}
		\label{subfig. duffing mpc EDMD}
	\end{subfigure}
	\begin{subfigure}{0.3\linewidth}
		\centering
		\includegraphics[width=0.95\linewidth]{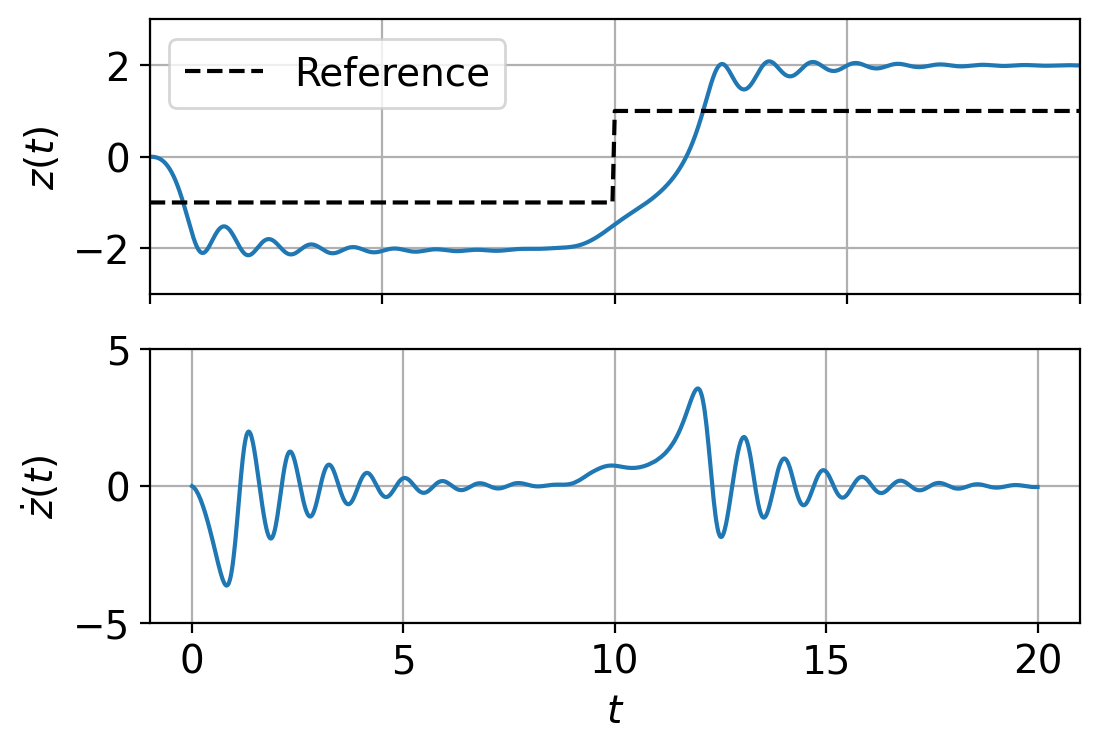}
		\caption{MPC (normal NN).}
		\label{subfig. duffing mpc normal}
	\end{subfigure}
	\begin{subfigure}{0.3\linewidth}
		\centering
		\includegraphics[width=0.95\linewidth]{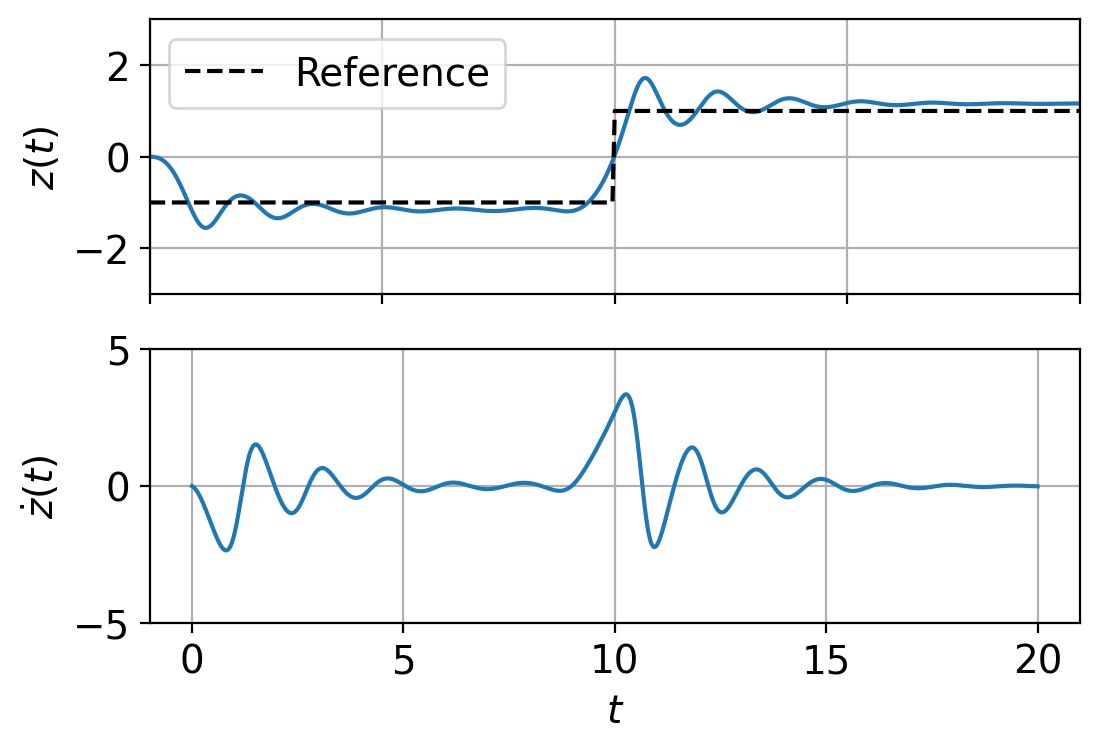}
		\caption{MPC (proposed).}
		\label{subfig. duffing mpc oblique}
	\end{subfigure}
	\caption{Results of the Duffing oscillator (control applications).}
	\label{fig. duffing control applications}
\end{figure}

The results of the four different tasks are shown in Figs. \ref{fig. duffing state prediction} and \ref{fig. duffing control applications}.
First, the EDMD model captures the behavior of the target dynamics quite accurately (Fig. \ref{subfig. duffing state prediction EDMD}).
This is because the Duffing oscillator only has a cubic nonlinearity as in \eqref{eq. duffing oscillator} and the embedded state \eqref{eq. EDMD model embedded state in numerical examples} also includes exactly the same nonlinear feature map $\chi_{(1)}^3$.
In this case, it is possible to achieve zero state prediction error\cite{control_aware_Koopman}.
A similar result is seen in the error contour plot (Fig. \ref{subfig. duffing error contour EDMD}), which shows one-step state prediction errors. 
On the other hand, the normal NN model fails to obtain accurate state prediction (Figs. \ref{subfig. duffing state prediction normal} and \ref{subfig. duffing error contour normal}).
This is considered as a result of neural network training terminated at an unsatisfactory local minimum for the given data, which showcases the difficulty of training a neural network-based model involving a high-dimensional non-convex optimization as mentioned in Section \ref{sec. on no convergence property of EDMD}.
Although the normal NN model could be more accurate, one needs to repeat the training process with different initialization of parameters until the learned model shows satisfactory results, which may be prohibitively time-consuming depending on the problem.
On the other hand, the proposed model shows quite accurate state predictive accuracy that is comparable to the EDMD model (Figs. \ref{subfig. duffing state prediction oblique} and \ref{subfig. duffing error contour oblique}).
In Section 3 of the Supplementary Material,
a comparison with other general nonlinear classes of models is provided.

The superiority of the proposed method is also observed when the models are deployed in the controller design tasks.
In the first application of stabilization by LQR, the EDMD and the normal NN models result in either having a steady state error (normal NN model in Fig. \ref{subfig. duffing cl normal}) or altering the closed-loop system unstable (EDMD model in Fig. \ref{subfig. duffing cl EDMD}).
On the other hand, the dynamics is successfully stabilized by a controller designed for the proposed model (Fig. \ref{subfig. duffing cl oblique}).
As more quantitative evaluations of stabilization by LQR, Figs. \ref{subfig. duffing basin of attraction EDMD}, \ref{subfig. duffing basin of attraction normal}, and \ref{subfig. duffing basin of attraction oblique} show estimation of basin of attraction, where 
the responses of the closed-loop dynamics starting from various initial conditions are overlaid in each single plot.
%various initial conditions are tested to see if they are stabilized by the controllers.
The proposed model stabilizes all the given initial conditions, whereas steady-state errors are present with the normal NN model and many simulations diverge with the EDMD model.

For the reference tracking task, the proposed and the normal NN models achieve the control objective (Figs. \ref{subfig. duffing cl servo normal} and \ref{subfig. duffing cl servo oblique}) and the EDMD model results in a divergent closed-loop dynamics (Fig. \ref{subfig. duffing cl servo EDMD}).
It is noted that the EDMD model struggles with both stabilization by LQR and the reference tracking although it shows quite accurate state prediction. It is owed to how the modeling error regarding the state prediction differently propagates from that regarding the closed-loop dynamics\cite{control_aware_Koopman}.

Since MPC is only concerned with relatively short-time predictions of the models controlled by the finite horizon of its optimization, none of the three models leads to divergent closed-loop dynamics.
Although all the cases do not perfectly track the reference signal, the proposed method has the least amount of steady-state error (Figs. \ref{subfig. duffing mpc EDMD}, \ref{subfig. duffing mpc normal}, and \ref{subfig. duffing mpc oblique}).

\subsection{Simple Pendulum}
As the second example, we consider the simple pendulum system:
\begin{align}
	\ddot{z}(t) = -\sin z(t) + u(t).
\end{align}

\begin{figure}[]
	\centering 
	\begin{subfigure}{0.3\linewidth}
		\centering
		\includegraphics[width=0.95\linewidth]{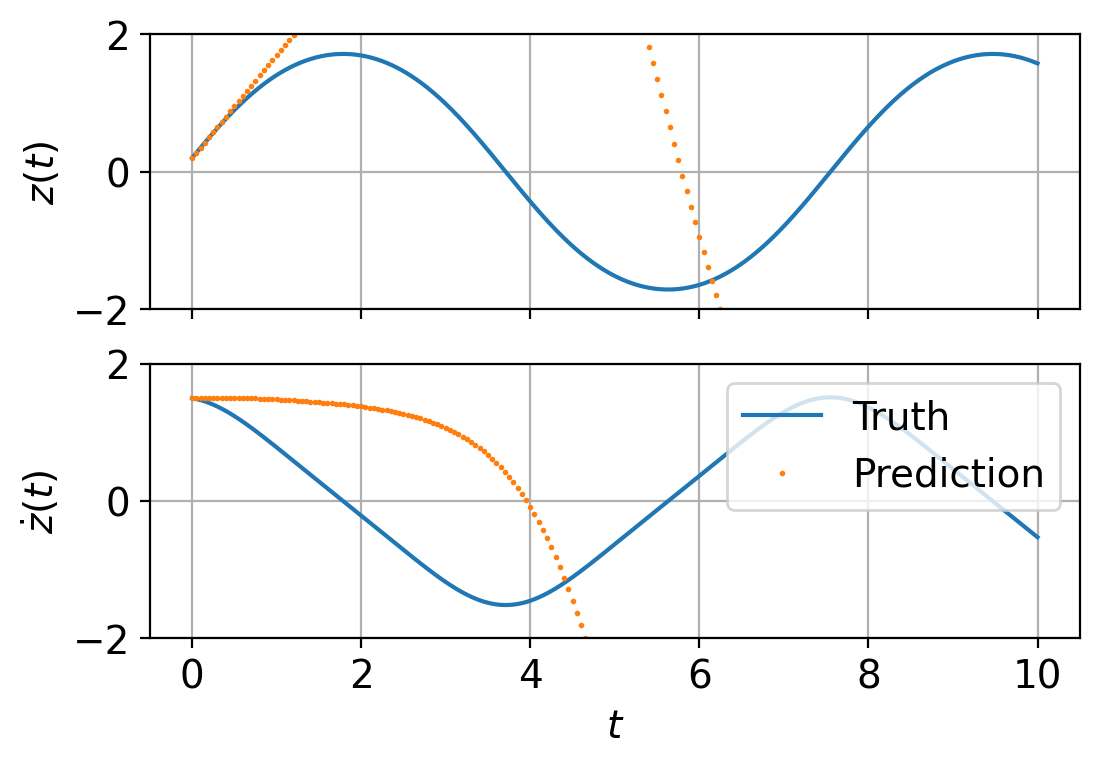}
		\caption{State prediction (EDMD).}
		\label{subfig. pendulum state prediction EDMD}
	\end{subfigure}
	\begin{subfigure}{0.3\linewidth}
		\centering
		\includegraphics[width=0.95\linewidth]{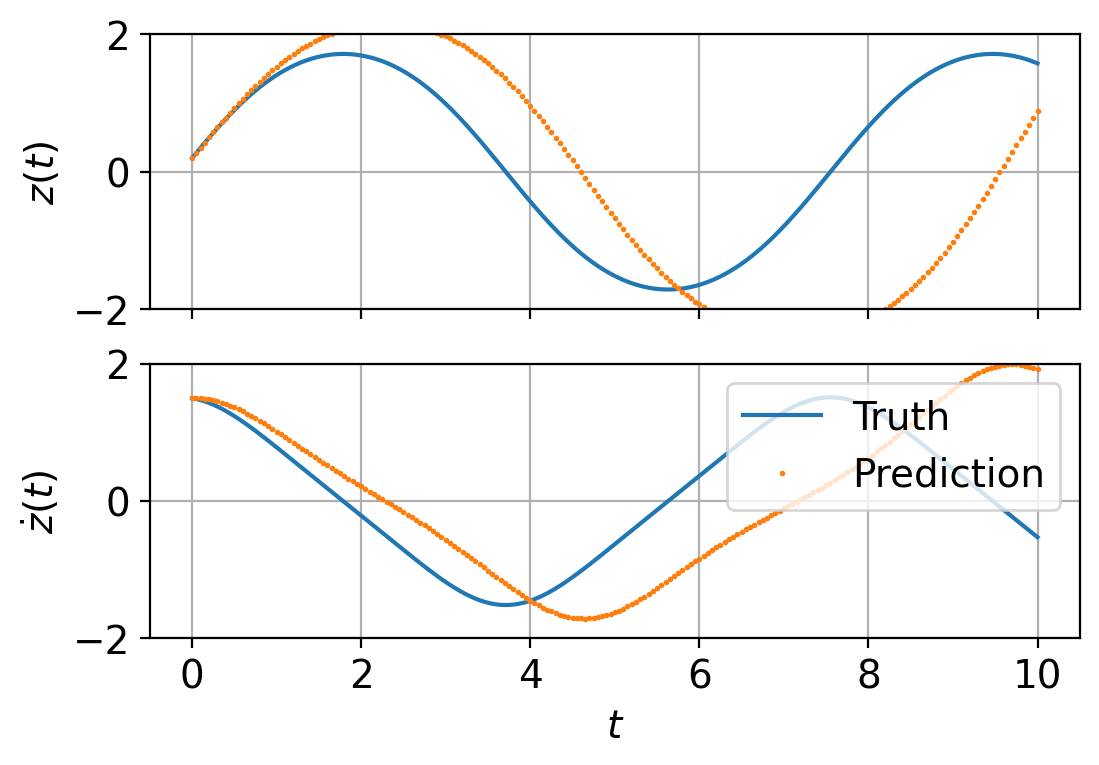}
		\caption{\scriptsize{State prediction (normal NN).}}
		\label{subfig. pendulum state prediction normal}
	\end{subfigure}
	\begin{subfigure}{0.3\linewidth}
		\centering
		\includegraphics[width=0.95\linewidth]{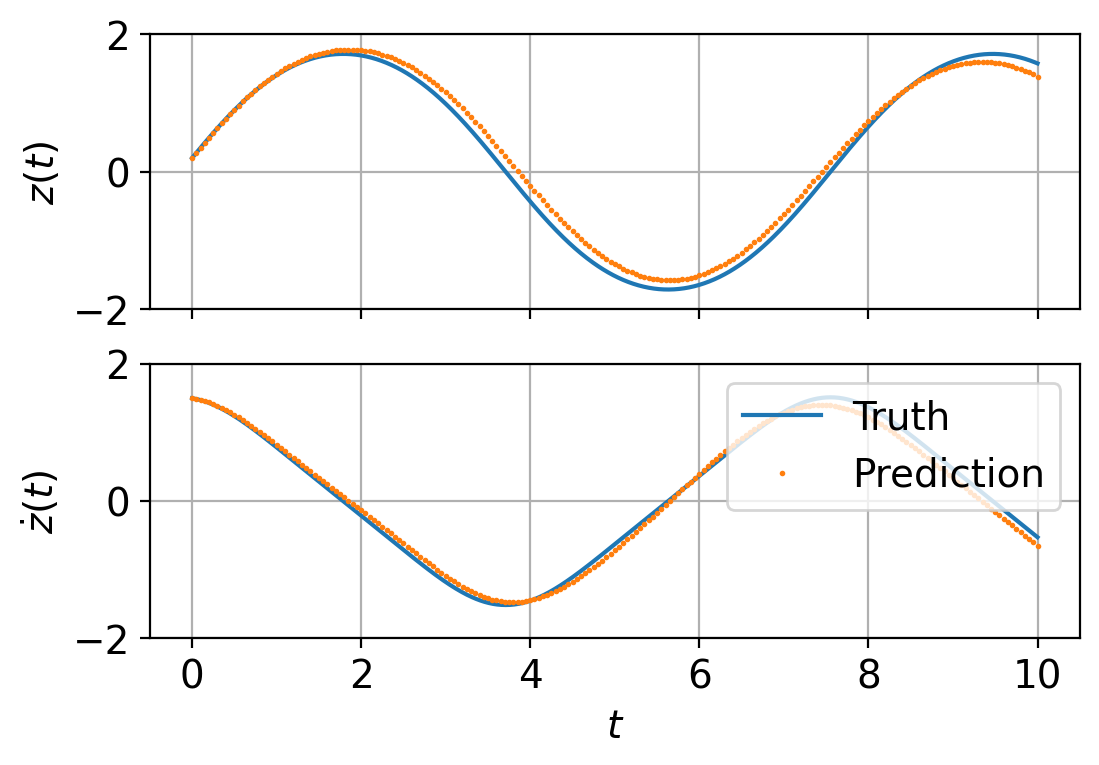}
		\caption{\footnotesize{State prediction (proposed).}}
		\label{subfig. pendulum state prediction oblique}
	\end{subfigure}
	\begin{subfigure}{0.3\linewidth}
		\centering
		\includegraphics[width=0.95\linewidth]{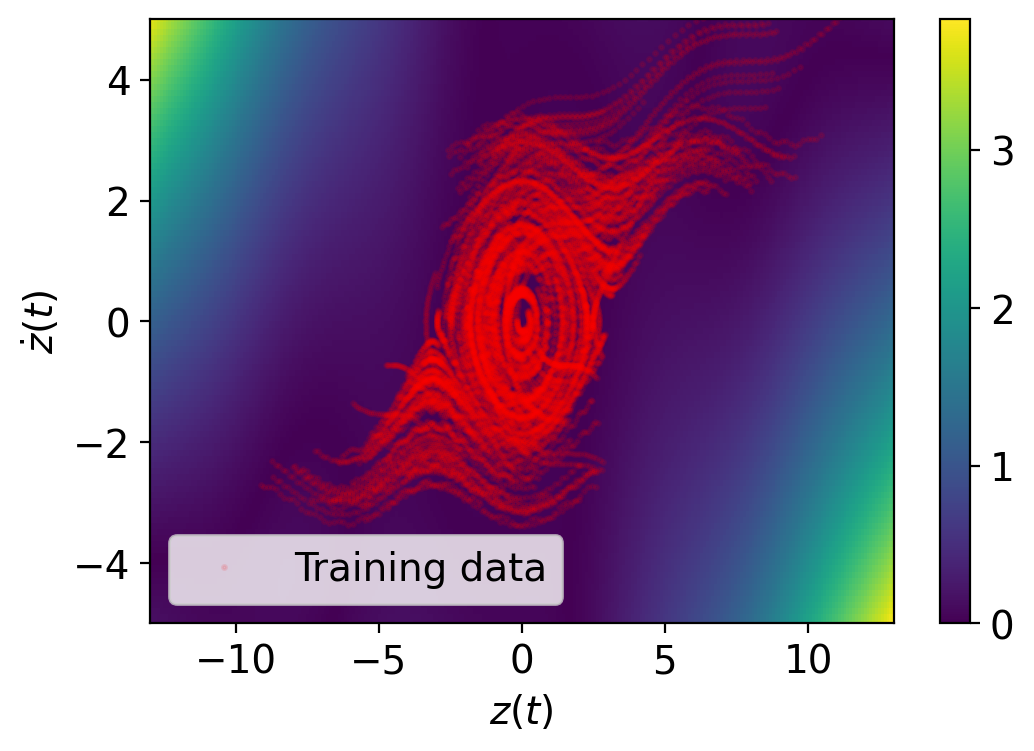}
		\caption{Error contour (EDMD).}
		\label{subfig. pendulum error contour EDMD}
	\end{subfigure}
	\begin{subfigure}{0.3\linewidth}
		\centering
		\includegraphics[width=0.95\linewidth]{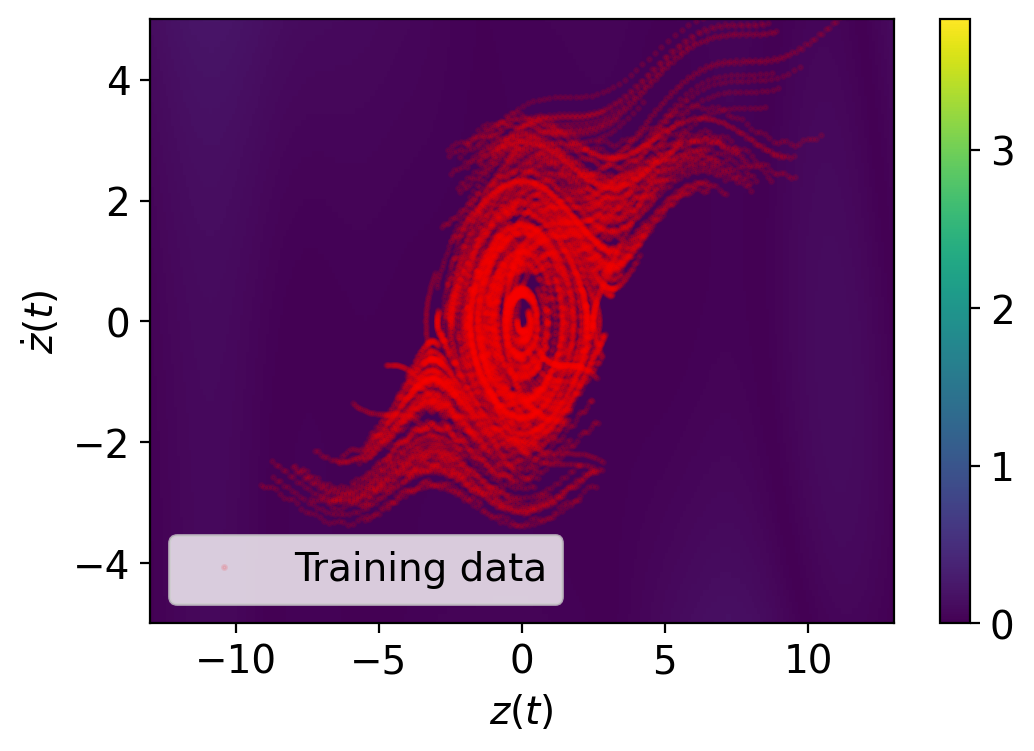}
		\caption{\footnotesize{Error contour (normal NN).}}
		\label{subfig. pendulum error contour normal}
	\end{subfigure}
	\begin{subfigure}{0.3\linewidth}
		\centering
		\includegraphics[width=0.95\linewidth]{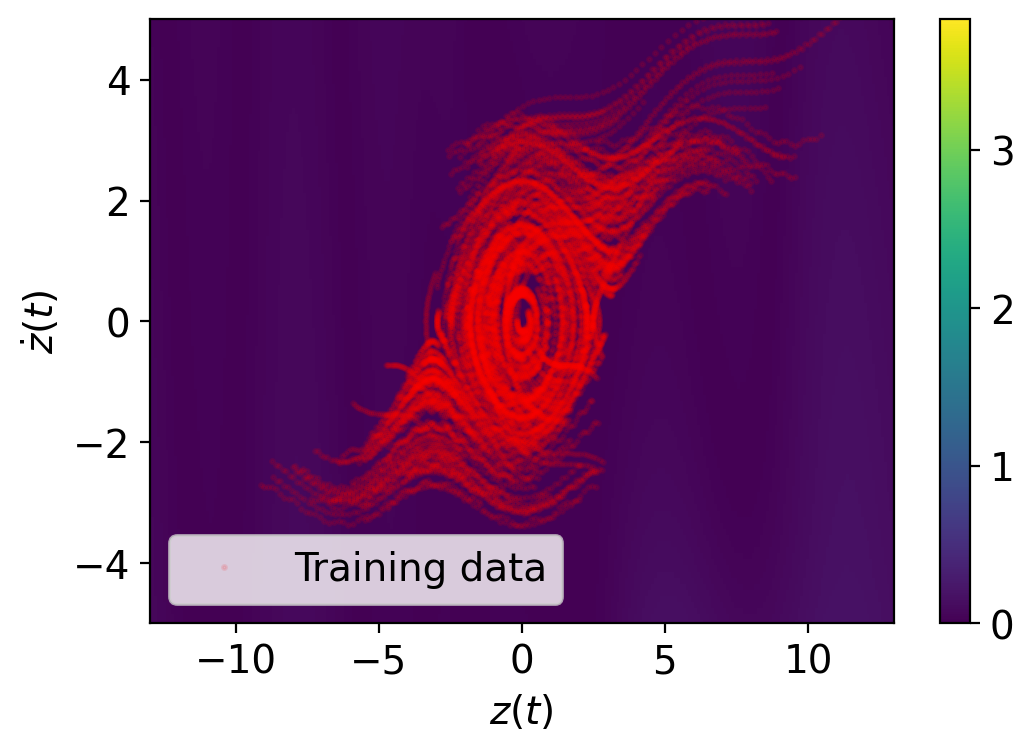}
		\caption{Error contour (proposed).}
		\label{subfig. pendulum error contour oblique}
	\end{subfigure}
	\caption{Results of the simple pendulum (state prediction).}
	\label{fig. pendulum state prediction}
\end{figure}

\begin{figure}[H]
	\centering 
	\begin{subfigure}{0.3\linewidth}
		\centering
		\includegraphics[width=0.95\linewidth]{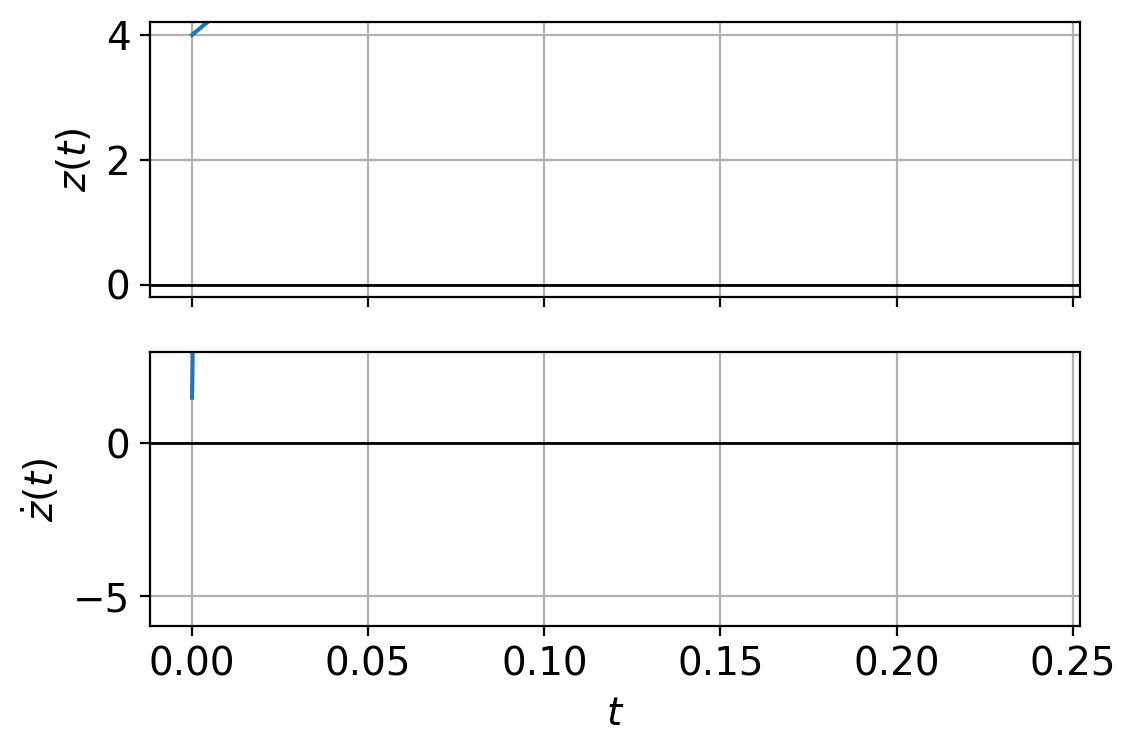}
		\caption{Stabilization by LQR (EDMD).}
		\label{subfig. pendulum cl EDMD}
	\end{subfigure}
	\begin{subfigure}{0.3\linewidth}
		\centering
		\includegraphics[width=0.95\linewidth]{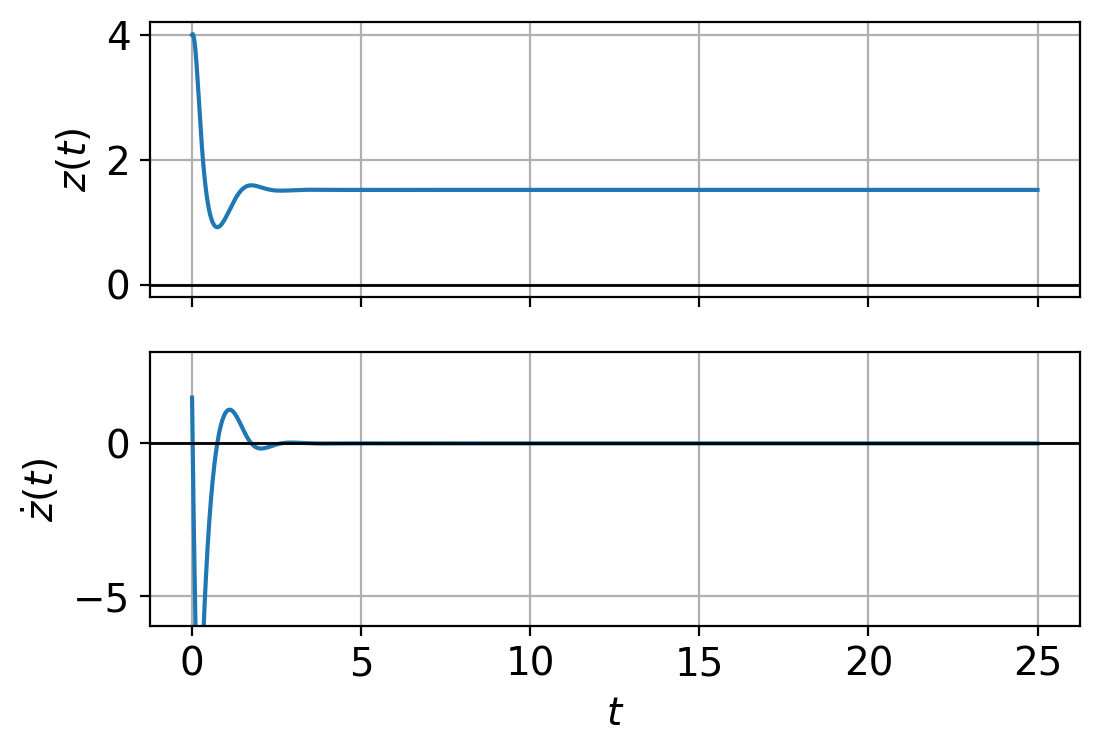}
		\caption{\footnotesize{Stabilization by LQR (normal NN).}}
		\label{subfig. pendulum cl normal}
	\end{subfigure}
	\begin{subfigure}{0.3\linewidth}
		\centering
		\includegraphics[width=0.95\linewidth]{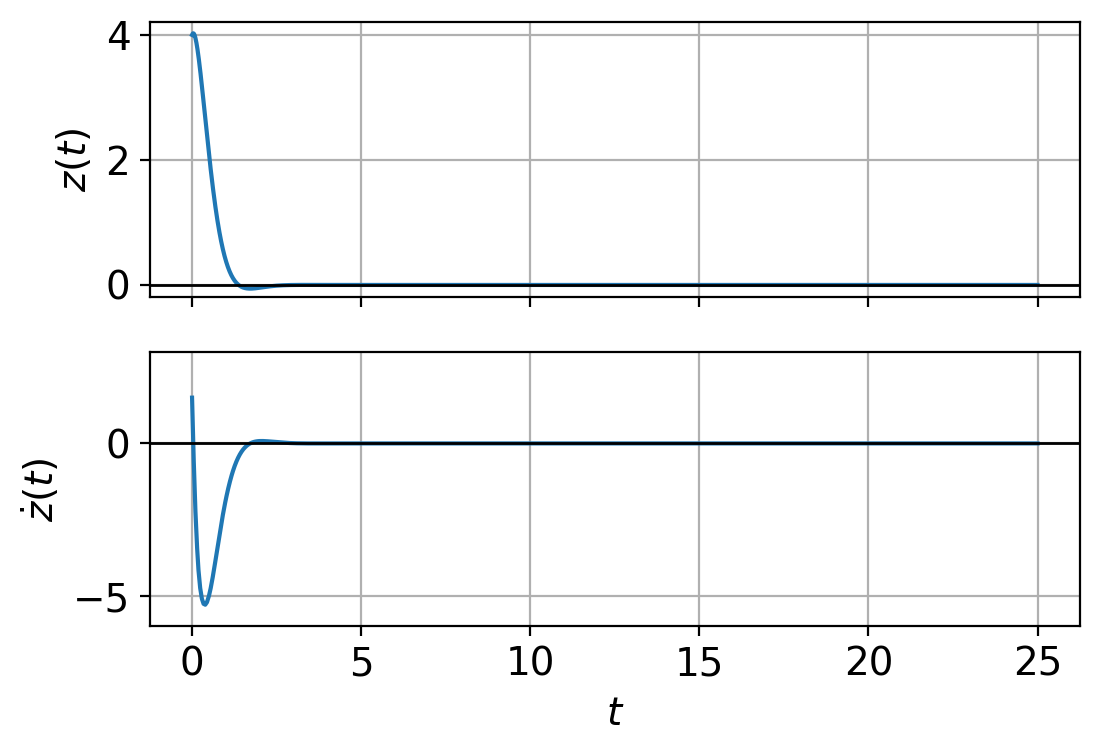}
		\caption{Stabilization by LQR (proposed).}
		\label{subfig. pendulum cl oblique}
	\end{subfigure}
	\begin{subfigure}{0.3\linewidth}
		\centering
		\includegraphics[width=0.95\linewidth]{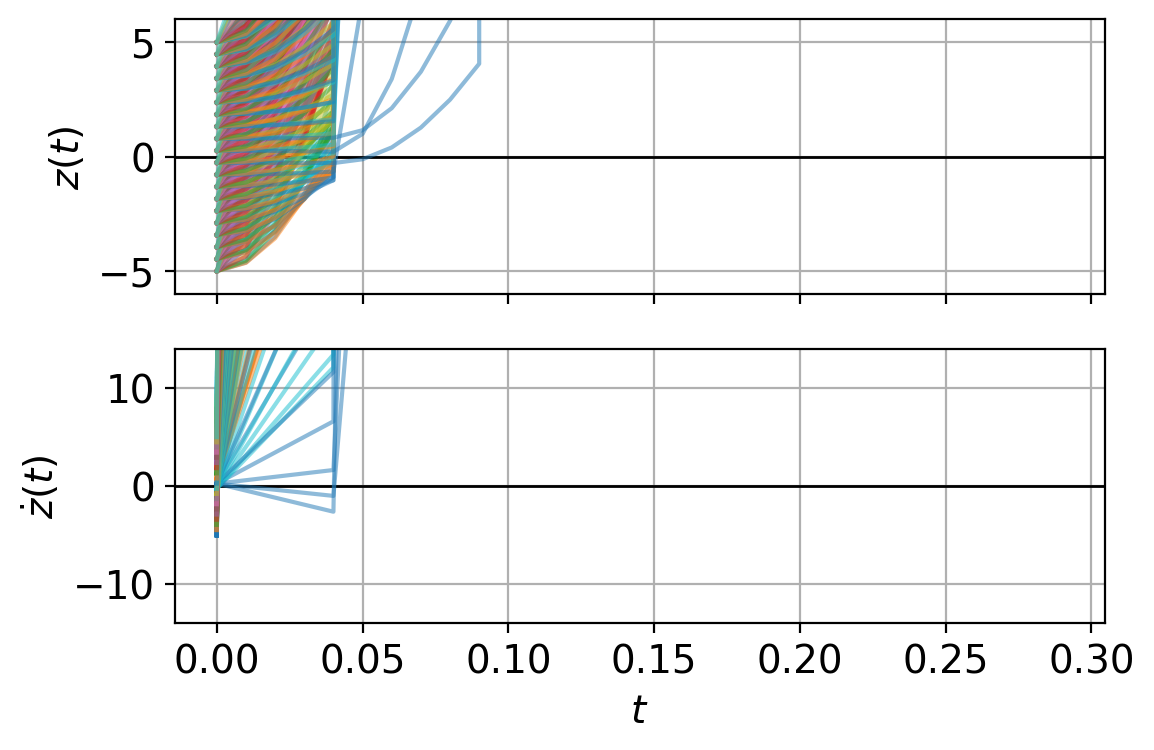}
		\caption{Estimation of basin of attraction (EDMD).}
		\label{subfig. pendulum basin of attraction EDMD}
	\end{subfigure}
	\begin{subfigure}{0.3\linewidth}
		\centering
		\includegraphics[width=0.95\linewidth]{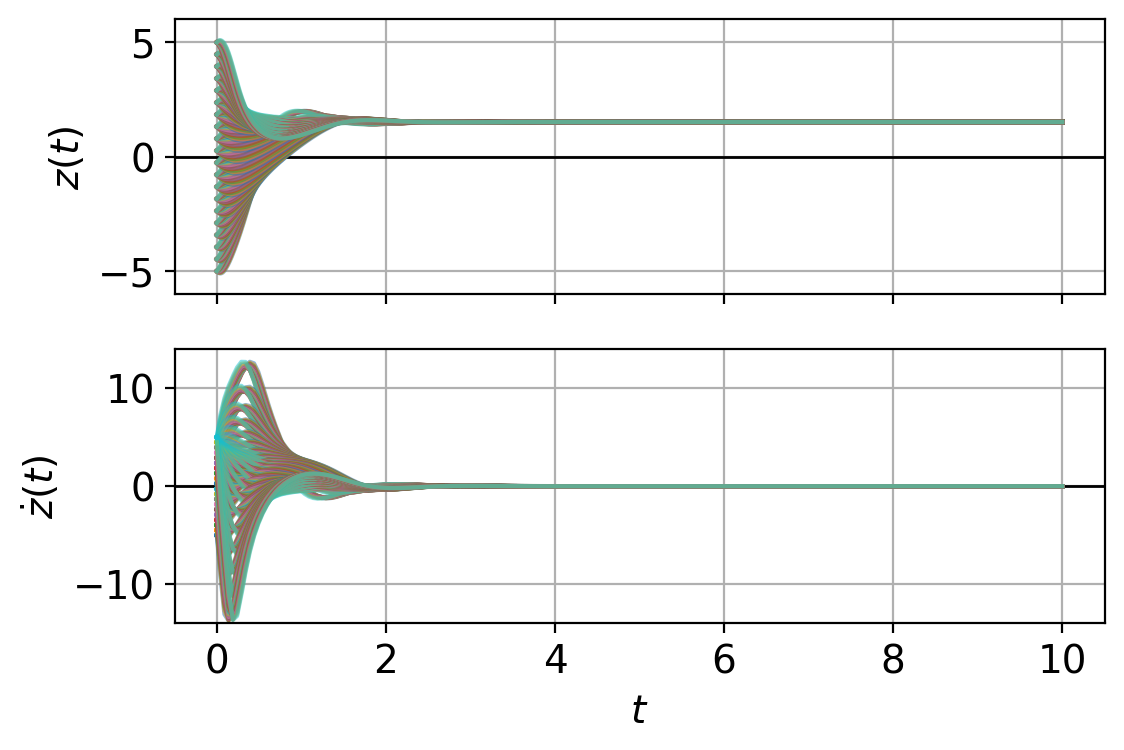}
		\caption{Estimation of basin of attraction (normal NN).}
		\label{subfig. pendulum basin of attraction normal}
	\end{subfigure}
	\begin{subfigure}{0.3\linewidth}
		\centering
		\includegraphics[width=0.95\linewidth]{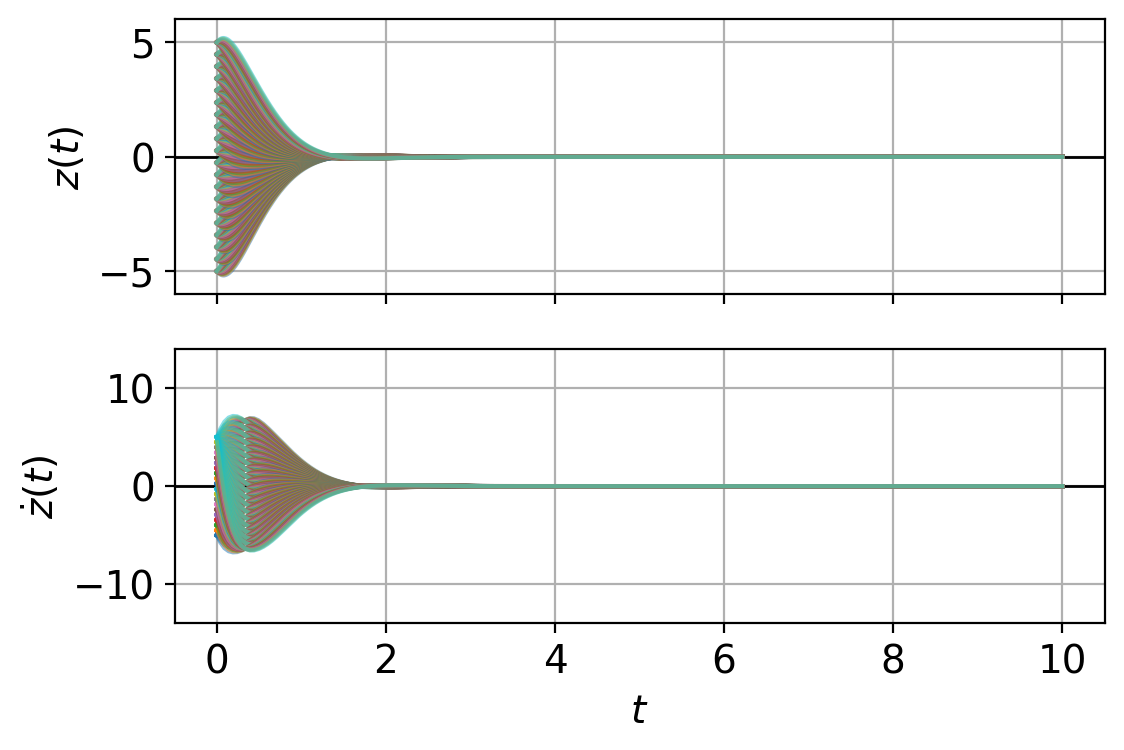}
		\caption{Estimation of basin of attraction (proposed).}
		\label{subfig. pendulum basin of attraction oblique}
	\end{subfigure}
	\begin{subfigure}{0.3\linewidth}
		\centering
		\includegraphics[width=0.95\linewidth]{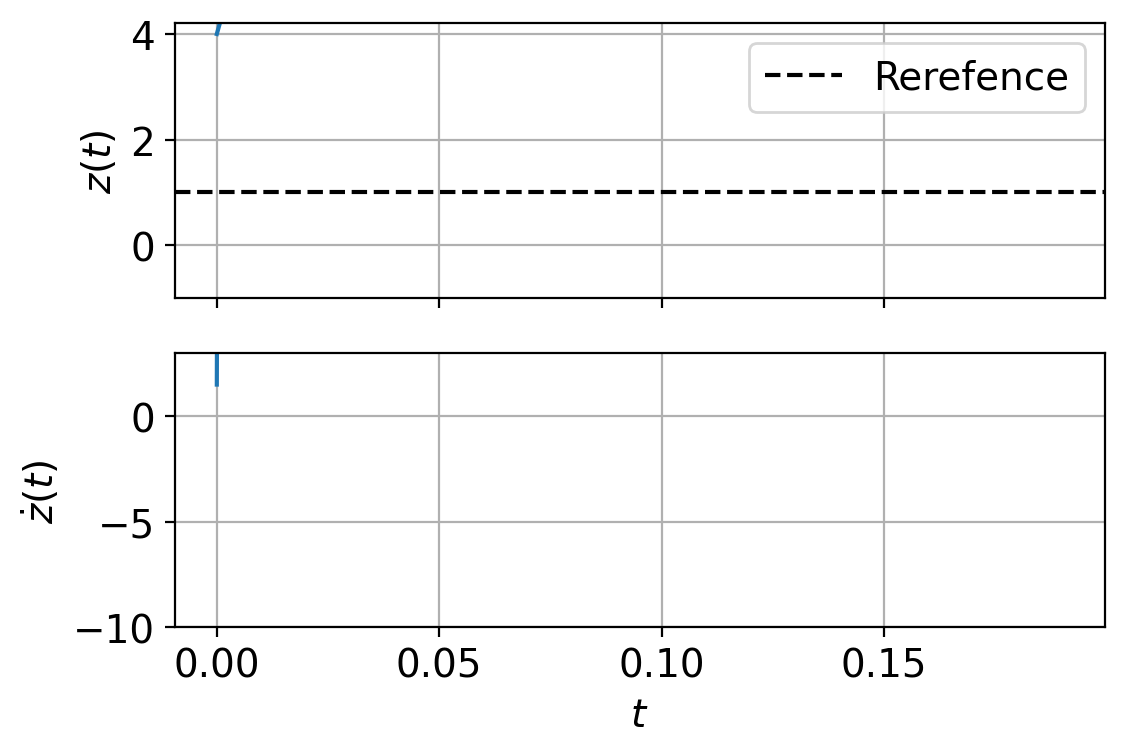}
		\caption{Reference tracking (EDMD).}
		\label{subfig. pendulum cl servo EDMD}
	\end{subfigure}
	\begin{subfigure}{0.3\linewidth}
		\centering
		\includegraphics[width=0.95\linewidth]{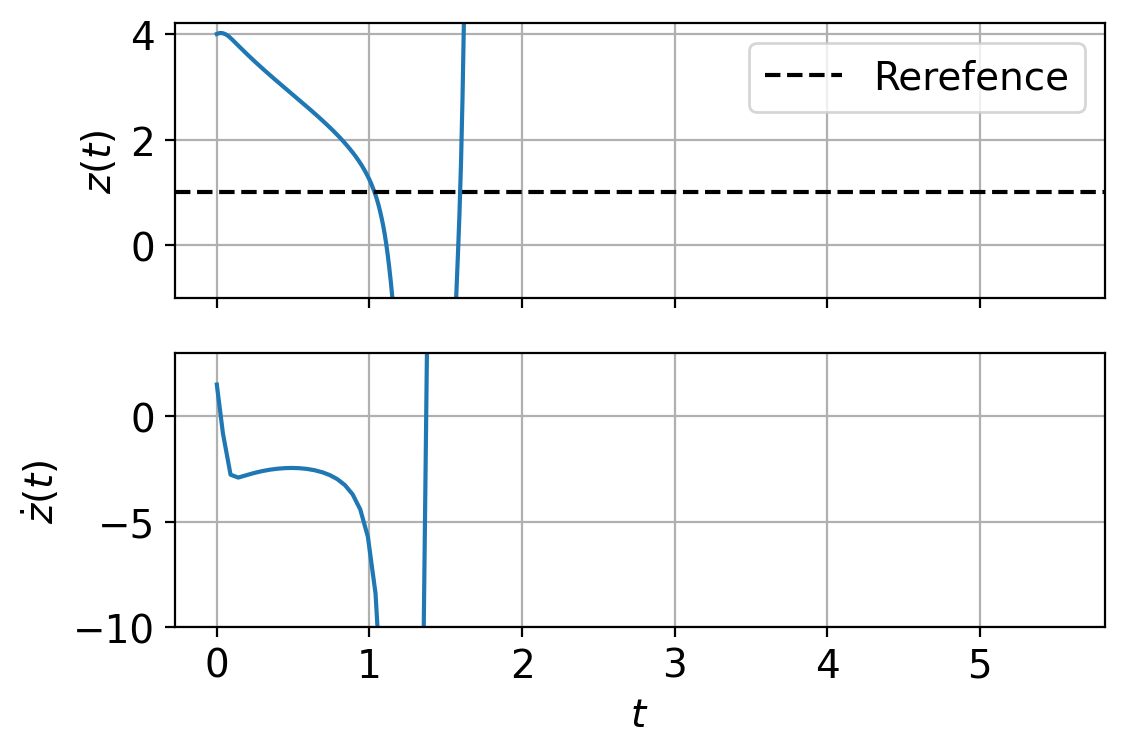}
		\caption{Reference tracking (normal NN).}
		\label{subfig. pendulum cl servo normal}
	\end{subfigure}
	\begin{subfigure}{0.3\linewidth}
		\centering
		\includegraphics[width=0.95\linewidth]{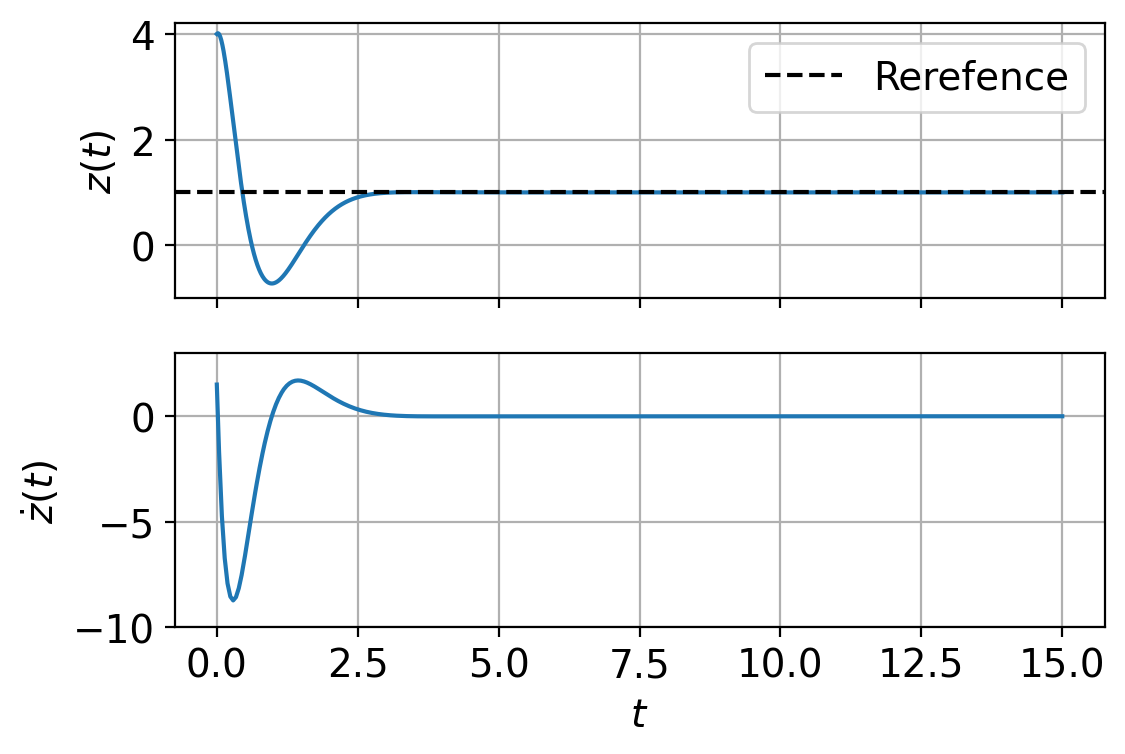}
		\caption{Reference tracking (proposed).}
		\label{subfig. pendulum cl servo oblique}
	\end{subfigure}
	\begin{subfigure}{0.3\linewidth}
		\centering
		\includegraphics[width=0.95\linewidth]{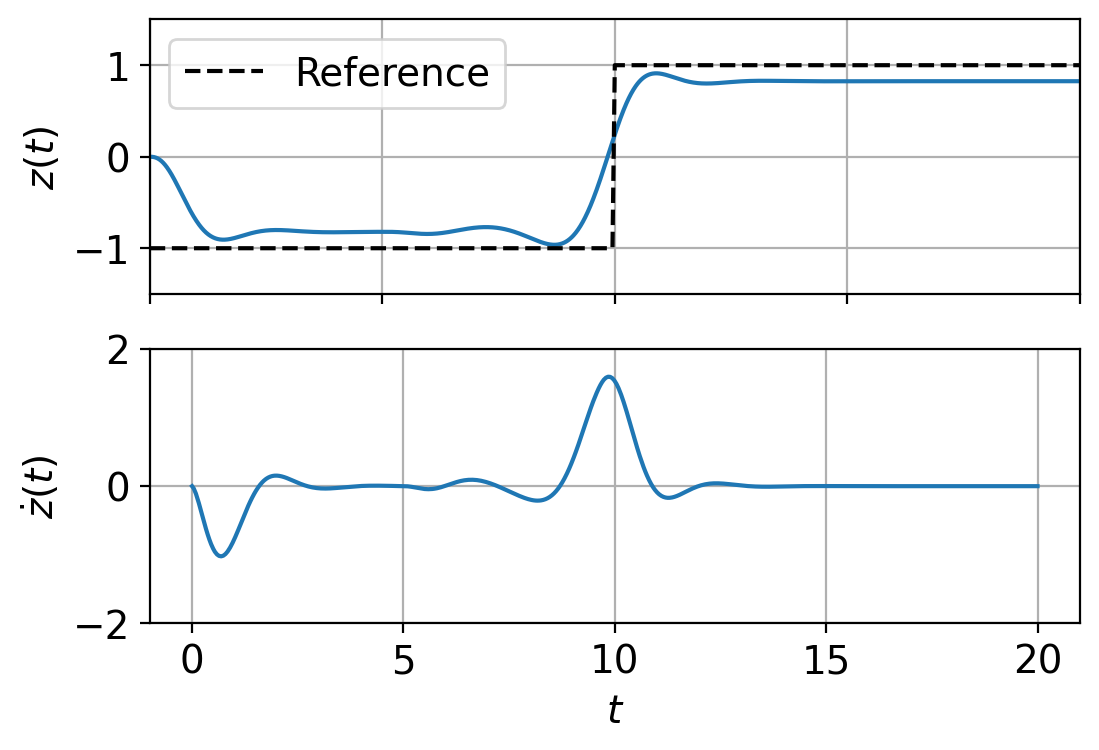}
		\caption{MPC (EDMD).}
		\label{subfig. pendulum mpc EDMD}
	\end{subfigure}
	\begin{subfigure}{0.3\linewidth}
		\centering
		\includegraphics[width=0.95\linewidth]{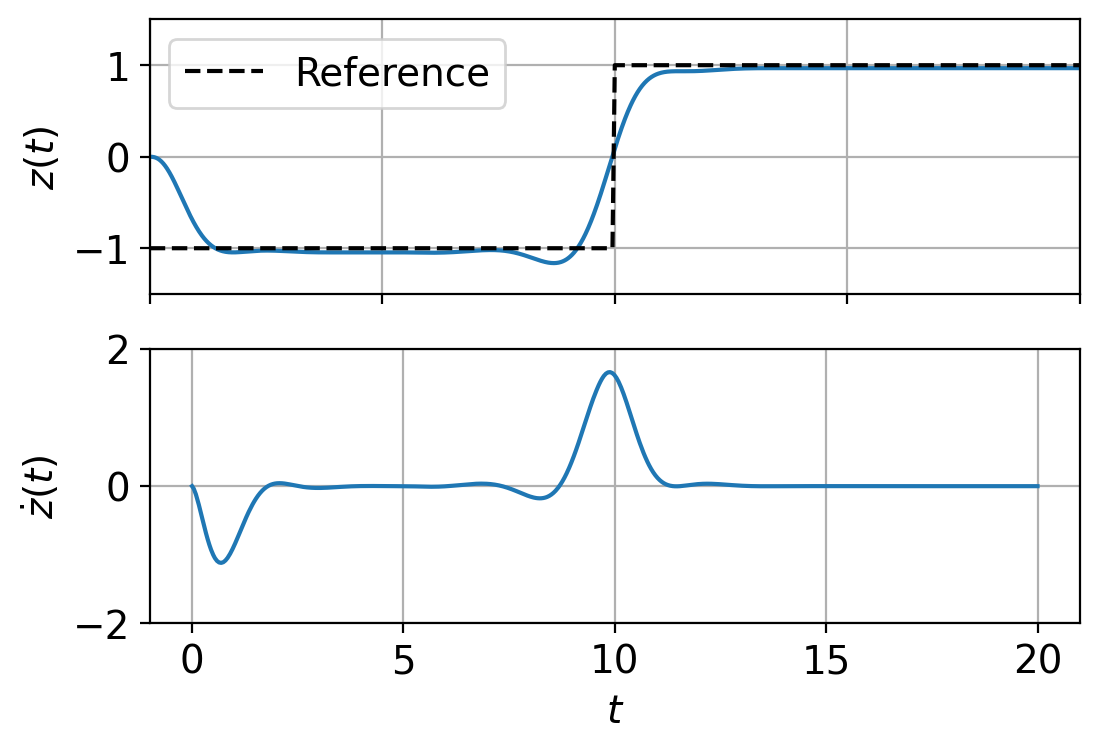}
		\caption{MPC (normal NN).}
		\label{subfig. pendulum mpc normal}
	\end{subfigure}
	\begin{subfigure}{0.3\linewidth}
		\centering
		\includegraphics[width=0.95\linewidth]{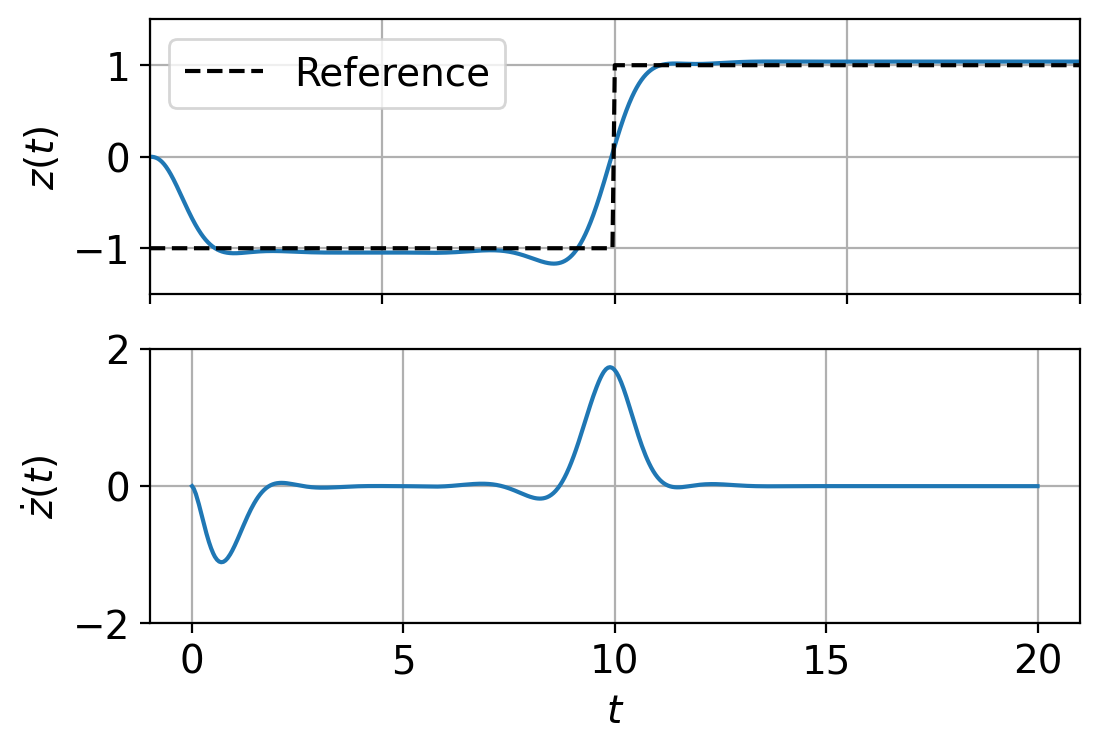}
		\caption{MPC (proposed).}
		\label{subfig. pendulum mpc oblique}
	\end{subfigure}
	\caption{Results of the simple pendulum (control applications).}
	\label{fig. pendulum control applications}
\end{figure}

The system has fixed points $z=k\pi$ ($k\in\mathbb{Z}$), which are stable for $k=0,\pm 2,\pm 4,\cdots,$ and unstable for $k=\pm1, \pm3,\cdots$.
We trained three data-driven models in the same way as the first example of the Duffing oscillator.
The conditions and setups for the controller design tasks are also the same as in the first example.
The results are shown in Figs. \ref{fig. pendulum state prediction} and \ref{fig. pendulum control applications}.
As opposed to the result in the Duffing oscillator example, the EDMD model fails to predict the state accurately for the simple pendulum (Fig. \ref{subfig. pendulum state prediction EDMD}).
This may be inferred from the observation that including more monomial features with higher orders is necessary to reproduce the dynamics near the origin of the state space. 
It is easily seen by the Taylor series expansion that the dynamics of the simple pendulum may be represented by an infinite number of monomials in the vicinity of local linear approximation at the origin as $\sin z(t)\approx z(t) - z^3(t)/3! + z^5(t)/5! - z^7(t)/7! + \cdots$.
Therefore, the given embedded state \eqref{eq. EDMD model embedded state in numerical examples} consisting of monomials up to the third order may not be enough to reproduce the original dynamics.
However, adding more features to this EDMD model will not be a good strategy for control.
Indeed, the embedded state \eqref{eq. EDMD model embedded state in numerical examples} already leads to divergent closed-loop systems in the stabilization by LQR and the reference tracking task (Figs. \ref{subfig. pendulum cl EDMD}, \ref{subfig. pendulum basin of attraction EDMD}, and \ref{subfig. pendulum cl servo EDMD}).
This implies that the effect of modeling error $\mathcal{E}(\chi,u)$ in \eqref{eq. def of modeling error} remains quite high leading to unreasonable approximation, and adding more feature maps in this situation may further increase $\|\mathcal{E}(\chi,u)\|$ as suggested in a numerical example in \cite{control_aware_Koopman}.

On the other hand, the other two neural network-based models have reasonable state predictive accuracy.
The state predictive error contours of both the normal NN and the proposed models are comparable to each other with quite low error profiles (Fig. \ref{subfig. pendulum error contour normal} and \ref{subfig. pendulum error contour oblique}).
Also, there are more successful control applications than the EDMD model (Figs. \ref{subfig. pendulum cl normal}-\ref{subfig. pendulum mpc oblique}).

In addition to the difference between the EDMD and the neural network-based models,
it is also apparent that the proposed model also outperforms the normal NN model in almost all tasks.
Whereas the one-step error contours for the normal NN and the proposed models seem quite similar (Figs. \ref{subfig. pendulum error contour normal} and \ref{subfig. pendulum error contour oblique}), multi-step state predictive accuracy of the proposed model is more accurate than the normal NN model (Figs. \ref{subfig. pendulum state prediction normal} and \ref{subfig. pendulum state prediction oblique}).
Also, the normal NN model only achieves the control objective in the MPC task (Fig. \ref{subfig. pendulum mpc normal}) and the other two control problems result in either a steady-state error (Figs. \ref{subfig. pendulum cl normal} and \ref{subfig. pendulum basin of attraction normal}) or a divergent simulation (Fig. \ref{subfig. pendulum cl servo normal}), whereas the proposed method achieves the control objectives in all tasks (Figs. \ref{subfig. pendulum cl oblique}, \ref{subfig. pendulum basin of attraction oblique}, \ref{subfig. pendulum cl servo oblique}, and \ref{subfig. pendulum mpc oblique}).

\subsection{Rotational/Translational Actuator (RTAC)}
As the third target system,
Rotational/Translational Actuator (RTAC) is considered, which is represented by the following equations\cite{benchmark_problems_nonlinear_control,RTAC_Hinf}:
\vspace{-4mm}
\begin{align}
	\left\{
	\begin{array}{l}
		\dot{z}_1(t) = z_2(t)
		\\
		\dot{z}_2(t) = 
		\cfrac{-z_1(t) + \epsilon z_4^2(t) \sin z_3(t)}{1 - \epsilon^2 \cos^2 z_3(t)}
		-
		\cfrac{\epsilon \cos z_3(t)}{1 - \epsilon^2 \cos^2 z_3(t)}u(t)
		\\
		\dot{z}_3(t) = z_4(t)
		\\
		\dot{z}_4(t) = 
		\cfrac{z_1(t) - \epsilon z_4^2(t) \sin z_3(t)}{1 - \epsilon^2 \cos^2 z_3(t)}
		+
		\cfrac{1}{1 - \epsilon^2 \cos^2 z_3(t)}u(t)
	\end{array}
	\right.,
\end{align}
where $\epsilon:=me/\sqrt{(J+me^2)(M+m)}$.
The variables $z_1(t)$, $z_2(t)$, $z_3(t)$, and $z_4(t)$ correspond to the first, second, third, and fourth components of $\chi\in \mathbb{R}^4$, respectively.
The system has stable fixed points $[z_1\ z_2\ z_3\ z_4]^\tr = [0\ 0\ \alpha\ 0]^\tr$, where $\alpha\in \mathbb{R}$ is arbitrary.

In this example, an EDMD model was trained with monomial feature maps up to the second order, which yields the 14-dimensional embedded state:
\vspace{-2mm}
\begin{align}
	[g_1(\chi)\cdots g_{14}(\chi)]^\tr=
	[\chi^\tr\ \ \chi_{(1)}^2\ \ \chi_{(2)}^2\ \ \chi_{(3)}^2\ \ \chi_{(4)}^2\cdots
	]^\tr,
	\label{eq. EDMD model embedded state in numerical examples RTAC}
\end{align}
where $[\chi_{(1)}\ \chi_{(2)}\ \chi_{(3)}\ \chi_{(4)}]^\tr:=\chi$.

For the normal NN and the proposed models, we set the embedded state as follows.
Three feature maps are included so that the embedded state is of the form $[\chi\ g_5(\chi)\ g_6(\chi)\ g_7(\chi)]^\tr$, where $[g_5(\chi)\ g_6(\chi)\ g_7(\chi)]^\tr$ is characterized by a fully connected feed-forward neural network with a single hidden layer consisting of 25 neurons.
For the test functions of the proposed method, we adopted a specific structure
\vspace{-2mm}
\begin{align}
	[\varphi_1(\chi,u)\cdots \varphi_6(\chi,u)]^\tr=[\chi^\tr\ u\ \varphi_6(\chi)]^\tr,
	\label{eq. specific test functions for RTAC}
\end{align}
with $\varphi_6(\chi)$ be a neural network consisting of a single hidden layer with 10 neurons.
	We found that imposing the above  structure on the design of test functions in \eqref{eq. specific test functions for RTAC} leads to not only better performance of the model, but also robust neural network learning as reported in Section 1 of the supplementary material.
	It suggests that flexible modeling is allowed by having more control over the learning process itself, compared to jointly learning $[\bm{A}\ \bm{B}]$ and $g_i$ as in \cite{control_aware_Koopman}, where it is difficult to directly regulate the neural network training.

In the LQR design, we set the weight matrices as in \eqref{eq. LQR weights duffing} with $Q_\text{state}$ modified as 
\vspace{-2mm}
\begin{align}
    Q_\text{state}=
    \left[
        \begin{array}{cccc}
             100 & 0 & 0 & 0   \\
             0 & 1 & 0 & 0 \\
             0 & 0 & 100 & 0 \\
             0 & 0 & 0 & 1
        \end{array}
    \right].
\end{align}

\vspace{-2mm}
In the reference tracking problem, we reset $\bm{C}=[0\ 0\ 1\ 0]$ with $r=1$.
In the MPC design, we modified the objective function in \eqref{eq. MPC objectives duffing} as $l(\xi,u_k,k)=(\xi_{k,(3)} - r_\text{MPC}(k))^2$ so that the third component of the state will track the reference signal.

\begin{figure}[]
	\centering 
	\begin{subfigure}{0.3\linewidth}
		\centering
		\includegraphics[width=0.95\linewidth]{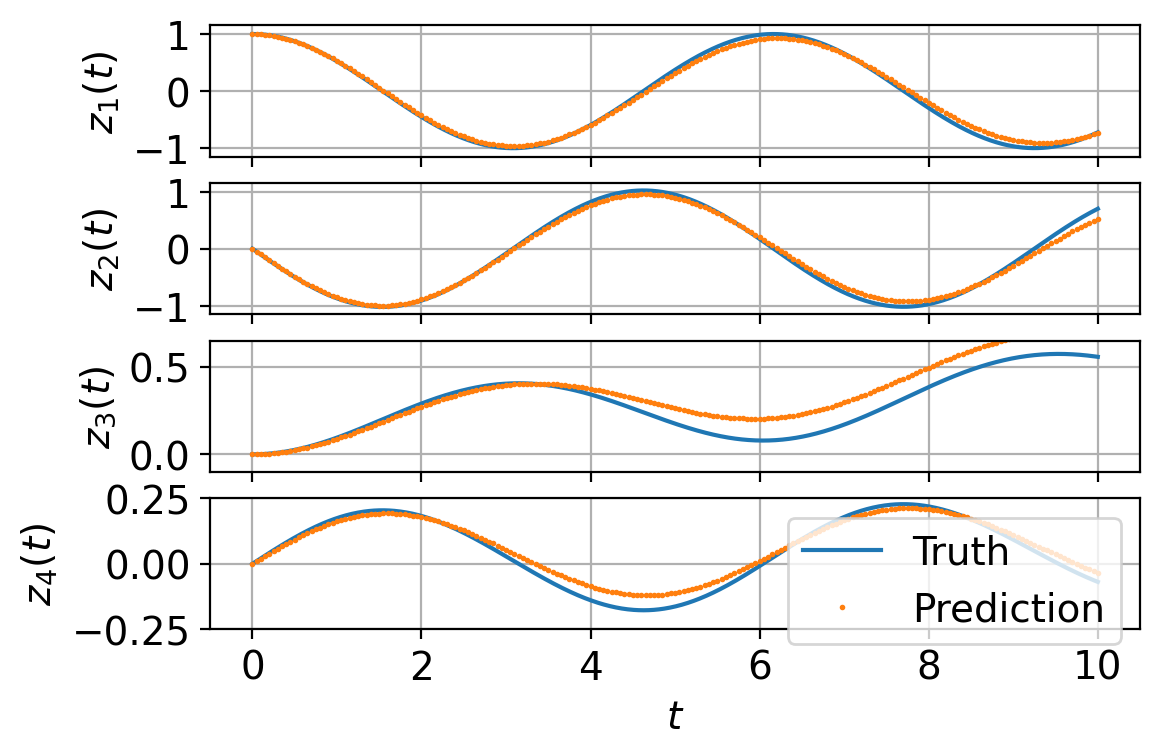}
		\caption{State prediction (EDMD).}
		\label{subfig. RTAC state prediction EDMD}
	\end{subfigure}
	\begin{subfigure}{0.3\linewidth}
		\centering
		\includegraphics[width=0.95\linewidth]{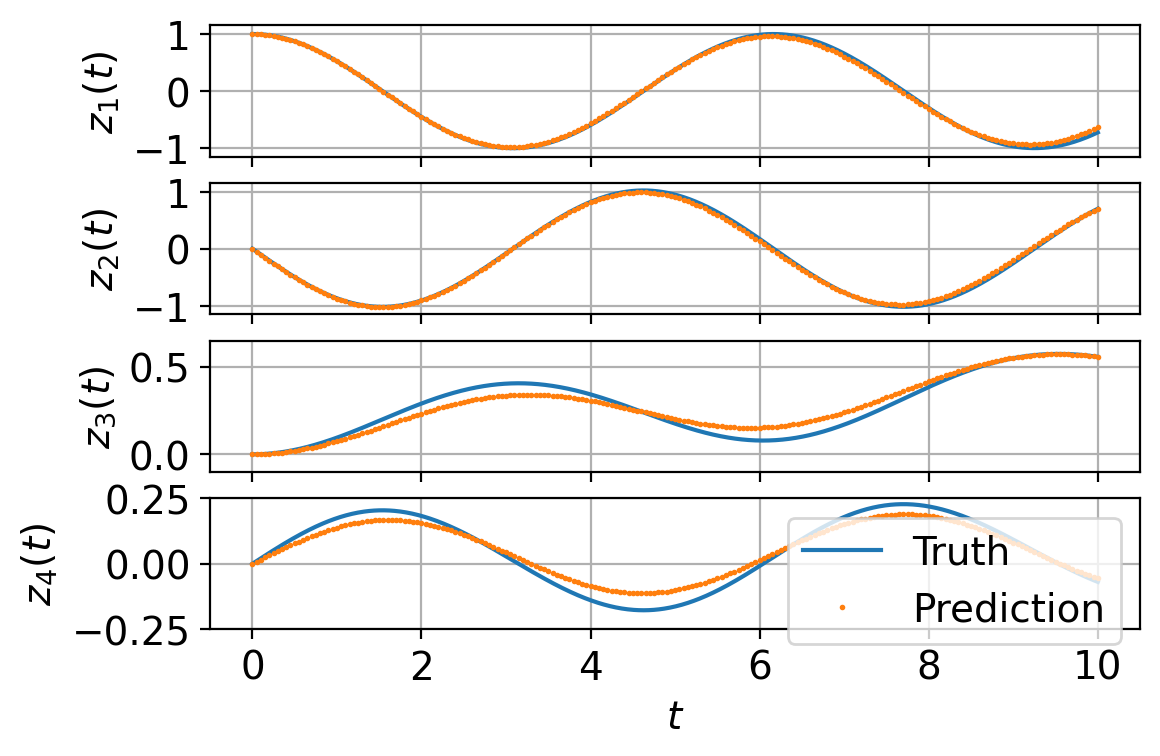}
		\caption{\scriptsize{State prediction (normal NN).}}
		\label{subfig. RTAC state prediction normal}
	\end{subfigure}
	\begin{subfigure}{0.3\linewidth}
		\centering
		\includegraphics[width=0.95\linewidth]{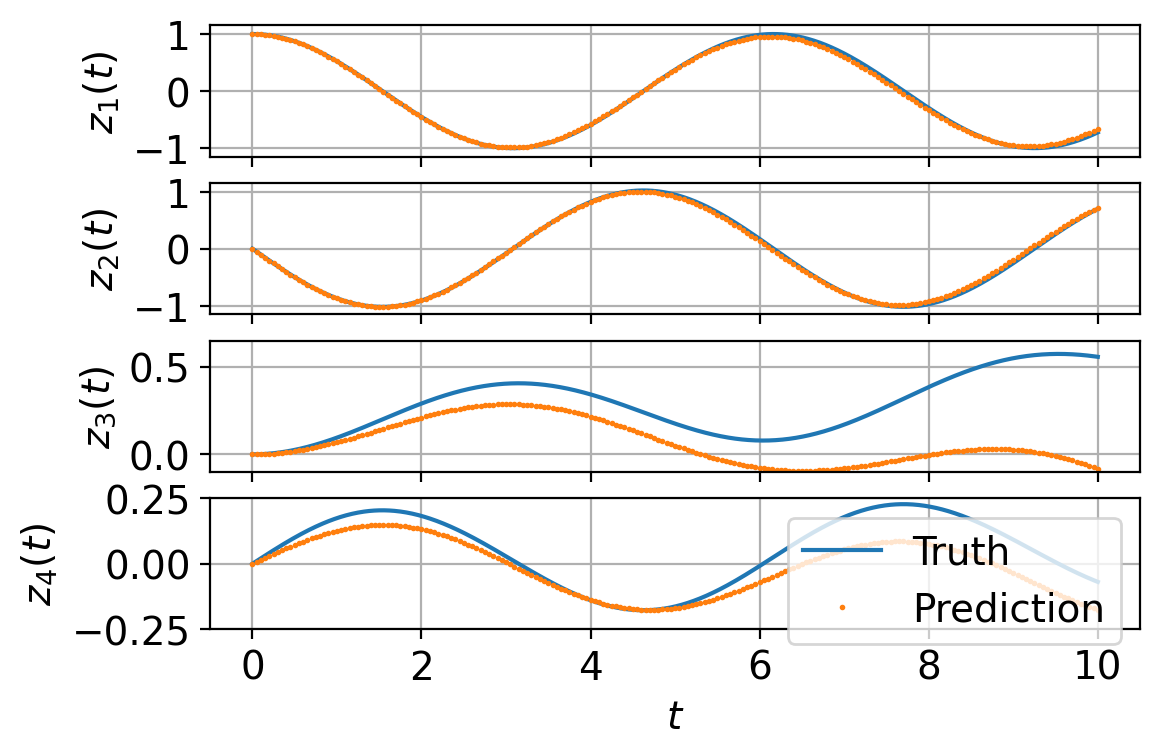}
		\caption{\footnotesize{State prediction (proposed).}}
		\label{subfig. RTAC state prediction oblique}
	\end{subfigure}
	\begin{subfigure}{0.3\linewidth}
		\centering
		\includegraphics[width=0.95\linewidth]{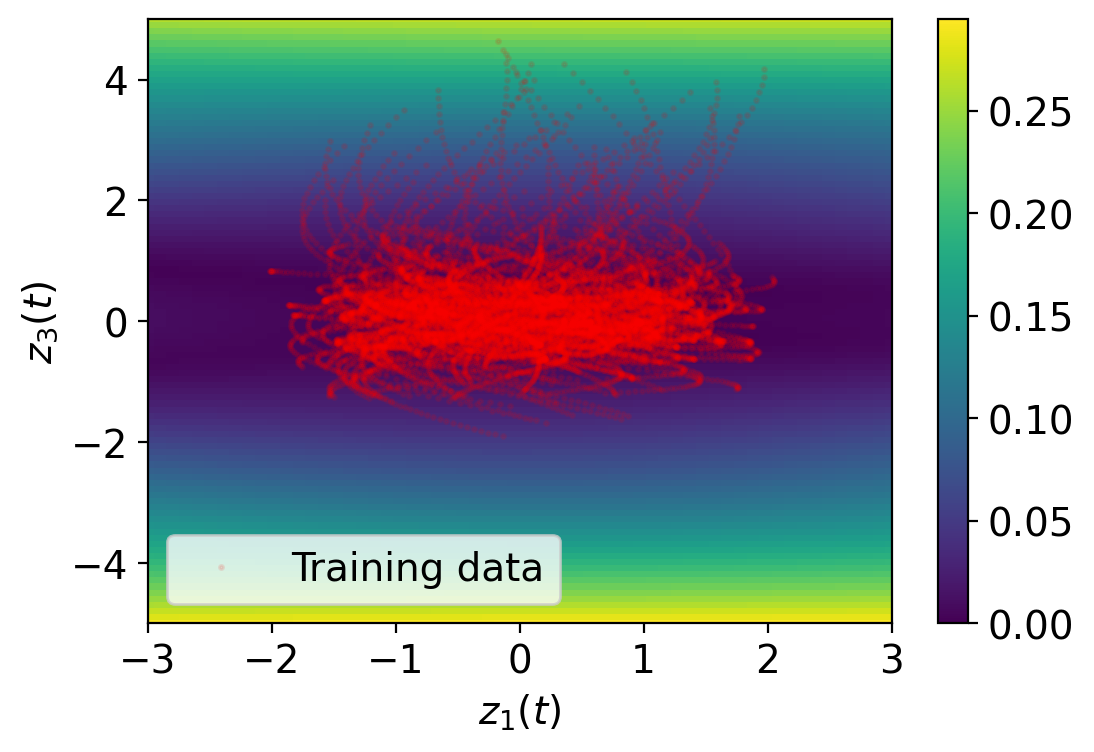}
		\caption{Error contour (EDMD).}
		\label{subfig. RTAC error contour EDMD}
	\end{subfigure}
	\begin{subfigure}{0.3\linewidth}
		\centering
		\includegraphics[width=0.95\linewidth]{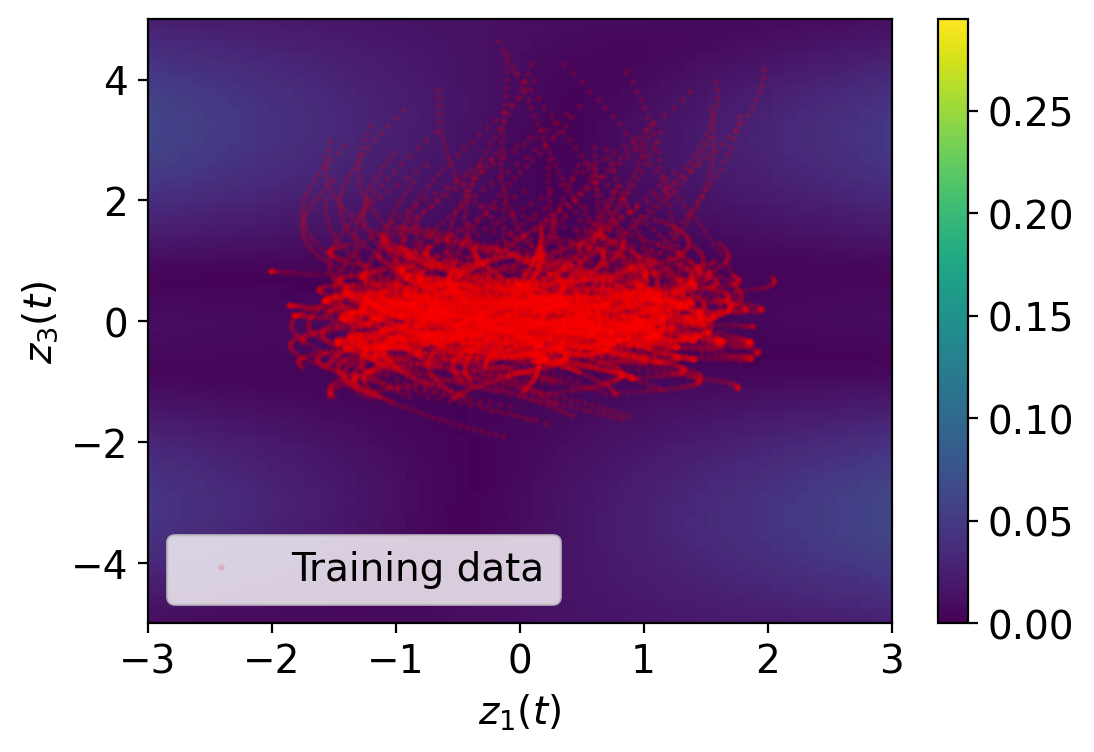}
		\caption{\footnotesize{Error contour (normal NN).}}
		\label{subfig. RTAC error contour normal}
	\end{subfigure}
	\begin{subfigure}{0.3\linewidth}
		\centering
		\includegraphics[width=0.95\linewidth]{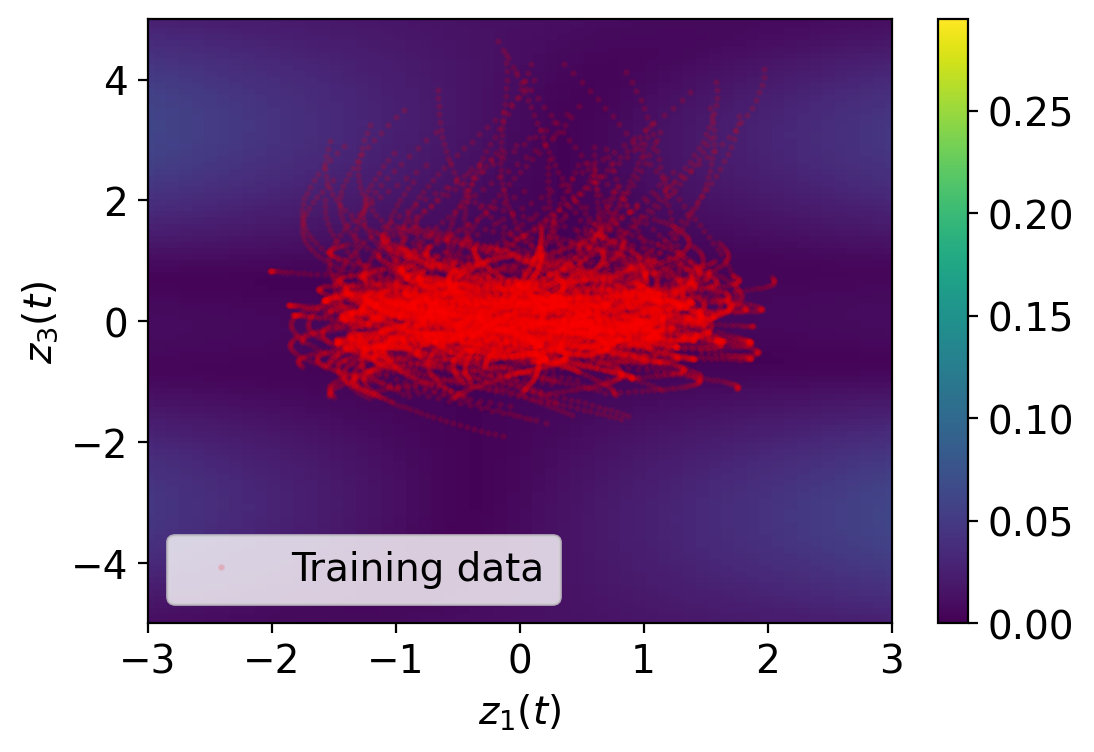}
		\caption{Error contour (proposed).}
		\label{subfig. RTAC error contour oblique}
	\end{subfigure}
	\caption{Results of RTAC (state prediction). In error contour plots, both $z_2$ and $z_4$ are fixed to 0.}
	\label{fig. RTAC state prediction}
\end{figure}

\begin{figure}[H]
	\centering 
	\begin{subfigure}{0.3\linewidth}
		\centering
		\includegraphics[width=0.95\linewidth]{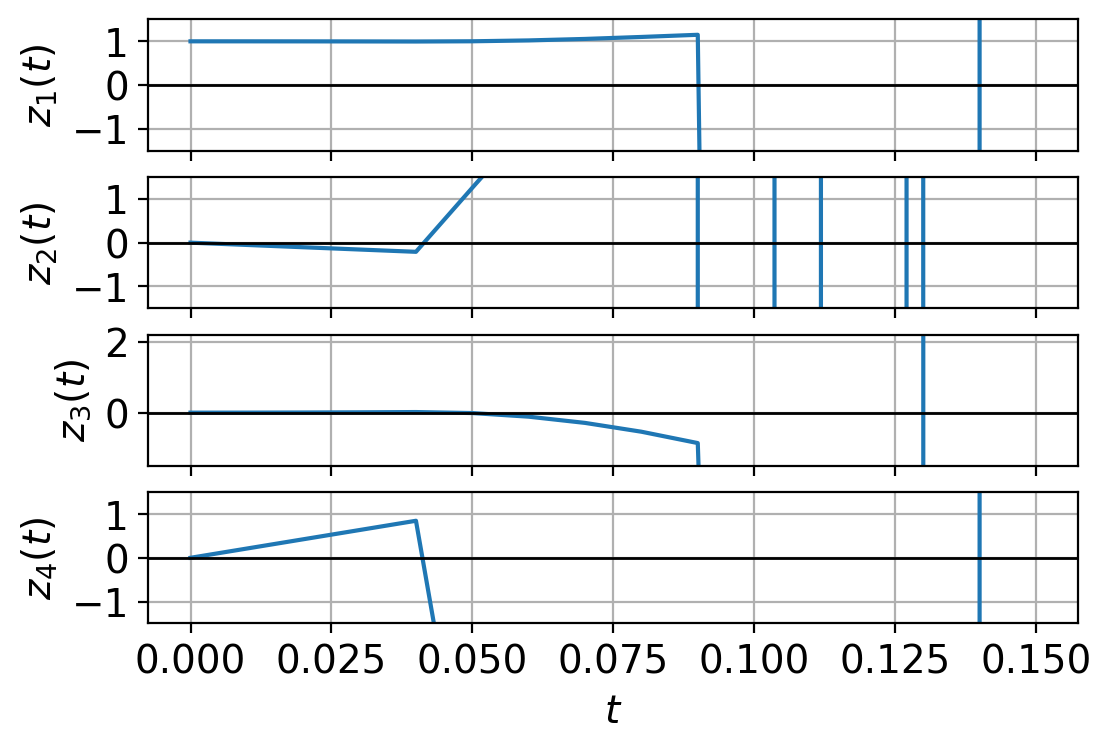}
		\caption{Stabilization by LQR (EDMD).}
		\label{subfig. RTAC cl EDMD}
	\end{subfigure}
	\begin{subfigure}{0.3\linewidth}
		\centering
		\includegraphics[width=0.95\linewidth]{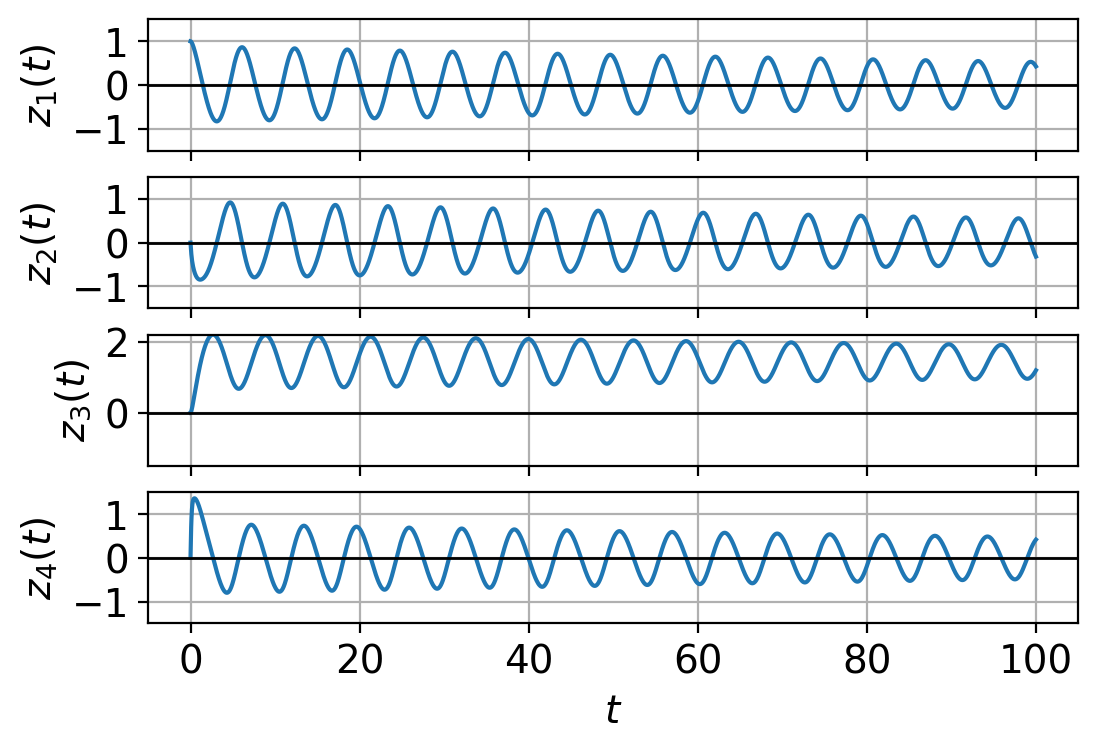}
		\caption{\footnotesize{Stabilization by LQR (normal NN).}}
		\label{subfig. RTAC cl normal}
	\end{subfigure}
	\begin{subfigure}{0.3\linewidth}
		\centering
		\includegraphics[width=0.95\linewidth]{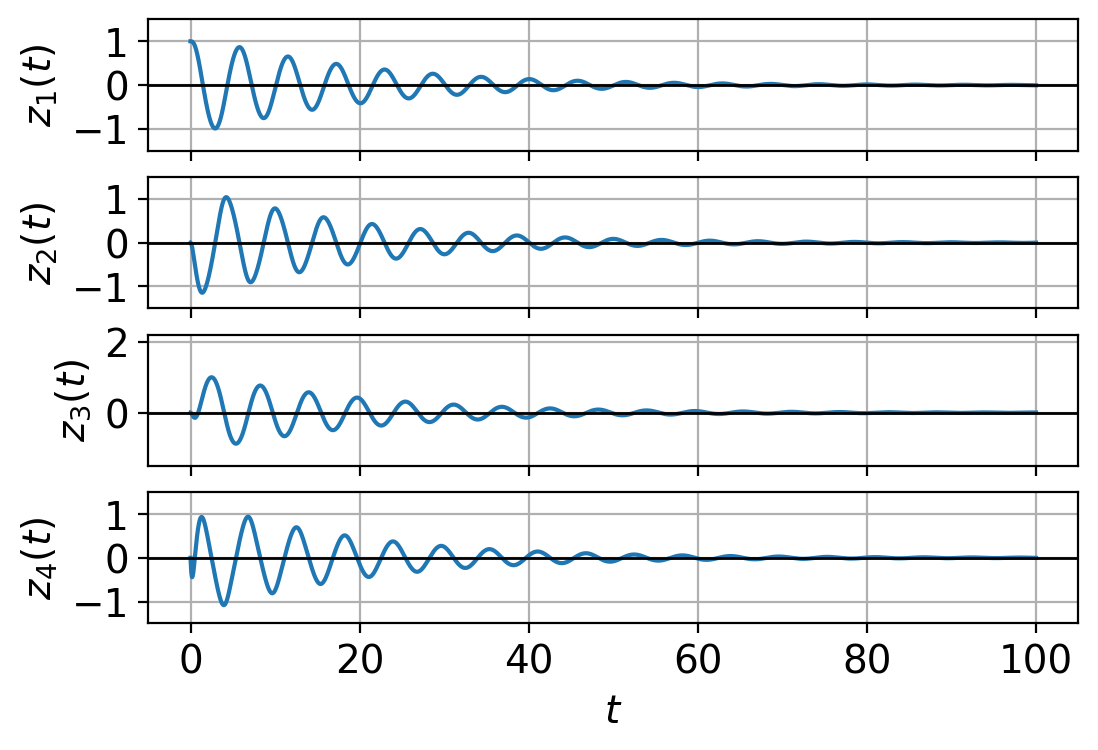}
		\caption{Stabilization by LQR (proposed).}
		\label{subfig. RTAC cl oblique}
	\end{subfigure}
	\begin{subfigure}{0.3\linewidth}
		\centering
		\includegraphics[width=0.95\linewidth]{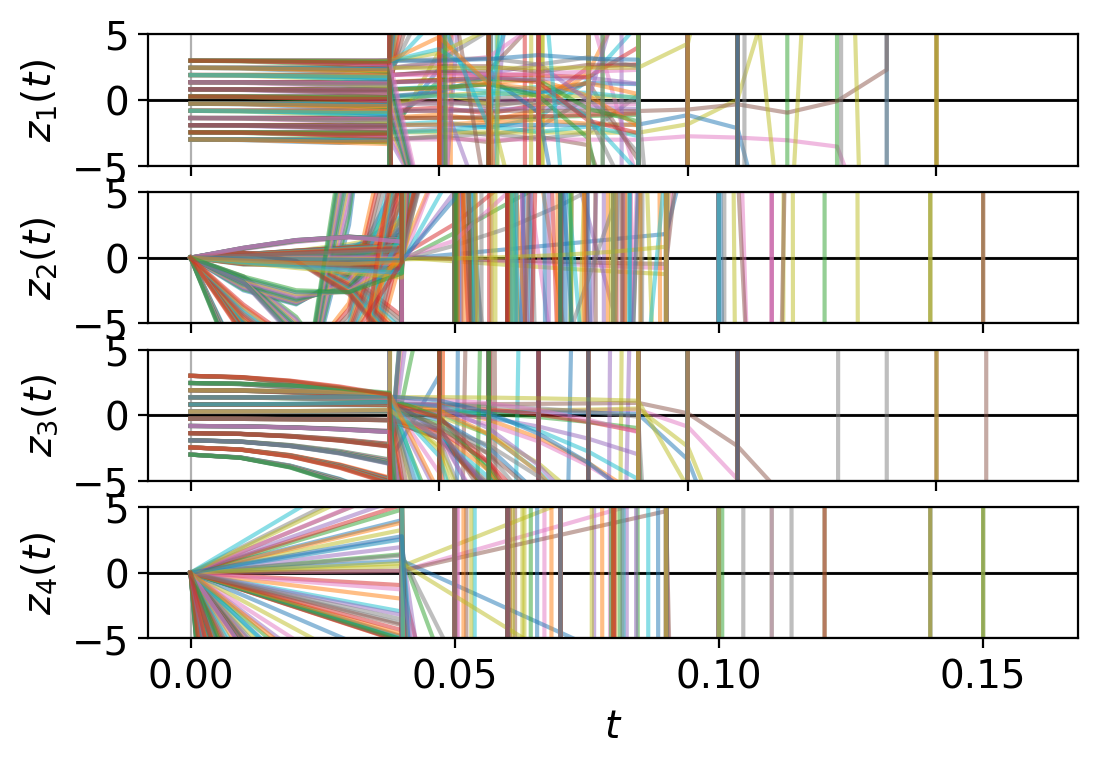}
		\caption{Estimation of basin of attraction (EDMD).}
		\label{subfig. RTAC basin of attraction EDMD}
	\end{subfigure}
	\begin{subfigure}{0.3\linewidth}
		\centering
		\includegraphics[width=0.95\linewidth]{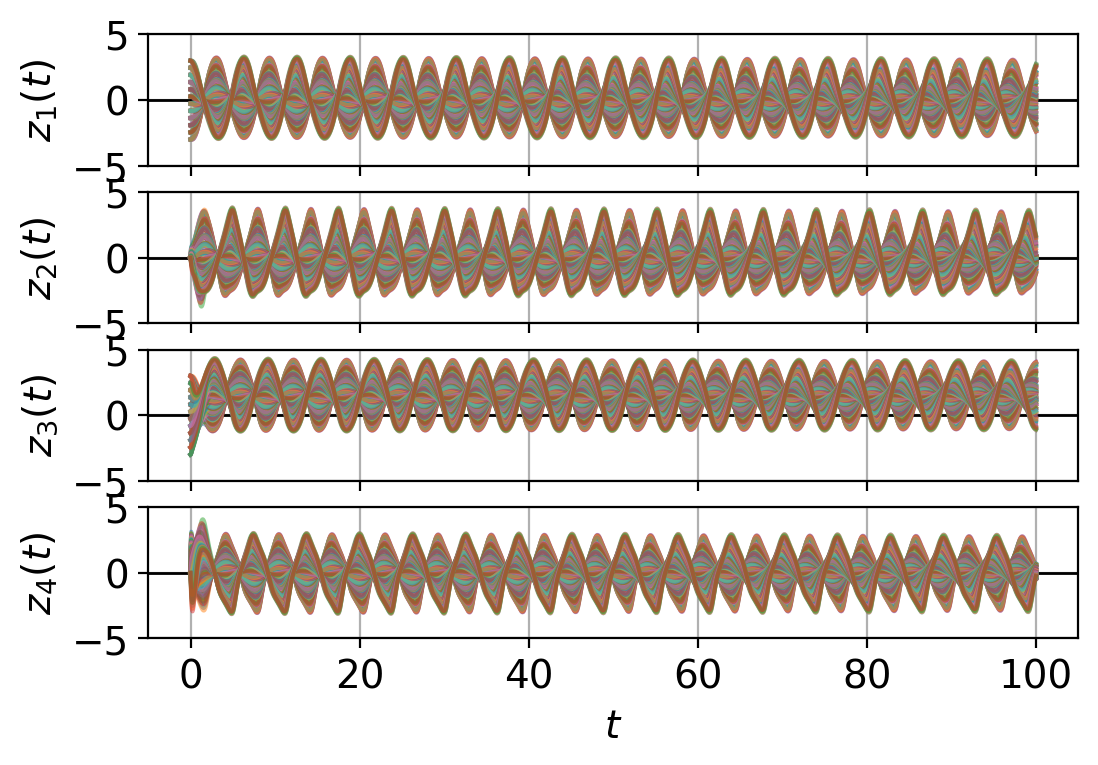}
		\caption{Estimation of basin of attraction (normal NN).}
		\label{subfig. RTAC basin of attraction normal}
	\end{subfigure}
	\begin{subfigure}{0.3\linewidth}
		\centering
		\includegraphics[width=0.95\linewidth]{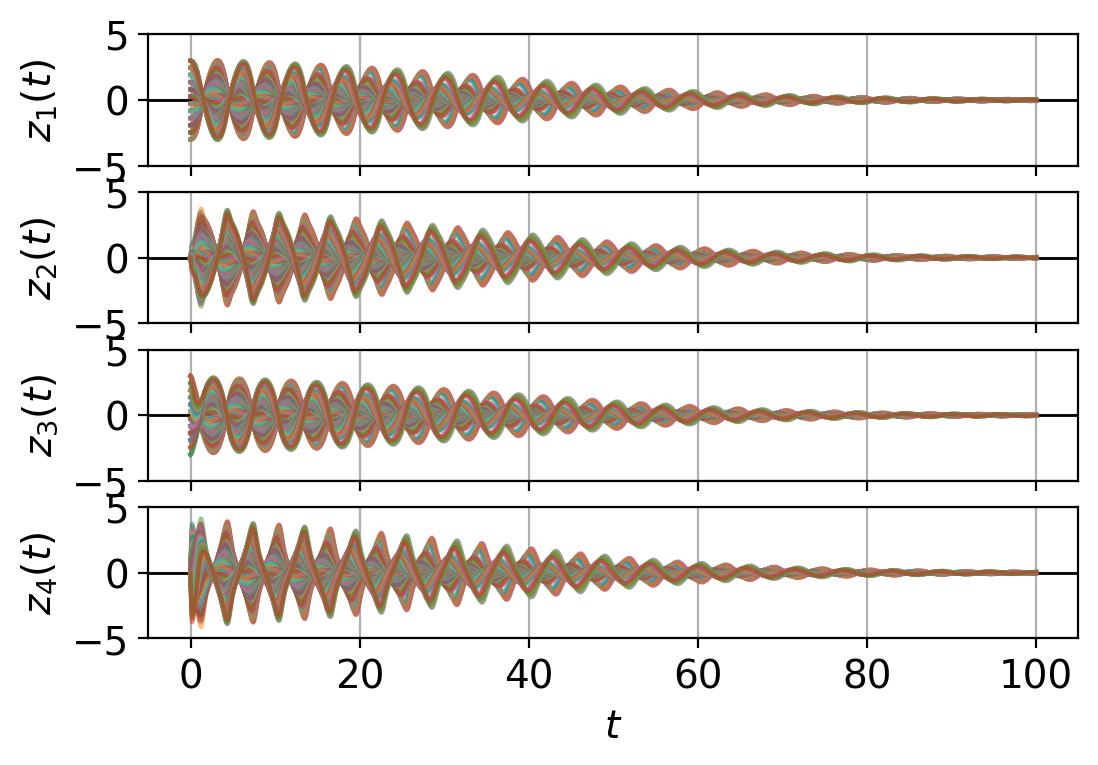}
		\caption{Estimation of basin of attraction (proposed).}
		\label{subfig. RTAC basin of attraction oblique}
	\end{subfigure}
	\begin{subfigure}{0.3\linewidth}
		\centering
		\includegraphics[width=0.95\linewidth]{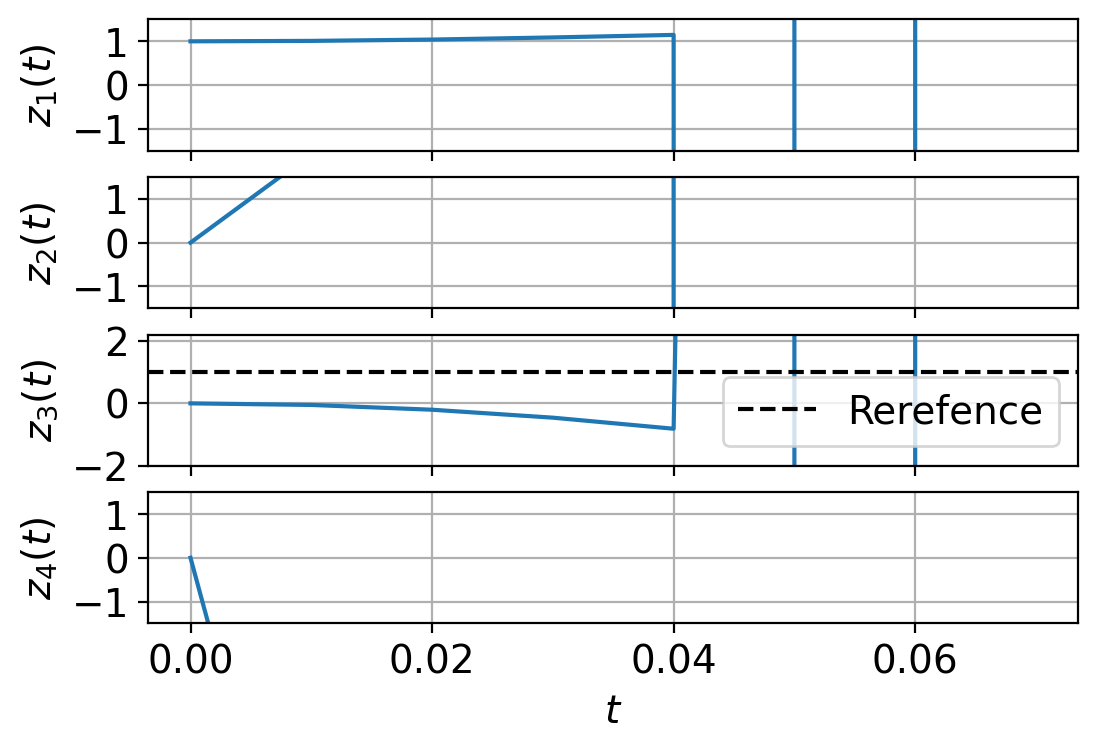}
		\caption{Reference tracking (EDMD).}
		\label{subfig. RTAC cl servo EDMD}
	\end{subfigure}
	\begin{subfigure}{0.3\linewidth}
		\centering
		\includegraphics[width=0.95\linewidth]{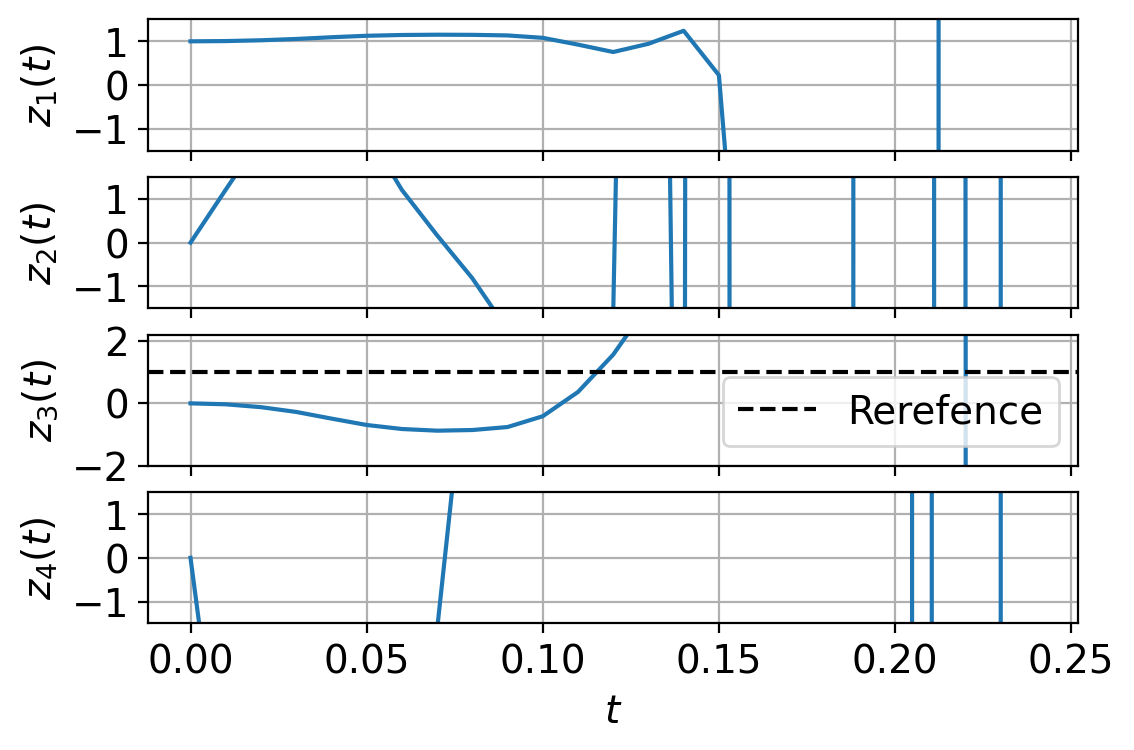}
		\caption{Reference tracking (normal NN).}
		\label{subfig. RTAC cl servo normal}
	\end{subfigure}
	\begin{subfigure}{0.3\linewidth}
		\centering
		\includegraphics[width=0.95\linewidth]{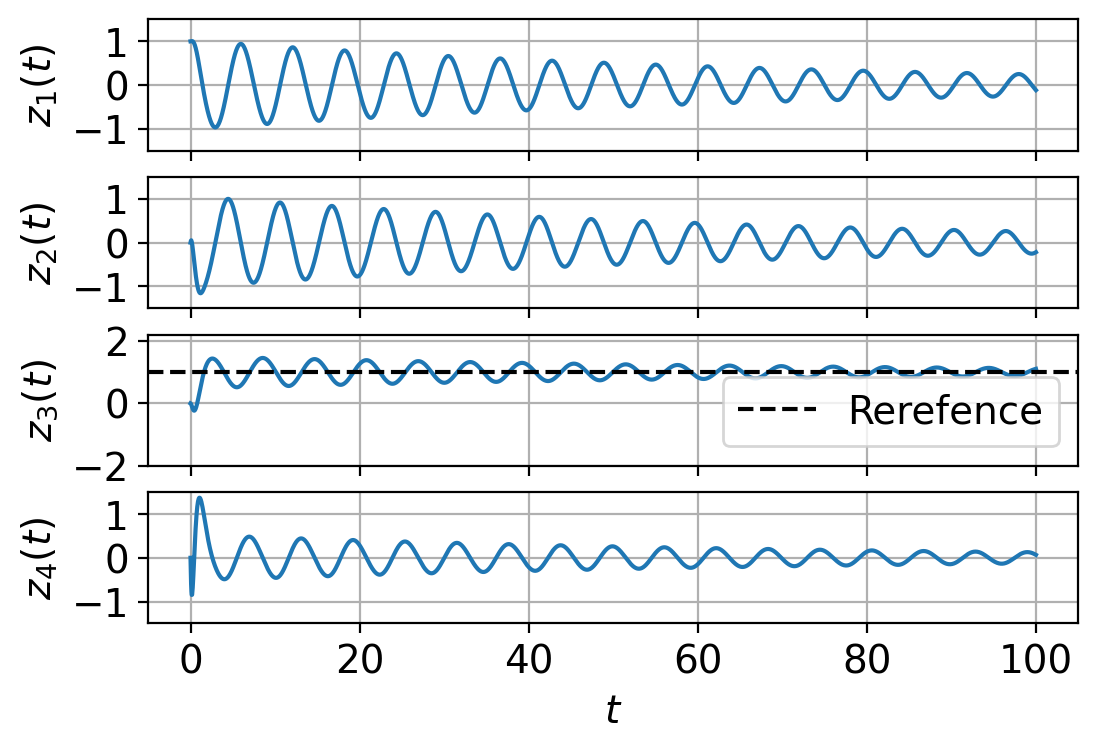}
		\caption{Reference tracking (proposed).}
		\label{subfig. RTAC cl servo oblique}
	\end{subfigure}
	\begin{subfigure}{0.3\linewidth}
		\centering
		\includegraphics[width=0.95\linewidth]{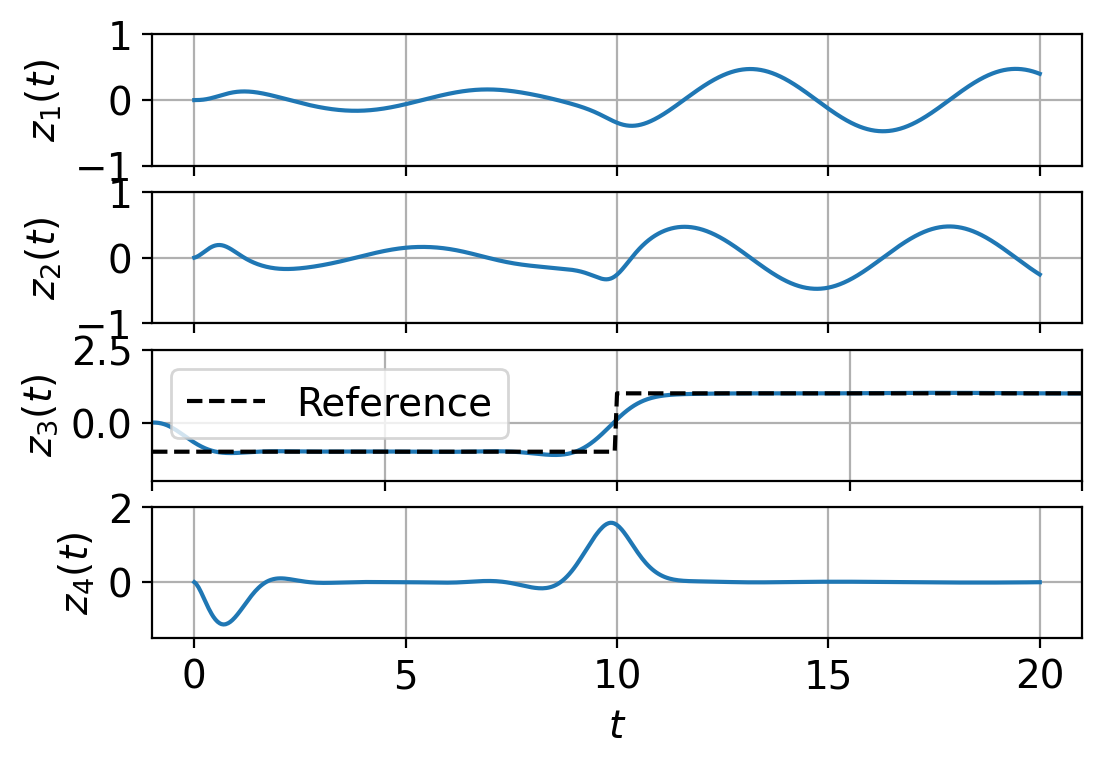}
		\caption{MPC (EDMD).}
		\label{subfig. RTAC mpc EDMD}
	\end{subfigure}
	\begin{subfigure}{0.3\linewidth}
		\centering
		\includegraphics[width=0.95\linewidth]{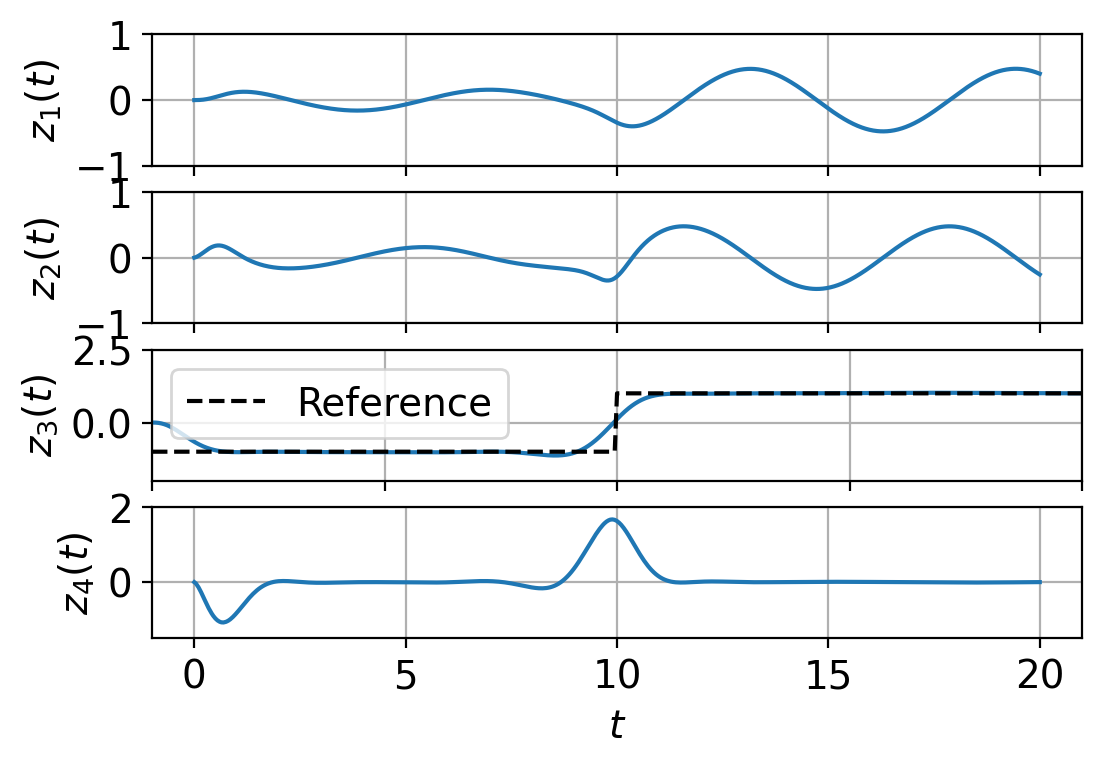}
		\caption{MPC (normal NN).}
		\label{subfig. RTAC mpc normal}
	\end{subfigure}
	\begin{subfigure}{0.3\linewidth}
		\centering
		\includegraphics[width=0.95\linewidth]{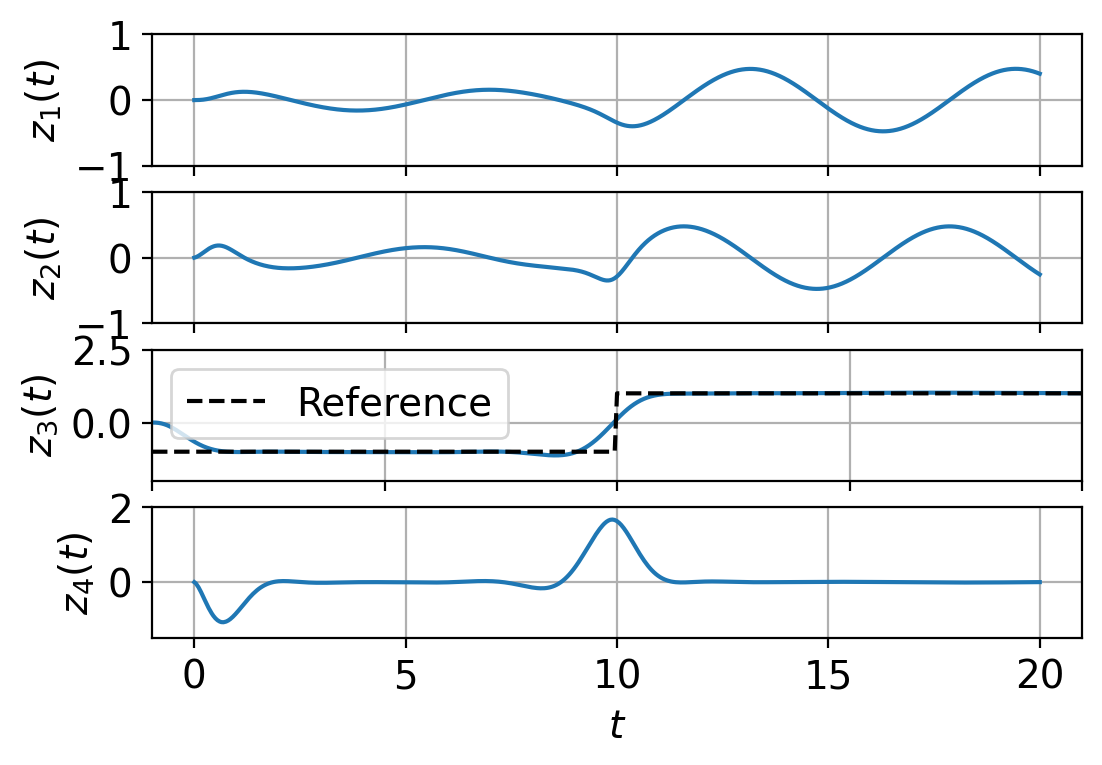}
		\caption{MPC (proposed).}
		\label{subfig. RTAC mpc oblique}
	\end{subfigure}
	\caption{Results of RTAC (control applications).}
	\label{fig. RTAC control applications}
\end{figure}

The results are shown in Figs. \ref{fig. RTAC state prediction} and \ref{fig. RTAC control applications}.
First, the error contour of the EDMD model shows that its one-step state predictive accuracy is not good far from the training data regime (Fig. \ref{subfig. RTAC error contour EDMD}), whereas the normal NN and the proposed models have lower error profiles across a wide range of operating points (Figs. \ref{subfig. RTAC error contour normal} and \ref{subfig. RTAC error contour oblique}).
Note that both $z_2(t)$ and $z_4(t)$ are fixed to 0 to visualize the error profiles in two dimensions.
On the other hand, multi-step time-series state predictions show different results, where the predictions of $z_3(t)$ and $z_4(t)$ of the proposed model deviate from the true values towards the end of the simulation (Fig. \ref{subfig. RTAC state prediction oblique}), whereas the other two models output more reasonable state predictions for all the four variables (Figs. \ref{subfig. RTAC state prediction EDMD} and \ref{subfig. RTAC state prediction normal}). 
Noticing that the multi-step predictions start from an initial condition in the training data regime, it implies that the proposed model especially has lower state predictive accuracy on the data points. This is the only result that the performance of the proposed model could not surpass the other two models in the evaluations of this paper.

As for control, the proposed method outperforms the other models in all applications.
Except for the MPC task, the EDMD and the normal NN models result in either divergent closed-loop dynamics or undesirable oscillatory behaviors (Figs. \ref{subfig. RTAC cl EDMD}, \ref{subfig. RTAC basin of attraction EDMD}, \ref{subfig. RTAC cl servo EDMD}, \ref{subfig. RTAC cl normal}, \ref{subfig. RTAC basin of attraction normal}, and \ref{subfig. RTAC cl servo normal}).
On the other hand, the proposed method enables successful control application in every task (Figs. \ref{subfig. RTAC cl oblique}, \ref{subfig. RTAC basin of attraction oblique}, \ref{subfig. RTAC cl servo oblique}, and \ref{subfig. RTAC mpc oblique}), which shows its generalizability to different types of tasks.

%\vspace{-2mm}
\section{Conclusion}
We propose a new data-driven modeling method for nonlinear,\\ non-autonomous dynamics based on the concept of linear embedding with oblique projection in a Hilbert space. 
Linear embedding models, which naturally arise as a practical class of models in the Koopman operator framework, have the advantage that linear systems theories can be applied to controller synthesis problems even if the target dynamics is nonlinear. However, there are fundamental limitations regarding the accuracy of the models.
In addition to convergence issues associated with Extended Dynamic Mode Decomposition (EDMD) in the non-autonomous setting, we provided a necessary condition for a linear embedding model to achieve zero modeling error as well as subsequent analyses that suggest a fundamental difficulty of obtaining a model that is accurate for a wide range of state/input values.
This condition reveals a trade-off relation between the expressivity of the model introduced by nonlinear feature maps and the unique model structure of being strictly linear w.r.t. the input to allow the use of linear systems theories for nonlinear dynamics.

To achieve good performance despite these fundamental limitations,
we proposed a neural network-based modeling method combined with two-staged learning, which is derived from a weak formulation of projection-based linear operator learning. After initializing a model with orthogonal projection only, which ensures a least square type of optimality w.r.t. given data points, the main training process follows in which test functions that characterize oblique projection are optimized on different data points along with the feature maps aiming to improve the generalizability of the model so that it can be applicable to many different tasks, including state prediction and various controller design problems.
The use of oblique projections is expected to alter the property of the initialized model so that the model accuracy improves for certain dynamics and a simple example was investigated to showcase that the proposed method is effective when the target dynamics possesses the property of non-normality.

To evaluate the proposed method, comprehensive  studies were conducted, where
four different tasks: state prediction, stabilization by Linear Quadratic Regulator (LQR), reference tracking, and linear Model Predictive Control (MPC) were considered. In each task, we compared the proposed method with other data-driven linear embedding models, targeting three nonlinear dynamics: the Duffing oscillator, the simple pendulum, and the Rotational/Translational Actuator (RTAC).
The superiority and effectiveness of the proposed model were confirmed through these various examples and it was shown that it has successfully gained sufficient generalizability.

As a future direction of the research,
relaxing the strictly linear structure of the model w.r.t. the input may lead to better accuracy of the model. 
For instance, bilinear model structures can be advantageous for a wider range of nonlinear dynamics\cite{advantages_bilinearization}.
Also, exploiting prior knowledge about the unknown dynamics may allow for a more rigorous basis for the method from a theoretical perspective.
If one restricts the attention to control-affine systems, several theoretical analyses become available, e.g., controllability analysis\cite{Koopman_bilinearization} and finite data error bounds\cite{towards_reliable_EDMD_control_affine,finite_data_error_bounds}.
For the proposed method developed in this paper, it may be also possible to supplement the modeling framework with similar mathematical analyses based on these works.
	This work considered situations where the dimension of the embedded space of the model is higher than the target dynamics. It may be also useful to explore whether the proposed method can be effectively applied to dynamical systems with large dimensions, in which case the problem is reduced to finding a latent space that can capture the essential properties of the higher-dimensional dynamics.
	Since the proposed method is only based on the single-step error it may be difficult for the model to learn long-time behaviors. Taking multi-step errors into account may help construct accurate models for these more complex dynamics.

\section*{Acknowledgments}
This work was funded by APRA-E under the project \textit{SAFARI: Secure Automation for Advanced Reactor Innovation}. 

\bibliographystyle{siamplain}
\bibliography{Koopman_oblique_projection_arxiv}  

\appendix
%\appendixpage
\begin{comment}
\section{Proof of Proposition \ref{prop. optimality of orthogonal projection in a Hilbert space}}
\label{appendix proof of optimality of orthogonal projection in a Hilbert space}
Suppose $E=W^\perp$ so that $P_{W,W^\perp}$ in \eqref{eq. def of projection operator} is the orthogonal projection.
Let $g,h\in W$ be arbitrary.
\begin{align}
    \| \mathcal{L}_{|W}g - h \|^2
    =&
    \| \mathcal{L}_{|W}g - Lg + Lg - h \|^2
    &\nonumber 
\\
    =&
    \| e + Lg - h \|^2
    \ \ 
    (\because \eqref{eq. direct decomposition each element} )
    &\nonumber
\\
    =&
    \| e \|^2 + \| Lg - h \|^2
    \ \ 
    \left(
        \because 
        \langle e,Lg-h \rangle=0 \text{ since } e\in W^\perp, Lg-h\in W
    \right)
    &\nonumber 
\\
    \geq&
    \| e \|^2.
    &\nonumber
\end{align}

The equality holds when $h=Lg$ since 
%\begin{align}
$    
    \| \mathcal{L}_{|W}g-h \|^2
    =
    \| \mathcal{L}_{|W}g-Lg \|^2
    =
    \| e \|^2,
$
%\end{align}
which implies
\begin{align}
    Lg = \underset{h\in W}{\text{argmin}}\ 
	\| \mathcal{L}_{|W}g - h  \|.
\end{align}
\end{comment}

\section{Proof of Proposition \ref{prop. optimality of orthogonal projection in the finite data setting}}
\label{appendix proof of prop on optimality of orthogonal projection}
If we define
\vspace{-2mm}
\begin{displaymath}
	\bm{X}:=
	\left[\hspace{-1mm}
		\begin{array}{ccc}
			\bm{\psi} (x_1)&\cdots&\bm{\psi} (x_M)
		\end{array}
	\hspace{-1mm}\right]\in \mathbb{R}^{N\times M},
	\ \ 
	\bm{Y}:=
	\left[\hspace{-1mm}
	\begin{array}{ccc}
		(\mathcal{L}_{|W}\bm{\psi}) (x_1)&\cdots&(\mathcal{L}_{|W}\bm{\psi}) (x_M)
	\end{array}
	\hspace{-1mm}\right]\in \mathbb{R}^{N\times M},
\end{displaymath}
the least-square problem \eqref{eq. least square problem} can be written as
\vspace{-2mm}
\begin{align}
	\underset{\hat{\bm{L}}\in \mathbb{R}^{N\times N}}{\text{min}}
	\sum_{l=1}^{M}
	\| (\mathcal{L}_{|W}\bm{\psi})(x_l) - \hat{\bm{L}}\bm{\psi} (x_l)\|_2^2
	&=
	\underset{\hat{\bm{L}}\in \mathbb{R}^{N\times N}}{\text{min}}
	\|
	\bm{Y} - \hat{\bm{L}}\bm{X}
	\|_F^2
	=
	\underset{\hat{\bm{L}}\in \mathbb{R}^{N\times N}}{\text{min}}
	\|
	\bm{Y}^\tr - \bm{X}^\tr\hat{\bm{L}}^\tr
	\|_F^2,
	&\nonumber
\end{align}
which corresponds to solving an over-determined system of equations $\bm{X}^\tr\hat{\bm{L}}^\tr=\bm{Y}^\tr$ with the assumption $N<M$.
The unique solution that achieves the least squared error can be obtained by solving the normal equation:
\vspace{-2mm}
\begin{align}
	&
	\bm{X}\bm{X}^\tr\hat{\bm{L}}^\tr=\bm{X}\bm{Y}^\tr
	%\ \
        &\nonumber
\\
	\Leftrightarrow
	\ \ 
        &
	\hat{\bm{L}}\bm{X}\bm{X}^\tr = \bm{Y}\bm{X}^\tr
	&\nonumber
\\
	\Leftrightarrow
	\ \ 
	&
	\hat{\bm{L}}
	\left\{
		\frac{1}{M}
		\sum_{l=1}^{M}
		\bm{\psi}(x_l)\bm{\psi}(x_l)^\tr 
	\right\}
	=
	\frac{1}{M} \sum_{i=1}^{M}
	(\mathcal{L}_{|W}\bm{\psi})(x_l) \bm{\psi}(x_l)^\tr,
	&
\end{align}
which gives \eqref{eq. def of L_m} with $\psi_i=\varphi_i$ on the assumption that $\bm{X}\bm{X}^\tr$ is invertible, or equivalently $\bm{X}$ has full rank.

\end{document}

\maketitle

\section{Sensitivity Analysis of the Modeling Methods}
\label{appendix sensitivity analysis}
In a data-driven modeling procedure, one may need to repeat the model training with or without new data if the learned model does not show satisfactory performance for the given task.
For neural network-based models, the training process is often terminated at an undesirable local minimum due to the nature of high-dimensional non-convex optimization.
Thus, it may be required to repeat the training process many times to obtain a satisfactory result and the computational cost of model training may be an issue in practice.
In this section, the sensitivities of the learning methods are evaluated.

Specifically, in the same settings as in Section 5, we repeated collecting data and training models 10 times as follows: 
\vspace{-2mm}
\begin{algorithm}
	\begin{algorithmic}
		\For {$k=1,\cdots,10$}
		\State Collect data to form a data set $\mathcal{D}_k:=\{ (y_i,\chi_i,u_i)\mid y_i=F(\chi_i,u_i) \}$.
		\State Train an EDMD, a normal NN, and a proposed models using the same data set $\mathcal{D}_k$.
		\EndFor 
	\end{algorithmic}
\end{algorithm}
\vspace{-3mm}

\begin{figure}[H]
	\centering
	\begin{subfigure}{0.3\linewidth}
		\centering
		\includegraphics[width=0.95\linewidth]{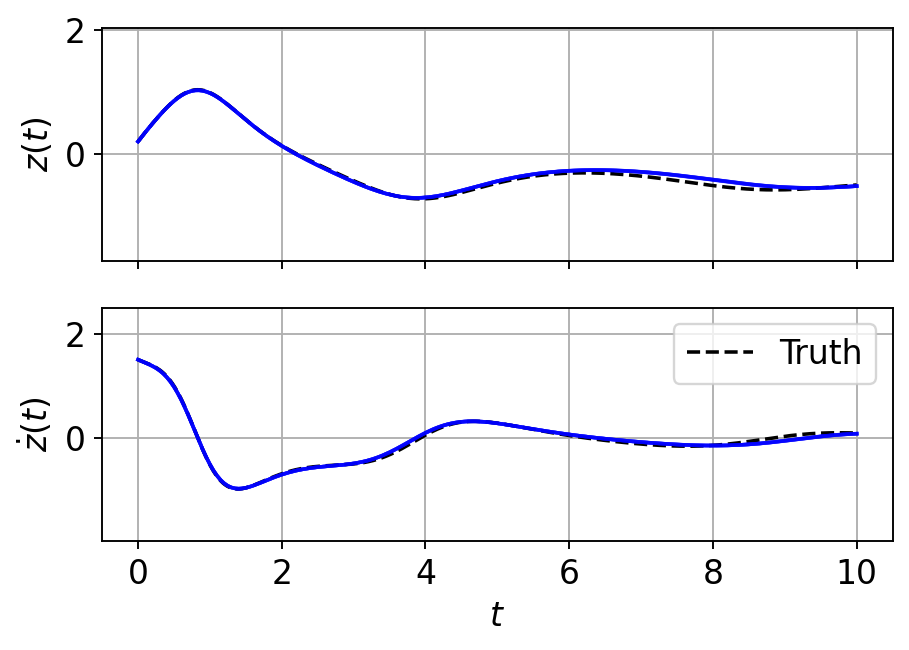}
		\caption{State prediction (EDMD).}
		\label{subfig. duffing single_plot state prediction EDMD}
	\end{subfigure}
	\begin{subfigure}{0.3\linewidth}
		\centering
		\includegraphics[width=0.95\linewidth]{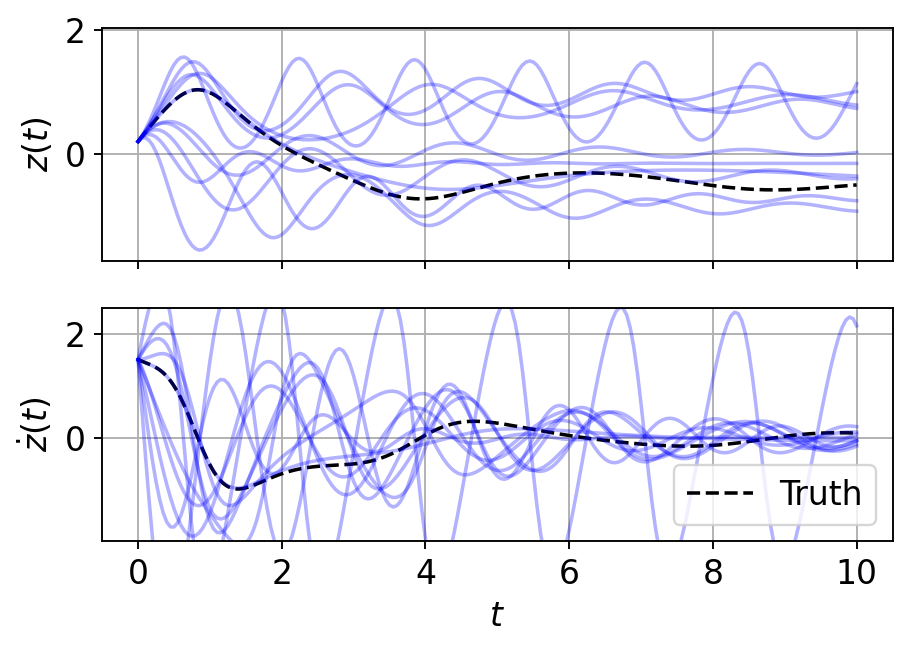}
		\caption{\scriptsize{State prediction (normal NN).}}
		\label{subfig. duffing single_plot state prediction normal}
	\end{subfigure}
	\begin{subfigure}{0.3\linewidth}
		\centering
		\includegraphics[width=0.95\linewidth]{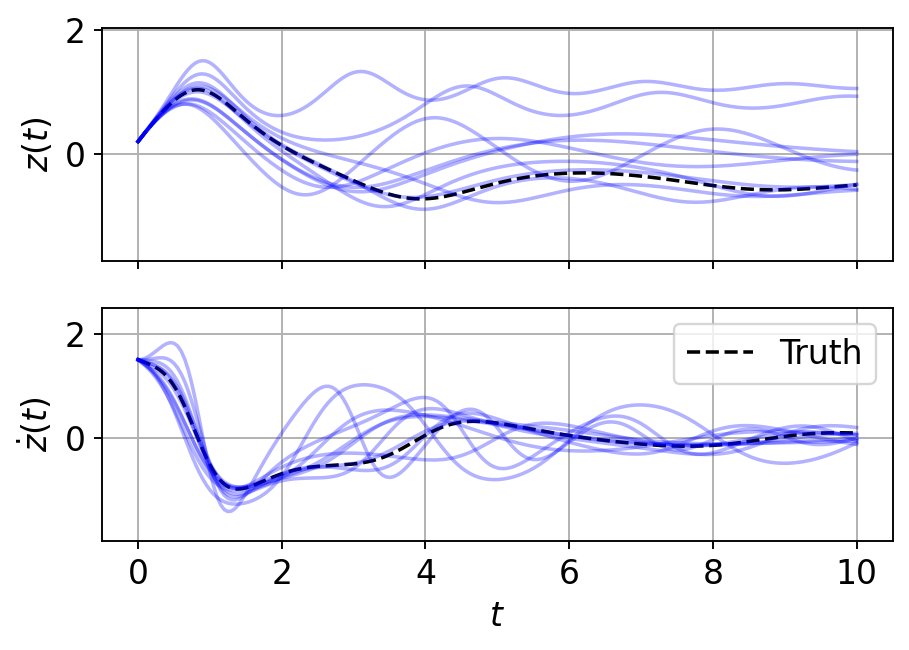}
		\caption{\footnotesize{State prediction (proposed).}}
		\label{subfig. duffing single_plot state prediction oblique}
	\end{subfigure}
	\begin{subfigure}{0.3\linewidth}
		\centering
		\includegraphics[width=0.95\linewidth]{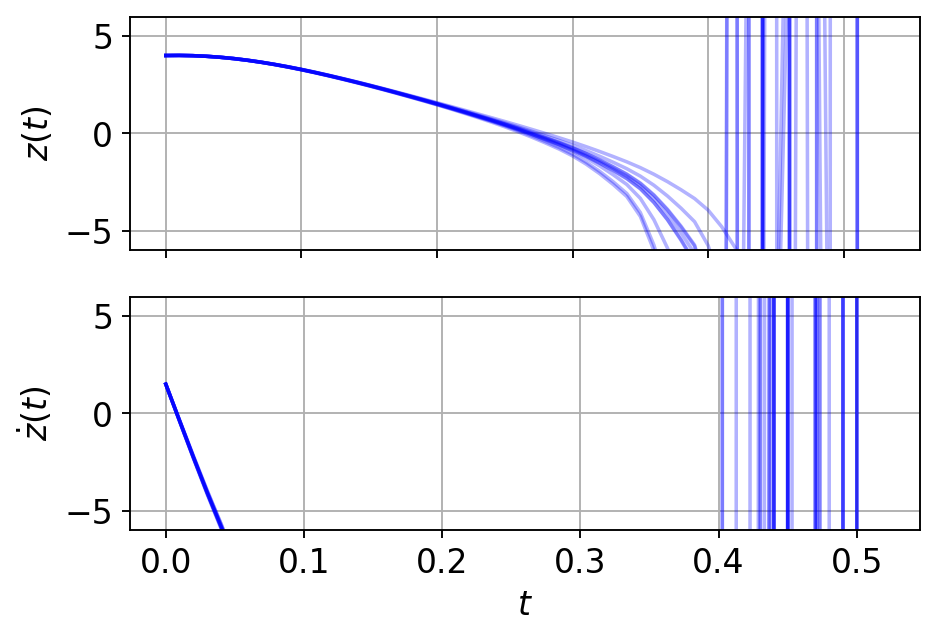}
		\caption{Stabilization by LQR (EDMD).}
		\label{subfig. duffing single_plot cl EDMD}
	\end{subfigure}
	\begin{subfigure}{0.3\linewidth}
		\centering
		\includegraphics[width=0.95\linewidth]{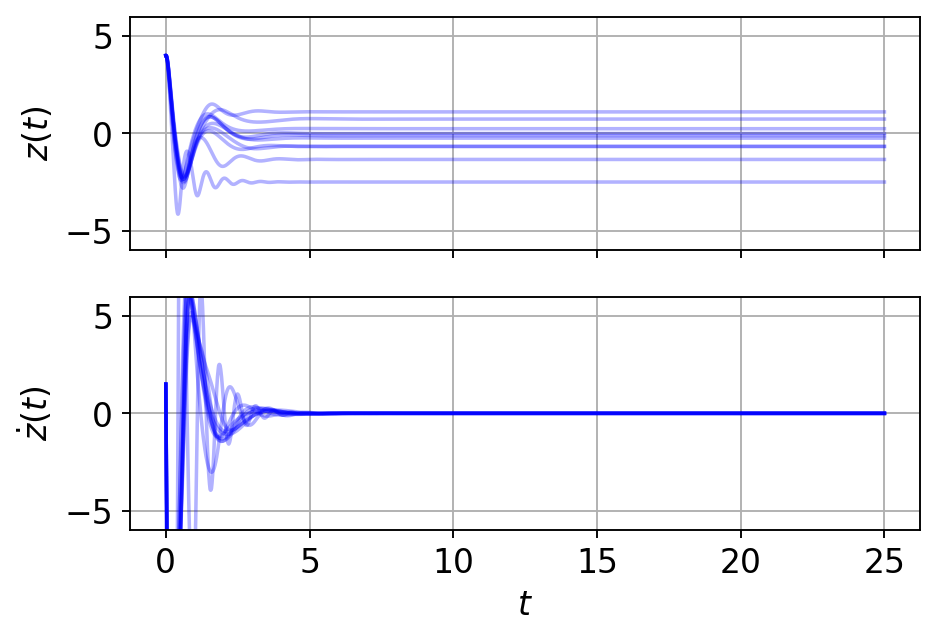}
		\caption{\footnotesize{Stabilization by LQR (normal NN).}}
		\label{subfig. duffing single_plot cl normal}
	\end{subfigure}
	\begin{subfigure}{0.3\linewidth}
		\centering
		\includegraphics[width=0.95\linewidth]{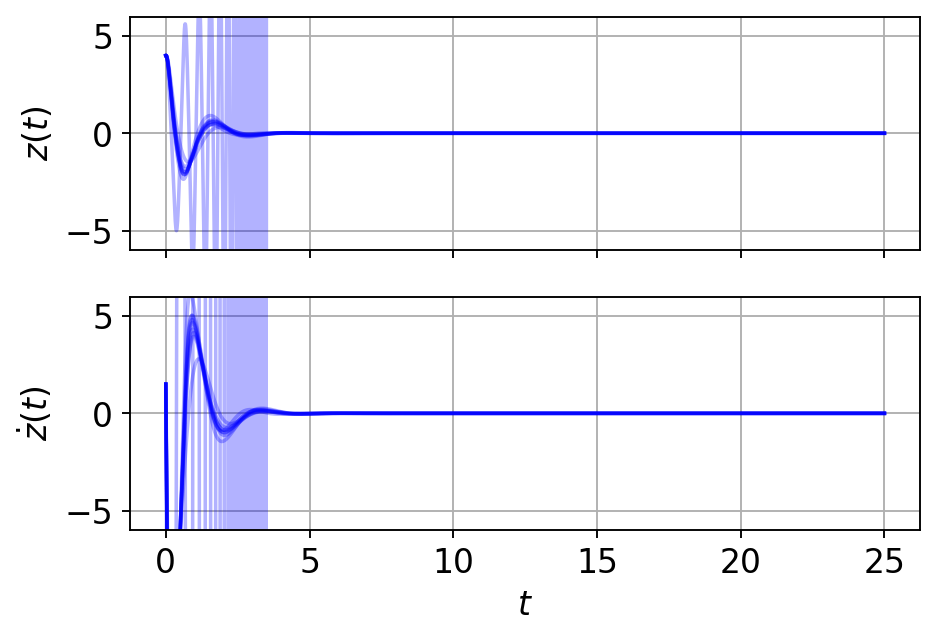}
		\caption{Stabilization by LQR (proposed).}
		\label{subfig. duffing single_plot cl oblique}
	\end{subfigure}
	\begin{subfigure}{0.3\linewidth}
		\centering
		\includegraphics[width=0.95\linewidth]{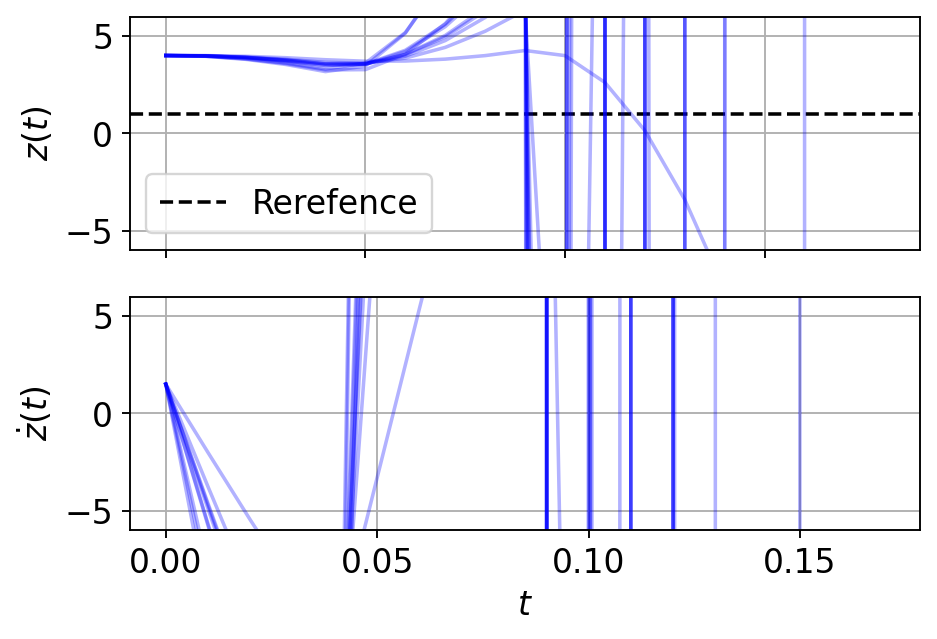}
		\caption{Reference tracking (EDMD).}
		\label{subfig. duffing single_plot cl servo EDMD}
	\end{subfigure}
	\begin{subfigure}{0.3\linewidth}
		\centering
		\includegraphics[width=0.95\linewidth]{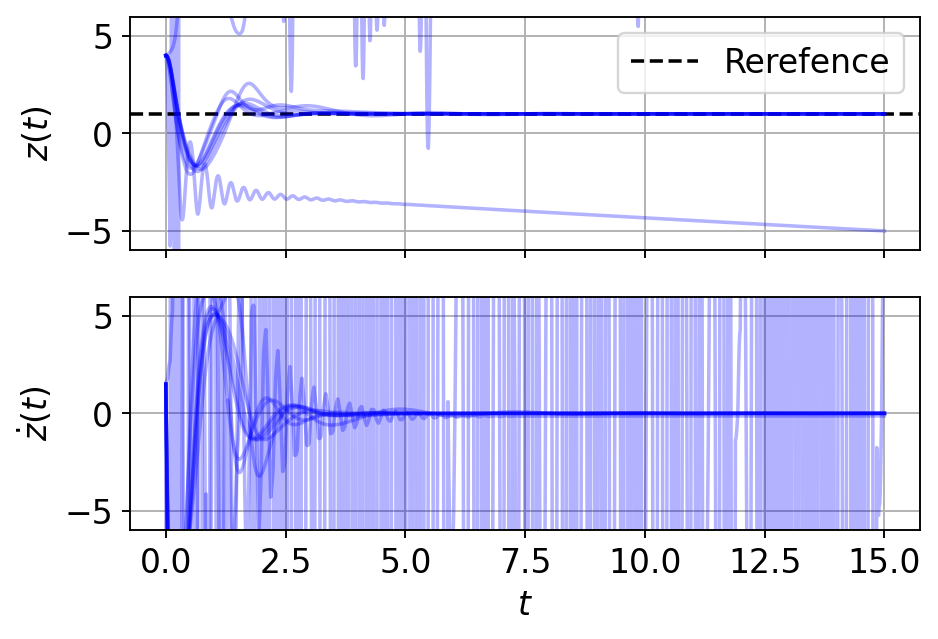}
		\caption{Reference tracking (normal NN).}
		\label{subfig. duffing single_plot cl servo normal}
	\end{subfigure}
	\begin{subfigure}{0.3\linewidth}
		\centering
		\includegraphics[width=0.95\linewidth]{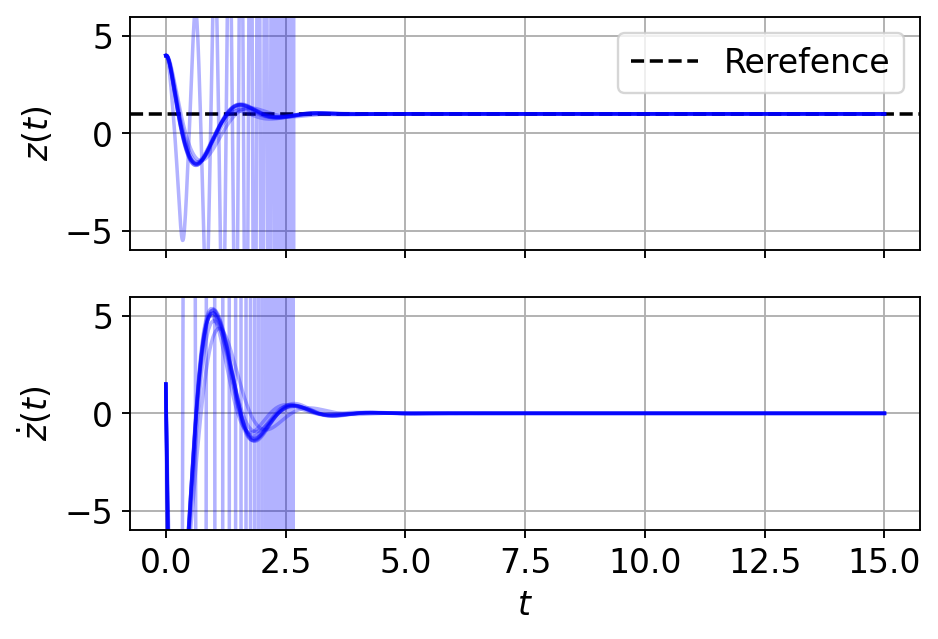}
		\caption{Reference tracking (proposed).}
		\label{subfig. duffing single_plot cl servo oblique}
	\end{subfigure}
	\begin{subfigure}{0.3\linewidth}
		\centering
		\includegraphics[width=0.95\linewidth]{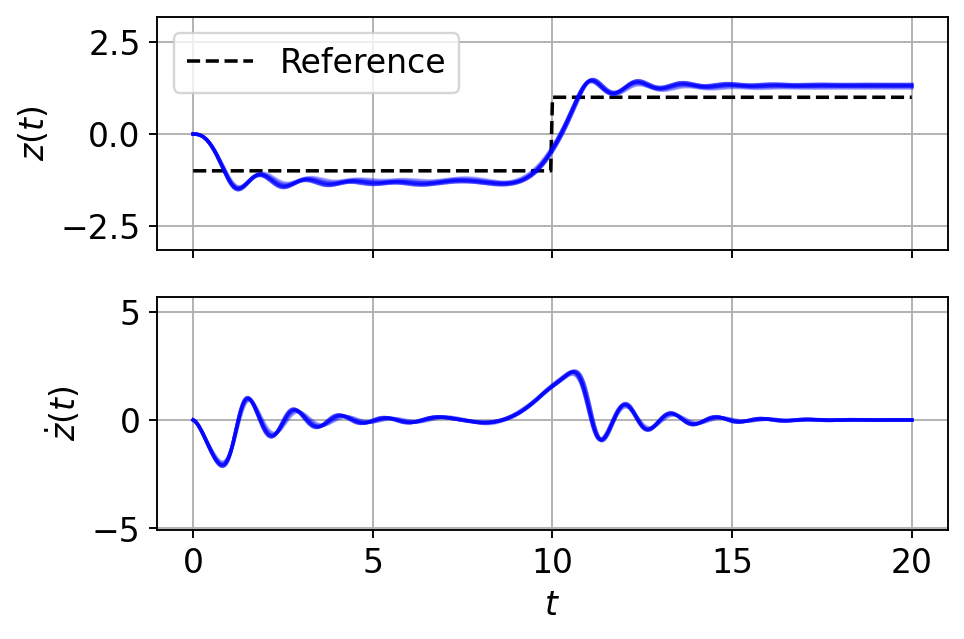}
		\caption{MPC (EDMD).}
		\label{subfig. duffing single_plot mpc EDMD}
	\end{subfigure}
	\begin{subfigure}{0.3\linewidth}
		\centering
		\includegraphics[width=0.95\linewidth]{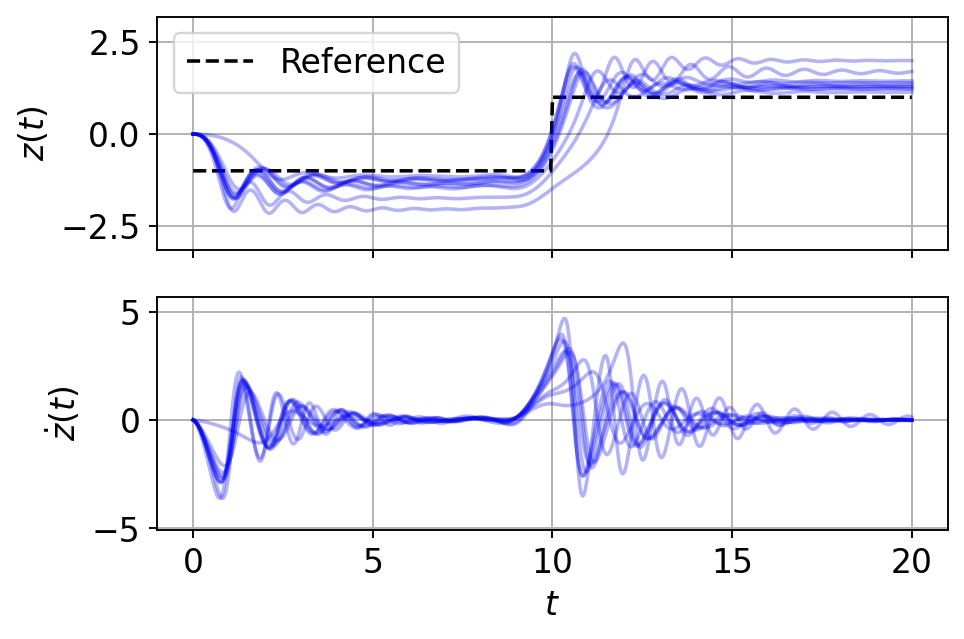}
		\caption{MPC (normal NN).}
		\label{subfig. duffing single_plot mpc normal}
	\end{subfigure}
	\begin{subfigure}{0.3\linewidth}
		\centering
		\includegraphics[width=0.95\linewidth]{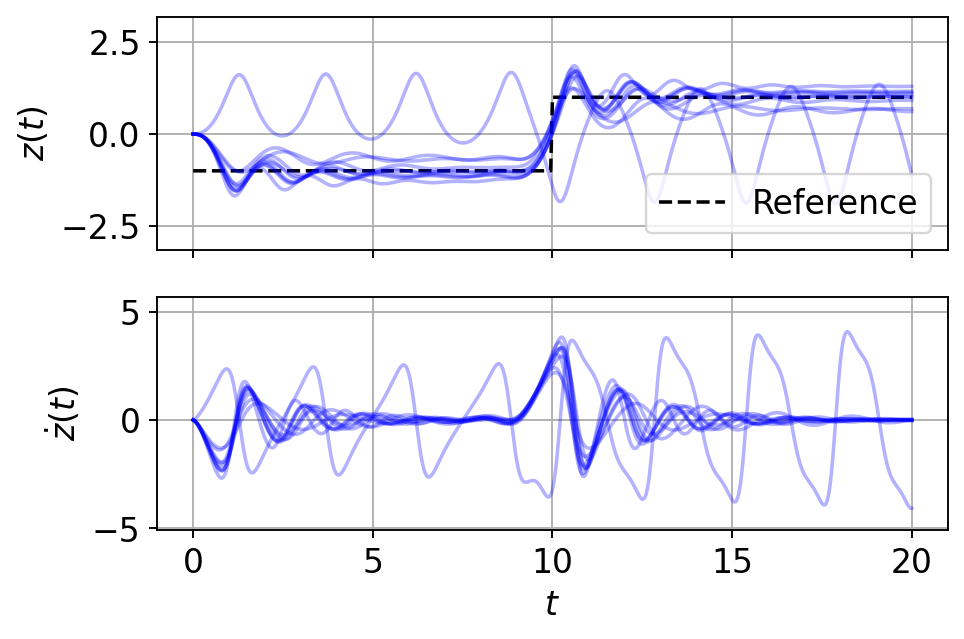}
		\caption{MPC (proposed).}
		\label{subfig. duffing single_plot mpc oblique}
	\end{subfigure}
	\caption{Sensitivity analysis (Duffing oscillator).}
	\label{fig. duffing uncertainty quantification}
\end{figure}

\begin{figure}[H]
	\centering
	\begin{subfigure}{0.3\linewidth}
		\centering
		\includegraphics[width=0.95\linewidth]{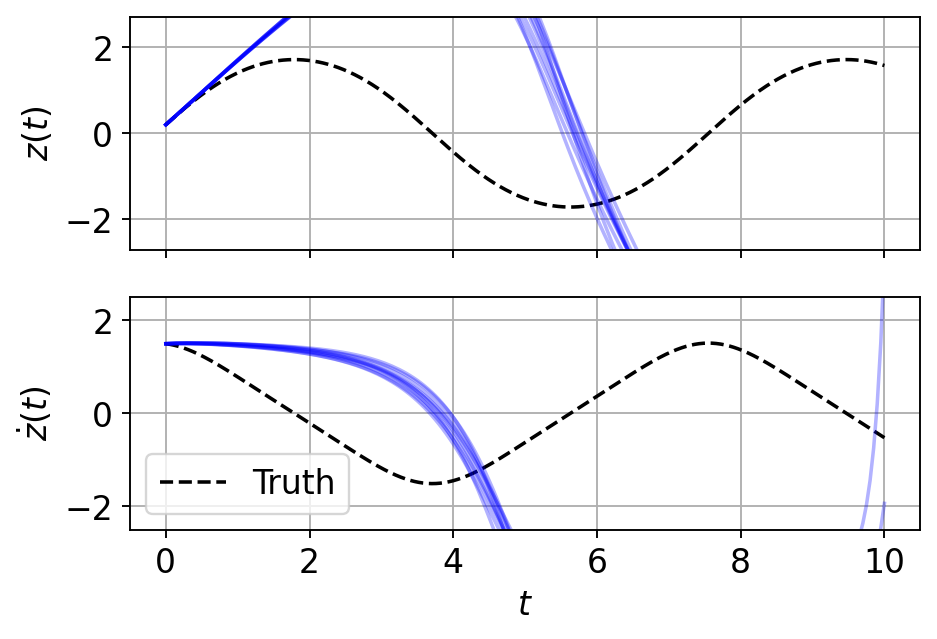}
		\caption{State prediction (EDMD).}
		\label{subfig. pendulum single_plot state prediction EDMD}
	\end{subfigure}
	\begin{subfigure}{0.3\linewidth}
		\centering
		\includegraphics[width=0.95\linewidth]{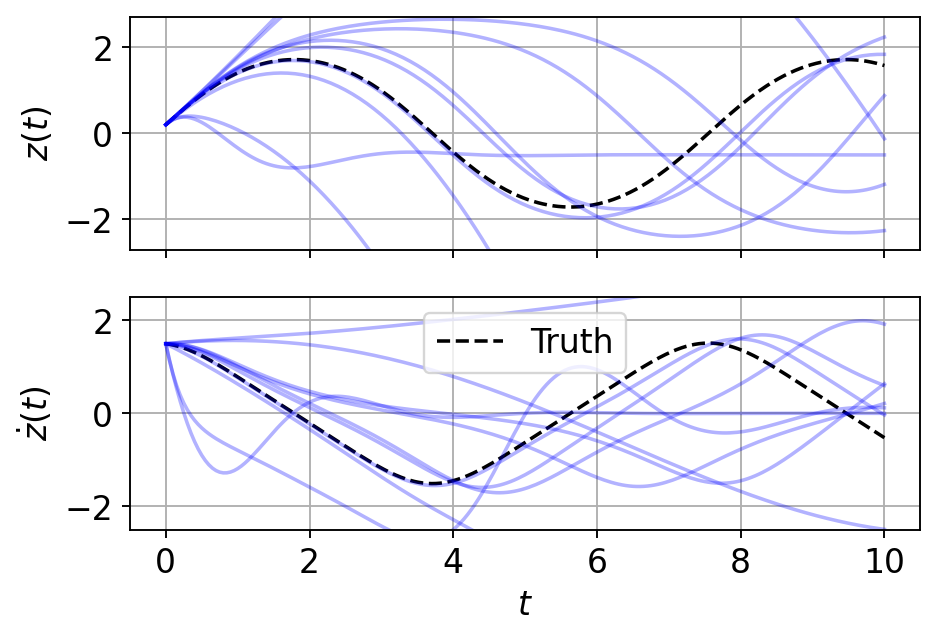}
		\caption{\scriptsize{State prediction (normal NN).}}
		\label{subfig. pendulum single_plot state prediction normal}
	\end{subfigure}
	\begin{subfigure}{0.3\linewidth}
		\centering
		\includegraphics[width=0.95\linewidth]{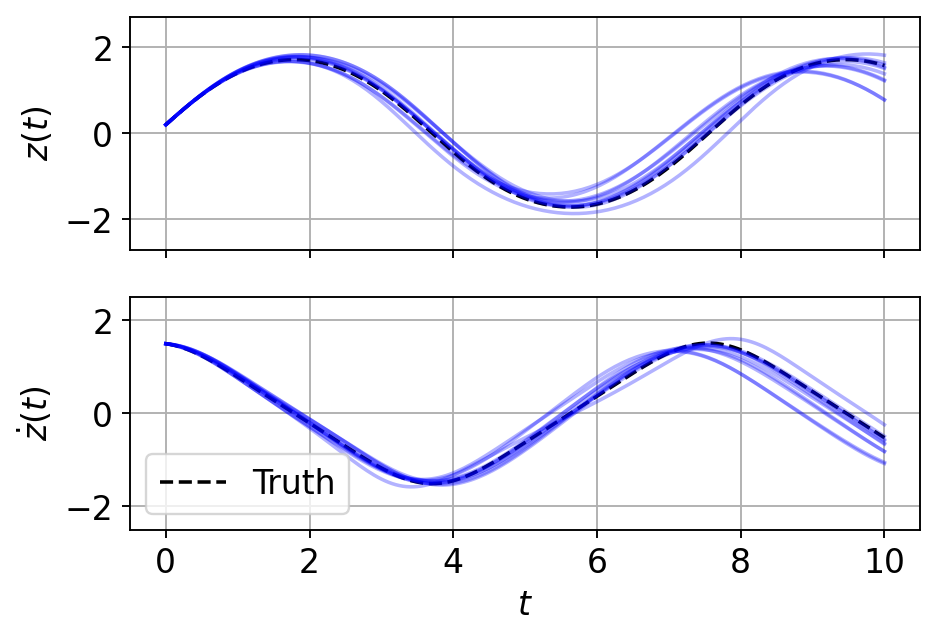}
		\caption{\footnotesize{State prediction (proposed).}}
		\label{subfig. pendulum single_plot state prediction oblique}
	\end{subfigure}
	\begin{subfigure}{0.3\linewidth}
		\centering
		\includegraphics[width=0.95\linewidth]{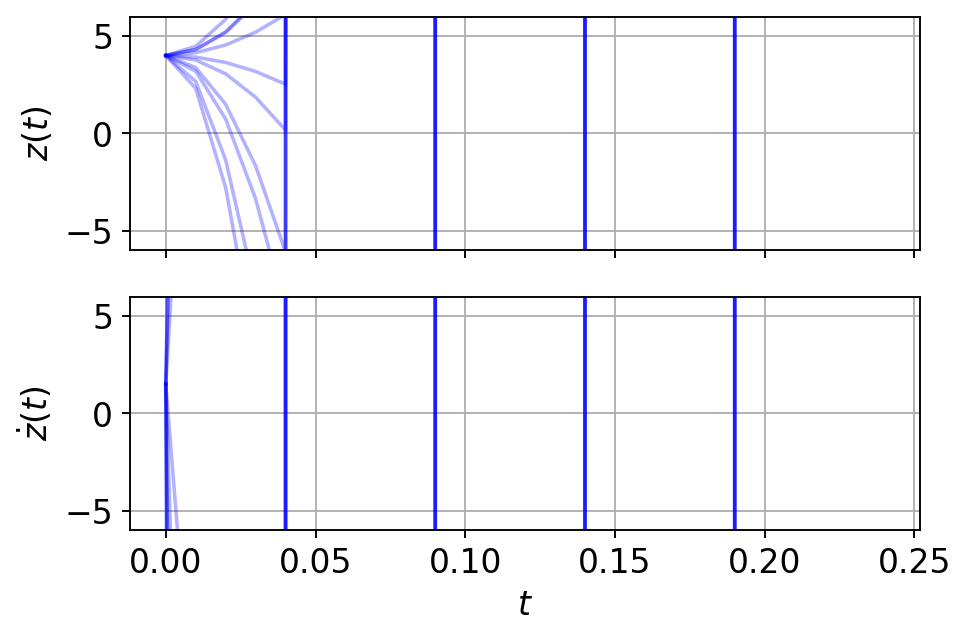}
		\caption{Stabilization by LQR (EDMD).}
		\label{subfig. pendulum single_plot cl EDMD}
	\end{subfigure}
	\begin{subfigure}{0.3\linewidth}
		\centering
		\includegraphics[width=0.95\linewidth]{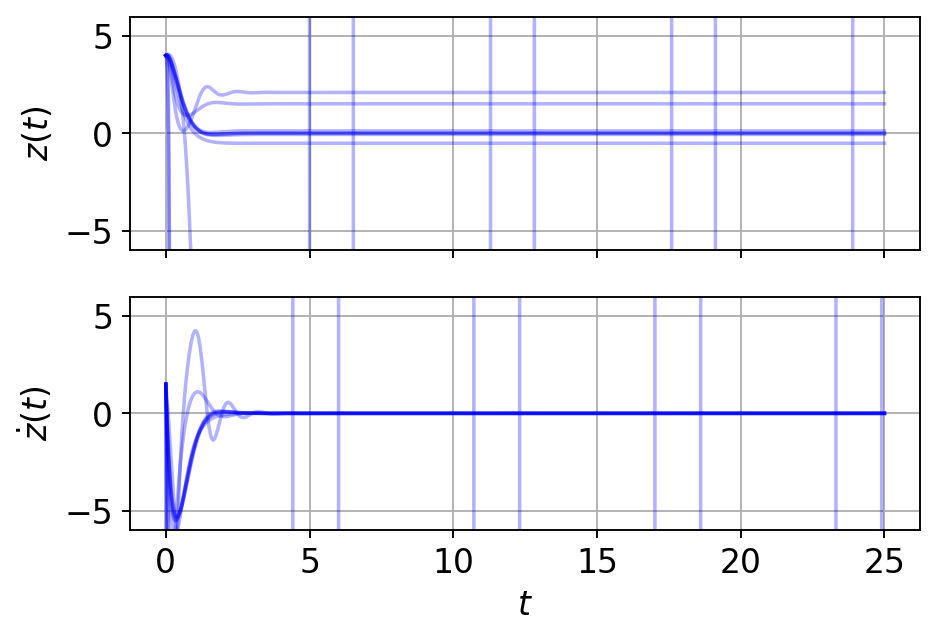}
		\caption{\footnotesize{Stabilization by LQR (normal NN).}}
		\label{subfig. pendulum single_plot cl normal}
	\end{subfigure}
	\begin{subfigure}{0.3\linewidth}
		\centering
		\includegraphics[width=0.95\linewidth]{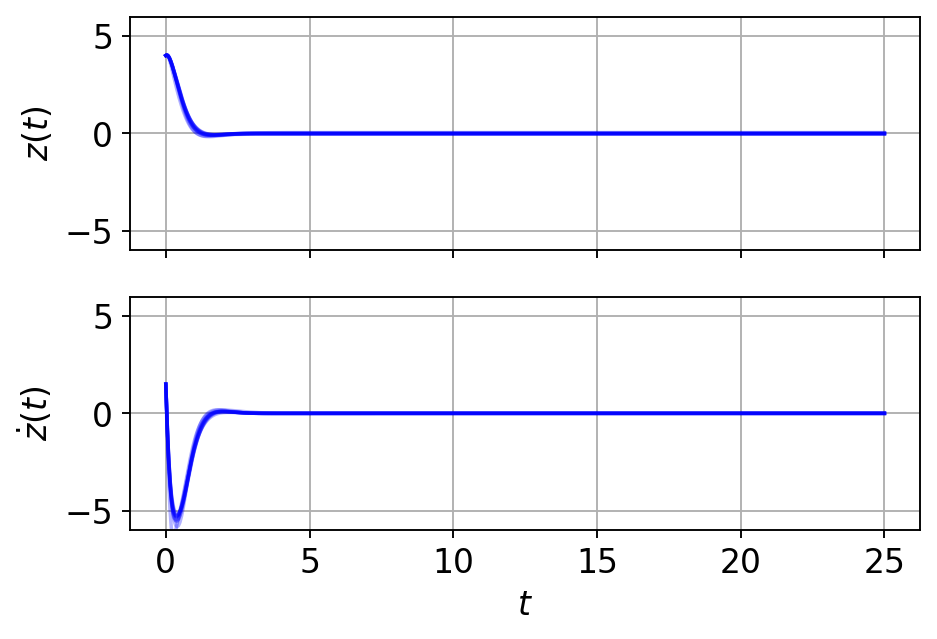}
		\caption{Stabilization by LQR (proposed).}
		\label{subfig. pendulum single_plot cl oblique}
	\end{subfigure}
	\begin{subfigure}{0.3\linewidth}
		\centering
		\includegraphics[width=0.95\linewidth]{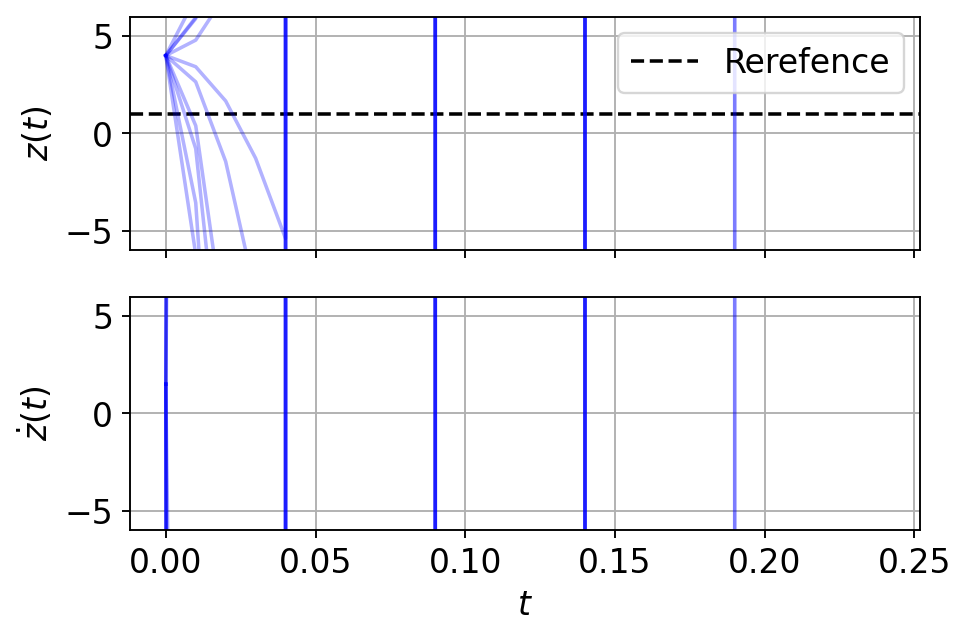}
		\caption{Reference tracking (EDMD).}
		\label{subfig. pendulum single_plot cl servo EDMD}
	\end{subfigure}
	\begin{subfigure}{0.3\linewidth}
		\centering
		\includegraphics[width=0.95\linewidth]{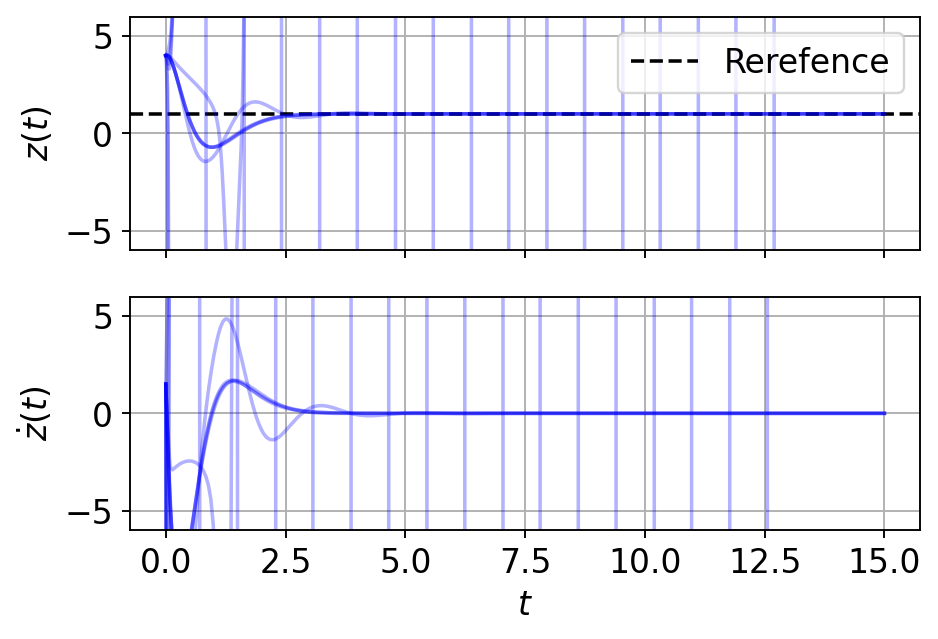}
		\caption{Reference tracking (normal NN).}
		\label{subfig. pendulum single_plot cl servo normal}
	\end{subfigure}
	\begin{subfigure}{0.3\linewidth}
		\centering
		\includegraphics[width=0.95\linewidth]{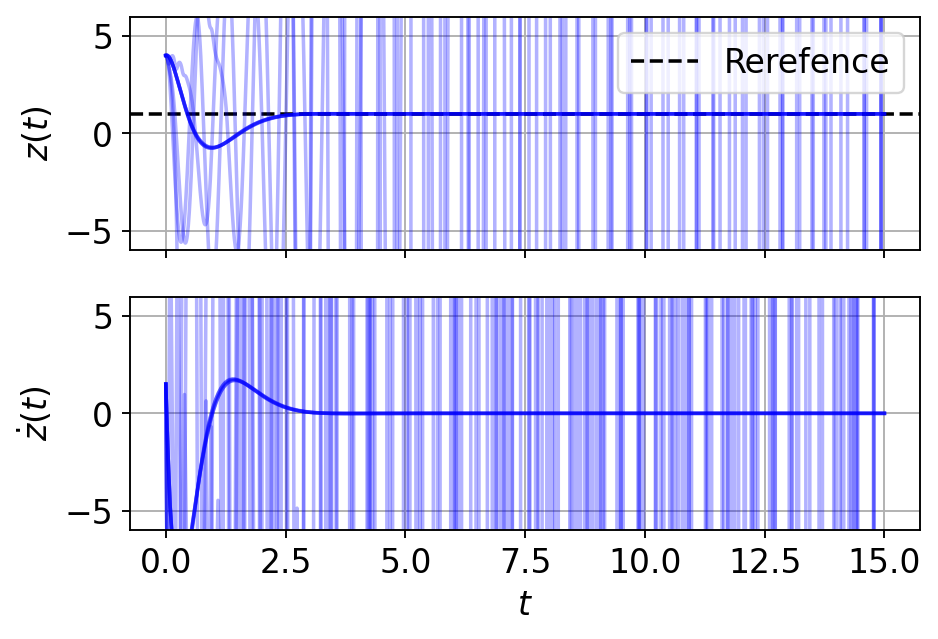}
		\caption{Reference tracking (proposed).}
		\label{subfig. pendulum single_plot cl servo oblique}
	\end{subfigure}
	\begin{subfigure}{0.3\linewidth}
		\centering
		\includegraphics[width=0.95\linewidth]{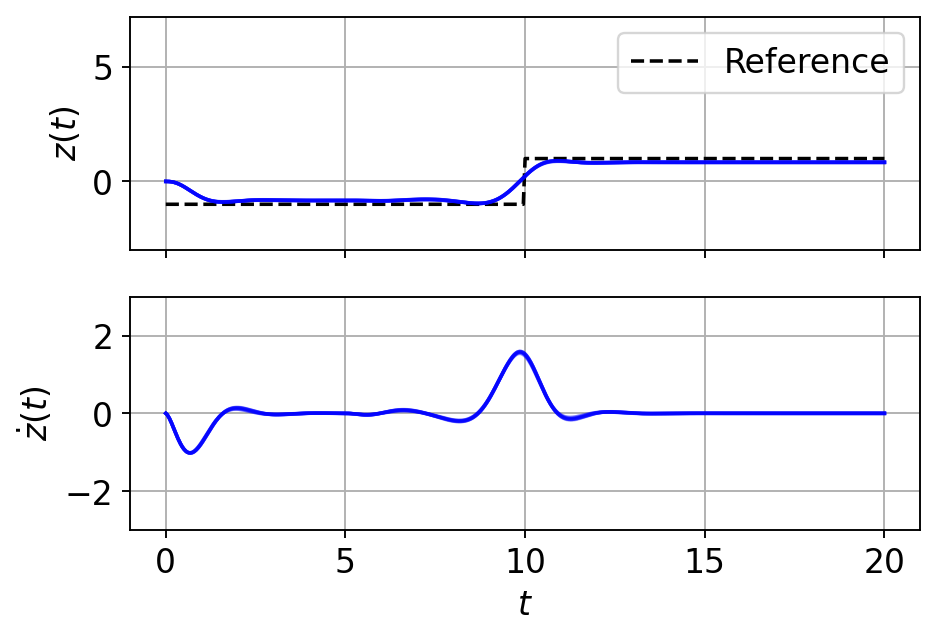}
		\caption{MPC (EDMD).}
		\label{subfig. pendulum single_plot mpc EDMD}
	\end{subfigure}
	\begin{subfigure}{0.3\linewidth}
		\centering
		\includegraphics[width=0.95\linewidth]{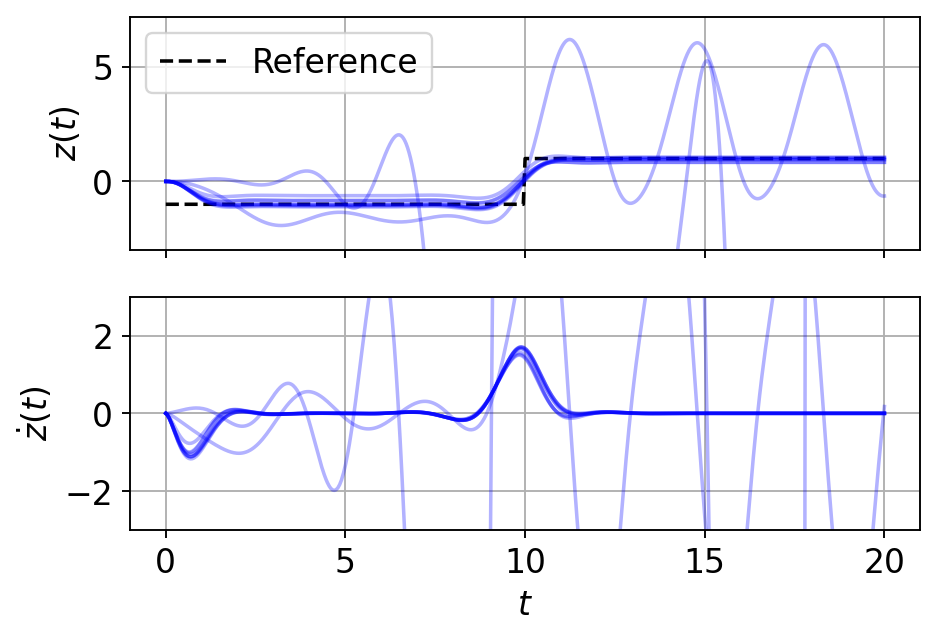}
		\caption{MPC (normal NN).}
		\label{subfig. pendulum single_plot mpc normal}
	\end{subfigure}
	\begin{subfigure}{0.3\linewidth}
		\centering
		\includegraphics[width=0.95\linewidth]{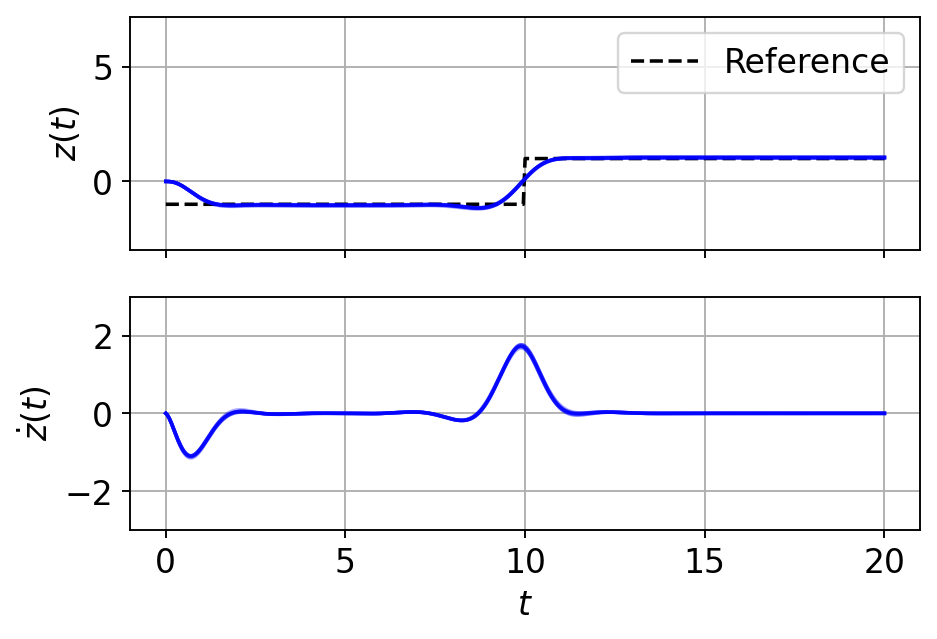}
		\caption{MPC (proposed).}
		\label{subfig. pendulum single_plot mpc oblique}
	\end{subfigure}
	\caption{Sensitivity analysis (simple pendulum).}
	\label{fig. pendulum uncertainty quantification}
\end{figure}

\begin{figure}[H]
	\centering
	\begin{subfigure}{0.3\linewidth}
		\centering
		\includegraphics[width=0.95\linewidth]{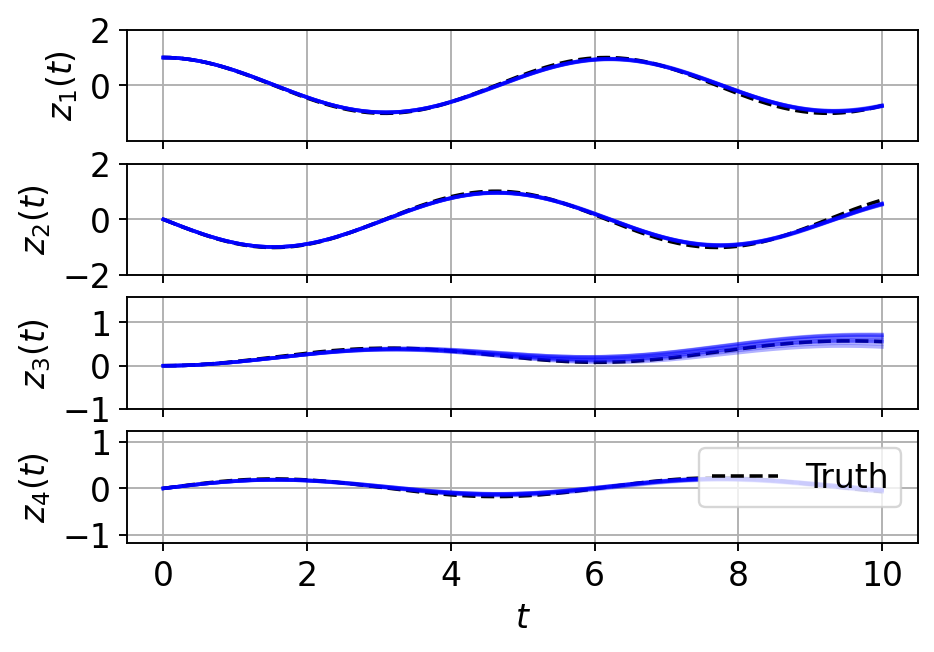}
		\caption{State prediction (EDMD).}
		\label{subfig. peRTAC single_plot state prediction EDMD}
	\end{subfigure}
	\begin{subfigure}{0.3\linewidth}
		\centering
		\includegraphics[width=0.95\linewidth]{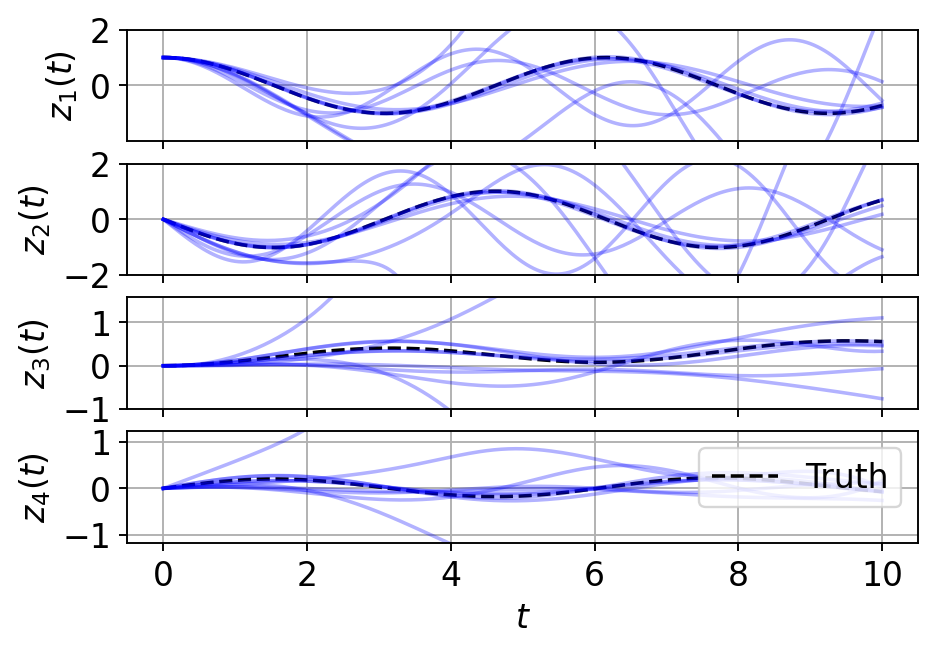}
		\caption{\scriptsize{State prediction (normal NN).}}
		\label{subfig. RTAC single_plot state prediction normal}
	\end{subfigure}
	\begin{subfigure}{0.3\linewidth}
		\centering
		\includegraphics[width=0.95\linewidth]{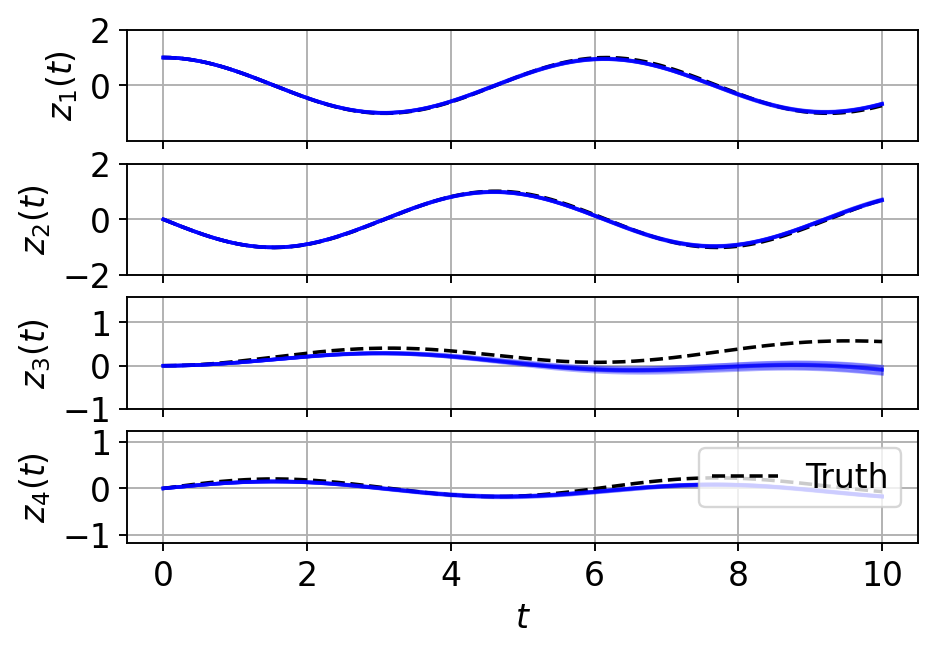}
		\caption{\footnotesize{State prediction (proposed).}}
		\label{subfig. RTAC single_plot state prediction oblique}
	\end{subfigure}
	\begin{subfigure}{0.3\linewidth}
		\centering
		\includegraphics[width=0.95\linewidth]{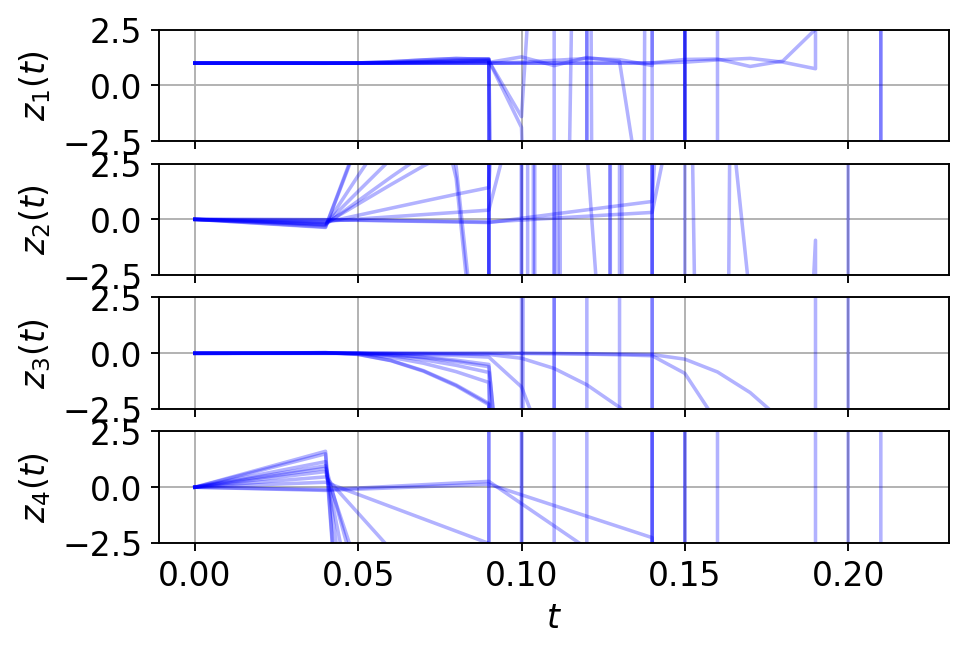}
		\caption{Stabilization by LQR (EDMD).}
		\label{subfig. RTAC single_plot cl EDMD}
	\end{subfigure}
	\begin{subfigure}{0.3\linewidth}
		\centering
		\includegraphics[width=0.95\linewidth]{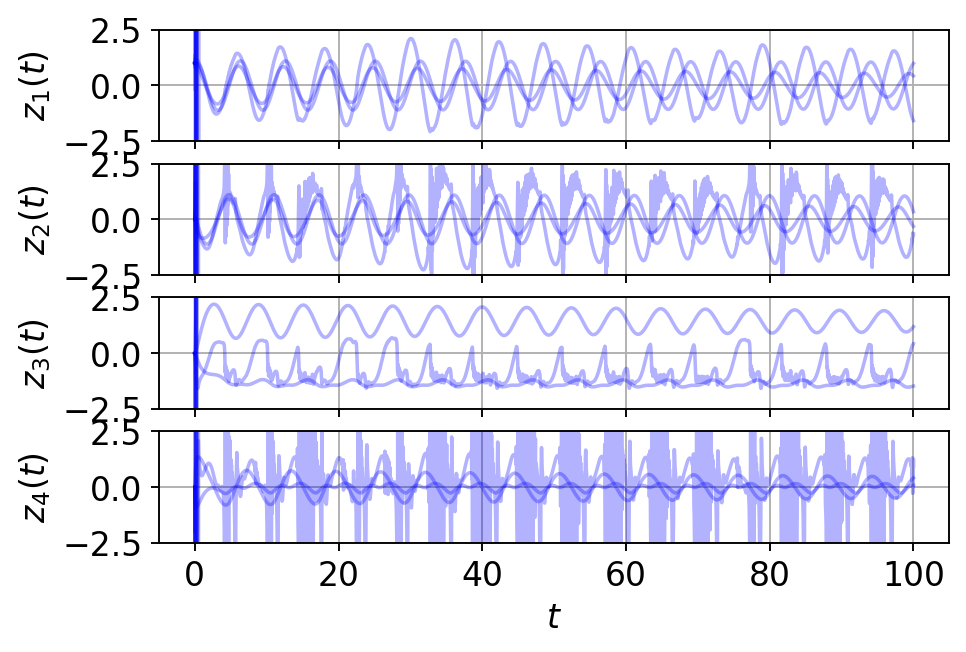}
		\caption{\footnotesize{Stabilization by LQR (normal NN).}}
		\label{subfig. RTAC single_plot cl normal}
	\end{subfigure}
	\begin{subfigure}{0.3\linewidth}
		\centering
		\includegraphics[width=0.95\linewidth]{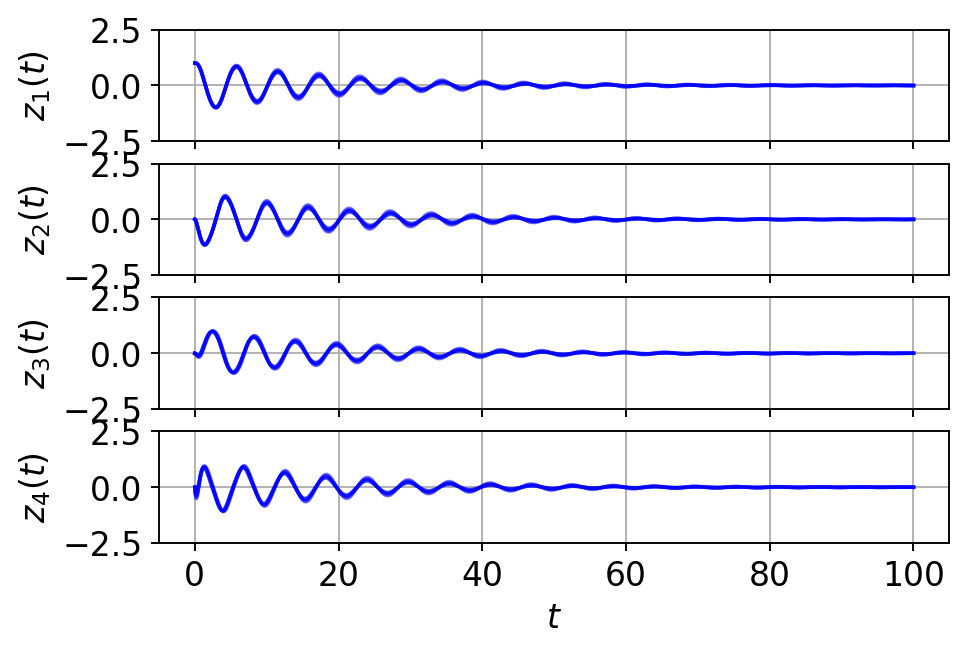}
		\caption{Stabilization by LQR (proposed).}
		\label{subfig. RTAC single_plot cl oblique}
	\end{subfigure}
	\begin{subfigure}{0.3\linewidth}
		\centering
		\includegraphics[width=0.95\linewidth]{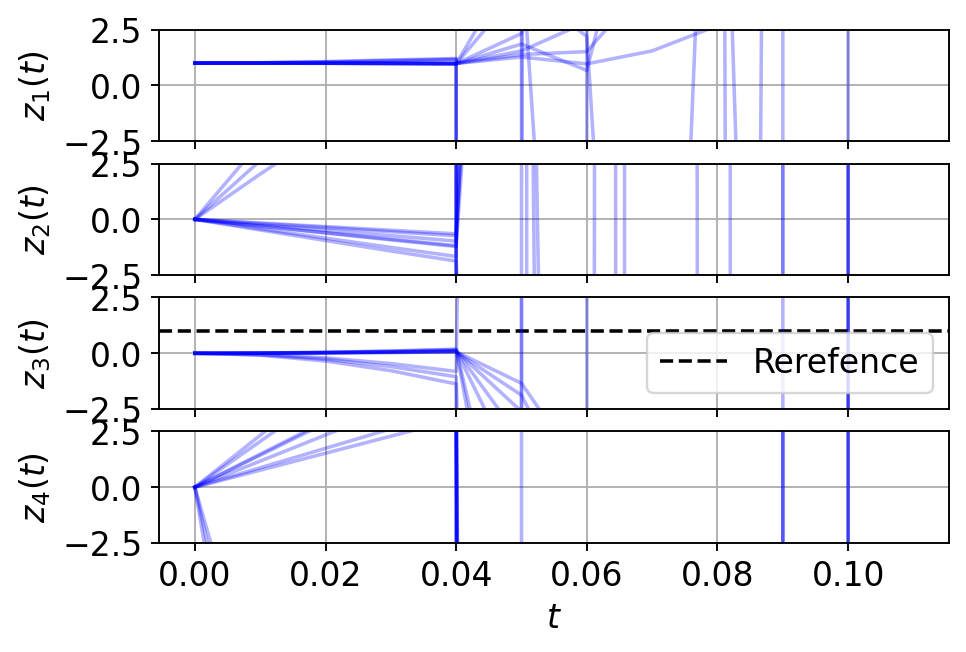}
		\caption{Reference tracking (EDMD).}
		\label{subfig. RTAC single_plot cl servo EDMD}
	\end{subfigure}
	\begin{subfigure}{0.3\linewidth}
		\centering
		\includegraphics[width=0.95\linewidth]{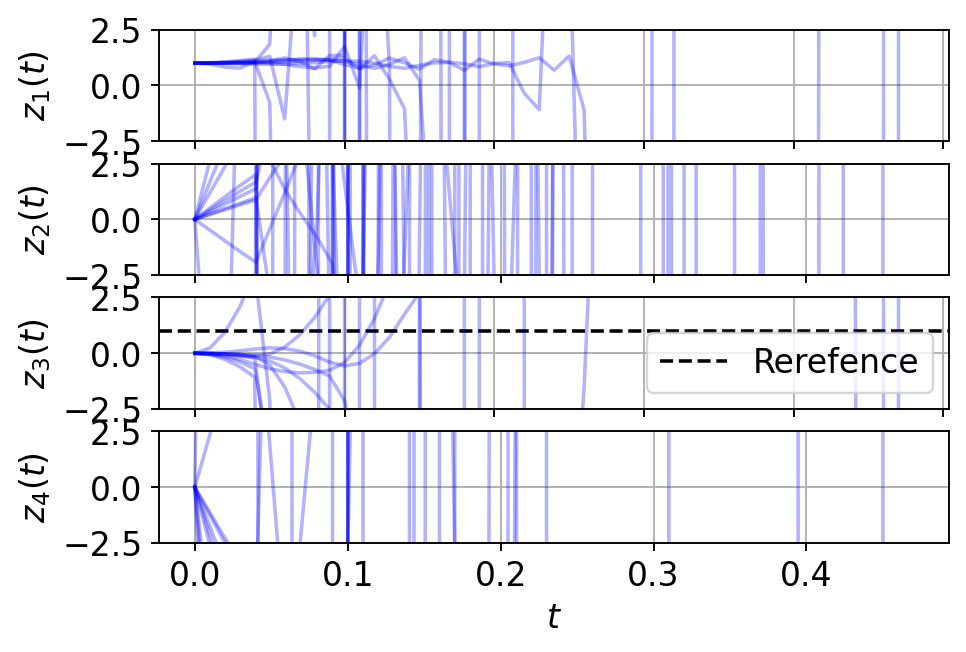}
		\caption{Reference tracking (normal NN).}
		\label{subfig. RTAC single_plot cl servo normal}
	\end{subfigure}
	\begin{subfigure}{0.3\linewidth}
		\centering
		\includegraphics[width=0.95\linewidth]{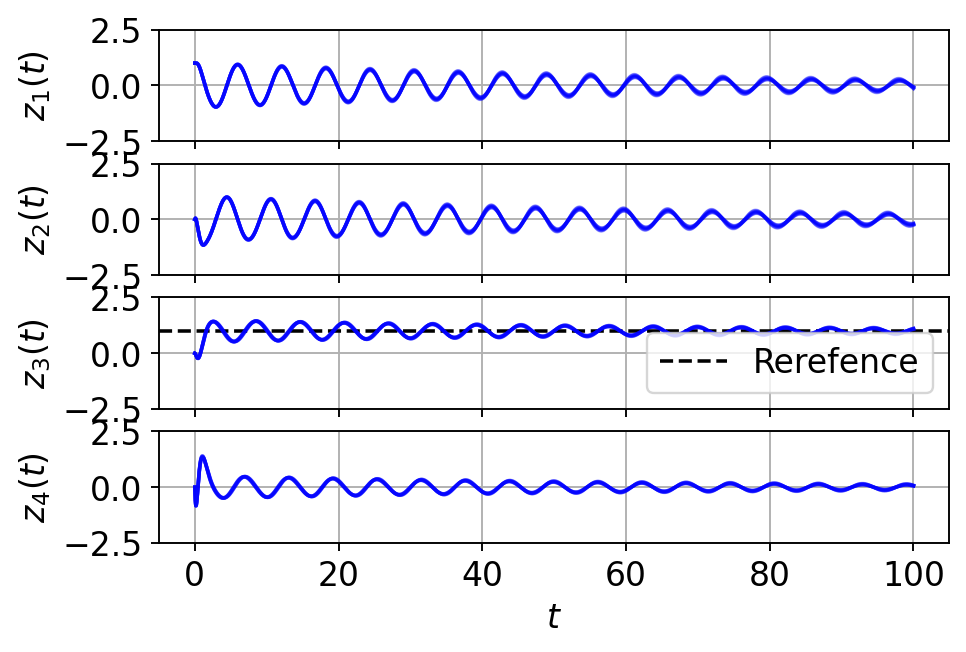}
		\caption{Reference tracking (proposed).}
		\label{subfig. RTAC single_plot cl servo oblique}
	\end{subfigure}
	\begin{subfigure}{0.3\linewidth}
		\centering
		\includegraphics[width=0.95\linewidth]{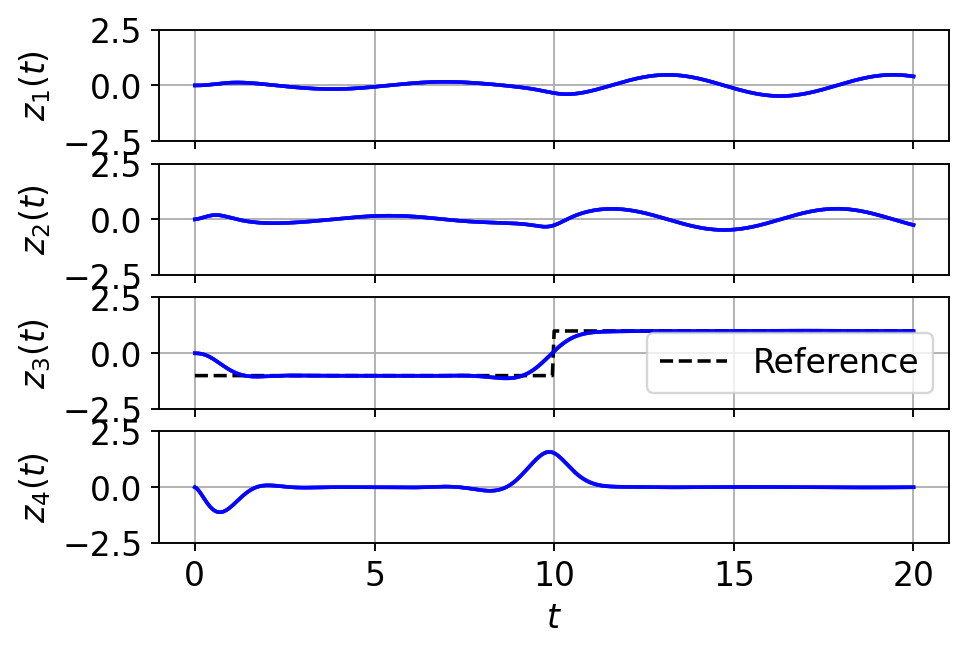}
		\caption{MPC (EDMD).}
		\label{subfig. RTAC single_plot mpc EDMD}
	\end{subfigure}
	\begin{subfigure}{0.3\linewidth}
		\centering
		\includegraphics[width=0.95\linewidth]{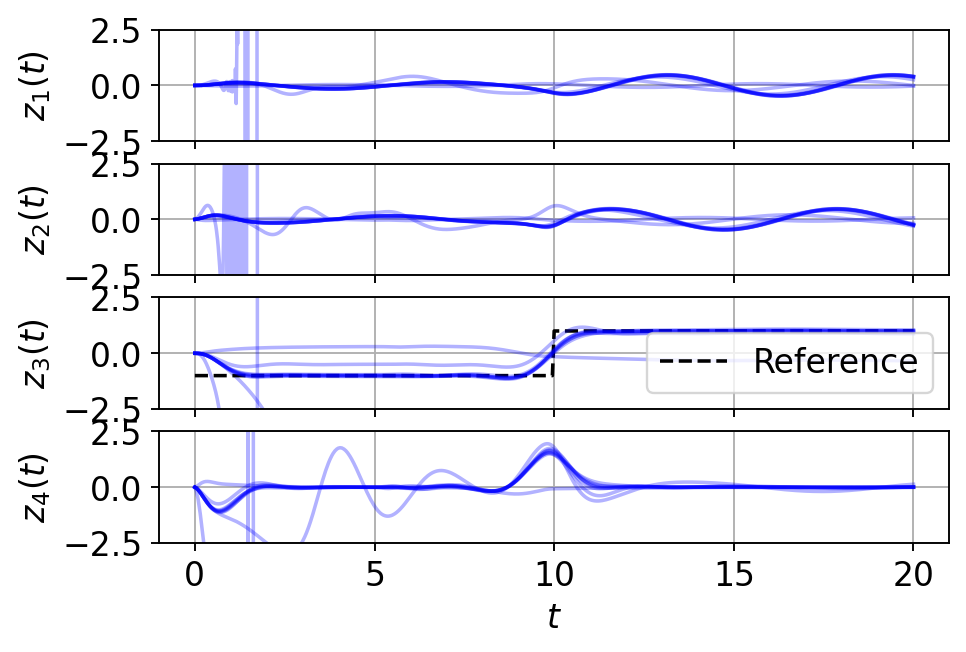}
		\caption{MPC (normal NN).}
		\label{subfig. RTAC single_plot mpc normal}
	\end{subfigure}
	\begin{subfigure}{0.3\linewidth}
		\centering
		\includegraphics[width=0.95\linewidth]{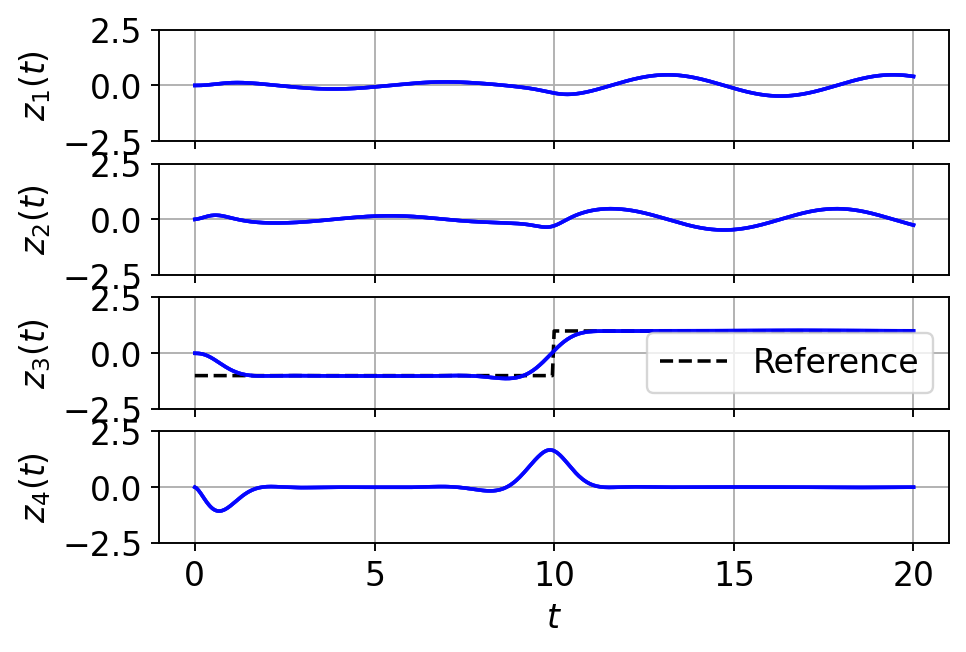}
		\caption{MPC (proposed).}
		\label{subfig. RTAC single_plot mpc oblique}
	\end{subfigure}
	\caption{Sensitivity analysis (RTAC).}
	\label{fig. RTAC uncertainty quantification}
\end{figure}

From the algorithmic perspective, a modeling method is considered robust and reliable if many trials result in high model performance, which implies that one can obtain a model with desirable performance in a small number of trials.
Figures \ref{fig. duffing uncertainty quantification}, \ref{fig. pendulum uncertainty quantification}, and \ref{fig. RTAC uncertainty quantification} show the results, where
the same four tasks as Section 5 are tested on the three dynamical systems and individual results from the 10 trials are overlaid in a single plot for each task.

First, since the EDMD model is obtained by a unique analytical solution to a linear regression problem as in (2.19), its learning results are quite robust, being only influenced by the difference of data points among the 10 trials.
However, the result also implies that it is not able to produce desirable results if the choice of feature maps or the quality of data is not appropriate, regardless of the number of trials.

On the other hand, neural network-based models have more chances to obtain better and preferable learning results with the same number of trials, at the expense of more sensitive learning procedures.
The results of the normal NN and the proposed models show that some trials can yield successful models even if other learning results completely fail in the given task.
Especially, neural network training initializes its parameter values randomly and this often leads to quite different learning results.

The superiority of the proposed method is also observed in this evaluation.
Specifically, in the RTAC case (Fig. \ref{fig. RTAC uncertainty quantification}), 
only the proposed method is both robust and reliable.
The robustness of the proposed method in this example may be attributed to the specific constraint (5.13) imposed on the test functions, which reduces the number of neural network parameters.

\newpage
\section{Invariance Proximity of the Models}
\label{appendix forecasting errors}
	While the modeling error (2.20) was introduced in Section 2 for analyzing the model accuracy and it represents the invariance proximity of the linear embedding model, i.e., how close the subspace spanned by the feature maps is to being invariant under the action of the Koopman operator, it is not a suitable metric in practice to measure the model performance with respect to multiple tasks.

\begin{figure}[]
	\centering
	\begin{subfigure}{0.45\linewidth}
		\centering
		\includegraphics[width=0.95\linewidth]{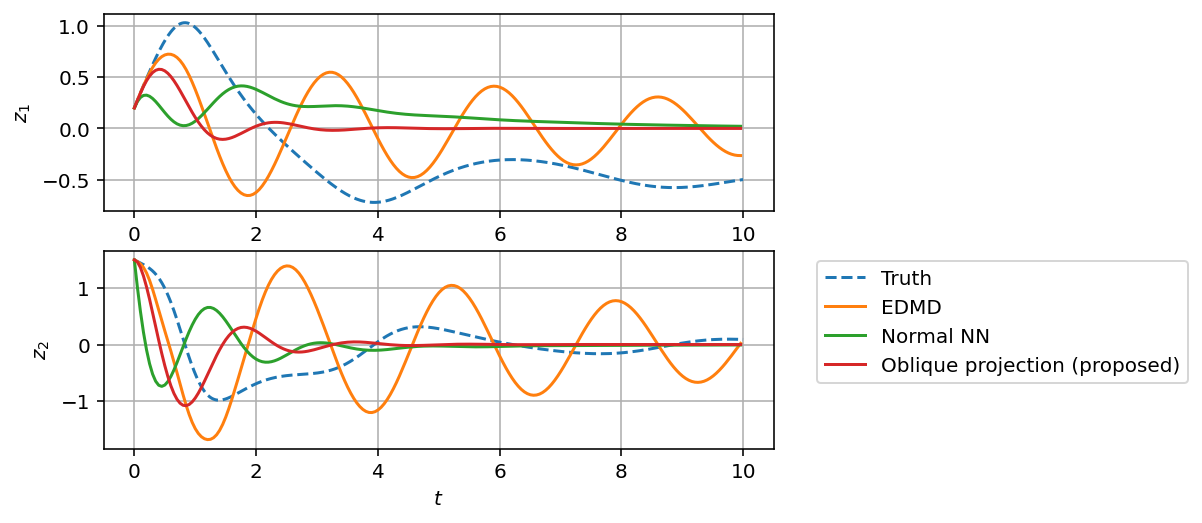}
		\caption{Duffing oscillator.}
		\label{subfig. appendix forecasting duffing}
	\end{subfigure}
	\begin{subfigure}{0.45\linewidth}
		\centering
		\includegraphics[width=0.95\linewidth]{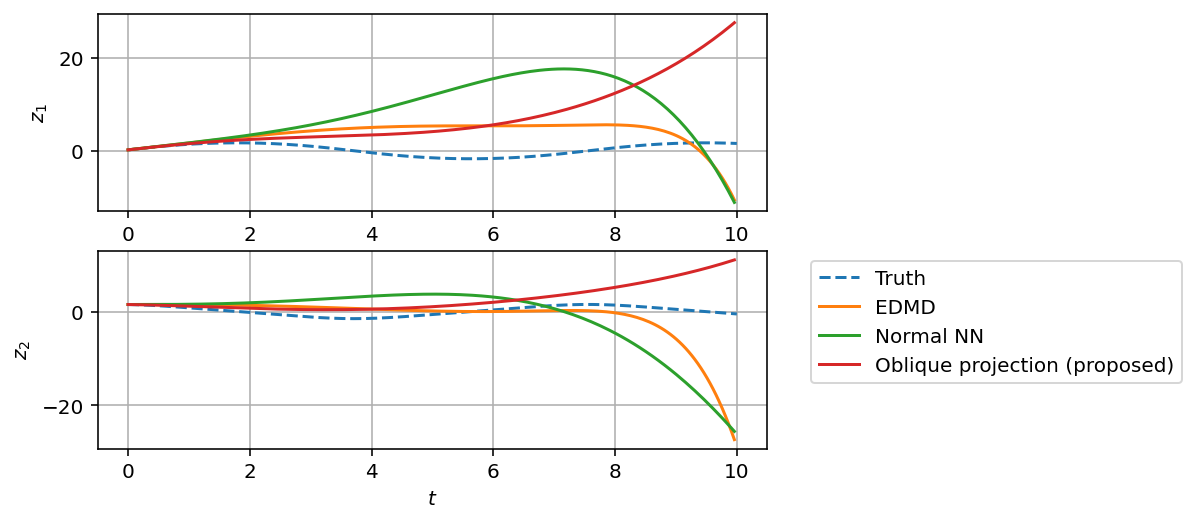}
		\caption{Simple pendulum.}
		\label{subfig. appendix forecasting pendulum}
	\end{subfigure}
	\begin{subfigure}{0.45\linewidth}
		\centering
		\includegraphics[width=0.95\linewidth]{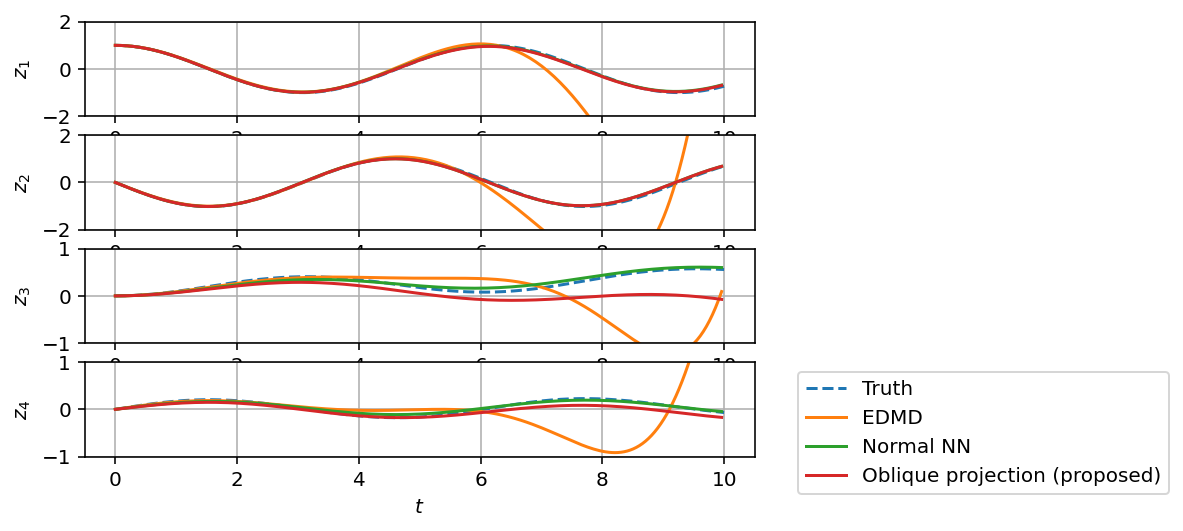}
		\caption{RTAC.}
		\label{subfig. appendix forecasting RTAC}
	\end{subfigure}
	\caption{Forecasting by Algorithm \ref{alg. forecasting}.}
	\label{fig. forecasting errors}
\end{figure}

	Given a state trajectory $\chi_0,\chi_1,\cdots$ of the true dynamics, the invariance proximity of the model can be evaluated by simply
	evolving a linear embedding model (2.3) as in Algorithm \ref{alg. forecasting}.
	\begin{algorithm}[H]
		\caption{Forecasting by simply evolving the linear embedding model}
		\label{alg. forecasting}
		\begin{algorithmic}
			\Require Initial condition $\chi_0\in \mathcal{X}$ and input sequence $(u_0,u_1,\cdots)\in l(\mathcal{U})$
			\Ensure Predicted embedded states $\bm{g}^\text{est}_0, \bm{g}^\text{est}_1,\cdots$
			\State $\bm{g}^\text{est}_0\leftarrow 
			\left[
				\begin{array}{c}
					\chi_0
				\\
					g_{n+1}(\chi_0)
				\\
					\vdots 
				\\
					g_{N_x}(\chi_0)
				\end{array}
			\right]	
			$
			\For {$k=1,2,\cdots$}
			\State 
			$
			\bm{g}^\text{est}_k \leftarrow 
			\bm{A}\bm{g}^\text{est}_{k-1} + \bm{B}u_{k-1}
			$
			\EndFor 
		\end{algorithmic}
	\end{algorithm}
	Figure \ref{fig. forecasting errors} shows the first $n$ components of the predicted embedded states $\bm{g}_k^\text{est}$ of Algorithm \ref{alg. forecasting}.
	It is apparent from the figure that the predictions by Algorithm \ref{alg. forecasting} differ from the actual state prediction results (Figs. 5, 7, and 9), where the feature maps are re-evaluated at every time step instead of simply evolving the model from the initial condition. 
    Especially, the EDMD and the proposed models show poor forecasting in Fig. \ref{subfig. appendix forecasting duffing} whereas both have satisfactory performance in the original state prediction task as in Fig. 6.
	Similarly, Algorithm \ref{alg. forecasting} yields unstable forecasting for the simple pendulum example (Fig. \ref{subfig. appendix forecasting pendulum}), whereas only the EDMD model results in the same behavior in the original prediction task (Fig. 8).

	Also, the results of the control applications in Section 5 show a diverse range of model performance depending on which controller design and linear embedding model are adopted (Figs. 6, 8, and 10). Therefore, it is unable to measure or infer the performance of individual methods w.r.t. multiple control tasks targeting multiple dynamics by looking at the invariance proximity computed by Algorithm \ref{alg. forecasting}.
	While the invariance proximity is useful to analyze the model accuracy analytically, it is difficult and not stable in practice to measure the generalizability, i.e., the model performance w.r.t. multiple tasks, by Algorithm \ref{alg. forecasting} or only considering the modeling error (2.20).

\section{Comparison of State Predictive Accuracy with Nonlinear Models}
\begin{figure}
    \centering
    \begin{subfigure}{0.45\linewidth}
        \centering
        \includegraphics[width=0.95\linewidth]{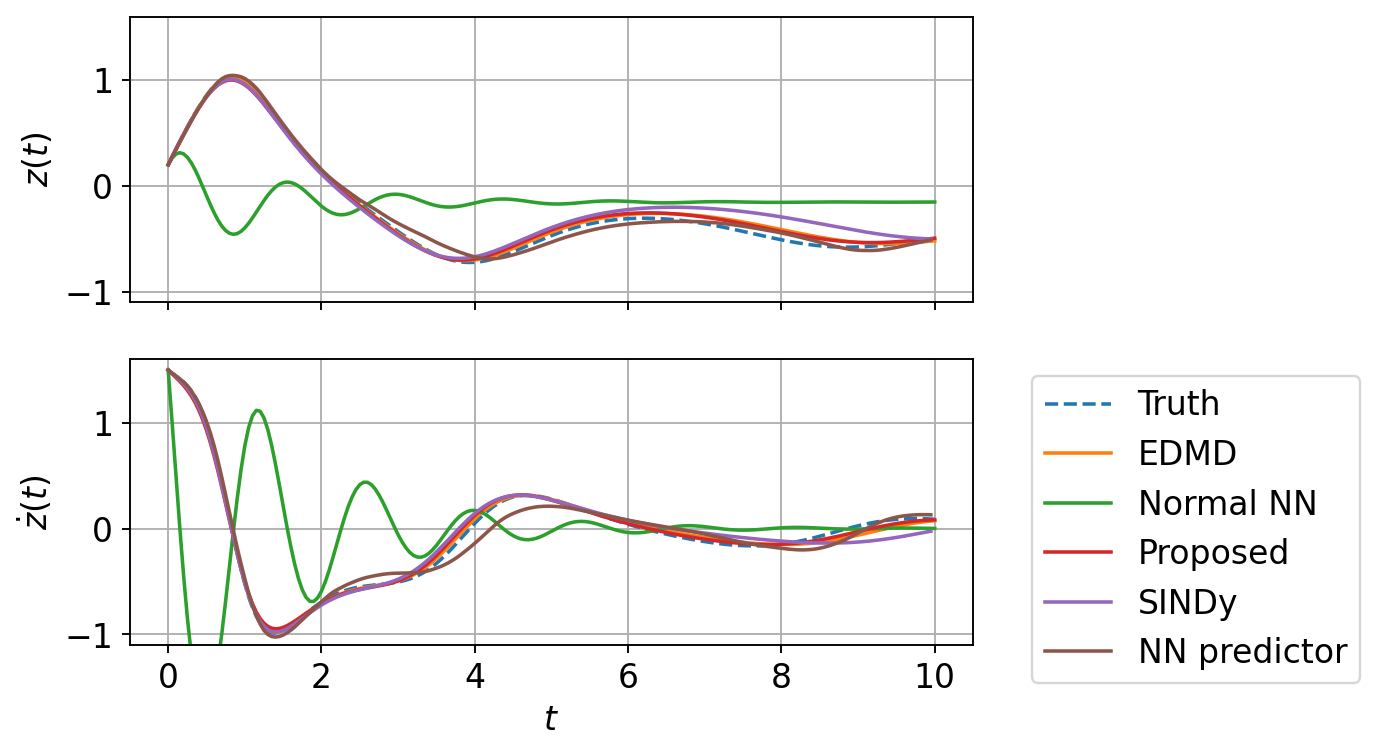}
        \caption{Duffing oscillator.}
        \label{subfig. state predictive accuracy duffing}
    \end{subfigure}
    \begin{subfigure}{0.45\linewidth}
        \centering
        \includegraphics[width=0.95\linewidth]{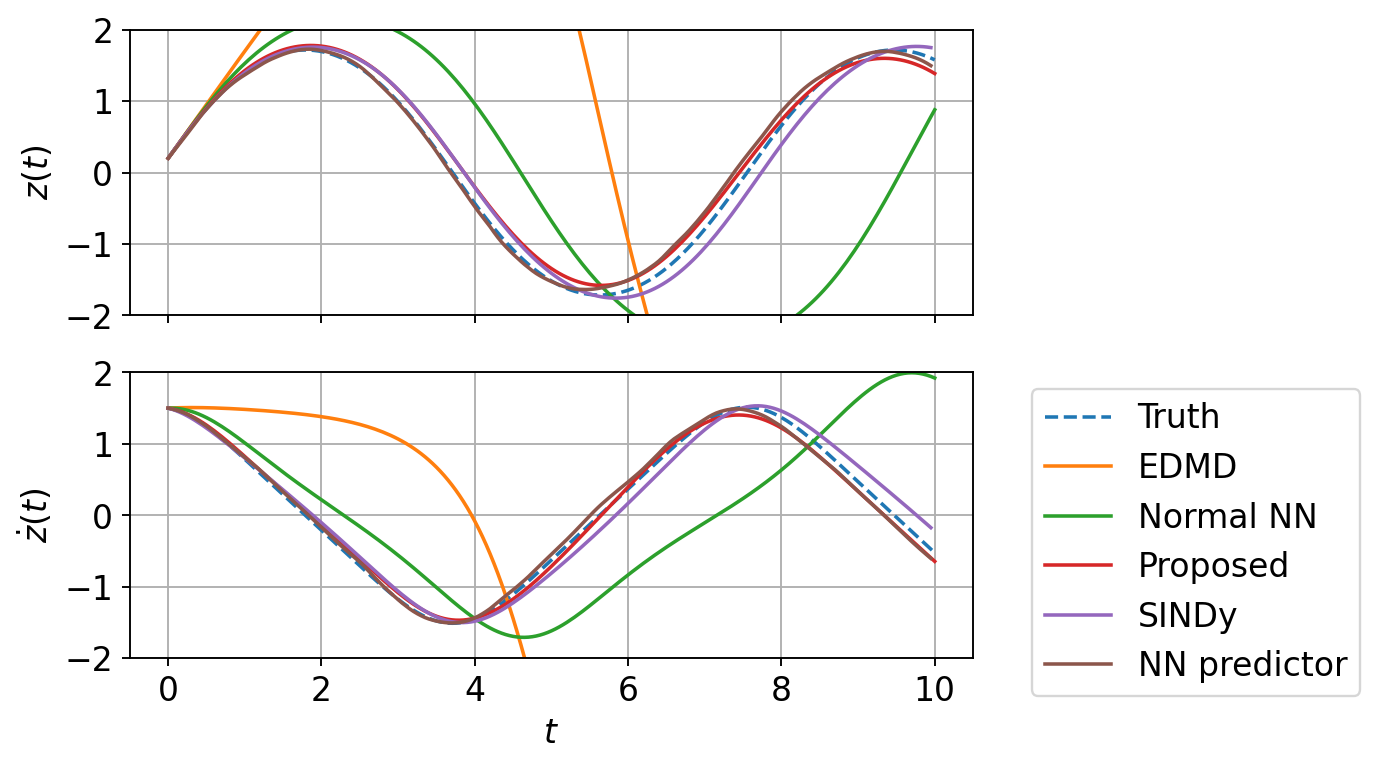}
        \caption{Simple pendulum.}
        \label{subfig. state predictive accuracy pendulum}
    \end{subfigure}
    \begin{subfigure}{0.45\linewidth}
        \centering
        \includegraphics[width=0.95\linewidth]{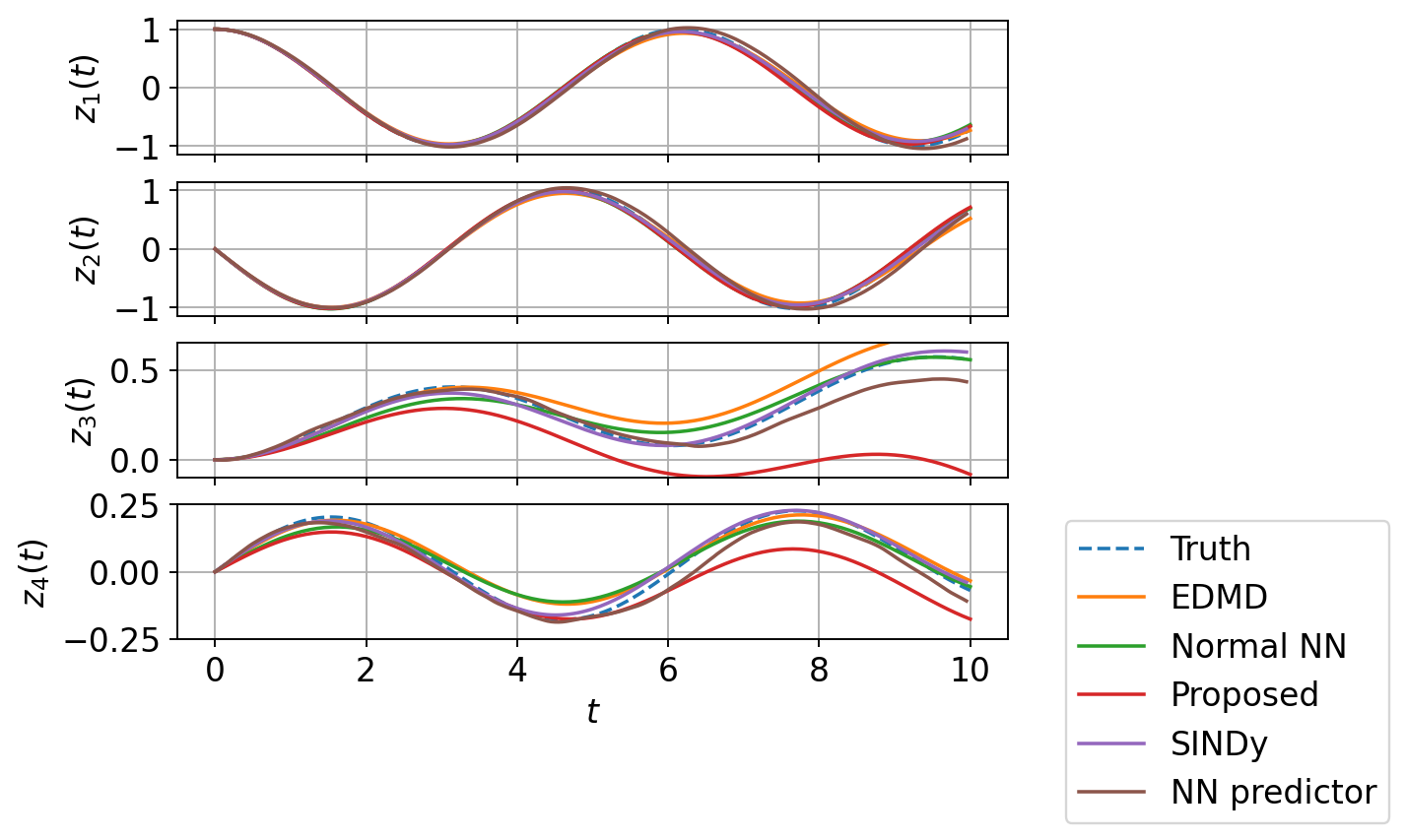}
        \caption{RTAC.}
        \label{subfig. state predictive accuracy RTAC}
    \end{subfigure}
    \caption{Comparison of state predictive accuracy.}
    \label{fig. comparison of state predictive accuracy with nonlinear models}
\end{figure}

In this section, a more detailed evaluation of the state predictive accuracy of linear embedding models is provided. 
Specifically, the state predictions of the linear embedding models trained in Section 5 are compared with nonlinear predictors. 
The linear embedding models considered in this paper are concerned with the state predictive accuracy by the second term $\lambda_2 \mathcal{E}^2_{\text{state}}(\chi_i,u_i)$ of the loss functions and the invariance proximity as an approximation of the Koopman operator, which corresponds to their first term $\lambda_1 \mathcal{E}^2(\chi_i,u_i)$. 
On the other hand, general nonlinear predictor models are likely to possess higher accuracy if they are trained with loss functions that only minimize the state prediction error. 
To compare the state predictive accuracy of different types of models, two nonlinear models are considered in this section. 
Sparse Identification of Nonlinear Dynamical Systems (SINDy)\cite{SINDy,SINDYc} is a modeling method to learn nonlinear predictors that identify functional forms of the target dynamics based on sparse regression. The resulting model equations are linear combinations of candidate dictionaries of feature maps that are pre-specified by the user. 
As another nonlinear model, a neural network-based predictor is considered, which is a feedforward neural network $F_{NN}:\mathcal{X}\times \mathcal{U}\rightarrow \mathcal{X}$ that directly predicts the state at the next time step given a pair of state and input.
Both nonlinear models are trained on exactly the same data set as the linear embedding models in Section 5.
For SINDy models, polynomial features are used and a feedforward neural network consisting of two hidden layers with the Rectified Linear Unit (ReLU) activation is used to train the neural netowrk predictor $F_{NN}$.
\\ \ \\
Figure \ref{fig. comparison of state predictive accuracy with nonlinear models} shows the comparison of the state predictive accuracy of both the linear embedding models from Section 5 and the two nonlinear predictors. 
Both the SINDy model and the neural network predictor show good state predictive accuracy in all three examples, which is considered due to the fact that they are general nonlinear predictors trained to only minimize the state prediction error.
Also, the proposed model is the only linear embedding model that also performs well in all examples.
Normal NN model fails the prediction task for the Duffing oscillator and the simple pendulum (Figs. \ref{subfig. state predictive accuracy duffing} and \ref{subfig. state predictive accuracy pendulum}) and the EDMD model results in an unstable prediction for the simple pendulum (Fig. \ref{subfig. state predictive accuracy pendulum}).
It is observed that the prediction of the proposed model deviates from the truth in the third component $z_3(t)$ in the example of RTAC (Fig. \ref{subfig. state predictive accuracy RTAC}). 
It is considered as a result of the trade-off relation between the robustness or reliability of the training procedure, which is stated in Section \ref{appendix sensitivity analysis} and Fig. \ref{fig. RTAC uncertainty quantification}, and the possible deterioration of model performance. 
In this example, a strong constraint (5.13) is imposed on the test functions and this is considered to have resulted in the less expressive model whereas the model training becomes robust while keeping satisfactory model performance for other controller design tasks. 

\bibliographystyle{siamplain}
\bibliography{supplementary_material_arxiv}